%% file: journal_version.tex
\documentclass[mnsc,nonblindrev]{informs3_mod}

\OneAndAHalfSpacedXI

\newcommand{\Ascr}{\ensuremath{\mathcal A}}

\newcommand{\Dscr}{\ensuremath{\mathcal D}}
\newcommand{\Escr}{\ensuremath{\mathcal E}}
\newcommand{\Fscr}{\ensuremath{\mathcal F}}

\newcommand{\Hscr}{\ensuremath{\mathcal H}}

\newcommand{\Mscr}{\ensuremath{\mathcal M}}
\newcommand{\Nscr}{\ensuremath{\mathcal N}}

\newcommand{\Pscr}{\ensuremath{\mathcal P}}

\newcommand{\Vscr}{\ensuremath{\mathcal V}}

\newcommand{\Xscr}{\ensuremath{\mathcal X}}
\newcommand{\Yscr}{\ensuremath{\mathcal Y}}

\def\abs#1{\vert #1 \vert}
\def\E{\mathbb{E}}
\def\R{\mathbb{R}}
\def\P{\mathbb{P}}
\def\Q{\mathbb{Q}}
\usepackage{threeparttable,booktabs,caption, graphicx, subfig, makecell, multicol, multirow}
\usepackage{soul}
\usepackage{bm}

\usepackage[shortlabels]{enumitem}

\usepackage{natbib}
 \bibpunct[, ]{(}{)}{,}{a}{}{,}%
 \def\bibfont{\small}%
\TheoremsNumberedThrough     
\ECRepeatTheorems

\EquationsNumberedThrough    

\MANUSCRIPTNO{}

\def\gi#1{\textcolor{black}{#1}}

\def\PDRO{\texttt{P-DRO}}
\def\PERM{\texttt{P-ERM}}

\begin{document}


\RUNAUTHOR{Garud, Henry, and Tianyu}

\RUNTITLE{Hedging Complexity in Generalization via a Parametric DRO Framework}

\TITLE{Hedging Complexity in Generalization via a Parametric Distributionally Robust Optimization Framework}

\ARTICLEAUTHORS{%
\AUTHOR{Garud Iyengar, Henry Lam, Tianyu Wang}
\AFF{Department of Industrial Engineering and Operations Research, Columbia University, New York, NY 10027,
\EMAIL{garud@ieor.columbia.edu, henry.lam@columbia.edu, tianyu.wang@columbia.edu}} 
} 

\ABSTRACT{%
Empirical risk minimization (ERM) and distributionally robust optimization
(DRO) are popular approaches for solving stochastic optimization
problems that appear in operations management and machine learning. Existing generalization error bounds for these methods depend on either the complexity of the cost function or dimension of the random perturbations. Consequently, the performance of these methods can be poor for high-dimensional problems with complex objective functions.
We propose a simple approach in which the distribution of random perturbations is approximated using a
parametric family of distributions. This mitigates both sources of
complexity; however, it introduces a model misspecification error. We show that this new source of error can be controlled by suitable DRO formulations. Our proposed parametric DRO approach has significantly improved generalization bounds over existing ERM and DRO methods and parametric ERM for a wide variety of settings. Our method is particularly effective under distribution shifts and works broadly in contextual optimization. We also illustrate the superior performance of our approach on both synthetic and real-data portfolio optimization and regression tasks.  
}%



\maketitle
%


\section{Introduction}
The goal of data-driven stochastic optimization is to solve
\begin{equation}\label{eq:true-obj}
    \min_{x \in \Xscr}\left\{Z(x):=\E_{\xi \sim \P^*}[h(x;\xi)]\right\},
\end{equation}
where $x \in \Xscr$ is the decision, 
$\xi$ is a random perturbation in the sample space $\Xi$ distributed according to $\P^\ast$, and
$h: \Xscr \times \Xi \to \R$ 
is the cost function. 
Typically, 
$\P^\ast$ is unknown and one only has access to
i.i.d. samples $\hat{\xi}_i \sim \P^*$, $i = 1, \ldots, n$. 
Problems 
of this nature 
arise 
in many different settings
from machine learning to 
decision making~\citep{shapiro2014lectures,birge2011introduction}. 

A standard approach to approximately
solve~\eqref{eq:true-obj} is Empirical Risk Minimization~(ERM), where one replaces the 
unknown $\P^\ast$ with the empirical measure 
$\hat{\P}_n := \frac{1}{n}\sum_{i = 1}^n\delta_{\hat{\xi}_i}$, 
leading to the 
problem~\citep{hastie2009elements}    
\begin{equation}\label{eq:EO-obj}
\min_{x \in \Xscr}\left\{\hat{Z}^{ERM}(x):= \E_{\hat{\P}_n}[h(x;\xi)] =
  \frac{1}{n}\sum_{i = 1}^n h(x;\hat{\xi}_i)\right\}.
\end{equation}
ERM 
is conceptually natural;   
see \cite{vapnik1999nature,shalev2014understanding,shapiro2014lectures} for 
comprehensive surveys on its statistical 
guarantees.
A second approach that 
is gaining 
popularity in recent years 
is Distributionally Robust
Optimization (DRO), 
where 
the unknown $\P^\ast$ is replaced by the worst-case distribution over a
so-called ambiguity set $\Ascr$, giving rise to the problem
\begin{equation}\label{eq:DRO-obj}
\min_{x \in \Xscr}\left\{\hat{Z}^{DRO}(x):= \max_{\P \in \Ascr}\E_{\P}[h(x;\xi)]\right\}.
\end{equation}
Here, $\Ascr$ is constructed using the data and, 
if 
$\Ascr$ is 
chosen so that, at least intuitively speaking, it
covers the 
unknown
$\P^\ast$ with high
confidence, then \eqref{eq:DRO-obj} outputs a solution with a worst-case
performance guarantee. 
In order to 
guarantee a statistically consistent solution, it is common to set $\Ascr=
\{\P | d(\P, \hat{\P}_n)\leq \varepsilon\}$ for some statistical distance
$d$, and  $\varepsilon$ shrinking to zero as $n$ increases. This approach
has been studied with $d$ set to the Wasserstein
distance~\citep{esfahani2018data},
$f$-divergence~\citep{ben2013robust}, kernel
distance~\citep{staib2019distributionally} and other variants. Compared to ERM, DRO offers a worst-case protection against model shifts and is especially useful to handle problems with only partial information. Moreover, its solutions possess different statistical behaviors from ERM that are advantageous for certain situations (as we will review in the sequel). DRO
has been successfully applied to 
many applications in machine learning and statistics including 
linear
regression~\citep{chen2018distributionally,shafieezadeh2019regularization}, neural
networks~\citep{Sagawa2020Distributionally,duchi2020distributionally} and
transfer learning~\citep{volpi2018generalizing}; in operations
including newsvendor problems~\citep{hanasusanto2015distributionally,natarajan2018asymmetry},
portfolio 
optimization~\citep{blanchet2022distributionally,doan2015robustness} and energy
systems~\citep{wang2018risk}. See, e.g., the  
survey~\citep{rahimian2019distributionally} for a recent overview.

Despite their wide usages, both ERM and DRO 
could have poor performances in high-dimensional complex-structured problems. This can be explained via their generalization errors. To this end, let $\hat{x}$ denote an approximate solution
for~\eqref{eq:true-obj}. The generalization error of $\hat x$, measured by the excess risk, or equivalently the optimality gap or regret, is defined as:
\begin{equation}
\mathcal{E}(\hat{x}) :=
Z(\hat{x}) - Z(x^\ast),\label{generalization error}
\end{equation}
where 
$x^\ast \in \argmin_{x \in
  \Xscr}Z(x)$ denotes an optimal solution of~\eqref{eq:true-obj}.  
The expression \eqref{generalization error} quantifies the performance of $\hat x$ relative to the oracle best benchmark $x^*$ in terms of the attained \emph{true} objective value. 
In the literature, bounds on $\mathcal{E}(\hat{x})$ are 
typically of the form
\begin{equation}\label{eq:general-bd}
    \mathcal{E}(\hat{x}) 
\leq \frac{B}{n^{\alpha}},
\end{equation}
where $B$, $\alpha$ 
are positive quantities dependent on the method used to obtain~$\hat{x}$. For 
ERM, $\alpha = \frac{1}{2}$ and $B$ depends on the complexity $\text{Comp}({\Hscr})$ of the
function class $\Hscr = \{h(x,\cdot):x\in\mathcal X\}$. 
This complexity measures how rich is the function class and the proneness to overfitting, exemplified by
well-known measures such as the Vapnik–Chervonenki (VC) dimension
\citep{vapnik1999nature,bartlett2002rademacher} and 
local Rademacher complexity \citep{bartlett2005local, xu2020towards}. 
On the other hand, approaches to analyze the generalization error of DRO 
generally fall into 
two 
categories.
The first approach 
views DRO as a  regularization of ERM, where the regularizer depends on the choice of
$d$~\citep{duchi2019variance,gotoh2021calibration,lam2019recovering,blanchet2019robust,gao2022finite,gupta2019near,blanchet2022confidence},
which results in $(\alpha, B)$ similar to ERM.
In the second approach, 
the ambiguity set $\Ascr$ is
constructed as a non-parametric confidence region for $\P^\ast$,  resulting
in a
worst-case performance bound on 
$\hat{x}$~\citep{esfahani2018data,bertsimas2018robust,delage2010distributionally,wiesemann2014distributionally,goh2010distributionally}.  
This can be 
converted
into a bound for $\mathcal{E}(\hat{x})$ where $B$ depends only on the true loss
$h(x^\ast,\cdot)$ instead of the complexity of the hypothesis
class~\citep{zeng2021generalization}, but now $\alpha$  
typically degrades as $1/D_{\xi}$ where $D_{\xi}$
denotes the dimension of the randomness $\xi$. 
In other words, in all the existing bounds for ERM and DRO, the
generalization error $\mathcal{E}(\hat{x})$ depends on \emph{either the
complexity of the cost function class or the distribution dimension}. 
Moreover, these bounds are in a sense tight. 
Hence, for a high-dimensional problem with a complex cost function, both ERM and DRO 
could incur poor performances.
Given the above challenge, our goal is to create a new approach with better guarantees for high-dimensional complex problems, 
by removing the dependence of %
the generalization error
bound 
on \emph{both} the complexity $\text{Comp}(\Hscr)$ in $B$ 
and the distribution dimension~$D_{\xi}$ in $\alpha$ from~\eqref{eq:general-bd}. 
Our approach is conceptually simple:
We center the ambiguity set $\Ascr$ 
at a suitable parametric distribution, instead of the empirical
distribution~$\hat{\P}_n$. Therefore,  
we call our approach
\emph{parametric DRO} (\PDRO). %
Intuitively speaking, we follow the second analysis 
approach for DRO described above  to
obtain
$B$ depending only on $h(x^*,\cdot)$ instead of the cost function complexity. Moreover, 
because we center $\Ascr$ at a parametric distribution, 
$\alpha $ no longer depends on $D_{\xi}$. 
These thus remove both $\text{Comp}({\Hscr})$ and $D_{\xi}$. However, 
they come
with the price of the model misspecification error 
associated with the chosen family of 
parametric distributions. 
Our 
main insight 
is that by choosing the size $\varepsilon$ of the ambiguity set $\Ascr$
appropriately, the worst-case nature of \PDRO\ can exert 
control on the impact of model misspecification, and ultimately exhibit a
desirable trade-off between this model error and the simultaneous removal of
complexity and dimension dependence. 

Our framework broadly generalizes to two important settings. First is distribution shifts, i.e., when training and testing data are different.
While previous literature has argued 
that
DRO 
is able to protect
against unexpected distribution shifts, the
arguments are based on a worst-case bound applied 
to the 
objective,
i.e.,  $Z(\hat{x})\leq \hat{Z}^{DRO}(\hat{x})$ with high probability. This protection guarantee 
does not directly
indicate whether the excess risk performance of the  
DRO
solution 
is superior to 
other possibilities. In contrast, we show how \PDRO\ solutions exhibit a better trade-off among distribution shift, parametric model misspecification errors, and complexity and dimension compared to other alternatives in terms of the excess risk of the shifted testing data.
Second, we generalize \PDRO\ 
to the contextual
optimization 
setting,
where the distribution of $\xi$ depends on observable exogenous features. In this setting, \PDRO\ is arguably even more dominant since non-parametric methods are arguably difficult to implement: As the decision is now a map from the feature, using the empirical distribution alone is unable to generalize the decision to feature values that are not observed before in a continuous feature space, and this necessitates the use of more sophisticated kernel, tree-based or other smoothing approaches. 





The rest of this paper is organized as follows. 
Section~\ref{sec:background} introduces 
existing 
generalization error
bounds of ERM and DRO, and uses them to motivate our investigation.
Section~\ref{sec:main} 
presents \PDRO\ and its foundational theory. Section~\ref{subsec:sample} incorporates computation errors incurred by the Monte Carlo sampling into the generalization bounds. Section~\ref{sec:extensions} extends our results to  
distribution shifts settings and 
contextual optimization. Section~\ref{sec:discussion} discusses the overall trade-off among different errors and compares \PDRO\ with other alternatives in detail. 
Section~\ref{sec:numerics} presents
numerical results on both synthetic and real data examples. Additional discussions and experimental results, and proofs of all theorems are deferred to the Appendix. 
\section{Background} \label{sec:background}
We briefly discuss how existing bounds for the excess risk 
$\mathcal{E}(\hat x)$ for
ERM and DRO are constructed, and 
set the stage for our new \PDRO\ bounds. 
Let $\hat Z(\cdot)$
denote the estimated objective function via a particular
approximation scheme, e.g., ERM uses 
the sample average objective $\hat Z^{ERM}(\cdot)$ 
in~\eqref{eq:EO-obj} and DRO uses the 
worst-case objective $\hat Z^{DRO}(\cdot)$ in~\eqref{eq:DRO-obj}, $\hat{x}$ denote the
corresponding approximate solution, and $x^{\ast}$ denote the optimal
solution of the stochastic optimization problem~\eqref{eq:true-obj}. Then
we have that 
\begin{align}
  \mathcal{E}(\hat{x})
  &= [Z(\hat{x}) - \hat{Z}(\hat{x})] +
    [\hat{Z}(\hat{x}) - \hat{Z}(x^*)] + [\hat{Z}(x^*) -
    Z(x^*)] \nonumber\\ 
  &\leq [Z(\hat{x}) - \hat{Z}(\hat{x})] +
    [\hat{Z}(x^*) - Z(x^*)],  \label{eq:error-decomp} 
\end{align}
where~\eqref{eq:error-decomp} 
follows from the definition
$\hat{x} \in \argmin_{x \in \Xscr} \hat{Z}(x)$. 
Next, we bound the two terms in 
\eqref{eq:error-decomp}, 
and the optimal 
overall generalization
bound relies on a balance between these two terms. 

\paragraph{ERM.} 
Here, $\hat{Z}(\cdot)$ is taken as $\hat Z^{ERM}(\cdot)$ and $\hat x$ is the ERM solution.
  A bound for the   
second term $\hat Z(x^*) - Z(x^*)$ in \eqref{eq:error-decomp}  
follows from standard bounds for the difference between expectation and the
sample
mean. On the other hand, the first term $Z(\hat{x}) - \hat{Z}(\hat{x})$ in
\eqref{eq:error-decomp} depends on 
the (random) solution $\hat x$, and is bounded by its supremum
$\sup_{x \in \Xscr}|Z(x) - \hat{Z}(x)|$ (or a localized version). Using
tools from empirical process theory~\citep{van1996weak}, 
one can obtain 
a bound of the form  $\mathcal{E}(\hat x) \leq O\big(\sqrt{\frac{M
    Z(x^*)\text{Comp}(\Hscr)\log
    n}{n}}\big)$~\citep{vapnik1999nature,boucheron2005theory}, 
where $M=\sup_{x\in 
  \Xscr}\|h(x;\cdot)\|_{\infty}$ and $\text{Comp}(\Hscr)$ is some
complexity measure, such as the VC dimension, metric entropy, or the  Rademacher complexity, of the
hypothesis class $\Hscr=\{h(x;\cdot)|x\in\mathcal X\}$. We use the metric entropy to represent $\text{Comp}(\Hscr)$ throughout the paper unless otherwise stated. Note that some regularized-ERM approaches can attain better generalization performance than the standard ERM and we defer this discussion to Section~\ref{sec:discussion}.

\paragraph{DRO bound using the regularization perspective.} 
Here $\hat Z(\cdot)$ is taken as $\hat Z^{DRO}(\cdot)$ and $\hat x$ is the DRO solution. Recent works \citep{lam2016robust,duchi2019variance,gao2022wasserstein} show
that, for a small enough 
$\varepsilon$, 
\begin{equation}\label{eq:dro-reg-erm}
\hat{Z}^{DRO}(x) = \hat{Z}^{ERM}(x) + \Vscr_d(x)\sqrt{\varepsilon} +
O(\varepsilon),\ \forall x \in \Xscr. 
\end{equation} 
Here $\Vscr_d(x)$ is a variability measure of the cost function $h$ that
depends on
the statistical distance $d$ used to define the ambiguity
set. For example, $\Vscr_d(x)$ is the Lipschitz norm of $h(x;\cdot)$
when $d$ is 1-Wasserstein distance
\citep{blanchet2019robust,gao2022wasserstein}, gradient norm when $d$ is $p$-Wasserstein distance \citep{gao2022wasserstein} and 
$\sqrt{\text{Var}_{\P^*}[h(x,\xi)]}$ when $d$ is an $f$-divergence  
\citep{lam2016robust,lam2018sensitivity,duchi2019variance,duchi2021statistics,gotoh2018robust,gotoh2021calibration,bennouna2021learning}.
The expansion~\eqref{eq:dro-reg-erm} can be used to bound the second term
$\hat{Z}(x^*) - 
Z(x^*)$ in \eqref{eq:error-decomp} by 
combining with ERM bounds. Moreover, 
for appropriately chosen
$\varepsilon$ 
(
depending on the
hypothesis class
complexity $\text{Comp}(\Hscr)$; see Examples~\ref{ex:w1-dro-eps} and~\ref{ex:chi-dro-eps} in Appendix~\ref{app:erm-dro-derive}), \eqref{eq:dro-reg-erm} 
together with the empirical Bernstein inequality~\citep{maurer2009empirical}
implies the bound 
\begin{equation}\label{eq:true-empirical-dro}
      Z(x) \leq \hat{Z}^{DRO}(x) + O\left(\frac{1}{n}\right),\ \forall x \in
  \Xscr,
\end{equation}
which can then be used to 
bound the first term $Z(\hat{x}) - \hat{Z}(\hat{x})$ \citep{duchi2019variance,gao2022finite}.
Putting these together one arrives at the 
bound 
$\mathcal E(\hat x) \leq
O\left(\Vscr_d(x^*)\sqrt{\frac{\text{Comp}(\Hscr)}{n}}\right)$. Compared 
to the bound for
ERM, in 
the DRO bound the term $MZ(x^*)$  
is replaced by 
$\Vscr_d(x^*)$; 
however, 
both bounds involve $\text{Comp}(\mathcal H)$.  

\paragraph{DRO bound from a robustness perspective.}
Again $\hat Z(\cdot)$ is taken as $\hat Z^{DRO}(\cdot)$ and $\hat x$ is the DRO solution. 
Suppose
$\varepsilon$
is chosen large enough 
so
that
\begin{equation} 
\P[d(\P^*,\hat{\P}_n) \leq \varepsilon] \geq 1-\delta, \label{robust bound}
\end{equation}
i.e., $\mathcal \Ascr$ covers the ground-truth $\P^*$ with 
probability $1-\delta$. Note that this choice of
$\varepsilon$ does not depend on the cost function $h$. The first term $Z(\hat{x}) -
\hat{Z}(\hat{x})$ in \eqref{eq:error-decomp} is non-positive with
probability at least $1-\delta$ 
\citep{ben2013robust,bertsimas2018robust}, 
while the second term $\hat{Z}(x^*) - Z(x^*)$ depends on $\varepsilon$ and
$h(x^\ast;\xi)$, but not 
on 
$\text{Comp}(\Hscr)$~\citep{zeng2021generalization}. 
However, for \eqref{robust bound} to hold,
we typically need to choose $\varepsilon= 
O(n^{-1/D_{\xi}})$, which in turn 
degrades the bound for the
second term. This is the case 
for the
Wasserstein distance due to its  concentration for the empirical distribution (\cite{fournier2015rate}). This is also the case for $f$-divergences since they are defined only for 
absolutely 
continuous distributions~\citep{jiang2018risk} and so the associated DRO ball center requires smoothing
$\hat\P_n$ 
appropriately, e.g., the kernel density
estimator \citep{zhao2015data,chen2022distributionally} requires $\varepsilon =
O((nh_n^{D_{\xi}})^{-\frac{1}{2}}\vee h_n^2)$ to ensure $\P^* \in \Ascr$ with high confidence. 
The only exception is the
maximum mean discrepancy 
(MMD) 
where one can set 
$\varepsilon = O(1/\sqrt n)$, but 
in this case, to bound the 
second term 
one needs to assume that $h(x^\ast, \cdot)$ belongs to  a reproducing kernel Hilbert
space (RKHS)~\citep{zeng2021generalization}, or otherwise the resulting bound degrades again.

\paragraph{New bounds based on parametric distributions.} The bounds
discussed above are shown in 
the 
rows marked ``standard'' in Table~\ref{tab:general-table} for the basic setup. 
As noted above, these bounds 
either depend on the hypothesis class complexity $\text{Comp}(\mathcal H)$
or the distributional dimension $D_{\xi}$.
Our approach \PDRO\ 
uses 
an appropriately chosen
parametric model as the center of ambiguity set $\Ascr$. This choice
results in a bound 
replacing both
$\text{Comp}(\mathcal H)$ and $D_{\xi}$ with a potentially much smaller
\emph{parametric complexity} $\text{Comp}(\Theta)$. However, in doing so,
we incur a model misspecification term $\mathcal E_{apx}$. The trade-off
between $\text{Comp}(\Theta)$ and $\mathcal E_{apx}$ is shown in the
row marked ``parametric'' in  Table~\ref{tab:general-table}. 
When 
sample size $n$ is 
not too large, 
the gain in $\text{Comp}(\Theta)$ over
$\text{Comp}(\mathcal H)$ 
can
be significant enough 
that outweighs 
the increase in errors due to 
$\mathcal E_{apx}$. Moreover, if we simply apply the same parametric model
in ERM, we obtain a bound that depends less desirably on $\mathcal
E_{apx}$ by having an additional $M \Escr_{apx}^{\frac{3}{4}}$ error (shown at the left entry of the second row). Overall, when a problem has large $\text{Comp}(\Hscr)$ and $D_\xi$, but small $\text{Comp}(\Theta)$ and $\Escr_{apx}$, \PDRO\ has a better generalization error bound than the existing approaches for relatively small sample size $n$. In this sense, \PDRO\ provides an efficient mechanism to take advantage of parametric distributional structures.

\begin{table*}[h]
\footnotesize
    \centering
    \begin{tabular}{c|c|cc}
    \toprule
         & ERM & \multicolumn{2}{c}{DRO with metric $d$} \\
         \midrule
        Standard & $\sqrt{\frac{MZ(x^*) {\text{Comp}(\Hscr)}}{n}}$ & (regularization) $\Vscr_d(x^*)\sqrt{\frac{\text{Comp}(\Hscr)}{n}}$  & (robustness) $\Vscr_d(x^*)\cdot\frac{1}{n^{1/D_{\xi}}}$ \\
        \midrule
        Parametric & $\Vscr_d(x^*)\left(\sqrt{\frac{\text{Comp}(\Theta)}{n}}+ \mathcal{E}_{apx}\right) + M\mathcal{E}_{apx}^{\frac{3}{4}}$ &  \multicolumn{2}{c}{$\mathbf{\Vscr}_d(x^*)\left(\sqrt{\frac{\text{Comp}(\Theta)}{n}}+ \mathcal{E}_{apx}\right)$}\\
        \bottomrule
    \end{tabular}
    \caption{Generalization errors of different methods under the basic setup. Each error is represented by a $(1-\delta)$-probability upper bound, where we ignore numerical constants and $\log (1/\delta)$ and $\log n$ terms in the numerator. 
    }
    \label{tab:general-table}
\end{table*}

For interested readers, we provide further details on the existing generalization error bounds that we have discussed earlier in Appendix~\ref{app:erm-dro-derive}.

\section{Parametric-DRO: Main Results}\label{sec:main}
Given i.i.d. sample $\{\hat{\xi}_i\}_{i = 1}^n$ and a class of parametric
distributions $\Pscr_{\Theta}= \{\P_{\theta}: \theta \in \Theta\}$, \PDRO\
solves \eqref{eq:DRO-obj} with the ambiguity set 
\[
\Ascr = \{\P | d(\P,  \hat{\Q})\leq \varepsilon\},
\]
where $\hat{\Q}$ is an appropriately chosen distribution in 
$\Pscr_{\Theta}$.
Parametric-ERM (\PERM) is the special case obtained by setting $\varepsilon
= 0$. 


We consider two main types of metrics $d$.
\begin{enumerate}[(a)]
\item \emph{Integral Probability Metric (IPM).} The IPM $d (\P,
  \Q)$~\citep{muller1997integral} is defined as 
  \begin{equation*}\label{eq:def-IPM}
    d (\P, \Q):= \sup_{\{f: \Vscr_d(f) \leq 1\}}\Big|\E_{\P}[f] - \E_{\Q}[f]\Big|,
  \end{equation*}
  where $\Vscr_d(f)$ is an appropriately defined variability
  measure 
  with $\Vscr_d(\alpha f) = \alpha \Vscr_d(f)$ for $\alpha \geq 0$~\citep{zhao2015data}. Special cases include the 1-Wasserstein distance ($\Vscr_d(f) = \|f\|_{\text{Lip}}$), total
  variation (TV) distance ($\Vscr_d(f) =  2\|f\|_{\infty}$) and MMD
  ($\Vscr_d(f) 
  = \|f\|_{\Hscr}$). For convenience, we also abbreviate $\Vscr_d(x) = \Vscr_d(h(x;\cdot))$ when no confusion arises, which leads to the $\Vscr_d(x^*)$ in Table \ref{tab:general-table} earlier.
\item 
  \emph{$f$-divergence lower bounded by the TV-distance.} Let $\P$ and $\Q$ be two
  distributions and $\P$ is absolutely continuous w.r.t. $\Q$. For a
  convex function $f:[0,\infty) \to (-\infty, \infty]$ such that $f(x)$ is
  finite $\forall x > 0, f(1) = 0$, the $f$-divergence of $\P$ from $\Q$ is
  defined as:
  \[
    d_f(\P, \Q) = \int f\left(\frac{d\P}{d\Q}\right)d\Q =
    \E_{\Q}\left[f\left(\frac{d\P}{d\Q}\right)\right].
  \]
  We consider $f$-divergences that can be lower bounded by the TV-distance as follows
  \begin{equation}\label{eq:general-f-tv-ineq}
    d_{TV}(\P, \Q) \leq C_{f}\sqrt{d_f(\P, \Q)},
  \end{equation}
  where $C_f > 0$ is a constant. The $\chi^2$-divergence ($f(t) = (t - 1)^2/2$), Kullback-Leibler (KL) divergence ($f(t) = t \log t - (t - 1)$) and
  squared Hellinger ($H^2$) distance ($f(t) = (\sqrt{t} - 1)^2$) are examples of $f$-divergences that satisfy
  the lower bound condition. 
\end{enumerate}
In order to analyze \PDRO, we first make the following general
assumption.
\begin{assumption}[Oracle estimator]\label{asp:oracle-param-est}
  Let $\text{Comp}(\Theta)$ be the complexity of $\Pscr_{\Theta}$, and
  $\mathcal{E}_{apx}(\P^*,\mathcal P_\Theta)$ (abbreviated to
  $\mathcal{E}_{apx}$) is a non-negative function 
  such that $\mathcal{E}_{apx}(\P^*,\mathcal P_\Theta)
  = 0$ if $\P^* \in \mathcal{P}_{\Theta}$. Then, for all $\delta \in
  (0,1)$, there exists $\alpha >
  0$
  such that 
  the center 
  $\hat{\Q} \in
  \mathcal{P}_{\Theta}$ 
  of the ambiguity set $\Ascr$ satisfies
  \begin{equation}
    \label{eq:param-est-ineq}
    d(\P^*,\hat{\Q}) \leq \mathcal{E}_{apx}(\P^*, \mathcal{P}_{\Theta}) +
    \left(\frac{\text{Comp}(\Theta)}{n}\right)^{\alpha} \log(1/\delta) =:\Delta(\delta,
    \Theta), 
  \end{equation}
  with probability $1-\delta$.
\end{assumption}

Assumption \ref{asp:oracle-param-est} holds under a wide range of
parametric models and estimation procedures, 
although the detailed verification of $\text{Comp}(\Theta)$ and $\Escr_{apx}$ must be done on a case-by-case basis.
Here we
discuss two important examples.
\begin{example}\label{ex:oracle-estimator-was}
  Suppose
  $d$ is given by the  $1$-Wasserstein distance, the set of parametric distributions
  $\Pscr_{\Theta} = \{\Nscr(\mu, \Sigma)|  \mu \in 
  \R^{D_{\xi}}\}$ with known $\Sigma$, and $\xi \sim \Q^\ast$ with  
  sub-Gaussian marginal distribution 
  with parameter $\sigma$,
  i.e., $\E[\exp(v^{\top}(\xi - \E[\xi]))] \leq \exp\left(\|v\|^2
    \sigma^2/2\right), \forall v \in \R^{D_{\xi}}$.
  Then 
  Assumption~\ref{asp:oracle-param-est} holds for
  $\hat{\Q}=\Nscr(\frac{1}{n}\sum_{i = 
    1}^n\hat{\xi}_i, \Sigma)$, 
  $\mathcal{E}_{apx}
  = W_1(\P^{\ast}, \Q^\ast)$ with $\Q^\ast = \Nscr(\E[\xi],\Sigma)$, 
  $\alpha = \frac{1}{2}$ and $\text{Comp}(\Theta) =
  D_{\xi}\sigma^2$. 

  This result is established as follows. By the triangle inequality, 
  $W_1(\P^*,
  \hat{\Q})\leq W_1(\P^*, \Q^*) + W_1(\Q^*, \hat{\Q})$. Next,  
  bound
  $W_1(\Q^*, \hat{\Q})\leq W_2(\Q^*, \hat{\Q})
  =\sqrt{\sum_{j =1}^{D_{\xi}}\vert \frac{1}{n}\sum_{i = 1}^n (\hat{\xi}_{i})_j -
  \E[\xi]_j|^2}$, where the equality follows from the fact that $\hat{\Q}$
and $\Q^\ast$ are both Gaussian~\citep{dowson1982frechet}. 
Next,
we apply the sub-Gaussian
concentration inequality \citep{wainwright2019high} to all $D_{\xi}$
components and obtain $W_2(\Q^*, \hat{\Q}) \leq \sigma \sqrt{\frac{D_{\xi}
    \log(1/\delta)}{n}}$.
\end{example}
 \begin{example}[Theorem 13 in
   \cite{liang2021well}]\label{ex:oracle-estimator-kl}
   Suppose 
   $d$ is given by the KL-divergence and  the parametric class $\mathcal P_\Theta$ is the class of all
   distributions of the random variable $g_\theta(Z)$, where $Z$ is a
   fixed random variable and 
   $g_\theta$ 
   is given by a feed-forward neural network parametrized by 
   $\theta\in\Theta$. Let 
   \[
     \hat{\theta}_n \in \argmin_{\theta: \theta \in \Theta}\max_{\omega:f_{\omega}\in \Fscr,\atop
       \|f_{\omega}\|_{\infty}\leq B}\left\{\E_{Z}f_{\omega}(g_{\theta}(Z)) -
    \E_{\hat{\P}_n} f_{\omega}(\xi)\right\},
   \]
   denote the Generative Adversarial Network (GAN) estimator with the
   discriminator class $\Fscr=  \{f_{\omega}(x): \R^{D_{\xi}} \to \R\}$,
   realized by a neural network with weight parameter $\omega$.  
   Then Assumption~\ref{asp:oracle-param-est} holds for $\hat{\Q}$ set to  the
   distribution of $g_{\hat{\theta}_n}(Z)$, $\alpha = \frac{1}{2}$, and 
   \begin{align*} 
     \mathcal{E}_{apx} &= \sup_{\theta}\inf_{\omega}\left\|\log
                         \frac{p_{*}}{p_{{\theta}}} - f_{\omega}\right\|_{\infty} + B
                         \inf_{\theta}\left\|\log
                         \frac{p_{\theta}}{p_{*}}\right\|_{\infty}^{\frac{1}{2}},
     \\ 
     \text{Comp}(\Theta) &=\texttt{Pdim}(\Fscr), 
   \end{align*}
   where $p_*$ (resp. $p_{{\theta}}$) denotes the density of $\P^\ast$
   (resp. $g_{\theta}(Z)$), 
   and $\texttt{Pdim}(\Fscr)$ is the pseudo dimension of $\Fscr$. Here, 
   $\Escr_{apx}$ reflects the expressiveness of the generator 
   and
   $\text{Comp}(\Theta)$ describes the statistical complexity of the
   discriminator. 
 \end{example}
 
 Note that $\alpha$ is dimension-independent in both
 examples above, and this is also 
 generally the case for other interesting metrics; see Appendix~\ref{app:asp-param-estimator} for more examples that satisfy Assumption~\ref{asp:oracle-param-est}. These examples include situations where: 1) $\hat \Q$ is estimated through GAN and $d$ is the $H^2$-distance, and 2) $\hat \Q$ is estimated through a Gaussian mixture model and $d$ is the 1-Wasserstein distance. 
 On the other hand, $\text{Comp}(\Theta)$ 
 scales with the number of parameters. 
 Connecting to operational practice, parametric modeling has been widely used, e.g., normal distributions in assets returns \citep{demiguel2009optimal} and exponential distributions in product demands \citep{liyanage2005practical}.
 Classical maximum likelihood and moment methods 
 possess
 finite-sample guarantees 
 \citep{spokoiny2012parametric,boucheron2013concentration},
 and so are modern generative models such as GAN
 \citep{xu2019modeling,goodfellow2020generative}, 
 e.g., Example~\ref{ex:oracle-estimator-kl} above and
\cite{zhang2017discrimination,liang2021well}. These all lead to the satisfaction of Assumption~\ref{asp:oracle-param-est} under a wide array of settings.
Our first main result is as follows. 
\begin{theorem}[Generalization bounds for \PDRO]\label{thm:general-param-dro-ipm}
  Let $x^{\PDRO}$ denote the solution to $\PDRO$. 
  Suppose Assumption~\ref{asp:oracle-param-est} holds and the size of the
  ambiguity set $\varepsilon \geq \Delta(\delta, \Theta)$ defined in
  \eqref{eq:param-est-ineq}. Then, with 
  probability at least $1-\delta$, the generalization error $\mathcal{E}(x^{\PDRO}) $ of \PDRO\ satisfies the following: 
  \begin{enumerate}[(a)]
  \item 
    When $d$ is an IPM,
    \begin{equation}\label{eq:general-param-dro-ipm}
      \mathcal E(x^{\PDRO}) \leq 2 \Vscr_d(x^*)\varepsilon.
    \end{equation}
  \item 
    When $d$ is 
    a non-IPM metric
    satisfying
    \eqref{eq:general-f-tv-ineq}, e.g.,  
     $\chi^2, KL, H^2$,
    \begin{equation}\label{eq:general-param-dro-2}
      \mathcal E(x^{\PDRO}) \leq 4 C_{d}\|h(x^*;\cdot)\|_{\infty}\sqrt{\varepsilon}. 
    \end{equation}
  \item 
    When $d $ is the 
    $\chi^2$-divergence, the above bound can be improved to  
    \begin{equation*}\label{eq:general-param-dro-chi2}
      \mathcal E(x^{\PDRO}) \leq
      2\sqrt{\varepsilon\text{Var}_{\P^*}[h(x^*;\xi)]} +
      2\varepsilon^{\frac{3}{4}}\|h(x^*;\cdot)\|_{\infty}. 
\end{equation*}
\end{enumerate}
\end{theorem}

Theorem~\ref{thm:general-param-dro-ipm} immediately gives the bounds on
$\mathcal E(x^{\PDRO})$ (excluding constant factors) with probability at
least $1-\delta$ for the following examples: 
\begin{enumerate}[(1)]
\item \emph{1-Wasserstein distance in Example~\ref{ex:oracle-estimator-was}:} 
  $\mathcal{E}(x^{\PDRO}) \leq 2\|h(x^*;\cdot)\|_{\text{Lip}}\Delta(\delta, \Theta)$.
\item \emph{KL-divergence in Example~\ref{ex:oracle-estimator-kl}:} 
  $\mathcal E(x^{\PDRO}) \leq
  2\|h(x^*;\cdot)\|_{\infty}\sqrt{\Delta(\delta, \Theta)}$. 
  \end{enumerate}

In Appendix~\ref{app:improve-f}, 
we further show that the improved bound for the $\chi^2$-divergence in Theorem~\ref{thm:general-param-dro-ipm}(c) can be extended to general $f$-divergence such as the KL divergence and $H^2$-distance, under additional 
conditions (conditions of Examples~\ref{ex:pdro-kl} and~\ref{ex:pdro-h2}). These extensions state that there exist constants $c_1$ and $c_2$ such
that 
$\mathcal{E}(x^{\PDRO})  \leq c_1 \sqrt{\varepsilon\text{Var}_{\P^*}[h(x^*;\xi)]} + c_2
\varepsilon^{\frac{3}{4}}\|h(x^*;\cdot)\|_{\infty}$ with probability at least $1-\delta$.

  Next, for comparison, we establish bounds for the generalization error of $\PERM$.

\begin{theorem}[Generalization bounds for \PERM]\label{thm:general-param-erm}
  Suppose Assumption~\ref{asp:oracle-param-est} holds and let 
  $M = \sup_{x, \xi}|h(x;\xi)| < \infty$. Then with probability at
  least $1-\delta$, the generalization error $\mathcal{E}(x^{\PERM})$ of \PERM\ satisfies the
  following: 
  \begin{enumerate}[(a)]
      \item When $d$ is an IPM,
      \begin{equation*}\label{eq:p-erm} 
        \mathcal E(x^{\PERM}) \leq 2\left(\sup_{x \in
            \Xscr}\Vscr_{d}(x)\right) \Delta(\delta, \Theta). 
      \end{equation*} 
      \item When $d$ is a non-IPM satisfying \eqref{eq:general-f-tv-ineq},
        e.g., $\chi^2, KL, H^2$, 
        \begin{equation*}\label{eq:p-erm-nonipm} 
          \mathcal E(x^{\PERM}) \leq 4 C_{d} M\sqrt{\Delta(\delta, \Theta)}.
        \end{equation*}
      \item When $d$ is the $\chi^2$-divergence, the above bound
        can be 
        improved 
        to
      \begin{equation}\label{eq:p-erm-2}
        \sqrt{2\Delta(\delta,\Theta)}
        \sqrt{\text{Var}_{\P^*}[h(x^*;\xi)]}+
        2M\left(\Delta(\delta,\Theta)\right)^{\frac{3}{4}}.
      \end{equation}
    \end{enumerate}
\end{theorem}

Note that 
the bound for
$\mathcal{E}(x^{\PERM})$ in Theorem \ref{thm:general-param-erm} involves 
a worst-case term of the form 
$\sup_{x\in \Xscr}\Vscr_d(x)$ or $M = \sup_{x,\xi} \abs{h(x;\xi)}$, which
can be 
significantly
larger than the terms of the form $\Vscr_d(x^*)$ or
$\|h(x^*;\cdot)\|_{\infty}$ that appear in the bound for 
$\mathcal E(x^{\PDRO})$ in Theorem~\ref{thm:general-param-dro-ipm}. That is, the bound for
$\mathcal{E}(x^{\PERM})$
amplifies the model error $\mathcal{E}_{apx}$ when using \PERM\ and, conversely, it demonstrates the power of \PDRO\ in curbing the impact of
model error.

The main results of this section are summarized in the ``Parametric'' row of Table~\ref{tab:general-table}, with each entry representing the best result for each method from Theorems~\ref{thm:general-param-dro-ipm} and~\ref{thm:general-param-erm} respectively. In particular, the bound in the \PDRO\ entry is attained with the 1-Wasserstein distance, and the bound in the \PERM\ entry is attained with the $\chi^2$-divergence. The following show explicitly these bounds presented in the table:
\begin{corollary}[More explicit generalization bounds for \PDRO\ and \PERM]\label{coro:specific}

Ignoring the appearance of $\log(1/\delta)$ and numerical constants, 
if we set $\varepsilon = C \cdot \Delta(\delta, \Theta)$ with $C \geq 1$, then when $d$ is taken as the 1-Wasserstein distance, we have:
\begin{align*}
   \Escr(x^{\texttt{P-DRO}}) &\leq \|h(x^*;\cdot)\|_{Lip}\left(\left(\frac{\text{Comp}(\Theta)}{n}\right)^{\alpha} + \Escr_{apx}\right) \\
   \Escr(x^{\texttt{P-ERM}}) &\leq \sup_{x \in \Xscr}\|h(x;\cdot)\|_{Lip}\left(\left(\frac{\text{Comp}(\Theta)}{n}\right)^{\alpha} + \Escr_{apx}\right).
\end{align*}

When $d$ is taken as the $\chi^2$-divergence, we have:
\begin{align*}
    \Escr(x^{\PDRO}) &\leq \sqrt{\text{Var}_{\P^*}[h(x^*;\xi)]}\left(\left(\frac{\text{Comp}(\Theta)}{n}\right)^{\frac{\alpha}{2}} + \Escr_{apx}\right) +\|h(x^*;\cdot)\|_{\infty}\Escr_{apx}^{\frac{3}{4}}\\
    \Escr(x^{\texttt{P-ERM}}) &\leq \sqrt{\text{Var}_{\P^*}[h(x^*;\xi)]}\left(\left(\frac{\text{Comp}(\Theta)}{n}\right)^{\frac{\alpha}{2}} + \Escr_{apx}\right) + M \Escr_{apx}^{\frac{3}{4}}.
\end{align*}
\end{corollary}

Note that $\alpha$ depends on the metric $d$, where $\alpha = \frac{1}{2}$ 
when $d$ is the 1-Wasserstein
distance  and $\alpha = 1$ when $d$ is the $\chi^2$-divergence. We can readily see that the bound in the \PDRO\ entry in Table~\ref{tab:general-table} is attained by 1-Wasserstein distance, and the bound in the \PERM\ entry is attained by $\chi^2$-divergence. For the other two bounds in Corollary~\ref{coro:specific} not shown in Table~\ref{tab:general-table}, note that \PDRO\ with $\chi^2$-divergence still improves its \PERM\ counterpart by turning the uniform quantity $M$ into $\|h(x^*;\cdot)\|_\infty$, a quantity that depends on $h$ evaluated only at $x^*$, and likewise, \PDRO\ with 1-Wasserstein turns the uniform quantity $\sup_{x\in\mathcal X}\|h(x;\cdot)\|_{Lip}$ in \PERM\ into $\|h(x^*;\cdot)\|_{Lip}$. Corollary~\ref{coro:specific} follows by observing that the choice of the ambiguity size $\varepsilon$ in \PDRO\ satisfies the conditions in Theorem~\ref{thm:general-param-dro-ipm}, and we plug in the concrete expression of $\Delta(\delta, \Theta)$ in Assumption~\ref{asp:oracle-param-est} to obtain the bounds for \texttt{P-DRO} and \texttt{P-ERM}.

We close this section by explaining the main proof ideas in establishing Theorems~\ref{thm:general-param-dro-ipm} and \ref{thm:general-param-erm}. We first discuss
Theorem~\ref{thm:general-param-dro-ipm}. Here $\hat{Z}(\cdot)$ is taken as $\hat{Z}^{P-DRO}(x) := \sup_{d(\P, \hat{\Q})\leq \varepsilon}[h(x;\xi)]$. Intuitively, we would like to set the ambiguity size $\varepsilon$ to achieve the best trade-off between the coverage probability $\P(d(\P^*, \hat{\Q}) \leq \varepsilon)$, which controls the probability that the first term $Z(\hat{x}) -
  \hat{Z}^{DRO}(\hat{x})$ in \eqref{eq:error-decomp} is non-positive, and the size of the ambiguity set, which determines the magnitude of the second term $\hat{Z}(x^*) - Z(x^*)$ in \eqref{eq:error-decomp}. Specifically,
  with the parametric distribution in Assumption~\ref{asp:oracle-param-est} set as the ball center, we can analyze this trade-off as follows:
\begin{enumerate}[(i)]
\item 
Assumption~\ref{asp:oracle-param-est} and our choice of $\varepsilon$ ensures
$\P[d(\P^*, \hat{\Q}) \leq \varepsilon] \geq 1-\delta$.
\item 
  In the event $d(\P^*, \hat{\Q}) \leq \varepsilon$, we have $\P^* \in \Ascr$ and
  $\E_{\P^*}[g(\xi)] \leq \sup_{\P \in \Ascr} \E_{\P}[g(\xi)]$ for any measurable
  function $g$. 
 Therefore, 
  the first term $Z(\hat{x}) -
  \hat{Z}^{P-DRO}(\hat{x})$ in \eqref{eq:error-decomp} is non-positive with
  probability at least $1-\delta$. This observation holds for all three
  cases in Theorem~\ref{thm:general-param-dro-ipm}.
\item 
  When $d$ is an IPM, the second term 
  \begin{align*}
    \hat{Z}^{DRO}(x^*) - Z(x^*)
    &\leq \max_{\P:d(\P, \hat{\Q}) \leq
      \varepsilon} \abs{\E_{\P}[h(x^*;\xi)] -
      \E_{\P^*}[h(x^*;\xi)]}\\
    &\leq \max_{\P: d(\P, \hat{\Q}) \leq
      \varepsilon} \Big\{\sup_{f: \Vscr_d(f) \leq \Vscr_d(h(x^*,\cdot))}\abs{\E_{\P}[f] -
      \E_{\P^*}[f]} \Big\}\\
    &\leq \Vscr_d(x^*)\max_{\P: d(\P, \hat{\Q}) \leq \varepsilon} d(\P, \P^*)\\
                          &\leq \Vscr_d(x^*)(d(\P, \hat{\Q}) + d(\hat{\Q},
                            \P^*)) \leq 2\Vscr_d(x^*) \varepsilon. 
  \end{align*}
\item When $d$ is an $f$-divergence
  satisfying~\eqref{eq:general-f-tv-ineq}, we have
  \[
    \hat{Z}^{DRO}(x^*) \leq \max_{\P: d_{TV}(\P, \hat{\Q})\leq
      C_{d}\sqrt{\varepsilon}} \E_{\P}[h(x^*;\xi)],
  \]
  which allows us to reduce to the previous case. 
\item For the 
  $\chi^2$-divergence more specifically, 
  the Cauchy-Schwarz
  inequality 
  implies that
  \begin{eqnarray*} \label{eq:chi2-pseudo-ipm}
    \lefteqn{\sup_{\P: \chi^2(\P, \hat{\Q})
    \leq
      \varepsilon} |\E_{\P}[h(x^*;\xi)] -
    \E_{\P^*}[h(x^*;\xi)]|}\\
    & \leq & \sup_{\P: \chi^2(\P, \hat{\Q})
             \leq
             \varepsilon} |\E_{\P}[h(x^*;\xi)] -
             \E_{\hat{\Q}}[h(x^*;\xi)]| + |\E_{\hat{\Q}}[h(x^*;\xi)] -
             \E_{\P^*}[h(x^*;\xi)]| \\
    & \leq & 2\sqrt{2\varepsilon \text{Var}_{\hat{\Q}}[h(x^*;\xi)]}.
  \end{eqnarray*}
  The final bound for this case follows by writing $\text{Var}_{\hat{\Q}}[h(x^*;\xi)]$ in terms
  of $\text{Var}_{\P^*}[h(x^*;\xi)]$ and $\chi^2(\P^*, \hat{\Q})$.
\end{enumerate}


Finally, for Theorem~\ref{thm:general-param-erm}, consider the decomposition~\eqref{eq:error-decomp}, without the worst-case machinery of DRO here, we bound the two terms by $\sup_{x\in \Xscr}|Z(x) - \hat{Z}(x)|$, which leads to the appearance of $\sup_{x \in \Xscr}\Vscr_d(x)$. The improved $\chi^2$ result follows by 
replacing the uniform bound $\sup_{x \in \Xscr}
\sqrt{\text{Var}_{\P^*}[h(x;\xi)]}$ with an alternative bound 
$\sqrt{\text{Var}_{\P^*}[h(x^{\PERM};\xi)]}
  \leq
    \sqrt{\text{Var}_{\P^*}[h(x^\ast;\xi)]}
    + 2 M (\chi^2(\P^*,
    \hat{\Q}))^{\frac{3}{4}}$ that is valid for the 
solution $x^{\PERM}$.

\section{Parametric-DRO with Monte Carlo Approximation}\label{subsec:sample} 
When the 
center $\hat{\Q}$ is continuous, the inner
maximization $\sup_{\P \in \Ascr}\E_{\P}[h(x;\xi)]$ in 
\eqref{eq:DRO-obj} 
can be computationally challenging.
For example, for the 1-Wasserstein distance
\citep{esfahani2018data} and $f$-divergence \citep{bayraksan2015data},
the dual problem of the inner maximization is reformulated as follows:
\begin{align*}
    \sup_{\P \in \Ascr}\E_{\P}[h(x;\xi)]& = \inf_{\lambda \geq 0}\{\lambda \varepsilon + \E_{\xi_0 \sim \hat{\Q}}[\sup_{\xi \in \Xi}\{h(x;\xi) - \lambda \|\xi_0 - \xi\|\}]\}\tag{1-Wasserstein Distance }\\
     \sup_{\P \in \Ascr}\E_{\P}[h(x;\xi)] &=\inf_{\lambda \geq 0, \mu \in \R}\{\mu + \lambda \varepsilon + \E_{\xi \sim \hat{\Q}}[(\lambda f)^*(h(x;\xi) - \mu)]\}, \tag{$f$-divergence}
\end{align*}
where the objective involves a high-dimensional integral over $\hat{\Q}$
instead of the empirical distribution. This makes the problem harder
to evaluate and optimize than DRO that 
uses
the empirical
distribution as the ball center. There are two approaches to handle this issue: sample average approximation~(SAA) that reduces the 
problem to a structure resembling the standard empirical DRO, 
and
stochastic approximation. 
We discuss the first approach in this section.
The second approach is discussed in 
Appendix~\ref{app:sa}. 


\subsection{Monte Carlo sampling to solve the inner problem in
  \PDRO }
Let 
$\tilde{\xi}_i \sim
\hat{\Q}$, $i = 1,\ldots, m$ denote $m$ i.i.d. samples.   
We approximate $\hat{\Q}$ with the Monte Carlo estimate
$\hat{\Q}_m:=
\frac{1}{m}\sum_{i = 1}^m \delta_{\tilde{\xi}_i}$ 
and then define
$\hat{\Ascr} =\{\P|d(\P,\hat{\Q}_m)\leq\varepsilon\}$. Let $x^{\PDRO_m}$
denote the corresponding solution. 
We show 
how to bound the 
generalization error $\Escr(x^{\PDRO_m})$ 
and investigate how to compute a 
sample
size $m$ 
that ensures that 
$\mathcal{E}(x^{{\PDRO}_m})\approx \mathcal{E}(x^{\PDRO})$ 
within a constant error. 
More precisely, note that the approximate 
solution $x^{\PDRO_m}$ 
incurs both the statistical error 
arising from finite training
data and the
Monte Carlo sampling error from finite $m$, and our goal is to
choose $m$ to ensure 
that the Monte Carlo sampling error is smaller than the statistical error.
In this section, we assume $h(x;\xi) \in [0, M]$, for all $x, \xi$.

We first have the following bound for Wasserstein DRO that can be straightforwardly derived:
\begin{theorem}[Generalization bounds for Wasserstein \PDRO\ with Monte Carlo errors]\label{coro:p-dro-sample-w1}
Suppose Assumption~\ref{asp:oracle-param-est} holds with $\E_{\Q}[\exp(\|\xi\|^a)] < \infty$ for some $a > 1, \forall \Q \in \Pscr_{\Theta}$. Suppose also that the size of the
ambiguity set satisfies $\varepsilon\geq 2\Delta(\delta,
\Theta)$, and the distance metric $d$ defining $\Ascr$ is the
1-Wasserstein distance. Then $m = O((2/\varepsilon)^{D_{\xi}})$
ensures that 
\[
  \mathcal{E}(x^{{\PDRO}_m}) \leq 2 \varepsilon \|h(x^*;\cdot)\|_{\text{Lip}}.
\]
with probability $1-\delta$.
\end{theorem}
Suppose $\mathcal{E}_{apx} \approx 0$, i.e., the family $\mathcal{P}_{\Theta}$
approximates $\P^*$ well. Then 
$\varepsilon = 2\Delta(\delta, \Theta)$ implies that the required Monte
Carlo sample size $m = O( n^{\alpha D_{\xi}})$. This means that $m$ does not depend on $\text{Comp}(\Hscr)$. However, it depends exponentially on $D_{\xi}$. Moreover, a key in proving
Theorem~\ref{coro:p-dro-sample-w1} is that we can maintain the bound
$\P[d(\P^*, \hat{\Q}_m) \leq \varepsilon] \geq 1-\delta$ since $W_1(\P^*,
\hat{\Q}_m) \leq W_1(\P^*, \hat{\Q}) + W_1(\hat{\Q},\hat{\Q}_m) \leq
\varepsilon$ for large $m$. 
This argument does not hold
more generally, because $d(\P^*, \hat{\Q}_m)$ can be infinite  for any
$m$ 
when
$\P^*$ is continuous and $d$ is an $f$-divergence. These challenges motivate us to derive a more general result that leverages 
the equivalence between DRO and regularization. 
\begin{theorem}[Generalization bounds for general \PDRO\ with Monte Carlo
  errors]\label{thm:p-dro-sample} 
Suppose Assumption~\ref{asp:oracle-param-est} holds and the size of the
ambiguity set $\varepsilon \geq \Delta(\delta, \Theta)$. When $d$ is $\chi^2$-divergence or 1-Wasserstein distance, if the Monte Carlo size satisfies $m
\geq C \left(\frac{M}{\Vscr_d(x^*)
    \varepsilon}\right)^6\text{Comp}(\Hscr)\log m$ for some constant $C$,   
then with probability at least $1-\delta$,  $\mathcal{E}(x^{\texttt{P-DRO}_m})
\leq 3\mathcal{E}_{\texttt{P-DRO}}$, where $\mathcal{E}_\texttt{P-DRO}$ is the corresponding
generalization error upper bound in
Theorem~\ref{thm:general-param-dro-ipm}.  
\end{theorem}
  
Suppose $\Escr_{apx}\approx 0$. Then $m
\approx \text{Comp}(\Hscr) n^{6\alpha}$. Thus, we remove the exponential dependence on $D_{\xi}$ for the Monte Carlo size $m$ at the cost of a linear dependence on the hypothesis class complexity. In addition to the proof techniques
  in Theorem~\ref{thm:general-param-dro-ipm}, the key idea here is to 
use  the variability regularization property of DRO to show the bounded Monte Carlo sampling error: 
\begin{equation}\label{eq:dro-sample-reg}
\bigg|\sup_{d(\P, \hat{\Q})\leq \varepsilon}\E_{\P}[h(x;\xi)] -\sup_{d(\P,
  \hat{\Q}_m)\leq \varepsilon} \E_{\P}[h(x;\xi)]\bigg| \leq \mathcal{E}_d, \forall x \in \Xscr.
\end{equation}
To prove \eqref{eq:dro-sample-reg}, the idea is to first apply the following variability regularization: 
\begin{equation}\label{eq:dro-sample-explicit}
    \sup_{d(\P, \hat{\Q})\leq \varepsilon}\E_{\P}[h(x;\xi)] =  \E_{\hat{\Q}}[h(x;\xi)] + \Vscr_d(x) \sqrt{\varepsilon} + O(\varepsilon), \forall x \in \Xscr,
\end{equation}
such that we can decompose the left-hand side of~\eqref{eq:dro-sample-reg} to be a combination of
$|\E_{\hat{\Q}}[h(x;\xi)] - \E_{\hat{\Q}_m}[h(x;\xi)]|$ and the difference in $\Vscr_d(x)$ between $\hat{\Q}$ and $\hat{\Q}_m$ (e.g., 
$\sqrt{\text{Var}_{\hat{\Q}}[h(x;\xi)]} -
\sqrt{\text{Var}_{\hat{\Q}_m}[h(x;\xi)]}$ when $d$ is $\chi^2$-divergence). Then we apply uniform concentration inequalities for these terms over the hypothesis class
$\text{Comp}(\Hscr)$ to show \eqref{eq:dro-sample-reg} holds under large $m$.
However, we need to apply \eqref{eq:dro-sample-explicit} carefully under finite samples. When $d$ is $\chi^2$-divergence, we do not have $\sup_{\P: \chi^2(\P, \hat{\Q})\leq \varepsilon}\E_{\P}[h(x;\xi)] \geq \E_{\P}[h(x;\xi)] + \Vscr_d(x) \sqrt{\varepsilon}, \forall x \in \Xscr$ when $\text{Var}_{\hat{\Q}}[h(x;\xi)]$ is small. Therefore, we split $\Xscr$ depending on the size of $\text{Var}_{\hat{\Q}}[h(x;\xi)]$ into several regions and investigate the variability regularization effect in each region to give a tighter bound of the Monte Carlo sampling error, i.e., the left-hand side in \eqref{eq:dro-sample-reg}.

In Appendix~\ref{subsec:pdro-thm}, we present Theorem~\ref{thm:p-dro-sample-chi2} for $\chi^2$-divergence and Theorem~\ref{thm:p-dro-sample-w1} for 1-Wasserstein distance which give generalization bounds on $\mathcal{E}(x^{\texttt{P-DRO}_m})$, from which Theorem~\ref{thm:p-dro-sample} follows. Moreover, the metric $d$ in Theorem~\ref{thm:p-dro-sample} can be extended to the $p$-Wasserstein distance for $p \in [1, 2]$ with $\Vscr_d(x)$ being the gradient norm, and the minimum required Monte Carlo size $m$ also does not depend exponentially in $D_{\xi}$. This result is in Corollary~\ref{coro:p-dro-sample-wp} in Appendix~\ref{app:p-was-general}. 

Define $x^{{\PERM}_m}$ as the corresponding \PERM\ solution using Monte Carlo approximation with $m$ samples. We end this section by providing a bound on the number of samples required
to ensure that 
$\mathcal{E}(x^{{\PERM}_m}) \approx \mathcal{E}(x^{\PERM})$. 
\begin{theorem}[Generalization bounds for \PERM\ with Monte Carlo
  errors]\label{thm:general-param-erm-sample} 
  Suppose Assumption \ref{asp:oracle-param-est} holds for the metric~$d$.
  Let $\mathcal{E}_{\texttt{P-ERM}}$ 
  denote the 
  generalization error upper bound in
  Theorem~\ref{thm:general-param-erm}. 
  Then, with probability at least $1-\delta$, $\mathcal{E}(x^{\texttt{P-ERM}_m})
  \leq 2\mathcal{E}_\texttt{P-ERM}$ 
  provided the Monte Carlo size $m$ satisfies:
  \[
    \frac{m}{M\text{Comp}(\Hscr) \log m} \geq\max\left\{\frac{1}{Z(x^*) +
        \mathcal{E}_\texttt{P-ERM}}, \frac{Z(x^*) +
        \mathcal{E}_\texttt{P-ERM}}{\mathcal{E}_\texttt{P-ERM}^2}\right\}.
    \]
\end{theorem}
This result utilizes 
arguments similar to those used to establish
Theorem~\ref{thm:p-dro-sample} to control the term $\sup_{x \in
  \Xscr}|\E_{\hat{\Q}}[h(x;\xi)] - \E_{\hat{\Q}_m}[h(x;\xi)]|$. However,
the analysis here is much simpler than \PDRO\ since 
one does not have to control
the variability regularization term. 
Note that, 
ignoring the model misspecification error
$\mathcal{E}_{apx}$, 
the required Monte Carlo sample size $m$ 
is proportional to the function complexity $M\text{Comp}(\Hscr)$ 
and $n^{2\alpha}$, which does not depend on $D_{\xi}$ exponentially. 


Considering the Monte Carlo error in this section, \PDRO\ can be viewed as translating the statistical
errors associated with non-parametric methods, 
entailed by the distribution dimension or function complexity, to 
model misspecification errors and
additional computational effort associated with Monte Carlo sampling in Theorems~\ref{coro:p-dro-sample-w1} and \ref{thm:p-dro-sample}.
The computational
cost is driven by two factors.
The first is the computational cost of 
sampling from the parametric model $\hat{\Q}$, 
which is likely to be low since most parametric models in practice are sampling-friendly, and therefore 
we can allow the number of Monte Carlo samples 
$m \gg n$,  
the number of data points available from~$\P^*$.
The second is the computational cost  
in solving the corresponding
DRO problem. Note that the optimization problem in $\PDRO_m$
is the same as that in non-parametric DRO, albeit with a larger number of
samples. 
One can leverage recently proposed procedures for large-scale DRO
with the
$f$-divergence~\citep{levy2020large,jin2021non} and Wasserstein
distance~\citep{SinhaND18}. 
We expect improvements in computing
large-scale DRO 
will lead to more efficient methods in 
solving \PDRO. 

\section{Extensions}\label{sec:extensions}
\subsection{Extension to Distribution Shifts}\label{subsec:dist-shift}
We extend our \PDRO\ 
framework to the setting with distribution shifts where the testing distribution $\P^{te}$ is different from the training distribution $\P^{tr}$ of the data used to
compute the solution $\hat{x}$. 
We are
interested in the generalization error $\Escr^{te}(\hat{x}) =
\E_{\P^{te}}[h(\hat{x};\xi)] - \min_{x \in \Xscr}\E_{\P^{te}}[h(x;\xi)]$.
\begin{assumption}[Oracle estimator under distribution shifts]\label{asp:oracle-dist-shift}
    Assumption~\ref{asp:oracle-param-est} holds when $\P^*$ there is replaced by $\P^{tr}$. Besides, there exist some numerical constants $c_1, c_2 > 0$ such that:
    \begin{equation}\label{eq:dist-shift}
        d(\P^{te}, \hat{\Q}) \leq c_1 d(\P^{te}, \P^{tr}) + c_2 d(\P^{tr}, \hat{\Q}).
    \end{equation}
\end{assumption}
The first part of Assumption~\ref{asp:oracle-dist-shift} holds naturally since we still construct $\hat{\Q}$ using i.i.d. samples $\{\hat{\xi}_i\}_{i = 1}^n$ from $\P^{tr}$. For the second part, when $d$ is an IPM, the triangle inequality immediately implies~\eqref{eq:dist-shift} with $c_1 = c_2 = 1$. When $d$ is not an IPM, \eqref{eq:dist-shift} may still hold when $\P^{te}$ is absolutely continuous with respect to $\P^{tr}$. For example, when $d$ is KL-divergence, we show that $c_1 = 1, c_2 = \|d\P^{te}/d\P^{tr}\|_{\infty} < \infty$ in Lemma~\ref{lemma:triangle-f-div} in Appendix~\ref{subsec:dist-shift-proof2}.

We have the following results:
\begin{corollary}[Generalization bounds for \PDRO\ under distribution shifts]\label{coro:param-dro-dist-shift}
Suppose Assumption~\ref{asp:oracle-dist-shift} holds, and the size
$\varepsilon \geq c_1 d(\P^{te}, \P^{tr}) + c_2 \Delta(\delta, \Theta)$. Then the
bounds in Theorem~\ref{thm:general-param-dro-ipm} hold when $\mathcal{E}(x^{\PDRO})$ there is replaced by $\mathcal{E}^{te}(x^{\PDRO})$. 
\end{corollary}
\begin{corollary}[Generalization bounds for \PERM\ under distribution
  shift]\label{coro:param-erm-dist-shift} 
  Suppose Assumption~\ref{asp:oracle-dist-shift} holds. Then the bounds in
  Theorem~\ref{thm:general-param-erm} hold when $\Escr(x^{\PERM})$ there is replaced by $\mathcal{E}^{te}(x^{\PERM})$ 
  if one 
  replaces $\Delta(\delta,
  \Theta)$ with $c_1 d(\P^{te}, \P^{tr}) + c_2 \Delta(\delta, \Theta).$ 
\end{corollary}
These 
two results 
can be established using techniques similar to those used to establish
Theorem~\ref{thm:general-param-erm}. 
They show the same insights and strengths of \PDRO\ compared with \PERM, but with an additional benefit under distribution shifts: The generalization error of
\PERM\ here suffers from the product of the distribution shift amount $d(\P^{tr}, \P^{te})$ and uniform quantity over
$\Xscr$ such as $M$ and $\sup_{x\in \Xscr}\Vscr_d(x)$. On the other hand, the additional term in the generalization error of \PDRO\ involves only $d(\P^{tr}, \P^{te}) \Vscr_d(x^*)$. For comparison, we also discuss generalization error bounds of standard nonparametric ERM and DRO under distribution shifts and highlight the additional terms they involve in Appendix~\ref{app:erm-dro-dist-shift}.

\subsubsection{Comparisons with Existing Approaches.}
We compare our bounds with those available in the literature on
distribution shifts. A commonly used evaluation metric under 
distribution shifts  is the  ``discrepancy metric'' 
$\textsf{disc}_L(\Hscr; \P^{te}, \P^{tr}) = \sup_{h_1, h_2 \in
  \Hscr}|\E_{\P^{te}}L (h_1, h_2) - \E_{\P^{tr}}L(h_1, h_2)|$ for a given
loss function $L: \Hscr \times \Hscr \to \R$
\citep{mansour2009domain,ben2010theory,zhang2019bridging}. This metric 
depends on the hypothesis class $\mathcal H$, and one requires samples
from both $\P^{tr}$ and $\P^{te}$ to evaluate the metric. In contrast, our
bound is in terms of 
the metric
$d(\P^{te}, \P^{tr})$ 
that 
does not consider the interactions between
$\Hscr$ and $\P^{tr}, \P^{te}$, which means our method could be less sensitive to large $\text{Comp}(\Hscr)$ under distribution shifts compared to the use of the ``discrepancy metric''.   
We also do not need 
samples from $\P^{te}$. However, 
we require the knowledge of an upper bound on
$d(\P^{te}, \P^{tr})$. \cite{lee2018minimax} 
establish 
a 
bound on the 
generalization error 
when $\hat{\Q}$
is 
the empirical distribution, and $d$ is the $p$-Wasserstein
distance. This approach suffers  
from the curse of dimensionality as is the case for standard
Wasserstein-DRO approaches, and this is also apparent in the  
numerical
results in Section~\ref{sec:numerics}. 
There are works that assume
specific types of distribution shifts, e.g.,  
group~\citep{Sagawa2020Distributionally}, latent covariate
shifts~\citep{duchi2020distributionally}, and conditional
shifts~\citep{sahoo2022learning}. 
Extending the
\PDRO\ 
approach to distributional shifts with specialized structures 
appears to be an interesting future direction.

\subsection{Extension to Contextual Optimization}\label{sec:context}
We extend our model formulation to contextual
optimization~\citep{ban2019big,bertsimas2020predictive} 
where 
$\xi \sim \P^*_{\xi|y}$ 
depends on the 
covariate $y \in \R^{D_{y}}$ 
that is observed 
before making the decision~$x$. 
In this setting the data $\Dscr_n := \{(\hat{y}_i,\hat{\xi}_i)\}_{i =
1}^n$, and the 
contextual DRO 
problem~\citep{bertsimas2021bootstrap,esteban2020distributionally} is given by
\begin{equation}\label{eq:context-dro}
 \min_{x \in \Xscr}\max_{d(\P_{\xi|y}, \hat{\Q}_{\xi|y}) \leq
   \varepsilon}\E_{\P_{\xi|y}}[h(x;\xi)], 
\end{equation}
where $\hat{\Q}_{\xi|y}$ is estimated from 
$\Dscr_n$ and depends on $y$, which is chosen from a class of parametric distributions $\Pscr_{\Theta, y}$. Here, we can extend the definition of distribution class $\Pscr_{\Theta}$ in Section~\ref{sec:main} to incorporate the presence of $y$ in the distribution class. For example, extending Example~\ref{ex:oracle-estimator-was}, $\Pscr_{\Theta, y} = \{\Nscr(B y,\Sigma)| B \in \R^{D_{\xi}\times D_y}\}$, and extending Example~\ref{ex:oracle-estimator-kl}, $\Pscr_{\Theta, y}$ is the class of all distributions of the random variable $g_{\theta}(Z, y)$ for a fixed random variable $Z$. Note that the size of the ambiguity set $\varepsilon$ is fixed
for all  
values of the covariate~$y$.

We 
show how the
\PDRO\ approach proposed in
Section~\ref{sec:main} can be generalized to contextualized
DRO~\eqref{eq:context-dro}.
We establish two types of generalization bounds: First is a point-wise bound that holds with high
probability for all $y \in \Yscr$, and second is a bound that holds on average over
the random covariate $Y$. We investigate both types for reasons that will be apparent as we discuss these bounds.    

\subsubsection{
  Point-wise Generalization Error Bound.
}
Here we are interested in bounding the error:
\begin{equation}\label{eq:context-error}
\Escr_y(\hat{x}) = \E_{\P^*_{\xi|y}}[h(\hat{x};\xi)] - \min_{x\in
  \Xscr}\E_{\P^*_{\xi|y}}[h(x;\xi)],
\end{equation}
and make the following assumption:
 
\begin{assumption}[Oracle conditional estimator]
  \label{asp:conditional-oracle-estimator}
  Let $\text{Comp}(\Theta, y)$ be the complexity of $\Pscr_{\Theta, y}$, and
$\mathcal{E}_{apx}(\P^*_{\xi|y},\mathcal P_{\Theta,y})$ (abbreviated to $\mathcal{E}_{apx}(y)$) is a non-negative function
such that $\mathcal{E}_{apx}(\P^*_{\xi|y},\mathcal P_{\Theta,y})
= 0$ if $\P_{\xi|y}^* \in \mathcal{P}_{\Theta, y}$.
Then given any $y \in \Yscr$, for all $\delta \in (0, 1)$, there exists $\alpha > 0$ such that the center $\hat{Q}_{\xi|y} \in \Pscr_{\Theta, y}$ of the ambiguity set satisfies
  \begin{equation}\label{eq:cond-param-est-ineq}
    \begin{aligned}
      d(\P^*_{\xi|y},\hat{\Q}_{\xi|y}) &\leq
      \mathcal{E}_{apx}(\P^*_{\xi|y}, \mathcal{P}_{\Theta, y}) +
      \left(\frac{\text{Comp}(\Theta, y)}{n}\right)^{\alpha} \log(1/\delta)=:\Delta(\delta, \Theta, y),
  \end{aligned}
\end{equation}
with probability $1-\delta$.
\end{assumption}


Note that the oracle estimation property in \eqref{eq:cond-param-est-ineq} holds with high probability for any given fixed $y$. Next, we  provide an example satisfying
Assumption~\ref{asp:conditional-oracle-estimator}.


\begin{example}\label{prop:conditional-oracle-estimator}

Suppose $d$ is given by the 1-Wasserstein distance, the set of parametric distributions $\Pscr_{\Theta, y} = \{N(f_{\theta}(y), \Sigma)|\theta \in \Theta, \Sigma \in \mathbb{S}_{++}^{D_{\xi}\times D_{\xi}}\}$ and the true distribution 
$\xi := f_{\theta^*}(y) + \eta$,
where $f_{\theta^*}(y)$ is a deterministic function of $y$ parametrized by the unknown $\theta^*$, and $\eta$ is
independent of $y$ 
with $\|\eta\|_2^2 \leq C_{\eta}$, 
$\E[\eta] = 0$ and $\E[\eta\eta^{\top}] =
\Sigma$.
Then under some mild conditions (i.e., Assumption~\ref{asp:new} in Appendix~\ref{app:context-proof}), 
Assumption~\ref{asp:conditional-oracle-estimator} holds for $\hat{\Q}_{\xi|y} = \Nscr(f_{\hat{\theta}}(y), \hat{\Sigma})$ 
$\forall y \in \Yscr$ where $\hat{\theta}
\in \argmin_{\theta\in \Theta}\sum_{i =
1}^n\|\hat{\xi}_i - f_{\theta}(\hat{y}_i)\|_2^2$, and 
$\hat{\Sigma} = \frac{1}{n}\sum_{i = 1}^n(\hat{\xi}_i -
f_{\hat{\theta}}(\hat{y}_i))(\hat{\xi}_i - f_{\hat{\theta}}(\hat{y}_i))^{\top}$ with $\alpha = \frac{1}{2}$, 
$\Escr_{apx}(y) = W_1(\P_{\xi|y}^*,
\Nscr(\E[\xi |y], \Sigma))$ and $\text{Comp}(\Theta, y) =
O(\text{Comp}(\Fscr) \vee \text{Tr}(\Sigma) \max\{C_{\eta}^2,
D_{\xi}\})$ where $\text{Comp}(\Fscr)$ is the function complexity of
$\Fscr(:= \{f_{\theta}(\cdot): \theta \in \Theta\})$.
\end{example}

The following result extends our guarantees for $\PDRO$ and $\PERM$ in Section \ref{sec:main} to contextual
optimization. 
\begin{corollary}[Generalization bounds for contextual \PDRO\ and \PERM]
\label{coro:param-dro-erm-context} 
Suppose Assumption~\ref{asp:conditional-oracle-estimator} holds and the
size $\varepsilon \geq \sup_y \Delta(\delta, \Theta, y)$. 
Then for any given $y \in \Yscr$, with probability
at least $1-\delta$, the bounds in Theorem~\ref{thm:general-param-dro-ipm} hold when $\Escr(x^{\PDRO})$ there are replaced by $\Escr_y(x^{\PDRO})$. And the bounds in Theorem~\ref{thm:general-param-erm} hold when $\Escr(x^{\PERM})$ there is replaced by $\Escr_y(x^{\PERM})$ if one replaces $\Delta(\delta, \Theta)$ 
with $\Delta(\delta, \Theta, y)$.
\end{corollary}


This result 
follows a similar argument as the non-contextual case by using Assumption~\ref{asp:conditional-oracle-estimator} that exerts a high-probability condition on $\Escr_y(\hat{x})$ for any
$y \in \Yscr$,
and the strengths of \PDRO\ and comparisons with \PERM\ presented previously all carry over to this contextual case. However, Assumption~\ref{asp:conditional-oracle-estimator} is arguably overly strong, as the covariate $Y$ is random and it could be difficult to ensure the high-probability condition holds for every single $y\in \Yscr$ and the satisfying $\varepsilon$ in Corollary~\ref{coro:param-dro-erm-context} can be overly large since the term $\text{Comp}(\Theta, y)$ in $\Delta(\delta, \Theta, y)$ can be extremely large under some $y$, as we will see in Example~\ref{prop:expect-conditional}. This motivates us to consider the average generalization error presented in the next subsection.

\subsubsection{Average Generalization Error Bound.}
The
\emph{average} generalization error $\Escr_{\Yscr}(\hat{x})$ of a contextual
decision $\hat{x}(y)$ is defined 
as
\begin{equation}\label{eq:expect-context-error}
    \Escr_{\Yscr}(\hat{x}):=\E_{\Dscr_n}\E_y\left(\E_{\P_{\xi|y}^*}[h(\hat{x}(y);\xi)] - \E_{\P_{\xi | y}^*}[h(x^*(y);\xi)]\right),
\end{equation}
where the first expectation is over the random 
dataset $\Dscr_n$, 
the second expectation is over the covariate distribution, and $x^*(y) \in
\argmin_{x\in \Xscr}\E_{\P_{\xi|y}^*}[h(x;\xi)]$. 
This average generalization error $\Escr_{\Yscr}(\hat{x})$ is the same as the
average regret introduced in \cite{hu2022fast}. 



Comparing~\eqref{eq:expect-context-error} with the point-wise generalization error in~\eqref{eq:context-error}, the difference is that we involve the expectation over the dataset $\Dscr_n$ and the new covariate distribution $Y$. The two expectations in \eqref{eq:expect-context-error} represent two sources of uncertainty: uncertainty 
in the historical data~$\Dscr_n$, and uncertainty 
in the value of the covariate $Y$. 
In order to balance these two sources of uncertainty, we assume that we
have access to an oracle that satisfies the following assumption.
\begin{assumption}[Average performance of oracle conditional
  estimator]\label{asp:expect-conditional-oracle-estimator}
  Denote $\hat d := \sup_{y \in \Yscr}\Escr_{apx}(y) + (\E_y[\text{Comp}(\Theta, y)]/n)^{\alpha}$ with $\Escr_{apx}(y)$ and $\text{Comp}(\Theta, y)$ defined in Assumption~\ref{asp:conditional-oracle-estimator}. Suppose the estimator $\hat{\Q}_{\xi|y}$
  satisfies
  \begin{equation}\label{eq:context-joint-ineq}
    \P\left(d(\P_{\xi|y}^*,
      \hat{\Q}_{\xi | y})  - \hat{d}
      \geq t\right) \leq c_1\exp(-c_2 a_n t^2),
  \end{equation}
  for some $c_1, c_2 \ge 0$, and all $t \ge 0$, $a_n$ with $n$ being the sample size, is an increasing deterministic sequence such that $a_n \to
  \infty$ as $n \to \infty$.
\end{assumption}
Next, we show an example of a setting where~\eqref{eq:context-joint-ineq}  
holds. 
\begin{example}\label{prop:expect-conditional} 
Suppose $d$ is given by 1-Wasserstein distance and the parametric class $\Pscr_{\Theta, y} = \{N(By, \Sigma)| B \in \R^{D_{\xi}\times D_y}\}$ with known $\Sigma$. Then for a general conditional distribution $\P_{\xi|y}^*$ (not necessarily $\E[\xi|y] = By, \forall y$), if Assumption~\ref{asp:bound-cov} in Appendix~\ref{app:context-proof} holds, and $\hat{\Q}_{\xi|y}:= N(\hat{B} y, \Sigma)$ for some fixed $\Sigma$ and
  $\hat{B} \in \R^{D_{\xi}\times D_{y}}$ where $\hat{B} = \argmin_{B} \sum_{i = 1}^n \|\hat{\xi}_i - B \hat{y}_i\|_2^2$, then:
  \begin{itemize}
      \item Assumption~\ref{asp:conditional-oracle-estimator} holds with 
  $\Escr_{apx}(y) = W_1(\P_{\xi|y}^*, \Nscr(B^*y, \Sigma))$, $\alpha = \frac{1}{2}$, $\text{Comp}(\Theta, y) = O (D_{\xi}D_y \|y\|_{\Sigma_y^{-1}})$, where $\Sigma_y = \E[y y^{\top}]$, and $B^* = \argmin_B \E_{(\xi,y)}[\|\xi - B y\|^2]$.
    \item Assumption~\ref{asp:expect-conditional-oracle-estimator} 
  holds with $a_n = n$, $\alpha = \frac{1}{2}, \text{Comp}(\Theta) = O\left(\frac{D_{\xi}D_y}{\lambda_{\min}(\Sigma_y)}\right)$;
  \end{itemize}
\end{example}
 
To derive Example \ref{prop:expect-conditional}, we observe $W_1(\P^*_{\xi|y}, \hat{\Q}_{\xi|y}) \leq
W_1(\P^*_{\xi|y}, \Nscr(B^*y, \Sigma)) + \|\hat{B} - B^*\|_{\Sigma_y}\|y\|_{\Sigma_y^{-1}}$ 
and apply standard concentration results under misspecified linear models (note that here we do not assume $\E[\xi|y] = By$) from \cite{hsu2012random}.
  
In Example \ref{prop:expect-conditional}, when the observed $\|y\|_{\Sigma_y^{-1}}$ is 
small, $\text{Comp}(\Theta, y)$ and $\Delta(\delta, \Theta, y)$ in \eqref{eq:cond-param-est-ineq} would be small. For these instances of $y$, choosing a small $\varepsilon \geq \Delta(\delta, \Theta, y)$ would lead to a generalization error bound $\Escr_y(x^{\PDRO}) \leq 2 \Vscr_d(x^*(y))\varepsilon$ with probability at least $1-\delta$. In contrast, Corollary~\ref{coro:param-dro-erm-context} demonstrates that for any $y$, if $\varepsilon \geq \sup_{y \in \Yscr}\Delta(\delta, \Theta, y)$, then $\Escr_{y}(x^{\PDRO}) \leq 2\Vscr_d(x^*(y))\varepsilon$ with probability at least $1-\delta$. Comparing these deductions hints that the generalization error bound in Corollary~\ref{coro:param-dro-erm-context} could be overly pessimistic in the choice of $\varepsilon$ because we can attain that bound for some $y \in \Yscr$ with small $\|y\|_{\Sigma^{-1}}$ by using a possibly much smaller $\varepsilon$, e.g., when $\Delta(\delta, \Theta, y) \leq  \varepsilon \ll \sup_{y \in \Yscr}\Delta(\delta, \Theta, y)$. 
The average generalization error and \eqref{eq:context-joint-ineq} are proposed to reduce this pessimism.






\begin{theorem}[Average generalization error bounds for contextual \PDRO\ and \PERM]\label{thm:expect-generalization-bd}
  Suppose Assumption~\ref{asp:expect-conditional-oracle-estimator}
  holds and $d$ is an IPM. 
  Letting the ambiguity size $\varepsilon \geq \hat{d}$,
  the average generalization error of the \PDRO\ solution
  $\hat{x}^{\PDRO}(y)$
  and \PERM\ solution
  $\hat{x}^{\PERM}(y)$ satisfy
  \begin{equation}
    \begin{array}{rl}
    \label{eq:p-dro-bound}
      \Escr_{\Yscr}[\hat{x}^{\PDRO}(y))]
      & \leq 2\varepsilon\E_y[\Vscr_d(x^*(y))] + M
        c_1 \exp\left(-c_2 a_n (\varepsilon - \hat d)^2\right)\left(\varepsilon + \frac{1}{\sqrt{c_2
          a_n}}\right),\\ 
      \Escr_{\Yscr}[\hat{x}^{\PERM}(y)]
      &\leq  2M  \left(\hat{d} + \frac{c_1}{\sqrt{c_2a_n}}\right),
    \end{array}
  \end{equation}
  where $M:=\sup_{x\in \Xscr}\Vscr_d(x)$.
  \end{theorem}
Compared with Corollary~\ref{coro:param-dro-erm-context}, the choice of $\varepsilon$ in Theorem~\ref{thm:expect-generalization-bd} only needs to be larger than $\E_y[\text{Comp}(\Theta, y)]$ instead of $\sup_y\text{Comp}(\Theta, y)$ in terms of the dependence of $\text{Comp}(\Theta, y)$, which reduces pessimism with a potentially much smaller $\varepsilon$. Theorem \ref{thm:expect-generalization-bd} reveals that 
the average generalization error
of \PDRO\ and \PERM\ both depend on the uniform term $M$. However, the second
term of \PDRO\ in~\eqref{eq:p-dro-bound} involving $M$ only occurs when $\P_{\xi|y}^*$ is not in the ambiguity set, thus bearing the exponentially decaying factor $\exp(-c_2 a_n (\varepsilon - \hat d)^2)$. If we plug in the expression of $\hat{d}$ in Assumption~\ref{asp:expect-conditional-oracle-estimator} into \eqref{eq:p-dro-bound}, $M \cdot\Escr_{apx}$ appears in the generalization bound of \PERM, while a lighter dependence of $\Escr_{apx}$ occurs in \PDRO\, where $\E[\Vscr_d(x^*(y))] \cdot \Escr_{apx}$ appears in the first term of \PDRO\ in~\eqref{eq:p-dro-bound} and $M \cdot \Escr_{apx}$ appears in the second term of \PDRO\ in~\eqref{eq:p-dro-bound} but multiplying with a term decaying to zero exponentially with $a_n$. Therefore, \PDRO\ advantageously alleviates the effects of model misspecification error, an insight in line with
the non-context case shown in Table~\ref{tab:general-table}.

Theorem \ref{thm:expect-generalization-bd} is established in the following steps.
In \texttt{P-ERM}, we first show $\Escr_{\Yscr}(\hat{x}) \leq 2M \E_{\Dscr_n}\E_y[d(\P^*_{\xi|y}, \hat{\Q}_{\xi|y})]$ following the same way as in Theorem~\ref{thm:general-param-erm} and then apply Assumption~\ref{asp:expect-conditional-oracle-estimator}. In \texttt{P-DRO}, 
we consider the event set 
$\Ascr_1 = \{(\Dscr_n, y): d(\P^*_{\xi|y},
\hat{\Q}_{\xi|y}) \leq \varepsilon\}$ and 
its complement. In $\Ascr_1$, $\P^*_{\xi|y}$ is contained in the ambiguity set,
and we can apply the same analysis as in Theorem~\ref{thm:general-param-dro-ipm} to obtain the first term in the generalization error bound in \eqref{eq:p-dro-bound}. In 
the complement event
$\Ascr_1^{c}$, we bound the generalization error term by
\[\E_{\P_{\xi|y}^*}[h(\hat{x}(y);\xi)] - \E_{\P_{\xi|y}^*}[h(x^*(y);\xi)] \leq M d(\P^*_{\xi|y}, \hat{\Q}_{\xi|y}) + 2\varepsilon \Vscr_d(x^*(y)).\] 
Then we can apply Assumption~\ref{asp:expect-conditional-oracle-estimator} to obtain the corresponding error bound.




\subsubsection{Comparison with Existing Literature.}
We compare the generalization error bounds in Corollary~\ref{coro:param-dro-erm-context} and Theorem~\ref{thm:expect-generalization-bd} with the existing contextual optimization literature. The model~\eqref{eq:context-dro} estimates $\hat{\Q}_{\xi|y}$ and then solves the optimization problem. 
To our best knowledge, 
$\hat{\Q}_{\xi|y}$ in the existing contextual DRO literature 
\gi{is obtained from the} empirical joint distribution between the covariate $y$ and the response variable $\xi$ through so-called probability trimming 
\citep{esteban2020distributionally, nguyen2021robustifying} 
or 
estimated using other 
non-parametric approaches that model the conditional distribution or
the noise~\citep{bertsimas2021bootstrap,wang2021distributionally,kannan2020residuals}. 
The generalization error of a
non-parametric 
estimator is 
at least $O(n^{-1/D_{\xi}})$ or $O\left(n^{-1/(D_{\xi} + D_y)}\right)$,
and therefore, these 
approaches 
are not statistically tractable when the 
distributional
dimension is large. In this case, \PDRO\ mitigates the exponential dependence on the dimension of the random variable and pays a controllable price of model misspecification in \eqref{eq:p-dro-bound}.
On the other hand, in end-to-end learning approaches without estimating $\hat{\Q}_{\xi|y}$~\citep{el2019generalization,hu2022fast}, their generalization errors still involve
the term $\text{Comp}(\Hscr)$. The bounds for \PDRO\ in Corollary~\ref{coro:param-dro-erm-context} and Theorem~\ref{thm:expect-generalization-bd} mitigate this complexity
through parametrizing the
distribution and robustification without requiring uniform bounds over the function class.


\section{Discussion of Results}\label{sec:discussion}
We close our theoretical study with several lines of discussion to highlight our strengths, trade-offs, and positioning relative to existing works.
\paragraph{Comparison with Other Related Works.}
Despite the popularity of parametric models in statistics and machine learning, 
such models have not been employed in the DRO literature until recently.
\cite{shapiro2021bayesian} formulate and derive asymptotic
results similar to variability regularization for the so-called 
Bayesian risk
optimization under parametric uncertainty. \cite{michel2021modeling,
  michel2022distributionally} propose ambiguity sets that only contain
parametric distributions,   
but do not provide  
generalization guarantee 
when 
the parametric distribution class is misspecified. Besides, they evaluate model performances
via some robust loss instead of the generalization error that we investigate.  
Relative to these works, 
we focus on the
  generalization error 
  defined as the excess risk 
  over
  the oracle solution. 
We also provide finite-sample theoretical guarantees and demonstrate their potential benefits under problem instances with small $\text{Comp}(\Theta)$ and large $\text{Comp}(\Hscr)$.
Moreover, our framework 
  accommodates most of the commonly used distance metrics including the
  Wasserstein distance and $f$-divergence. Finally, \cite{lam2020parametric} study the required Monte Carlo size to approximate DRO centered at a parametric distribution, but they consider chance-constrained problems that are very different from our current focus. 

\paragraph{Model Selection on Parametric Models.}~Note that the success of the 
\PDRO\ framework hinges on the availability of a 
parametric model 
with low 
$\Escr_{apx}$.
There is a 
rich literature 
on choosing good parametric models including
information-based model selection \citep{anderson2004model} 
and decision-driven parameter calibration \citep{ban2018machine}. 
We do not propose 
new methods for parametric model selection; rather, we leverage  
this 
existing literature. More precisely, one of our major contributions is to propose \texttt{P-DRO} 
to transform the
\texttt{P-ERM} solution obtained from directly using these parametric models 
into consistently better solutions 
by adding robustness. 

\begin{figure}[h] 
    \centering    
    \subfloat[Without Shift] 
    {
        \begin{minipage}[t]{0.45\textwidth}
            \centering          
            \includegraphics[width=0.9\textwidth, trim = 40 40 10 100]{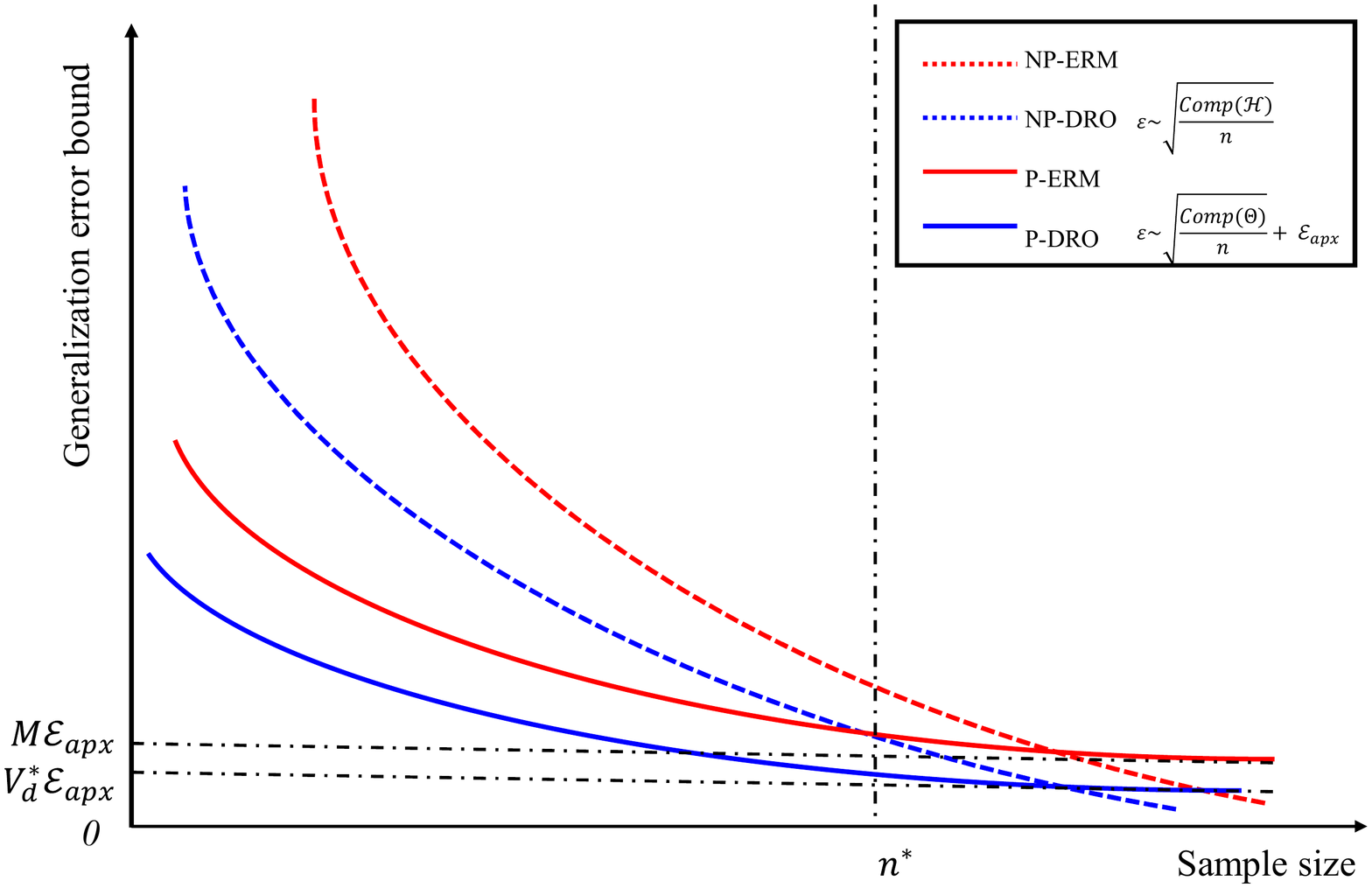}  
        \end{minipage}
    }
    \subfloat[With Shift] 
    {
        \begin{minipage}[t]{0.45\textwidth}
            \centering  
            \includegraphics[width=0.9\textwidth, trim = 40 40 10 100]{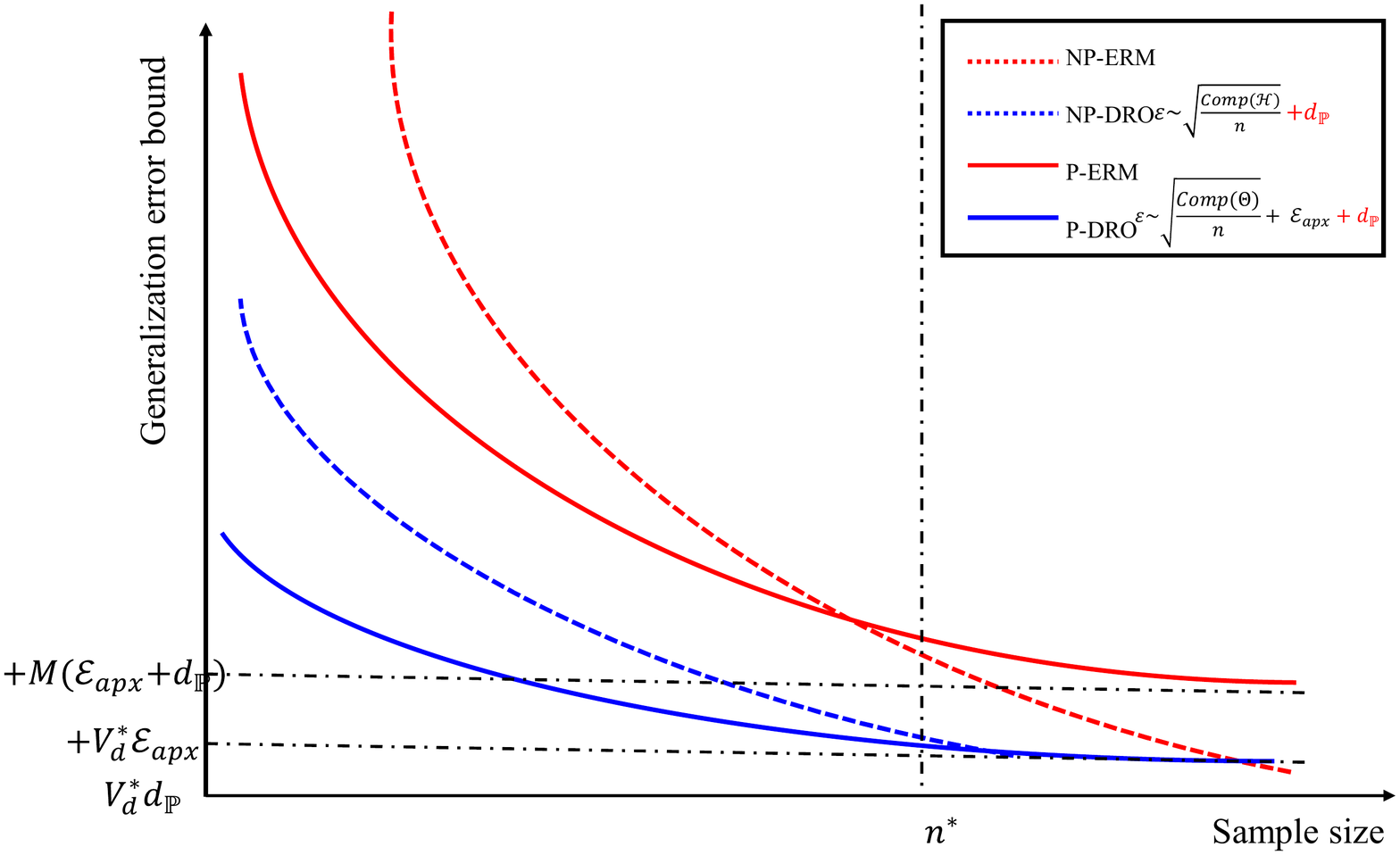} 
        \end{minipage}
    }%
    \caption{
      Generalization error as a function of sample size $n$ for the
      setting where 
      $\text{Comp}(\Hscr)\gg \text{Comp}(\Theta)$ and $d_{\P} =
      d(\P^{tr}, \P^{te}), V_d^* = \Vscr_d(x^*)$. The threshold sample size $n^*$ increases
      with $\text{Comp}(\Hscr)$ 
      and decreases with $\text{Comp}(\Theta)$ and $\mathcal{E}_{apx}$.} 
    \label{fig:concept}  
  \end{figure}

\paragraph{Generalization Error Trade-Off Compared to Existing Methods.}
 We discuss several implications of \PDRO\ regarding generalization.
In Figure~\ref{fig:concept} we plot the generalization errors of ERM and DRO centered at the empirical distribution (which we call \texttt{NP-ERM} and
\texttt{NP-DRO} respectively),
and \PERM\ and \PDRO, as functions of the sample size $n$,
with and without distribution shifts when $\text{Comp}(\Hscr) \gg \text{Comp}(\Theta)$. 
We see that 
when the sample size is not too large,
i.e. $n \leq n^*$, 
the generalization error of 
\PDRO\ 
is lower than the other competing methods. 
In addition, 
under distribution shifts, both \texttt{NP-ERM} and \PERM\ 
are
further negatively impacted by the amplification of the effect of function
class complexity. 
In general, the threshold
sample size $n^*$ increases 
with $\text{Comp}(\Hscr)$ 
while decreases with $\text{Comp}(\Theta)$ and
$\mathcal{E}_{apx}$, in which case \PDRO\ is more effective than existing approaches under a relatively larger sample size. That being said, in the limit $n \rightarrow
\infty$, the error of non-parametric approaches converges to $0$, while
the error of \PDRO\ is lower bounded by $\mathcal{E}_{apx}$. Besides, \PDRO\ is not likely to be very
competitive when 
$\text{Comp}(\Hscr) \approx
\text{Comp}(\Theta)$, and when the parametric class $\Pscr_{\Theta}$ provides
a poor approximation for $\P^\ast$, i.e. $\mathcal{E}_{apx}$ is large.

Overall,
when the true distribution is ``simple'' in that we have a good parametric class to represent the uncertainty, but the objective function is complex, then \PDRO\ yields better performance than existing approaches. Otherwise, it would be better to deploy other non-parametric approaches. In practice, this trade-off to decide whether to employ \PDRO\ can be made through cross-validation.

\paragraph{Generalization via Regularized ERM.} Besides DRO, regularized ERM approaches  with better generalization error bounds than standard ERM have been investigated, ranging from $\ell_1$-regularization \citep{koltchinskii2011oracle,bartlett2012ell} to generalized moment penalization \citep{foster2019orthogonal} including variance regularization \citep{maurer2009empirical,xu2020towards}. Some of these regularization approaches are equivalent to DRO \citep{shafieezadeh2019regularization,duchi2019variance,gao2022wasserstein}. Note that we may attain a faster generalization error rate with $\alpha = 1$ in~\eqref{eq:general-bd} under suitable curvature or margin conditions between $\Hscr$ and $\P^*$ for the obtained solution~\citep{liang2015learning,zhivotovskiy2018localization}. However, in general, we cannot remove the dependence of $\text{Comp}(\Hscr)$ on $B$ with $\alpha = \frac{1}{2}$ in~\eqref{eq:general-bd}. To see this, nearly all of the existing regularized ERM approaches build on the argument of uniform convergence in the hypothesis class. That is, with probability $1-\delta$,
\begin{equation}\label{eq:uniform-bd}
\E_{\P^*}[h(x;\xi)]  \leq \E_{\hat{\P}_n}[h(x;\xi)] + \sqrt{\frac{\text{Comp}(\Hscr_r) \log n + \log (1/\delta)}{n}} \hat{V}(x), \forall x \in \Xscr_{r}.
\end{equation}
for some problem-dependent measure $\hat{V}(x)$, e.g., $\hat{V}(x) = \sqrt{\text{Var}_{\hat{\P}_n}[h(x;\xi)]}$. Here $\Xscr_r$ is a subset of $\Xscr$ built on some localized arguments, for example localized Rademacher Complexity with the form $\Xscr_r = \{x \in \Xscr: \hat{V}(x) \leq r\}$~\citep{bartlett2005local} and $\Hscr_r = \{h(x;\cdot)| x \in \Xscr_r\}$. Then a regularized ERM formulation with $\hat{Z}(x) := \E_{\hat{\P}_n}[h(x;\xi)] + \lambda \hat{V}(x)$ (for some $\lambda > 0$) can achieve a generalization error bound with $O\left(\hat{V}(x^*) \sqrt{\frac{\text{Comp}(\Hscr_r)}{n}}\right)$ following~\eqref{eq:uniform-bd}. Although this bound is better than the standard ERM bound since it replaces the uniform term $M$ by $\hat{V}(x^*)$ from the problem-dependent regularization, the complexity term involving $\Hscr$ or a refined $\Hscr_r$ cannot be removed. This is because the key to the generalization error in these methods still involves a uniform convergence argument to bound $\sup_{x \in \Xscr_r}|\E_{\P^*}[h(x;\xi)]  - \E_{\hat{\P}_n}[h(x;\xi)]|$, as demonstrated in, e.g., \cite{xu2020towards}. Thus in practice, similar to DRO using the regularization perspective, regularized ERM approaches may not yield good performance when the hypothesis class complexity is high.

In contrast, the crux of our argument relies on the distribution coverage perspective that with probability $1-\delta$,
\begin{equation}\label{eq:uniform-dist-bd}
\E_{\P^*}[h(x;\xi)] \leq \sup_{d(\P, \hat{\Q})\leq \varepsilon}\E_{\P}[h(x;\xi)], \forall x \in \Xscr.
\end{equation}
Here, \eqref{eq:uniform-dist-bd} holds regardless of the complexity of $\Hscr$ as long as $d(\P^*, \hat{\Q}) \leq \varepsilon$. The choice of $\varepsilon$ is also independent of $\text{Comp}(\Hscr)$, leveraging Assumption~\ref{asp:oracle-param-est}. This avoids the uniform concentration argument applied to all $x \in \Xscr_r$ in~\eqref{eq:uniform-bd} and leads to better performance bounds of \PDRO\ when $\text{Comp}(\Hscr)$ is large.

\paragraph{Minimax versus Instance-Dependent Rate.} 
The hypothesis class complexity $\text{Comp}(\Hscr)$ cannot be improved in the minimax sense.
More formally, it is known from the machine learning literature (e.g., Chapter 19 in \cite{wainwright2019high} and Section 5.5 in \cite{boucheron2005theory}) that 
for any 
solution $\hat{x}$ which is a function of the data 
$\Dscr_n$ and a
\emph{large} class of 
distributions $\Pscr$, we have
\begin{equation}\label{eq:minimax}
  \lim\inf_{n \to \infty}\inf_{\hat{x}}\sup_{\P \in \Pscr}\sqrt{n}
  \E_{\Dscr_n}\left(\E_{\P}[h(\hat{x};\xi)] -
    \E_{\P}[h(x^*;\xi)]\right) = \sqrt{\text{Comp}(\Hscr)}. 
\end{equation}
This means that $\text{Comp}(\Theta)$ controls 
the minimax error decay when the data-driven solution $\hat{x}$ is only a function of the data $\Dscr_n$ and $h(x;\cdot)$. However, this rate is based on a large class of $\Pscr$, which can be too conservative and does not take into account the information of $\P^*$ in real-world instances.
Therefore, this minimax rate does not contradict our results that $\text{Comp}(\Hscr)$ can be improved to a smaller term $\text{Comp}(\Theta)$ in some cases (e.g., when $\Escr_{apx} \approx 0$), and in fact, improving $\text{Comp}(\Hscr)$ is exactly the motivation for our proposed \PDRO. 

On one hand, 
our results can be seen as an instance-dependent generalization error with a smaller distribution class $\Pscr$. For example, with $\Pscr = \Pscr_{\Theta}$ in Assumption~\ref{asp:oracle-param-est}, 
we attain the rate
$O(\sqrt{\text{Comp}(\Theta)/n})$ for $\hat{x}^{\PDRO\ }$ in Theorem \ref{thm:general-param-dro-ipm}, which can be much smaller than the minimax lower bound $O(\sqrt{\text{Comp}(\Hscr)/n})$ in~\eqref{eq:minimax}. On the other hand, 
the oracle estimator and associated parametric class in
Assumptions~\ref{asp:oracle-param-est}, ~\ref{asp:conditional-oracle-estimator}
and~\ref{asp:expect-conditional-oracle-estimator} reflect the decision
makers' belief regarding the distribution, and our data-driven solution $\hat{x}$ is a function of the data $\Dscr_n$, $h(x;\cdot)$ and $\Pscr_{\Theta}$ which could lead to better results under good $\Pscr_{\Theta}$.

\paragraph{Asymptotic Comparisons.}~While we have demonstrated the strengths of \PDRO\ against classical ERM
and DRO under finite sample, 
in the large-sample regime it is known that 
ERM is optimal 
\citep{lam2021impossibility},  and 
the generalization error of non-parametric DRO 
for a
proper
choice of $\varepsilon_n$ is $O_p(1/\sqrt{n})$
\citep{blanchet2019robust,duchi2021statistics}. 
Our \PDRO\ does not beat these latter methods asymptotically as the sample size grows large, which is also indicated from Figure~\ref{fig:concept}. This is because under large sample, 
the true distribution $\P^*$ can be learned 
well from data, and there is no need to parametrize 
distributions. 
Nonetheless, the situation 
can be very 
different 
when the sample size is small as we have shown. 


\section{Numerical Studies}\label{sec:numerics}


We 
compare the performances of 
both the non-contextual and contextual versions of \PDRO\
with non-parametric ERM (\texttt{NP-ERM}) and 
non-parametric DRO (\texttt{NP-DRO}), 
on 
synthetic and real-world datasets.  
We set the ambiguity
size $\varepsilon$ through cross-validation in DRO methods, and the
Monte
Carlo size $m = 50 n$ 
(unless noted otherwise). 


\subsection{Synthetic Example}\label{subsec:synthetic}
We consider the problem of minimizing the following objective: 
\begin{equation}
  \label{eq:high-order-simulation}
h_{\gamma}(x;\xi) = \left(\mu - \xi^{\top} x\right)_{+}^{\gamma} := \big|\min\{0, \xi^{\top}x - \mu\}\big|^{\gamma},   
\end{equation}
where $\gamma >0$, $\xi$ denotes the random  asset return, $\mu$ is a specified
deterministic  target return and the vector $x$ denotes the
allocation weights. The feasible set $\Xscr = \{\sum_{i} x_i = 1, x_i \geq
-\tau\}$ with $\tau > 0$. 
The objective $h_{\gamma}(x;\xi)$ is called the downside risk when $\gamma = 2$
~\citep{sortino2001managing}. Here
$\text{Comp}(\Hscr)$ grows with $\tau$ and $\gamma$, and is provably
large. We vary $\tau  \in \{2, 10\}$ and $\gamma \in \{1,2,4\}$, 
and define the parametric family of distributions
\begin{equation}
  \label{eq:P-beta-def}
  \Pscr_{\Theta} = \left\{\P: \xi  = (\xi_1,\ldots,\xi_{D_{\xi}}), \xi_i \sim
  2r\times \text{Beta}(\eta_i, 2) - r, \eta_i \in [1.5, 3],\ \forall i \right\}, 
\end{equation}
for a given fixed constants $r$ with unknown $\eta$.  
We use
$\chi^2$-divergence as the DRO metric here. 


Assumption~\ref{asp:oracle-param-est} holds for $\hat{\Q}$ when $\P^* \in
\Pscr_{\Theta}$ and 
$\eta$ in \eqref{eq:P-beta-def} is estimated 
using the moment method with
$\text{Comp}(\Theta) = O(D_{\xi})$ and $\text{Comp}(\Hscr) = O(D_{\xi}^{\gamma})$ when $n$ is large, and
$\alpha = 1, \Escr_{apx} = 0$ in \eqref{eq:param-est-ineq} (we provide more details in Appendix~\ref{subsec:comp-theta-h}). In addition, 
$M =\sup_{x \in \Xscr}\|h(x;\cdot)\|_{\infty} \leq (D_{\xi} \tau r +
\mu)^{\gamma}$ and $M^* = \sqrt{\text{Var}_{\P^*}[h(x^*;\cdot)]} \leq
\|h(x^*;\cdot)\|_{\infty} \leq (r \|x^*\|_{1} + \tau)^{\gamma}$. 
Generalization error  bounds for the competing methods are presented in
Table~\ref{tab:general-synthetic-bound}. 

\begin{table}[t]
  \centering
  \begin{center}
    \caption{Generalization error bounds for portfolio selection on
    synthetic data where \PERM\ and \PDRO\ are fitted using the Beta-family in \eqref{eq:P-beta-def}.}
    \label{tab:general-synthetic-bound}
    \begin{tabular}{c|c|c|c|c}
      \hline
      Method &\texttt{NP-ERM}
      &\texttt{P-ERM}&\texttt{NP-DRO}&\textbf{\texttt{P-DRO}}\\ 
      \hline 
     $\mathcal{E}(\hat{x})$&$M \sqrt{\frac{(D_{\xi}+\gamma)^{\gamma}\log
                             M}{n}}$&$M
                                      \sqrt{\frac{D_{\xi}}{n}}$&$M^*\sqrt{\frac{(D_{\xi}+\gamma)^{\gamma}
                                                                 \log
                                                                 M}{n}}$&$M^*\sqrt{\frac{D_{\xi}}{n}}$\\
      \hline
                                                                        
    \end{tabular}
  \end{center}
\end{table}

\begin{figure}[h] 
    \centering    
    \subfloat[Base Case] 
    {
        \begin{minipage}[t]{0.5\textwidth}
            \centering  
            \includegraphics[width=0.9\textwidth, trim = 100 10 40 60]{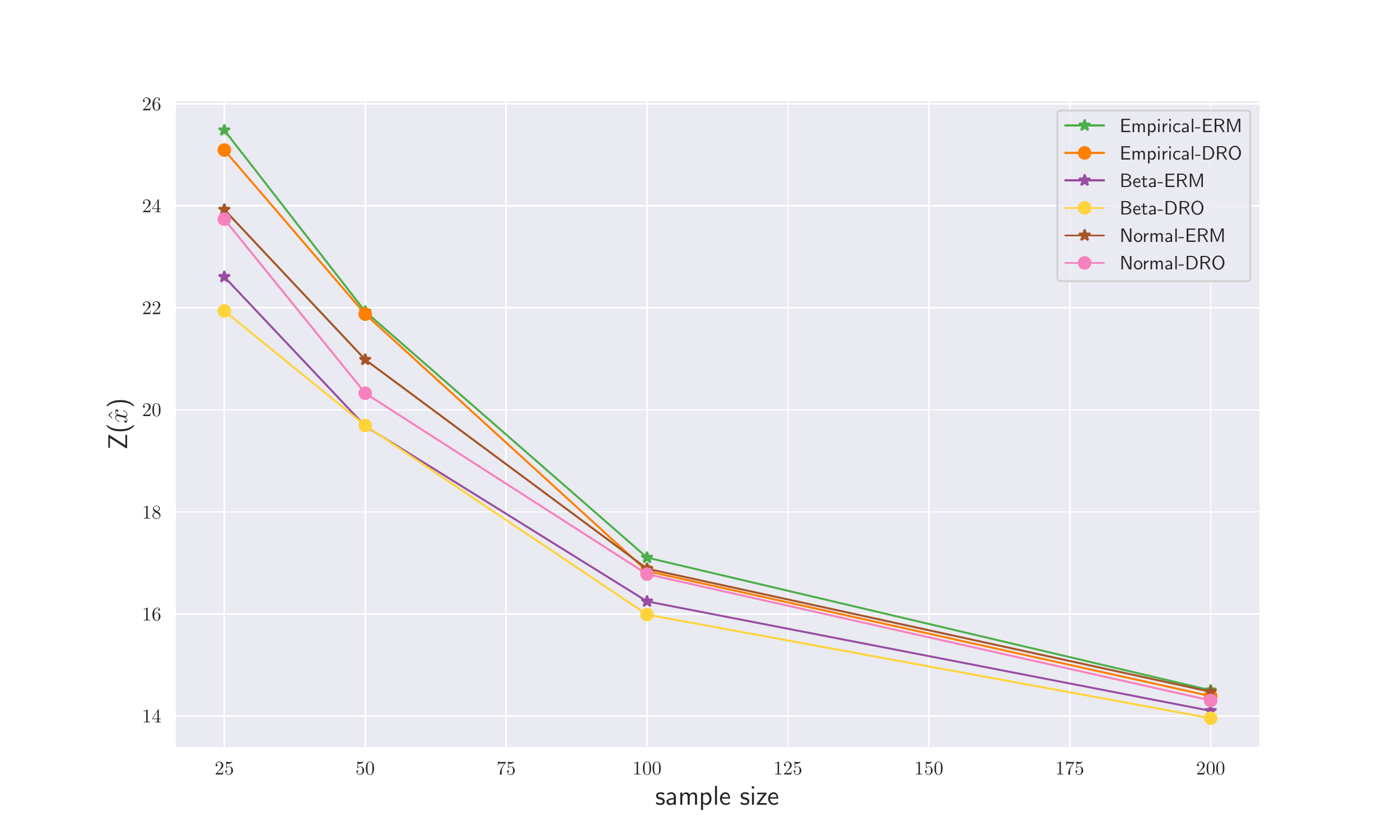} 
        \end{minipage}
    }%
    \subfloat[Distribution Shift]
    {
        \begin{minipage}[t]{0.5\textwidth}
            \centering      
            \includegraphics[width = 0.9\textwidth, trim = 100 10 40 60]{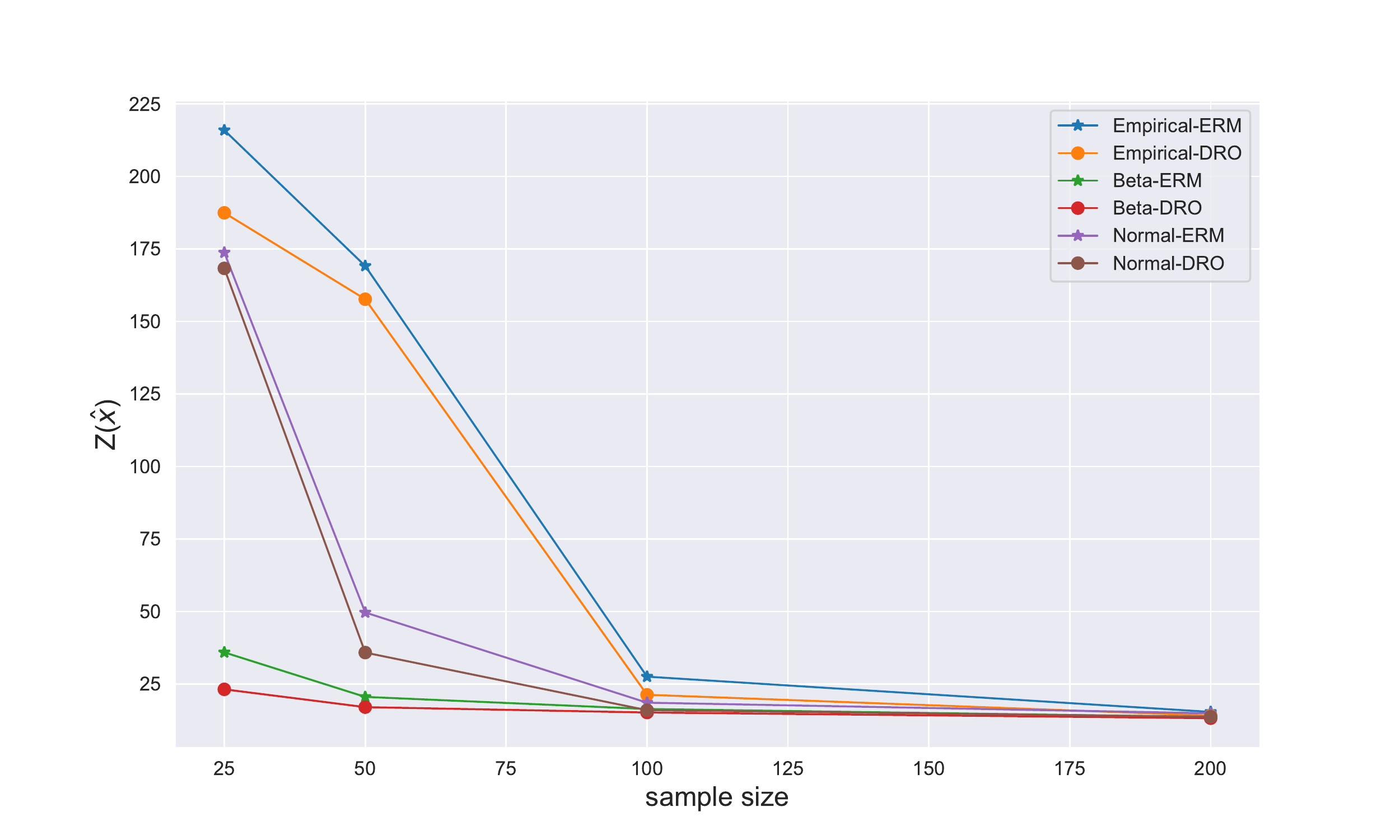}  
        \end{minipage}
    }%
    \caption{Average $Z(\hat{x})$ across different ERM-DRO models varying $n$ with $(\gamma, \tau) = (2,2)$.}
    \label{fig:simu-concept}  
  \end{figure}

In Figure~\ref{fig:simu-concept}~(a), the plots marked ``Empirical-*''
correspond to using $\hat{\Q} = \hat{\P}_n$, the ones marked ``Beta-*''
correspond to 
$\hat{\Q}$ fit using the Beta family defined in \eqref{eq:P-beta-def}, and  
the plots marked ``Normal-*'' correspond to fitting a multivariate Gaussian model to $\xi$. 
Since $\text{Comp}(\Hscr) \gg \text{Comp}(\Theta)$ when $\gamma$ is large for this example from the comparison in Table~\ref{tab:general-synthetic-bound} and Appendix~\ref{subsec:comp-theta-h}, 
we expect and do see that 
parametric models 
outperform their
non-parametric
counterparts. 
Although \texttt{P-DRO} statistically outperforms \texttt{P-ERM}
for all 
sample sizes, 
the absolute margins between \texttt{P-DRO} and \texttt{P-ERM} are not
obvious 
under
large sample sizes. 
In Figure~\ref{fig:simu-concept}~(b), we plot the results when 
distribution shifts  occur.
Here we find that \PDRO\ significantly 
outperforms all
other models. 
We present results for other values of 
$(\gamma, \tau)$ in 
Appendix~\ref{subsec:synthetic2}, which show that the  
performance gain of
\PDRO\ 
grows with
$(\gamma,
\tau)$. 

  \subsection{Synthetic Example Extension: Contextual
  Optimization}\label{subsec:synthetic-context} 


We report the results of our numerical experiments with
contextual DRO on synthetic data. 
We consider 
objective \eqref{eq:high-order-simulation} with $\gamma = 2$ and $\tau =
10$. However, the 
return distribution is now given by an exogenous variable $y$:
\[
  \xi|y := B y + g(y) + u,
\]
where $B \in \R^{D_{\xi}\times D_y}$, and each element $b_{ij}$ of the
matrix $B$ is drawn i.i.d.
$U(-0.5, 0.5)$ in the high signal-to-noise-ratio~(SNR) case, 
and i.i.d. $U(-0.1, 0.1)$
in the low SNR case, 
$u$  is a noise term independent of the covariate $y$, and $g(y) = 2
\sin(\|y\|_2)$. We parametrize the distribution by $\xi | y = By + u$ with normal noise $u$, and the 
term $g(y)$ represents a deterministic misspecification of the
distribution depending on the covariate $y$.
We use the Fama-French 3-factor model~\citep{FamaFrench}
for
the covariate $y$ 
and set the DRO metric
$d = \chi^2$-divergence.

Results of our numerical experiments are summarized in
Table~\ref{table:context-result}. In this table, 
``noncontext-p" denotes results for the Normal-*  models
used in Section~\ref{subsec:synthetic}  that  ignore the covariates, 
the columns labeled 
``kernel-*" correspond to 
$\hat{\Q}_{\xi|y}$ by the Nadaraya-Watson non-parametric estimator 
\citep[see][]{bertsimas2021bootstrap}, 
the columns labeled ``residual-p'' correspond to $\hat{\Q}_{\xi|y}$ by the empirical residual estimator 
\citep[see][]{kannan2020residuals}, 
and the columns 
labeled ``context-p''
correspond to $\hat{\Q}_{\xi|y} = N(\hat{B} y, \hat{\Sigma})$, where 
$\hat{B} \in \argmin_{B \in
  \R^{D_y\times D_{\xi}}}\sum_{i = 1}^n\|\hat{\xi}_i - B \hat{y}_i\|_2^2$
and 
$\hat{\Sigma} = \frac{1}{n}\sum_{i = 1}^n (\hat{\xi} -
\hat{B}\hat{y}_i)(\hat{\xi} - \hat{B}\hat{y}_i)^{\top}$. 
For each DRO model, we consider the same tuning
procedure for the ambiguity size as in Section~\ref{subsec:synthetic}.   
We show the result for the average cost function value for each model
across $50$ independent runs (each with a different random seed) 
and varying the
sample size $n$, high/low SNR, and including/excluding the misspecified term $g(y)$. 
We find that 
for all scenarios, except one, 
context-p DRO model 
has the best performance. 
Notably, the robustness in the context-p DRO model helps
mitigate the effects of model misspecification error in the regression and this leads to a significantly improved performance over the 
context-p ERM model. 
When the number of samples is limited, 
the context-p DRO model outperforms the 
residual-DRO and kernel-DRO models which are non-parametric, in line with our theoretical implications.
The results of this set of numerical experiments clearly illustrate the
power of contextual \PDRO.

\begin{table}[t]
  \footnotesize
 \caption{Comparison of avg $h$ under different models in each subcase, the first line representing average performance, and the second line standard deviation. Boldfaced values mean that the corresponding approach is the best in the considered setting.}
\label{table:context-result}
  \centering
  \begin{tabular}{ccc|cccc|cccc}
    \toprule
-&-&-&\multicolumn{4}{c}{ERM}&\multicolumn{4}{c}{DRO}\\
$n$ & $snr$ & mis &kernel&noncontext-p&residual&context-p&kernel&noncontext-p&residual&context-p\\
\midrule
50  & high & No   & 2490.6 & 34037.5 & 116.3  & 611.2  & 456.3 & 22060.1 & 36.2   & \textbf{27.4}   \\
    &      &       & 365.1  & 7027.0  & 50.8   & 309.1  & 63.2  & 3948.0  & 14.7   & 12.8   \\
50  & high & Yes & 2502.6 & 21907.4 & 122.1  & 450.4  & 936.7 & 12961.8 & 41.8   & \textbf{20.4}   \\
    &      &       & 373.1  & 2308.7  & 52.4   & 264.5  & 108.2 & 1980.9  & 19.1   & 7.7    \\
50  & low  & No   & 1194.1 & 8552.4  & 1544.1 & 1309.5 & 394.1 & 2726.7  & 318.5  & \textbf{152.7}  \\
    &      &       & 165.4  & 2358.7  & 254.5  & 257.9  & 66.7  & 863.2   & 101.2  & 28.4   \\
50  & low  & Yes & 1671.2 & 13437.1 & 2313.5 & 6656.4 & 614.2 & 4513.7  & 1067.6 & \textbf{543.4}  \\
    &      &       & 212.7  & 1531.5  & 312.3  & 1235.5 & 79.4  & 600.9   & 186.1  & 101.9  \\
100 & high & No   & 1597.5 & 17012.8 & 62.4   & 108.4  & 462.9 & 5813.7  & \textbf{1.6}    & 2.0    \\
    &      &       & 221.8  & 5311.1  & 52.1   & 93.8   & 59.0  & 1314.0  & 1.2    & 1.3    \\
100 & high & Yes & 1605.1 & 10231.8 & 63.6   & 68.7   & 493.8 & 4375.0  & 37.0   & \textbf{0.5}    \\
    &      &       & 228.2  & 1889.2  & 52.7   & 51.7   & 78.1  & 788.0   & 33.0   & 0.5    \\
100 & low  & No  & 1194.1 & 8552.4  & 1544.1 & 1309.5 & 394.1 & 2726.7  & 318.5  & \textbf{152.7} \\
    &      &       & 165.4  & 2358.7  & 254.5  & 257.9  & 66.7  & 863.2   & 101.2  & 28.4   \\
100 & low  & Yes & 1196.8 & 4105.6  & 1547.3 & 1750.7 & 527.1 & 1385.8  & 773.6  & \textbf{147.5}  \\
    &      &       & 164.8  & 640.5   & 252.1  & 359.1  & 50.0  & 172.4   & 156.3  & 32.6  \\
\bottomrule
  \end{tabular}
\end{table}
  
\subsection{Real Data I: Portfolio Optimization}\label{subsec:real1}
We report the results of our numerical experiment with
portfolio allocation on real data. We continue to
use objective~\eqref{eq:high-order-simulation} with $\gamma = 2$ and
$\tau\in \{2, 10\}$. We use the Fama-French data \citep{FamaFrench}
with $D_{\xi} \in \{6,10,25,30\}$. Note that the asset
returns are neither stationary nor generated from some simple parametric
families. Therefore, any approach would face the problem of distribution
shift and model misspecification. We compare \PDRO\ against 
benchmarks using the 
``rolling-sample" approach to
estimate the cost 
$\hat{h}=\frac{1}{N}\sum_{i =
  1}^N(\mu -\hat{r}_i)_+^2$ with $N$ out-of-sample returns $\{\hat{r}_i\}_{i
  = 1}^n$. We still fit parametric models with 
Beta 
and Normal distributions. We provide more details on our setup in Appendix~\ref{app:experiment-real-portfolio}.

\begin{table}[t]\small
\centering
\caption{
  Performances of 
  different models for
  the portfolio allocation problem with $\tau = 2$. The quantity 
  in the bracket is the
  the empirical cost $\hat{h}$ for $\tau =
  10$ as a multiple of the cost for $\tau = 2$. 
  Here $^+$ means the DRO model 
  outperforms the ERM counterpart, and  
  $^*$
  means \texttt{P-DRO} outperforms \texttt{NP-DRO} up to statistical
  significance of $p$-value $<0.001$. Boldfaced values mean that the corresponding approach is the best in the considered setting.
} 
\label{table:empirical-portfolio-1}
\begin{tabular}{c|cc|cc|cc}
\toprule
                 & \multicolumn{2}{c}{empirical} & \multicolumn{2}{c}{Beta} & \multicolumn{2}{c}{Normal}  \\
dataset / method & ERM        & DRO             & ERM     & DRO            & ERM     & DRO               \\
\midrule
10-Industry &36.26 (1.00)   & 33.35 (1.00)$^{+}$  & 31.27 (1.00) & \textbf{30.64 (1.00)} & 35.4 (1.00)   & 31.88 (1.00)$^{+}$   \\
6-FF&28.91 (1.02)  & 27.98 (1.01) & 35.93 (1.00) & 35.81 (1.00) & 28.75 (1.01) & \textbf{27.93 (1.00)}   \\
30-Industry & 210.07 (9.97) & 195.1 (9.58)$^{+}$ & 35.26 (1.00) & \textbf{34.33 (1.00)}$^{*}$ & 84.58 (1.03) & 62.06 (1.01)$^{+*}$  \\
25-FF& 60.86 (2.90)   & 53.39 (2.94)$^{+}$ & 37.62 (1.00) & \textbf{36.94 (1.00)}$^{*}$ & 48.58 (1.11) & 37.41 (1.04)$^{+*}$ \\
\bottomrule
\end{tabular}
\end{table}

The results summarized in
Table~\ref{table:empirical-portfolio-1}  show that the parametric models
(Beta, Normal) 
outperform the non-parametric methods, 
especially when 
$D_{\xi}$ is large. The performances of \texttt{NP-DRO} / \texttt{NP-ERM} are
very sensitive to the choice of $\tau$, and are significantly dominated
by parametric approaches in this regime. \texttt{P-DRO} can reduce the
problem of misspecification from \texttt{P-ERM}. 
We find that
the
Beta parametric models 
generate generally better decisions 
compared to the Normal parametric models. 


\subsection{Real Data II: Regression on LDW Data}\label{subsec:real2}
We report the result of using \PDRO\ on a 
regression problem.  
We work with the 
PSID 
data set \citep{dehejia1999causal} 
that contains 
$8$ features and $n = 2490$ samples, and the goal is to predict the household
earning using these features. 
The DRO problem in this setting is given by 
\begin{equation}\label{eq:ml-task}
\min_{h\in \Hscr} \max_{d(\P, \hat{\Q}) \leq \varepsilon}\E_{(x,y) \sim \P}[\ell_2(y; h(x))],   
\end{equation}
where $\Hscr$ 
is the set of all quadratic polynomials in $x$ and $\ell_2$ is the squared loss. For DRO models, we choose the distance
metric $d$ to be 2-Wasserstein distance with the definition deferred to Definition~\ref{def:p-was} in Appendix~\ref{app:p-was-general}.
We set $\Pscr_{\Theta}$ as the mixture of jointly Gaussian distributions 
$(x,y)$, where each component of the Gaussian mixture model represents individuals with one possible choice from a list of binary categories (\textsf{black}, \textsf{Hispanic}, \textsf{married} and \textsf{nondegree}) in $x$ and there are 16 components in all. For example, married white Hispanic individuals with degrees and unmarried white Hispanic individuals without degrees are two components.  Further details on our setup can be found in Appendix~\ref{app:experiment-real-reg}. In this case, $\text{Comp}(\Theta) = O(D_{\xi}^2)$ and $\text{Comp}(\Hscr) = O(D_{\xi}^4)$. We 
report
out-of-sample
$R^2$ (note that under the squared loss, $Z(\hat{x}) = K(1-R^2)$ for some
  $K > 0$) for different methods averaged over $50$ independent runs, each run with a random train-test-split using a different random seed, for each training sample size. 

\begin{figure}[h] 
    \centering    
    \subfloat[Without Shift]
    {
        \begin{minipage}[t]{0.5\textwidth}
            \centering          
            \includegraphics[width=0.9\textwidth, trim = 100 10 60 60]{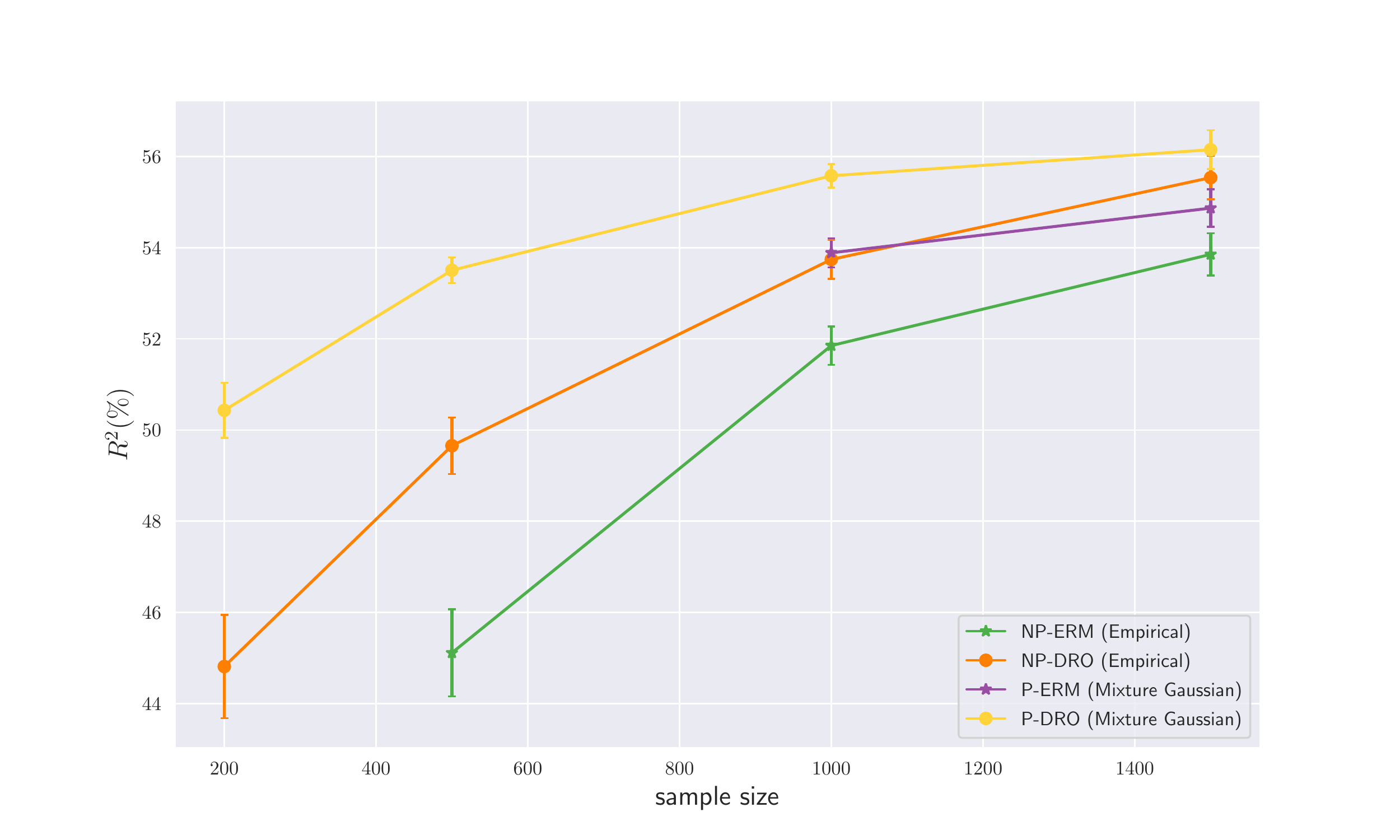}  
        \end{minipage}
    }
    \subfloat[With Age Shift] 
    {
        \begin{minipage}[t]{0.5\textwidth}
            \centering  
            \includegraphics[width=0.9\textwidth, trim = 100 10 60 60]{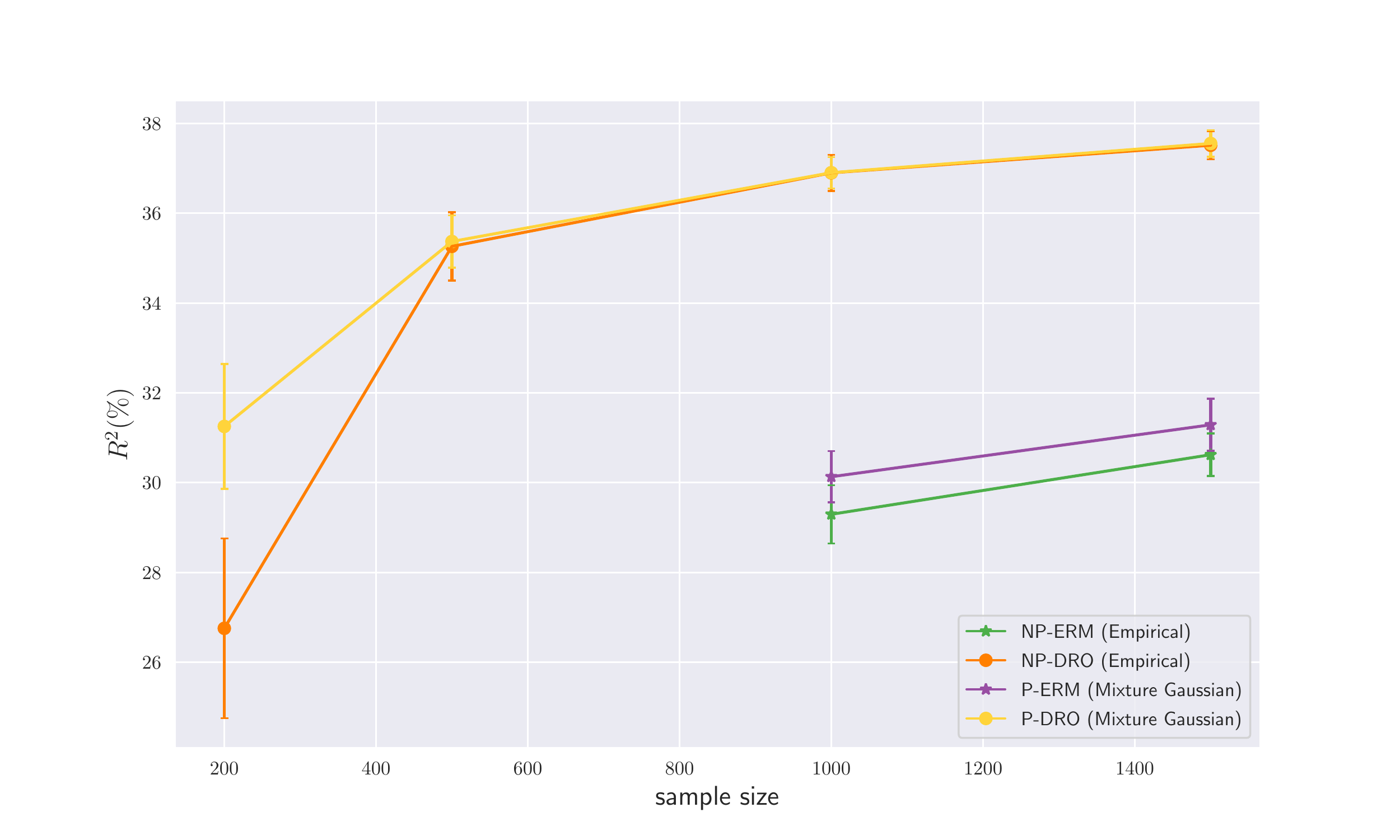} 
        \end{minipage}
    }%
    \caption{Comparison of average $R^2$ (\%) under different models in PSID Datasets. \texttt{P-DRO} is statistically more significant than \texttt{NP-DRO} ($p$-value $<0.001$) in all sample sizes without shifts and $n = 200$ under age shifts.}
    \label{fig:LDW}  
\end{figure}

Figure~\ref{fig:LDW} (a) 
shows that ERM models are dominated by DRO models, especially when the
sample size is not too large. Moreover, the performance of \texttt{P-DRO}
is superior to \texttt{NP-DRO} 
when the number of samples is small. 
In the case of distribution shift, we consider one type of marginal
distribution shifts on the feature vector in Figure~\ref{fig:LDW} (b). We model distribution shifts 
by training on individuals who 
are above $25$, but testing the model on individuals below $25$. We also tune
the ambiguity size $\varepsilon$ in \texttt{P-DRO} and \texttt{NP-DRO} from a small separate validation dataset sampled from the test distribution to
approximate the extent of distribution shifts. Under such case, the performance of \texttt{P-DRO} is slightly
better than \texttt{NP-DRO} but the difference is not statistically significant. Nonetheless, both of these models have significantly superior results than ERM models. We also present further results showing the persistence of the outperformance of \PDRO\ over \PERM\ under different parametric models at the end of Appendix~\ref{app:experiment-real-reg}.

Finally, we present an additional synthetic numerical example in Appendix~\ref{subsec:synthetic1} to further illustrate the benefits of \PDRO\ over \PERM\ in replacing the term $\sup_{x \in \Xscr}\Vscr_d(x)$ with $\Vscr_d(x^*)$ in the corresponding generalization error bounds.

\section{Conclusion and Future Directions}
In this paper, we studied a data-driven stochastic optimization framework named \PDRO, by setting the ambiguity set center of a DRO with a suitably fit parametric model. Our investigation is motivated by the challenges faced by existing ERM and DRO methods that have generalization performance degradation either when the problem dimension is high or when the cost function is complex. We showed how \PDRO\ exhibits better generalization performance under high-dimensional complex-cost problems, at the expense of a model misspecification error that is alleviated via the worst-case nature of DRO.
Our \PDRO\ hinges on and leverages the abundant literature on parametric model selection and fitting, and we showed how \PDRO\ can be generally solved with Monte Carlo sampling of the fit parametric distribution to reduce to conventional nonparametric DRO with similar generalization guarantees. Furthermore, we showed how \PDRO\ can be extended to distribution shifts and contextual optimization. In particular, we demonstrated the additional benefit of \PDRO\ in reducing the amplification of the model misspecification error from a factor that depends on the decision space size in nonparametric counterparts to one that only involves the cost function evaluated at the ground-truth solution. Future directions include further investigations on the interplay between nonparametric and parametric formulations in DRO and other data-driven optimization approaches, and the design of ambiguity sets incorporating other prior knowledge that has similar effects as parametric information.

\ACKNOWLEDGMENT{We gratefully acknowledge support from the InnoHK initiative, the Government of the HKSAR, and Laboratory for AI-Powered Financial Technologies.
}





\bibliographystyle{informs2014}
\bibliography{ref,response_ref,ref_cond}

\input{ECompanion}
\end{document}

%% file: ECompanion.tex
\ECSwitch

\ECHead{Proofs of Statements and Additional Experiments}


\section{Basic Lemmas for the Considered Distances}
We require the following standard measure concentration results of 1-Wasserstein distance and $f$-divergence in our analysis. 
\begin{lemma}[Measure concentration, from Theorem 2 in \cite{fournier2015rate}]\label{lemma:measure-concentration}
Suppose $\P^*$ is a light-tailed distribution such that $A:=\E_{\P^*}[\exp(\|\xi\|^a)] < \infty$ for some $a > 1$. Then there exist some constants $c_1, c_2$ only depending on $a, A$ and $D_{\xi}$ such that $\forall \delta \geq 0$, if $n \geq \frac{\log(c_1/\delta)}{c_2}$, then $W_1(\P^*, \hat{\P}_n) \leq \left(\frac{\log(c_1/\delta)}{c_2 n}\right)^{1/\max\{D_{\xi}, 2\}}.$
\end{lemma}


The following inequality shows that $\chi^2$-divergence, modified $\chi^2$-divergence ($f(t) = \frac{(t - 1)^2}{2t}$), $KL$-divergence and $H^2$-distance satisfy~\eqref{eq:general-f-tv-ineq}.
\begin{lemma}[Pinsker's inequality]\label{lemma:pinsker-ineq}
For any two distributions $\P, \Q$, under our definitions of $f$-divergences in the main body, we have:
$$d_{TV}(\P,\Q)\leq\sqrt{H^2(\P, \Q)} \leq\sqrt{\frac{1}{2}KL(\P, \Q)}\leq\frac{\sqrt{\chi^2(\P, \Q)}}{2}, d_{TV}(\P, \Q) \leq \frac{\sqrt{\chi^2(\Q, \P)}}{2}.$$
\end{lemma}
The following result shows that the standard and modified $\chi^2$-divergences can be represented in a similar form to IPMs in the main body with $\Vscr_d(x) = \sqrt{\text{Var}_{\P^*}[h(x;\xi)]}$ when $d$ is taken as $\chi^2$-divergence.
\begin{lemma}[Pseudo IPM property for $\chi^2$-divergence]\label{lemma:chi2-cauchy}
For distributions $\P, \Q$, under our definitions of $\chi^2$-divergence $\chi^2(\P, \Q)$ and modified $\chi^2$-divergence $\chi^2(\Q, \P)$, we have:
$$\bigg|\E_{\xi\sim\P}[g(\xi)] -  \E_{\xi\sim\Q}[g(\xi)]\bigg| \leq \sqrt{2\min\{\chi^2(\P, \Q)\text{Var}_{\xi \sim \P}[g(\xi)], \chi^2(\Q, \P)\text{Var}_{\xi \sim \Q}[g(\xi)]\}}.$$
\end{lemma}

\textit{Proof of Lemma~\ref{lemma:chi2-cauchy}.}~This result follows directly from the definition of $\chi^2$-divergence and the Cauchy-Schwarz inequality. Denote $M^* = \E_{\Q}[g(\xi)]$. Then we have:
\begin{equation*}
\begin{aligned}
  \E_{\P}[g(\xi)]- \E_{\Q}[g(\xi)] =  \E_{\Q}\left[\left(\frac{d \P}{d\Q} - 1\right)\left(g(\xi)-M^*\right)\right]&\leq \sqrt{\E_{\Q}\left(\frac{d \P}{d \Q} - 1\right)^2}\sqrt{\text{Var}_{{\Q}}[g(\xi)]}\\
  &= \sqrt{2\chi^2(\P, \Q) \text{Var}_{\Q}[g(\xi)]}.
\end{aligned}    
\end{equation*}
\begin{equation*}
\E_{\P}[g(\xi)]- \E_{\Q}[g(\xi)] = \E_{\Q}\left[\left(\frac{d \P}{d\Q} - 1\right)\left(g(\xi)-M^*\right)\right]\geq -\sqrt{2\chi^2(\P, \Q) \text{Var}_{\Q}[g(\xi)]}.
\end{equation*}
The other side follows from considering the term $\E_{\Q}[g(\xi)] - \E_{\P}[g(\xi)]$. $\hfill \square$

Following this result, all properties in our derived generalization error bounds exhibit the same behavior for $\chi^2$-divergence and modified $\chi^2$-divergence.
\vspace{0.2cm}

In the following, if not especially noted, all $C$ with different superscripts and subscripts are denoted as some constants independent of problem-dependent complexity terms.

\section{Details of Existing ERM and DRO Approaches}\label{app:erm-dro-derive}
In the following, we detail some key results and examples of how existing standard ERM and DRO approaches are derived in Table~\ref{tab:general-table} in Section~\ref{sec:background}. We ignore numerical constants in the bounds. 

\paragraph{Complexity of the Hypothesis Class.}We specify the term $\text{Comp}(\Hscr)$ here. Specifically, we use the logarithm of the covering number (i.e. metric entropy) to represent $\text{Comp}(\Hscr)$ that appears throughout the main body (i.e. Table~\ref{tab:general-table}, Theorem~\ref{thm:p-dro-sample}), and we utilize some established covering number arguments from Section 2.2.2 in \cite{maurer2009empirical}, \cite{duchi2019variance}, and Section 5.3.2 in \cite{shapiro2014lectures}. Nonetheless, the proof framework in this part can be generalized to other refined complexity measures such as localized Rademacher Complexity; see examples at the end of this subsection and \cite{bartlett2005local}.
\begin{definition}[Complexity of Hypothesis Class $\text{Comp}_n(\Hscr)$]\label{def:comp-class}
    Recall the hypothesis class $\Hscr = \{h(x;\cdot), x \in \Xscr\}$, we define $\text{Comp}_n(\Hscr) := \log N_{\infty}(\Hscr, \frac{1}{n}, n)$, where the \emph{empirical $\ell_{\infty}$ covering number} $N_{\infty}(\Hscr, \varepsilon, n)$ is defined to be:
\begin{equation}\label{eq:covering-num}
N_{\infty}(\Hscr, \varepsilon, n) := \sup_{\bm{\xi} \in \Xi^n}N(\Hscr_n(\bm{\xi}), \varepsilon, \|\cdot\|_{\infty}),   
\end{equation}
where $\Hscr_n(\bm{\xi}) = \{(h(\xi_1),\ldots, h(\xi_n)): h \in \Hscr\}\subseteq \R^n$ and for $A \subseteq \R^n$, the number $\Nscr(A, \varepsilon, \|\cdot\|_{\infty})$ is the \textit{covering number} denoting the smallest cardinality $|A^{\prime}|$ of a set $A^{\prime}\subseteq A$ such that $A \subset \cup_{x_0 \in A^{\prime}}\{x: \|x - x_0\|_{\infty}\leq \varepsilon\}$. Furthermore, we denote $\text{Comp}(\Hscr) = \log N_{\infty}(\Hscr, \tau, 1)$ for some constants $\tau > 0$.
\end{definition}

For the function class $\Hscr$ in practice, $\text{Comp}_n(\Hscr) = O((\log n)^{c_1}\text{Comp}(\Hscr))$ in $n$ for some constant $c_1$; see \cite{maurer2009empirical} and Chapter 5 of \cite{wainwright2019high} for more details. It is well-known that for the VC dimension, $\text{Comp}(\Hscr)\leq C\text{VC}(\Hscr)\log n$ for some numerical constant $C$. And we will detail in Appendix \ref{subsec:synthetic2} an example of estimating $\text{Comp}(\Hscr)$ for the numerical study shown in Section~\ref{subsec:synthetic}.

\paragraph{ERM.}
Denote $x^{\texttt{NP-ERM}}$ as the solution obtained by solving~\eqref{eq:EO-obj}. We have the following:

\begin{lemma}[Extracted from \cite{vapnik1999nature, boucheron2005theory}]\label{lemma:nonparam-erm}
Consider $x^{\texttt{NP-ERM}}$ as the minimizer of $\min_x \hat{Z}^{ERM}(x)$ in~\eqref{eq:EO-obj}. Denote $M:=\sup_{x \in \Xscr}\|h(x;\cdot)\|_{\infty}$. Then we have the following generalization error of $x^{\texttt{NP-ERM}}$ with probability at least $1-\delta$: 
\begin{equation}\label{eq:EO}
Z(x^{\texttt{NP-ERM}})-Z(x^*)\leq \log(1/\delta)\left[\sqrt{\frac{M Z(x^*)\text{Comp}(\Hscr) \log n}{n}} + \frac{\text{Comp}(\Hscr)M}{n}\right].
\end{equation}
\end{lemma}
$\text{Comp}(\Hscr)$ appears due to the bounding of $\sup_{x \in \Xscr}|Z(x) - \hat{Z}^{ERM}(x)|$ from~\eqref{eq:error-decomp}. This result is minimax optimal in terms of the function complexity $\text{Comp}(\Hscr)$, e.g., the case of $\text{VC}(\Hscr)$ shown in Section 5 of \cite{boucheron2005theory}.

\paragraph{DRO.} 
Denote the optimal solution to DRO~\eqref{eq:DRO-obj} as $x^{\texttt{NP-DRO}}$. To bound the excess risk, we let $\hat{x}:= x^{\texttt{NP-DRO}}$ and $ \hat{Z}^{DRO}(\cdot):=\hat{Z}(\cdot)$ in~\eqref{eq:error-decomp}, whose second term will lead to an error $O\left(\frac{\varepsilon\Vscr_d(x^*)}{\sqrt{n}}\right)$. The key lies in the first term, i.e. $Z(x^{\texttt{NP-DRO}}) - \hat{Z}^{DRO}(x^{\texttt{NP-DRO}})$. We restate the two bounding perspectives in Section~\ref{sec:background} here.
\begin{theorem}\label{thm:nonparam-dro}
The following generalization error bound of $x^{\texttt{NP-DRO}}$ holds with probability at least $1-\delta$ for some metrics $d$.

\textit{(Regularization Perspective)}~When the ambiguity size $\varepsilon =\Omega((\frac{\text{Comp}_n(\Hscr)}{n})^{\beta})$ with $\beta$ being a constant depending on different metrics $d$, we have: 
\begin{equation}\label{eq:dro-reg}
\begin{aligned}
Z(x^{\texttt{NP-DRO}})-Z(x^*)\leq & \log(1/\delta)\left[\varepsilon^{\beta}\Vscr_d(x^*) +\frac{\text{Comp}_n(\Hscr)(\sup_{x \in \Xscr}\Vscr_d(x))}{n}\right.\\
+& \left. \frac{\mathcal{E}_1(x^*)}{\sqrt{n}} + \mathcal{E}_2(x^*,\varepsilon)\right],     
\end{aligned}
\end{equation}
where $\mathcal{E}_1(x^*)$ only depends on $h(x^*;\xi)$ and $\P^*$,  $\mathcal{E}_2(x^*,\varepsilon) = O(\Vscr_d(x^*) \varepsilon^{2\beta})$, which is of order $\frac{1}{n}$. 

\textit{(Robustness Perspective)}~When the ambiguity size $\varepsilon = \Omega\left(n^{-1/g(D_{\xi})}\right)$, we have: 
\begin{equation}\label{eq:dro-pessim}
Z(x^{\texttt{NP-DRO}})-Z(x^*)\leq\frac{\Vscr_d(x^*)\log(1/\delta)}{n^{1/g(D_{\xi})}},
\end{equation}
where $g(D_{\xi})$ is a function of $D_{\xi}$.
\end{theorem}

Theorem~\ref{thm:nonparam-dro} unifies several streams of results in the DRO literature from the regularization and robustness perspectives. We present some examples below, where we denote $r_n^* \leq \frac{\text{VC}(\Hscr)\log\left(n/\text{VC}(\Hscr)\right)}{n}$ as the fixed point of some sub-root Rademacher Complexity in \cite{bartlett2005local}:
\begin{example}[\cite{esfahani2018data} and \cite{gao2022finite}]\label{ex:w1-dro-eps}
When $d$ is taken as 1-Wasserstein distance and $\Vscr_{W_1}(x) = \|h(x;\cdot)\|_{Lip}$:
\begin{itemize}
    \item \eqref{eq:dro-reg} holds with $\beta = \frac{1}{2}$ and $\varepsilon = \min\left\{\sqrt{\frac{\tau \log(N(\Hscr, \frac{1}{n}, n)/\delta)}{n}}, \sqrt{\frac{\tau\log(1/\delta)}{n}} + \sqrt{r_n^*} + \frac{1}{n\sqrt{r_n^*}}\right\}$ where $\tau$ is a constant only depending on $\P^*$. 
    \item \eqref{eq:dro-pessim} holds with $g(D_{\xi}) = D_{\xi}$.
\end{itemize}
\end{example}

\begin{example}[\cite{duchi2019variance}]\label{ex:chi-dro-eps}
When $d$ is taken as $\chi^2$-divergence, \eqref{eq:dro-reg} holds with $\beta = 1$, $\varepsilon = \frac{\log(N(\Hscr, \frac{1}{n}, n)/\delta)}{n}$ (or $\frac{M \log(1/\delta)}{n} + r_n^*)$. 
And $\Vscr_{\chi^2}(x) = \sqrt{\text{Var}_{\P^*}[h(x;\cdot)]}$. 
\end{example}

In the above examples, the bound~\eqref{eq:dro-reg} comes from a combination of the following results:
\begin{itemize}
    \item Variability regularization in the form:
    \begin{equation}\label{eq:var-reg}
     Z(x) \leq \hat{Z}^{DRO}(x) + \frac{\text{Comp}_n(\Hscr)\sup_{x \in  \Xscr}\Vscr_d(x)}{n}
    \end{equation}
    \item DRO expansions in the form with $\hat{Z}_n(\cdot):= \frac{1}{n}\sum_{i = 1}^n h(\cdot;\hat{\xi}_i)$:
    $$\hat{Z}^{DRO}(x^*) \leq \hat{Z}_n(x^*) + \varepsilon^{\beta}\Vscr_d(x^*) + \mathcal{E}_2(x^*,\varepsilon).$$
    \item Standard concentration bound for the empirical mean:
    $$|\hat{Z}_n(x^*) - Z(x^*)| \leq \frac{\mathcal{E}_1(x^*)}{\sqrt{n}}.$$
\end{itemize} 

The bound~\eqref{eq:dro-pessim} is achieved by using a ball size $\varepsilon$ large enough to cover the true $\P^*$ with probability at least $1-\delta$. In this case, we have:
\begin{equation}\label{eq:dro-pessim-cover}
Z(x^{\texttt{NP-DRO}}) - \hat{Z}^{DRO}(x^{\texttt{NP-DRO}}) \leq 0\\    
\end{equation}
\begin{equation}\label{eq:dro-ipm-bd}
    \hat{Z}^{DRO}(x^*) - Z(x^*) \leq \Vscr_d(x^*)\varepsilon.
\end{equation}
Typically, to ensure that the ball size is large enough to cover the true distribution with probability at least $1-\delta$, the ambiguity size needs to depend on $D_{\xi}$.


\section{Further Details and Proofs for Section~\ref{sec:main}}
\subsection{Further Examples of Parametric Estimators that Satisfy Assumption~\ref{asp:oracle-param-est}}\label{app:asp-param-estimator}
In the main body, we present two examples where $d$ is Wasserstein distance or KL-divergence. Here we give some additional estimators $\hat{\Q}$, and pairing distribution metrics $d$ that satisfy Assumption~\ref{asp:oracle-param-est}.
\begin{example}
 $d$ is $H^2$-distance, $\mathcal P_\Theta$ is the class of all distributions governing $g_\theta(Z)$ for some random variable $Z$ and function $g_\theta$ is given by a feed-forward neural network parametrized by $\theta\in\Theta$.
Then Assumption~\ref{asp:oracle-param-est} holds under the same estimation procedure as in Example~\ref{ex:oracle-estimator-kl}. Following the same notation there, we have:
\begin{align*}
&\mathcal{E}_{apx} = \sup_{\theta}\inf_{\omega}\left\|\frac{\sqrt{p_*} - \sqrt{p_{\theta}}}{\sqrt{p_*} + \sqrt{p_{\theta}}} - f_{\omega}\right\|_{\infty} + B \inf_{\theta}\left\|\frac{\sqrt{p_*} - \sqrt{p_{\theta}}}{\sqrt{p_*} + \sqrt{p_{\theta}}}\right\|_{\infty},\\
&\text{Comp}(\Theta) = \texttt{Pdim}(\Fscr), \alpha = \frac{1}{2}.
\end{align*}
\end{example}

\begin{example}\label{ex:add-ex2}
$d$ is $\chi^2$-divergence and $\P^* \in \Pscr_{\Theta}$ is a location variant of the Beta distribution. See Proposition~\ref{prop:chi2-oracle-estimator} for the specific result of $\text{Comp}(\Theta)$ and $\alpha$. This distribution class is considered in the numerical experiments in Sections~\ref{subsec:synthetic} and~\ref{subsec:synthetic-context}.
\end{example}

\begin{example}\label{ex:add-ex3}
$d$ is 1-Wasserstein distance. We consider the following two different mixture distribution classes:
\begin{itemize}
    \item Case I: Suppose that the true distribution comes from the Gaussian mixture model $\P^* \in \Pscr_{\Theta} (:= \{\frac{1}{2}\Nscr(\mu, \Sigma) + \frac{1}{2}\Nscr(-\mu, \Sigma)| \mu \in \R^{D_{\xi}}\})$ with known $\Sigma:= \sigma^2 I_{D_{\xi}\times D_{\xi}}$. 
    
    Then Assumption~\ref{asp:oracle-param-est} holds for $\hat{\Q}:=\frac{1}{2}\Nscr(\hat{\mu}, \Sigma) + \frac{1}{2}\Nscr(-\hat{\mu}, \Sigma)$, with $\hat{\mu}$ output by the EM algorithm. In addition, $\text{Comp}(\Theta) =D_{\xi} \sigma^2, \alpha = \frac{1}{2}$ and $\Escr_{apx} = 0$, which is implied by $W_1(\P^*, \hat{\Q}) \leq \|\hat{\mu} - \mu^*\|_2 = O\left(\sigma \sqrt{\frac{D_{\xi} \log (1/\delta)}{n}}\right)$ following from $\|\hat{\mu} - \hat{\mu}\|_2^2 = O\left(\frac{\sigma^2D_{\xi}\log(1/\delta)}{n}\right)$ in Theorem 6 of~\cite{xu2020towards} and Corollary 2 of~\cite{balakrishnan2017statistical} under some mild conditions.
    \item Case II: Suppose that the true distribution is represented by: $\P^*:=\sum_{k = 1}^K p_k^* \P_k^*$ for some unknown probability density $p_k$ and the distribution $\P_k^*$ for each group.  Define $\Pscr_{\Theta} = \{\sum_{k = 1}^K p_k \Nscr(\mu_k, \Sigma)|(p_1, \ldots, p_K)' \in \Delta_K, \mu_k \in \R^{D_{\xi}},\forall k\in [K]\}$ for some known $\Sigma$ where $\Delta_K$ represents the $K$-dimensional probability simplex. In addition, we are given the group labels $\{g_i\}_{i = 1}^n$ associated with $\{\hat{\xi}_i\}_{i = 1}^n$, where each $g_i \in [K]$. 
   
    Then Assumption~\ref{asp:oracle-param-est} holds for $\hat{\Q}:=\sum_{k\in [K]}\hat{p}_k\Nscr(\hat{\mu}_k, \Sigma)$, where $\hat{p}_k := \frac{\sum_{i = 1}^n \mathbb{I}_{\{g_i = k\}}}{n}, \hat{\mu}_k := \frac{\sum_{i = 1}^n \hat{\xi}_i \mathbb{I}_{\{g_i = k\}}}{n \hat{p}_k}, \forall k \in [K]$. $\mathcal{E}_{apx} = W_1(\P^*, \Q^*)$ with $\Q^*:=\sum_{k \in [K]}p_k^*\Nscr(\E_{\xi \sim \P_k^*}[\xi], \Sigma)$, $\alpha = \frac{1}{2}$, $\text{Comp}(\Theta) = CD_{\xi} \sigma^2 K^2$ with some constant $C$ depending on $\P^*$ (e.g., scaling with $\frac{1}{\min_{k\in [K]} p_k^*}$ and $\max_{i,j \in [K]}\|\E_{\xi\sim \P_i^*}[\xi] - \E_{\xi\sim \P_j^*}[\xi]\|$).
\end{itemize}
\end{example}
The result in Case II of Example~\ref{ex:add-ex3} follows:
$$W_1(\P^*, \hat{\Q}) \leq W_1(\P^*, \Q^*) + W_1(\Q^*, \tilde{\Q}) + W_1(\tilde{\Q}, \hat{\Q}),$$
where $\tilde{\Q}:=\sum_{k \in [K]}\hat{p}_k\Nscr(\E_{\xi \sim \P_k^*}[\xi], \Sigma)$.

\subsection{Proof of Theorem~\ref{thm:general-param-dro-ipm}} 
To distinguish from their nonparametric counterparts, we denote $\hat{Z}^{P-DRO}(x) := \sup_{\P: d(\P, \hat{\Q})\leq \varepsilon}\E_{\P}[h(x;\xi)]$.

In case (a) where $d$ is an IPM, we have:
\begin{equation}\label{eq:general-param-dro-ipm-derive}
        \begin{aligned}
      Z(x^{\PDRO}) - Z(x^*)&\leq \sup_{\P: d(\P, \hat{\Q})\leq \varepsilon}\E_{\P}[h(x^*;\xi)] - \E_{\P^*}[h(x^*;\xi)]\\
      &\leq \sup_{\P: d(\P, {\P}^*)\leq 2\varepsilon}\E_{\P}[h(x^*;\xi)] - \E_{\P^*}[h(x^*;\xi)]\\
      &\leq 2 \Vscr_d(x^*)\varepsilon, 
    \end{aligned}
\end{equation}
where the first inequality follows from the fact that when $\varepsilon \geq \Delta(\delta, \Theta)$, by Assumption~\ref{asp:oracle-param-est}, we have $\P^* \in \Ascr$ (i.e. $d(\P^*,\hat{\Q}) \leq \varepsilon$ with probability at least $1-\delta$). Therefore, the term $Z(x) - \hat{Z}^{P-DRO}(x)$ in~\eqref{eq:error-decomp}  is non-positive with probability at least $1-\delta$. Furthermore, the second inequality follows from the triangle inequality property of distance, $\forall \P \in \Ascr$, $d(\P, \P^*) \leq d(\P, \hat{\Q}) + d(\hat{\Q},\P^*)\leq 2\varepsilon$. Finally, the last inequality follows from the fact that $d$ is an IPM with $d(\P, \Q) = \sup_{f: V_{d}(f)\leq 1}\bigg|\E_{\P}[f] - \E_{\Q}[f]\bigg|$.

In case $(b)$ where $d$ satisfies the inequality~\eqref{eq:general-f-tv-ineq}, we have:
\begin{equation*}
\begin{aligned}
          \mathcal{E}(x^{\PDRO}) &\leq \sup_{\P: d(\P, \hat{\Q}) \leq  {\varepsilon}}\E_{\P}[h(x^*;\xi)]- \E_{\P^*}[h(x^*;\xi)]\\
          &\leq \sup_{\P: d_{TV}(\P, \hat{\Q})\leq C_{d}\sqrt{\varepsilon}}\E_{\P}[h(x^*;\xi)]- \E_{\P^*}[h(x^*;\xi)]\\
          &\leq 4C_{d}\sqrt{\varepsilon}\|h(x^*;\cdot)\|_{\infty},
\end{aligned}
\end{equation*}
where the first inequality follows from $\P^*\in \Ascr$ with probability at least $1-\delta$. The second inequality follows from the fact that $d_{TV}(\P, \hat{\Q}) \leq C_{d}\sqrt{d(\P, \hat{\Q})}$ such that $\{\P: d(\P, \hat{\Q})\leq \varepsilon\} \subseteq \{\P: d_{TV}(\P, \hat{\Q})\leq C_d\sqrt{\varepsilon}\}$. The remaining part follows the same argument as~\eqref{eq:general-param-dro-ipm-derive} above since TV distance is an IPM.

In the special case $(c)$ where $d$ is $\chi^2$-divergence, we have:
\begin{equation*}
\begin{aligned}
          \mathcal{E}(x^{\PDRO}) &\leq \sup_{\P: \chi^2(\P, \hat{\Q})\leq \varepsilon}\E_{\P}[h(x^*;\xi)]- \E_{\P^*}[h(x^*;\xi)]\\
          &\leq \E_{\hat{\Q}}[h(x^*;\xi)] + \sqrt{2\varepsilon \text{Var}_{\hat{\Q}}[h(x^*;\xi)]} - \E_{\P^*}[h(x^*;\xi)]\\
          &\leq 2\sqrt{\varepsilon \text{Var}_{\hat{\Q}}[h(x^*;\xi)]}\\
          &\leq 2\sqrt{\varepsilon\text{Var}_{\P^*}[h(x^*;\xi)]} + 2^{\frac{5}{4}}\varepsilon^{\frac{3}{4}}\left[(\text{Var}_{\P^*}[h^2(x^*;\xi)])^{\frac{1}{4}} + 2^{\frac{1}{4}}\|h(x^*;\cdot)\|_{\infty}^{\frac{1}{2}}(\text{Var}_{\P^*}[h(x^*;\xi)])^{\frac{1}{4}}\right]\\
          &\leq 2\sqrt{\varepsilon\text{Var}_{\P^*}[h(x^*;\xi)]}+ 4\varepsilon^{\frac{3}{4}}\|h(x^*;\cdot)\|_{\infty},
\end{aligned}    
\end{equation*}
where the first inequality follows from the fact that $\chi^2$-divergence satisfies Assumption~\ref{asp:oracle-param-est} with probability at least $1-\delta$. The second and third inequalities follow from Lemma~\ref{lemma:chi2-cauchy} for two pairs $(\P,\hat{\Q})$ and $(\P^*, \hat{\Q})$. And the fourth inequality follows from:
\begin{equation}\label{eq:var-decomp-chi2}
\begin{aligned}
    \text{Var}_{\hat{\Q}}[h] - \text{Var}_{\P^*}[h] &\leq \left|\E_{\hat{\Q}}[h^2] - \E_{\P^*}[h^2]\right| + 2\|h\|_{\infty}\left|\E_{\hat{\Q}}[h] - \E_{\P^*}[h]\right|\\
    &\leq \sqrt{2\chi^2(\P^*,\hat{\Q}) \text{Var}_{\P^*}[h^2]} + 2\|h\|_{\infty}\sqrt{2\chi^2(\P^*,\hat{\Q})\text{Var}_{\P^*}[h]}.\quad\hfill \square
\end{aligned}
\end{equation}
\subsection{Improved Results from Theorem~\ref{thm:general-param-dro-ipm} for General $f$-divergence}\label{app:improve-f}
The result in Theorem~\ref{thm:general-param-dro-ipm} can be improved for general $f$-divergence from $\|h(x^*;\cdot)\|_{\infty}$ to $\sqrt{\text{Var}_{\P^*}[h(x^*;\cdot)]}$ in terms of the dependence of $\varepsilon$, without requiring~\eqref{eq:general-f-tv-ineq} as long as some mild conditions hold for the cost function $h(x;\cdot)$ and sample size $n$. This is achieved by applying a general duality property of $f$-divergence DRO to bound $\hat{Z}^{P-DRO}(x^*) - Z(x^*)$. The result, which is presented in Theorem~\ref{thm:general-f-div-bd} below, holds for any $f$-divergence DRO regardless of the distribution center $\hat{\Q}$. 

\begin{theorem}\label{thm:general-f-div-bd}
When the sample size $n$ is large enough, for general metrics $d_f$ in $f$-divergence and $\|h(x^*;\cdot)\|_{\infty} < \infty$, we have:
\begin{equation}\label{eq:general-f-div-expansion-upper}
    \sup_{\P: d_f(\P, \hat{\Q})\leq \varepsilon}\E_{\P}[h(x;\xi)] \leq \E_{\hat{\Q}}[h(x;\xi)]+ C(f, \varepsilon)\sqrt{\text{Var}_{\hat{\Q}}[h(x;\xi)]\varepsilon},
\end{equation}
\begin{equation}\label{eq:general-f-div-expansion-lower}
    \inf_{\P: d_f(\P, \hat{\Q})\leq \varepsilon}\E_{\P}[h(x;\xi)] \geq \E_{\hat{\Q}}[h(x;\xi)]- C(f, \varepsilon)\sqrt{\text{Var}_{\hat{\Q}}[h(x;\xi)]\varepsilon},
\end{equation}
where $C(f, \varepsilon)$ only depends on the metric $d_f$ and $\varepsilon$ and is bounded by some numerical constant in classical $f$-divergences (See Examples~\ref{ex:pdro-kl} and~\ref{ex:pdro-h2} below).
\end{theorem}
With this, the result in case $(b)$ in Theorem~\ref{thm:general-param-dro-ipm} can be improved to:
\begin{equation*}
    \begin{aligned}
    \mathcal{E}(x^{\PDRO}) &\leq \sup_{\P: d_{f}(\P, \hat{\Q}) \leq \varepsilon}\E_{\P}[h(x^*;\xi)]- \E_{\P^*}[h(x^*;\xi)]\\
    &\leq \E_{\hat{\Q}}[h(x^*;\xi)] - \E_{\P^*}[h(x^*;\xi)] + C(f,\varepsilon)\sqrt{\text{Var}_{\hat{\Q}}[h(x^*;\xi)]\varepsilon}\\
    &\leq \E_{\hat{\Q}}[h(x^*;\xi)] - \inf_{\P: d_f(\P, \hat{\Q})\leq \varepsilon}\E_{\P}[h(x^*;\xi)] + C(f, \varepsilon)\sqrt{\text{Var}_{\hat{\Q}}[h(x^*;\xi)]\varepsilon}\\
    &\leq2C(f, \varepsilon)\sqrt{\text{Var}_{\hat{\Q}}[h(x^*;\xi)]\varepsilon},
    \end{aligned}
\end{equation*}
where the second and fourth inequalities follow from the result in Theorem~\ref{thm:general-f-div-bd}. The first and third inequalities follow from $\P[d_f(\P^*, \hat{\Q})\leq \varepsilon] \geq 1-\delta$ such that $\inf_{d_{f}(\P, \hat{\Q}) \leq \varepsilon}\E_{\P}[h(x^*;\xi)]\leq \E_{\P^*}[h(x^*;\xi)]\leq \sup_{\P: d_{f}(\P, \hat{\Q}) \leq \varepsilon}\E_{\P}[h(x^*;\xi)]$ with probability at least $1-\delta$. After that, we can use the same argument as before to bound $\sqrt{\text{Var}_{\hat{\Q}}[h(x^*;\xi)]}$. 

\textit{Proof of Theorem~\ref{thm:general-f-div-bd}.} We first show~\eqref{eq:general-f-div-expansion-upper}. We have:
\begin{equation}\label{eq:pdro-fdiv-expansion}
\begin{aligned}
    \sup_{d_f(\P^*, \hat{\Q})\leq \varepsilon}\E_{\P}[h(x^*;\xi)] &\leq \min_{\lambda \geq 0, \mu}\left\{\lambda \E_{\hat{\Q}}\left[f^*\left(\frac{h(x^*;\xi) - \mu }{\lambda}\right)\right] + \lambda \varepsilon + \mu\right\}\\
    &\leq \hat{\lambda}\E_{\hat{\Q}}\left[f^*\left(\frac{h(x;\xi) - \hat{\mu} }{\hat{\lambda}}\right)\right] + \sqrt{\frac{\text{Var}_{\hat{\Q}}[h(x^*;\xi)]\varepsilon}{f''(1)}} + \E_{\hat{\Q}}[h(x^*;\xi)]\\
    &\leq \E_{\hat{\Q}}[h(x^*;\xi)] + \left(\frac{1}{\sqrt{f^{''}(1)}} + \frac{\sqrt{f^{''}(1)} (f^*)^{''}(0) C_0(\varepsilon)}{2}\right)\sqrt{\text{Var}_{\hat{\Q}}[h(x^*;\xi)]\varepsilon},
\end{aligned}
\end{equation}
where the first inequality above is based on weak duality, i.e., Theorem 1 in \cite{ben2013robust} (Note that although strong duality holds generally in this problem, we only need weak duality in our proof). The second inequality above is given by $\hat{\lambda} = \sqrt{\frac{\text{Var}_{\hat{\Q}}[h(x^*;\xi)]}{f''(1)\varepsilon}}, \hat{\mu} = \E_{\hat{\Q}}[h(x^*;\xi)]$ as the feasible dual solution, and the third inequality follows from plugging in the values of $\hat{\lambda}$ and $\hat{\mu}$, and then taking the Taylor expansion up to the second order for $f^*$ with a Maclaurin remainder $C_0(\varepsilon)$ upper bounded by some constant and $C_0(0) = 1$:
\begin{align*}
    \hat{\lambda}\E_{\hat{\Q}}\left[f^*\left(\frac{h(x^*;\xi) - \hat{\mu} }{\hat{\lambda}}\right)\right] &\leq \hat{\lambda}\E_{\hat{\Q}}\left[f^*(0) + (f^*)'(0)\left(\frac{h(x^*;\xi) - \hat{\mu} }{\hat{\lambda}}\right) + \frac{(f^*)^{''}(0) C_0(\varepsilon)}{2}\left(\frac{h(x^*;\xi) - \hat{\mu} }{\lambda}\right)^2\right],\\
    &=\hat{\lambda}\E_{\hat{\Q}}\left[ \frac{(f^*)^{''}(0) C_0(\varepsilon)}{2}\left(\frac{h(x^*;\xi) - \hat{\mu}}{\hat{\lambda}}\right)^2\right] = \frac{(f^*)^{''}(0)C_0(\varepsilon)\text{Var}_{\hat{\Q}}[h(x^*;\xi)]}{2\hat{\lambda}},
\end{align*}
where the first equality follows from $f^*(0) = 0$ and $\hat{\mu} = \E_{\hat{\Q}}[h(x^*;\xi)]$. Then \eqref{eq:general-f-div-expansion-upper} follows from \eqref{eq:pdro-fdiv-expansion} if we denote $C(f, \varepsilon) = \frac{1}{\sqrt{f^{''}(1)}} + \frac{\sqrt{f^{''}(1)} (f^*)^{''}(0) C_0(\varepsilon)}{2}$. For~\eqref{eq:general-f-div-expansion-lower}, we only need to consider $-h(x;\cdot)$ and plug in the result of~\eqref{eq:general-f-div-expansion-upper}.
$\hfill \square$

We now show several common divergences satisfying~\eqref{eq:general-f-tv-ineq} and give some concrete values for $C(f, \varepsilon)$ above. To avoid redundancy, we only consider~\eqref{eq:general-f-div-expansion-upper} and ignore~\eqref{eq:general-f-div-expansion-lower}. 
\begin{example}[KL divergence]\label{ex:pdro-kl}
We take $f(t) = t\log t - (t - 1)$. Then $f^*(t) = e^t - 1$ with $f^{''}(1) = 1$. We use the inequality $e^t - 1 \leq t + t^2$ when $t \in (-1,1)$, i.e., we need $\left|\frac{h(x^*;\xi) - \hat{\mu}}{\hat{\lambda}}\right| \leq 1$, which implies when $\sqrt{\frac{\text{Var}_{\hat{\Q}}[h(x^*;\xi)]}{\varepsilon}} = \hat{\lambda} \geq 2 \|h(x^*;\cdot)\|_{\infty}$, i.e., $\varepsilon \leq \frac{\text{Var}_{\hat{\P}}[h(x^*;\xi)]}{4\|h(x^*;\cdot)\|_{\infty}^2}$. Therefore, if $\varepsilon \leq \frac{\text{Var}_{\hat{\P}}[h(x^*;\xi)]}{4\|h(x^*;\cdot)\|_{\infty}^2}$, we have:
$$\sup_{\P: KL(\P, \hat{\Q})\leq \varepsilon}\E_{\P}[h(x^*;\xi)] \leq \E_{\hat{\Q}}[h(x^*;\xi)]+ 3\sqrt{\text{Var}_{\hat{\Q}}[h(x^*;\xi)] \varepsilon}.$$
\end{example}

\begin{example}[$H^2$-distance]\label{ex:pdro-h2}
We take $f(t) = (\sqrt{t} - 1)^2$ and $f^{''}(1) = \frac{1}{2}$. Then for $t < 1$, $f^*(t) = \frac{t}{1 - t} = \frac{1}{1-t} - 1 \leq t + 2t^2$ when $t \in [-\frac{1}{2}, \frac{1}{2}]$. Therefore, if $\varepsilon \leq \frac{\text{Var}_{\hat{\Q}}[h(x^*;\xi)]}{2\|h(x^*;\cdot)\|_{\infty}^2}$, we have:
$$\sup_{\P: H^2(\P, \hat{\Q})\leq \varepsilon}\E_{\P}[h(x^*;\xi)] \leq \E_{\hat{\Q}}[h(x^*;\xi)] + (2+\sqrt{2})\sqrt{\text{Var}_{\hat{\Q}}[h(x^*;\xi)]\varepsilon}.$$
\end{example}
Therefore, following from the same argument as for case $(c)$ of Theorem~\ref{thm:general-param-dro-ipm}, if $\Delta(\delta, \Theta)\leq \varepsilon \leq \frac{\text{Var}_{\hat{\Q}}[h(x^*;\xi)]}{c_0 2\|h(x^*;\cdot)\|_{\infty}^2}$, $\mathcal{E}(x^{\PDRO})$ can be improved by $c_1 \sqrt{\text{Var}_{\P^*}[h(x^*;\xi)]} + c_2 \varepsilon^{\frac{3}{4}}\|h(x^*;\cdot)\|_{\infty}$ for KL divergence and Hellinger distance with probability at least $1-\delta$.

\subsection{Proof of Theorem~\ref{thm:general-param-erm}}

In case (a) where $d$ is an IPM, we have:
\begin{equation*}\label{eq:generalization-param-erm-derive}
\begin{aligned}
    \mathcal{E}(x^{\PERM}) &\leq |\E_{\hat{\Q}}[h(x^{\PERM};\xi)] - \E_{\P^*}[h(x^{\PERM};\xi)]| + |\E_{\hat{\Q}}[h(x^*;\xi)] - \E_{\P^*}[h(x^*;\xi)]|\\
    &\leq 2 \sup_{x \in \Xscr}\bigg|\E_{\P^*}[h(x;\xi)] - \E_{\hat{\Q}}[h(x;\xi)]\bigg|\\
    &\leq 2\sup_{x \in \Xscr}\Vscr_d(x) d(\P^*, \hat{\Q}). 
\end{aligned}
\end{equation*}
where the second inequality follows from the uniform bound $\forall x \in \Xscr$, and the last inequality follows from the fact that $d$ is an IPM such that $d(\P, \Q) = \sup_{\Vscr_d(f)\leq 1}|\E_{\P}[f] - \E_{\Q}[f]|.\hfill \square$

In case (b) where $d$ satisfies the inequality~\eqref{eq:general-f-tv-ineq}, similarly, we have:
\begin{align*}
\mathcal{E}(x^{\PERM}) &\leq |\E_{\hat{\Q}}[h(x^{\PERM};\xi)] - \E_{\P^*}[h(x^{\PERM};\xi)]| + |\E_{\hat{\Q}}[h(x^*;\xi)] - \E_{\P^*}[h(x^*;\xi)]|\\
    &{\leq} 2 \sup_{x \in \Xscr}\bigg|\E_{\P^*}[h(x;\xi)] - \E_{\hat{\Q}}[h(x;\xi)]\bigg|\\
    &{\leq} 4 M d_{TV}(\P^*, \hat{\Q}) \leq 4 C_d M \sqrt{d(\P^*, \hat{\Q})}.
\end{align*}

In the special case (c) where $d$ is (modified) $\chi^2$-divergence, following from the previous decomposition and Lemma~\ref{lemma:chi2-cauchy}, we have:
\begin{equation}\label{eq:generalization-param-erm-2-derive}
\begin{aligned}
    \mathcal{E}(x^{\PERM})&\leq |\E_{\hat{\Q}}[h(x^{\PERM};\xi)] - \E_{\P^*}[h(x^{\PERM};\xi)]| + |\E_{\hat{\Q}}[h(x^*;\xi)] - \E_{\P^*}[h(x^*;\xi)]|\\
    &\leq \sqrt{2\chi^2(\hat{\Q},\P^*)}\left(\sqrt{\text{Var}_{\P^*}[h(x^{\PERM};\xi)]} + \sqrt{\text{Var}_{\P^*}[h(x^*;\xi)]}\right)\\
    &\leq 2\sqrt{2\chi^2(\hat{\Q},\P^*)}\sqrt{\text{Var}_{\P^*}[h(x^*;\xi)]} + \sqrt{2\chi^2(\hat{\Q},\P^*)}\sqrt{\left|\E_{\P^*}[h^2(x^{\PERM};\xi)] - \E_{\P^*}[h^2(x^*;\xi)]\right|}\\
    &\leq 2\sqrt{2\chi^2(\hat{\Q},\P^*)}\sqrt{\text{Var}_{\P^*}[h(x^*;\xi)]} + \sqrt{2\chi^2(\hat{\Q},\P^*)}\sqrt{4M^2 d_{TV}(\hat{\Q},\P^*)}\\
    &\leq 2\sqrt{2\chi^2(\hat{\Q},\P^*)}\sqrt{\text{Var}_{\P^*}[h(x^*;\xi)]} + 2M(\chi^2(\hat{\Q},\P^*))^{\frac{3}{4}},  
\end{aligned}
\end{equation}
where the fourth inequality in~\eqref{eq:generalization-param-erm-2-derive} follows from:
\begin{equation*}
\begin{aligned}
\E_{\P^*}[h^2(x^{\PERM};\xi)] - \E_{\P^*}[h^2(x^*;\xi)]&\leq \E_{\P^*}\left[(h(x^{\PERM};\xi) + h(x^*;\xi))((h(x^{\PERM};\xi) - h(x^*;\xi))\right]\\
&\leq 2M\E_{\P^*}[h(x^{\PERM};\xi) - h(x^*;\xi)] \leq 4M^2 d_{TV}(\hat{\Q}, \P^*).
\end{aligned}
\end{equation*}
And the fifth inequality in~\eqref{eq:generalization-param-erm-2-derive} follows from Lemma~\ref{lemma:pinsker-ineq}.

For each case above, we then apply Assumption~\ref{asp:oracle-param-est} to obtain the result.$\hfill \square$

\section{Further Details and Proofs for Section~\ref{subsec:sample}}\label{sec:mc-approximation}

We denote $x^{\PDRO_m}\in \argmin_{x\in \Xscr}\max_{d(\P,\hat{\Q}_m) \leq \varepsilon}\E_{\P}[h(x;\xi)]$. In general, we want to investigate the required sample size $m$ such that $\mathcal{E}(x^{\PDRO_m}) \approx \mathcal{E}(x^{\PDRO})$ in Theorem~\ref{thm:general-param-dro-ipm}. The idea is to understand when the Monte Carlo sampling error is dominated by the generalization error in each \texttt{P-DRO} case.  

Across statements in Section~\ref{subsec:sample}, we ignore the polynomial dependence on $\log(1/\delta)$ for the required Monte Carlo size when we express the required sample sizes for different generalization error bounds. That is to say, when we write ``the required sample size $\cdots \geq \sqrt{\text{Comp}_m(\Hscr)}$", we mean ``the required sample size $\cdots \geq \sqrt{\text{Comp}_m(\Hscr) + c_1 \log(1/\delta)}$" for some constant $c_1$. In other words, we ignore the dependence of $\log(1/\delta)$ and the required sample sizes are at most polynomial in this ignored term. 


\subsection{Proof of Theorem~\ref{coro:p-dro-sample-w1}} 

Comparing the result with Theorem~\ref{thm:general-param-dro-ipm}, we only need to show that $\P^* \in \hat{\Ascr}$ with probability at least $1-\delta$. Then the other parts follow directly from the case (a) in Theorem~\ref{thm:general-param-dro-ipm}.
Following the triangle inequality, we have:
\begin{equation*}
\begin{aligned}
    W_1(\P^*,\hat{\Q}_m)&\leq W_1(\P^*,\hat{\Q}) +W_1(\hat{\Q},\hat{\Q}_m)\\
    &\leq \frac{\varepsilon}{2} + \left(\frac{C }{m}\right)^{\frac{1}{D_{\xi}}}\log(1/\delta)
    \leq \varepsilon,
\end{aligned}
\end{equation*}
where the second inequality follows from Lemma~\ref{lemma:measure-concentration} and third inequality follows from $\frac{\varepsilon}{2}\geq \Delta(\delta, \Theta)$ in Assumption~\ref{asp:oracle-param-est} and $\left(\frac{C }{m}\right)^{\frac{1}{D_{\xi}}}\log(1/\delta) \leq \frac{\varepsilon}{2}$, i.e., $m \geq C(\frac{2 \log(1/\delta)}{\varepsilon})^{D_{\xi}}$ for some constant $C$. Then the subsequent steps are analogous to the proof of case (a) in Theorem~\ref{thm:general-param-dro-ipm}, only replacing $\hat{\Q}$ with $\hat{\Q}_m$.
$\hfill \square$

\subsection{Proofs of Theorem~\ref{thm:p-dro-sample}}\label{subsec:pdro-thm}
We present results and proofs with respect to \textit{$\chi^2$-divergence} and \textit{1-Wasserstein distance} separately, with $m \approx \text{Comp}(\Hscr) n^{\alpha}$ and $\alpha$ being independent of $D_{\xi}$. Theorem~\ref{thm:p-dro-sample-chi2} and Theorem~\ref{thm:p-dro-sample-w1} are more general results from which Theorem~\ref{thm:p-dro-sample} follows.
\begin{theorem}[Generalization bounds for $\chi^2$ \texttt{P-DRO} with Monte Carlo errors]\label{thm:p-dro-sample-chi2}
Suppose Assumption~\ref{asp:oracle-param-est} holds and the cost function $h(x;\xi) \in [0, M],\forall x,\xi$ with $\text{Var}_{\P^*}[h(x^*;\cdot)] > 0$. The size of the ambiguity set $\varepsilon \geq \Delta(\delta, \Theta)$. If the Monte Carlo size $m \geq C_0\left(\frac{LM}{\sqrt{\text{Var}_{\P^*}[h(x^*;\cdot)]}\varepsilon}\right)^2{\text{Comp}_m(\Hscr)}$ for some numerical constant $C_0$, when $d$ is $\chi^2$-divergence, with probability at least $1-\delta$, we have:
$$\mathcal{E}(x^{\PDRO_m})\leq\begin{cases} 2\mathcal{E}_{\chi^2} + C_1 \sqrt{\frac{\varepsilon}{L}} M,&~\text{if}~\text{Var}_{\hat{\Q}}[h(x^{\PDRO_m};\xi)] \leq 2\varepsilon M^2\\ 2\mathcal{E}_{\chi^2},&~\text{otherwise}~\end{cases},$$
where $L \geq 1$ and $\mathcal{E}_{\chi^2}$ is the generalization error upper bound in the case $(c)$ of Theorem~\ref{thm:general-param-dro-ipm}.
\end{theorem}

Note that this result depends on another term $L$ due to ``incomplete" exact variance regularization of $\chi^2$-divergence. However, when $\text{Var}_{\P^*}[h(x;\xi)]$ is sufficiently large, as long as the required Monte Carlo size $m \geq C_0\left(\frac{M}{\sqrt{\text{Var}_{\P^*}[h(x^*;\cdot)}\varepsilon}\right)^2{\text{Comp}_m(\Hscr)}$, $\mathcal{E}(x^{\PDRO_m}) \leq 2\mathcal{E}_{\chi^2}$. On the other hand, even if the variance is not large enough, as long as $\sqrt{\frac{\varepsilon}{L}} M \leq \mathcal{E}_{\chi^2} \leq  \sqrt{\text{Var}_{\P^*}[h(x^*;\cdot)]\varepsilon}$, i.e., $L \geq \frac{M^2}{\text{Var}_{\P^*}[h(x^*;\cdot)]}$ and $m \geq \left(\frac{M}{\sqrt{\text{Var}_{\P^*}[h(x^*;\cdot)]}\varepsilon}\right)^6{\text{Comp}_m(\Hscr)}$ for some numerical constant $C_0$, we still have $\mathcal{E}(x^{\PDRO_m}) \leq 3\mathcal{E}_{\chi^2}$, which is the required sample size in Theorem~\ref{thm:p-dro-sample}.




We obtain a dimension-free required Monte Carlo sample size for the 1-Wasserstein case too. 
\begin{theorem}[Generalization bounds for 1-Wasserstein \texttt{P-DRO} with Monte Carlo errors]\label{thm:p-dro-sample-w1}
Suppose Assumption~\ref{asp:oracle-param-est} holds and the random quantity $h(x;\xi)$ is sub-Gaussian with parameter $M$ when $\xi$ follows any distribution $\Q \in \Pscr_{\Theta}, \forall x \in \Xscr$. Besides, $\Xi$ is unbounded and $h(x;\xi)$ is Lipschitz continuous and convex w.r.t. $\xi$. Furthermore, there exists $\xi_0 \in \Xi$ such that $\limsup_{\|\tilde{\xi} - \xi_0\| \to \infty}\frac{h(x;\xi) - h(x;\xi_0)}{\|\tilde{\xi} -\xi_0\|} = \|h(x;\cdot)\|_{\text{Lip}}, \forall x \in \Xscr$. The size of the ambiguity set $\varepsilon \geq \Delta(\delta, \Theta)$. 

If the Monte Carlo size $m \geq C_0\left(\frac{M}{\|h(x^*;\cdot)\|_{\text{Lip}}\varepsilon}\right)^2 \text{Comp}_m(\Hscr)$ for some numerical constant $C_0$, when $d$ is 1-Wasserstein distance, then with probability at least $1-\delta$, we have: $\mathcal{E}(x^{\PDRO_m}) \leq 4\|h(x^*;\cdot)\|_{\text{Lip}}\varepsilon.$
\end{theorem}

Before presenting the proofs, we introduce two uniform concentration inequalities for the empirical mean and variance over $\Hscr$.
\begin{definition}[Sub-Gaussian Random Variable]
    A random variable $g(\xi)$ over $\R$ is called sub-Gaussian with parameter $\sigma$ when $\xi \sim \Q$ if $\E_{\Q}[g(\xi)] < \infty$ and $\E_{\Q}[\exp(t (g(\xi) - \E_{\Q}[g(\xi)]))] \leq \exp(\sigma^2 t^2/2), \forall t \in \R$.
\end{definition}

\begin{lemma}[Uniform Hoeffding (Sub-Gaussian) Inequality]\label{lemma:uniform-hoeffding}
Suppose the random quantity $h(x;\xi)$ is sub-Gaussian with parameter $M$ when $\xi \sim \hat{\Q}, \forall x \in \Xscr$, then with probability at least $1-\delta$, we have:
\begin{equation}\label{eq:uniform-hoeffding}
    \begin{aligned}
    \E_{\hat{\Q}}[h(x;\xi)] -\E_{\hat{\Q}_m}[h(x;\xi)]&\leq C_1 M\sqrt{\frac{\text{Comp}_m(\Hscr)}{m}},
    \end{aligned}
\end{equation}
where $C_1$ is some numerical constant independent of the function complexity and sample size. 
\end{lemma}
This result is extracted from Theorem 6 in \cite{maurer2009empirical}. And a special case in Lemma~\ref{lemma:uniform-hoeffding} is when $0 \leq h(x;\xi) \leq M, \forall x, \xi$.

\begin{lemma}[Uniform Variance Concentration Inequality]\label{lemma:uniform-var}
When $0 \leq h(x;\xi) \leq M, \forall x, \xi$, with probability at least $1-\delta$, we have:
\begin{equation}\label{eq:uniform-var}
\sqrt{\text{Var}_{\hat{\Q}_m}[h(x;\xi)]} \geq \sqrt{1-\frac{1}{m}}\sqrt{\text{Var}_{\hat{\Q}}[h(x;\xi)]} - \frac{2M^2}{m}  - C_2M \sqrt{\frac{\text{Comp}_m(\Hscr)}{m}},   
\end{equation}
where $C_2$ is some numerical constant independent of the function complexity and sample size.
\end{lemma}
\proof{Proof of Lemma~\ref{lemma:uniform-var}.}
Consider the variance concentration inequality (extracted from Lemma A.1 in \cite{duchi2019variance}), $\forall x \in \Xscr$, when $m \geq 3$, with probability at least $1-\delta$:
\begin{equation}\label{eq:var-concentrate-lower}
\sqrt{\text{Var}_{\hat{\Q}_m}[h(x;\xi)]} \geq \sqrt{1-\frac{1}{m}}\sqrt{\text{Var}_{\hat{\Q}}[h(x;\xi)]} - \frac{2M^2}{m}  - M\sqrt{\frac{2\log(1/\delta)}{m}}.    
\end{equation}
\begin{equation}\label{eq:var-concentration-upper}
\sqrt{\text{Var}_{\hat{\Q}_m}[h(x;\xi)]} \leq \sqrt{1+\frac{1}{m}}\sqrt{\text{Var}_{\hat{\Q}}[h(x;\xi)]}+ M\sqrt{\frac{2\log(1/\delta)}{m}}.   
\end{equation}
From Definition~\ref{def:comp-class} with $\ell := N_{\infty}(\Hscr, \frac{1}{m}, m)$, we consider functions $h(x_1;\cdot), \ldots, h(x_{\ell};\cdot)$ such that $\Hscr_m(\bm \xi)$ is contained in the union of balls $D_k := \{(h(x;\xi_1),\ldots, h(x;\xi_m))|x\in \Xscr,\sup_{\bm \xi \in \Xi^m}\sup_{i \in [m]}|h(x;\xi_i) - h(x_k;\xi_i)| \leq 1/m\}, k \in [\ell]$. Then we apply the union bound to~\eqref{eq:var-concentrate-lower} such that with probability at least $1-\delta$:
$$\sqrt{\text{Var}_{\hat{\Q}_m}[h(x_i;\xi)]} \geq \sqrt{1-\frac{1}{m}}\sqrt{\text{Var}_{\hat{\Q}}[h(x_i;\xi)]}-\frac{2M^2}{m}- M\sqrt{\frac{2\log(N_{\infty}(\Hscr, \frac{1}{m}, m)/\delta)}{m}}, \forall i \in [\ell].$$  
Besides, for any other $h(x;\cdot)$ with $x \not\in \{x_1,\ldots, x_{\ell}\}$, we can always find one $h(x_k;\cdot)$ with $k \in [\ell]$ such that $\sup_{\bm\xi\in \Xi^m}\sup_{i \in [m]}|h(x;\xi_i) - h(x_k;\xi_i)| \leq \frac{1}{m}$ by the definition of the covering number $\ell$.

Therefore, we have:
\begin{align*}
    \sqrt{\text{Var}_{\hat{\Q}_m}[h(x;\cdot)]} &= \sqrt{\E_{\hat{\Q}_m}[h^2(x;\cdot)] - (\E_{\hat{\Q}_m}[h(x;\cdot)])^2}\\
    &\geq \sqrt{\E_{\hat{\Q}_m}[h^2(x_k;\cdot)] -(\E_{\hat{\Q}_m}[h(x_k;\cdot)])^2- 4M|\E_{\hat{\Q}_m}[h(x_k;\cdot) - h(x;\cdot)]|}\\
    &\geq \sqrt{\text{Var}_{\hat{\Q}_m}[h(x_k;\cdot)]} - 2\sqrt{M \frac{1}{m}}\\
    &\geq \sqrt{1-\frac{1}{m}}\sqrt{\text{Var}_{\hat{\Q}}[h(x_k;\cdot)]}-\frac{2M^2}{m}- M\sqrt{\frac{2\log(N_{\infty}(\Hscr, \frac{1}{m}, m)/\delta)}{m}} - 2\sqrt{\frac{M}{m}}\\
    &\geq \sqrt{1-\frac{1}{m}}\sqrt{\text{Var}_{\hat{\Q}}[h(x;\cdot)]}-\frac{2M^2}{m}- M\sqrt{\frac{2\log(N_{\infty}(\Hscr, \frac{1}{m}, m)/\delta)}{m}} - 4\sqrt{\frac{M}{m}},
\end{align*}
where the first inequality follows from $h(x;\cdot)$ being ``close'' to some $h(x_k;\cdot)$ from the covering number argument. And the second inequality follows from the ball size being $\frac{1}{m}$ and $\sqrt{a - b}\geq \sqrt{a} - \sqrt{b}$ when $\sqrt{a - b} \geq 0$. And the third inequality follows from the union bound above. 
The fourth inequality follows from $|\E_{\hat{\Q}}[h(x_k;\cdot) - h(x;\cdot)]| \leq \E_{\hat{\Q}}[h(x_k;\xi) - h(x;\xi)] \leq \frac{1}{m}$, i.e.:
\begin{align*}
\sqrt{\text{Var}_{\hat{\Q}}[h(x_k;\cdot)]} &\geq \sqrt{\text{Var}_{\hat{\Q}}[h(x;\cdot)] - 2M |\E_{\hat{\Q}}[h(x_k;\cdot) - h(x;\cdot)]|}\\
&\geq \sqrt{\text{Var}_{\hat{\Q}}[h(x;\cdot)} - \sqrt{\frac{2M}{m}},
\end{align*}
Then denote $\text{Comp}_{m}(\Hscr) =2 \log(N_{\infty}(\Hscr, \frac{1}{m}, m)/\delta)$, we obtain~\eqref{eq:uniform-var}. $\hfill \square$

\subsubsection{Proof of Theorem~\ref{thm:p-dro-sample-chi2}} 
We denote $\hat{Z}^{P-DRO}(x) := \sup_{\P: \chi^2(\P,\hat{\Q}) \leq \varepsilon}\E_{\P}[h(x;\xi)]$ and the discrete approximation $\hat{Z}_m^{P-DRO}(x): = \sup_{\P: \chi^2(\P,\hat{\Q}_m)\leq \varepsilon}\E_{\P}[h(x;\xi)]$ here. The key is to show that $\sup_{x\in \Xscr}|\hat{Z}^{P-DRO}(x) - \hat{Z}_m^{P-DRO}(x)|$ is small so that we can borrow results from Theorem~\ref{thm:general-param-dro-ipm}. In the beginning, we present the error decomposition between the empirical variance and the true variance under $\P^*$. That is, with probability at least $1-\delta$:
\begin{equation}\label{eq:var-decomp}
    \begin{aligned}
    \left|\text{Var}_{\hat{\Q}_m}[h(x;\xi)] - \text{Var}_{\P^*}[h(x;\xi)]\right| &= \left|\text{Var}_{\hat{\Q}_m}[h(x;\xi)] - \text{Var}_{\hat{\Q}}[h(x;\xi)]\right| + \left|\text{Var}_{\hat{\Q}}[h(x;\xi)] - \text{Var}_{\P^*}[h(x;\xi)]\right|\\
    &\leq M^2\left(C_1\sqrt{\frac{\text{Comp}_m(\Hscr)}{m}} +3\sqrt{2\varepsilon}\right), \forall x \in \Xscr,
    \end{aligned}
\end{equation}
where the first term in the inequality follows from the uniform Hoeffding inequality. And the second term in the inequality follows from the choice of $\varepsilon$ such that $\chi^2(\P^*,\hat{\Q})\leq \varepsilon$ and $\|h\|_{\infty} \leq M$ in~\eqref{eq:var-decomp-chi2}.

The main proof is divided into the following three steps. 

\noindent\textit{Step 1: Variance Regularization.} Following from Lemma~\ref{lemma:chi2-cauchy}, we have:
\begin{equation}\label{eq:chi2-cauchy-theory}
    \sup_{\P: \chi^2(\P,\hat{\Q})\leq \varepsilon}\E_{\P}[h(x;\xi)]\leq \E_{\hat{\Q}}[h(x;\xi) ] + \sqrt{2\varepsilon\text{Var}_{\hat{\Q}}[h(x;\xi)]},
\end{equation}
\begin{equation}\label{eq:chi2-cauchy-sample}
    \E_{\hat{\Q}_m}[h(x;\xi)]\leq\sup_{\P: \chi^2(\P,\hat{\Q}_m)\leq \varepsilon}\E_{\P}[h(x;\xi)]\leq \E_{\hat{\Q}_m}[h(x;\xi) ] + \sqrt{2\varepsilon\text{Var}_{\hat{\Q}_m}[h(x;\xi)]},
\end{equation}
We now choose the required size $m$ so that the equality condition of RHS holds in~\eqref{eq:chi2-cauchy-sample}. 

Note that for any $\varepsilon$, the objective value $\sup_{\P: \chi^2(\P,\hat{\Q}_m)\leq \varepsilon}\E_{\P}[h(x;\xi)]$ is equivalent to the optimal objective value of the following optimization problem:
$$\max_{p\in \R^m_{+}} \sum_{i = 1}^m  p_i h(x;\xi_i),~\text{s.t.}:~\sum_{i = 1}^m\left(p_i - \frac{1}{m}\right)^2  \leq \frac{2\varepsilon}{m}, \sum_{i = 1}^m p_i = 1.$$
The optimal objective value above is the same as RHS of~\eqref{eq:chi2-cauchy-sample} if $\sqrt{2\varepsilon}\frac{h(x;\xi) - \E_{\hat{\Q}_m}[h(x;\xi)]}{\sqrt{\text{Var}_{\hat{\Q}_m}[h(x;\xi)]}} 
\geq -1$. Since $h(x;\xi) \in [0, M],\forall x, \xi$, this condition holds if:
\begin{equation}\label{eq:cauchy-eq-condition}
\text{Var}_{\hat{\Q}_m}[h(x;\xi)] \geq 2\varepsilon M^2,\forall x \in \Xscr.    
\end{equation}
In general, we obtain the following \textit{variance-dependent} lower bound of $\sup_{\P: \chi^2(\P, \hat{\Q}_m)\leq \varepsilon}\E_{\P}[h(x;\xi)]$:
\begin{equation}\label{eq:chi2-cauchy-lower}
\sup_{\P: \chi^2(\P, \hat{\Q}_m)\leq \varepsilon}\E_{\P}[h(x;\xi)] \geq \E_{\hat{\Q}_m}[h(x;\xi)] + \sqrt{\Delta \text{Var}_{\hat{\Q}_m}[h(x;\xi)]}, \end{equation}
as long as $\text{Var}_{\hat{\Q}_m}[h(x;\xi)] \geq \Delta M^2$.

For any integer $L \geq 1$ and the output $\hat{\Q}$, we partition the decision space $\Xscr$ into the following regions $\Xscr_1 \cup \ldots \Xscr_{L + 1} \cup \Xscr_{L + 2}$, where $\Xscr_{L + 2} = \{x\in \Xscr: \text{Var}_{\hat{\Q}}[h(x;\xi)] \geq 2\varepsilon M^2\}$, and:  
$$\Xscr_{\ell} = \left\{x \in \Xscr: \text{Var}_{\hat{\Q}}[h(x;\xi)] \in \left[\frac{\ell - 1}{L}2\varepsilon M^2, \frac{\ell}{L}2\varepsilon M^2\right)\right\}, \ell \in \{1,\ldots, L + 1\}.$$

We first choose the Monte Carlo size such that:
\begin{equation}\label{eq:Monte-Carlo-size-requirement}
    M C_1 \sqrt{\frac{\text{Comp}_m(\Hscr)}{m}}=:\Delta_{\E} \leq \frac{2\varepsilon M}{L},   
\end{equation}
so that with probability at least $1-\delta$, following from~\eqref{eq:var-decomp}:
\begin{equation}\label{eq:var-reg-sample}
\bigg|\text{Var}_{\hat{\Q}}[h(x;\xi)] - \text{Var}_{\hat{\Q}_m}[h(x;\xi)]\bigg| \leq \frac{2\varepsilon M^2}{L}, \forall x \in \Xscr.    
\end{equation}

    


    

\noindent\textit{Step 2: Monte Carlo Error Decomposition.} We 
partition the problem of bounding $\hat{Z}(x) - \hat{Z}_m(x)$ into the following different regimes of the decision space. 

\textit{(2.1) $\forall x \in \Xscr_{L + 2}$}, due to~\eqref{eq:var-reg-sample}, \eqref{eq:cauchy-eq-condition} holds. Therefore by~\eqref{eq:chi2-cauchy-theory} and~\eqref{eq:chi2-cauchy-lower}, we obtain:
\begin{equation}\label{eq:differ-sample}
    \begin{aligned}
        \hat{Z}^{P-DRO}(x) -\hat{Z}_m^{P-DRO}(x) &\leq \E_{\hat{\Q}}[h(x;\xi)] + \sqrt{2\varepsilon\text{Var}_{\hat{\Q}}[h(x;\xi)]} - \E_{\hat{\Q}_m}[h(x;\xi)] - \sqrt{2\varepsilon\text{Var}_{\hat{\Q}_m}[h(x;\xi)]}\\
        &=(\E_{\hat{\Q}}[h(x;\xi)] - \E_{\hat{\Q}_m}[h(x;\xi)]) + \sqrt{2\varepsilon}(\sqrt{\text{Var}_{\hat{\Q}}[h(x;\xi)]} - \sqrt{\text{Var}_{\hat{\Q}_m}[h(x;\xi)]}).
    \end{aligned}
\end{equation}
For the first term of RHS in~\eqref{eq:differ-sample}, we choose the Monte Carlo size $m$ such that $C_1 M\sqrt{\frac{\text{Comp}_m(\Hscr)}{m}}:=\Delta_{\E} \leq \frac{2\varepsilon M}{L}$ in~\eqref{eq:uniform-hoeffding} in Lemma~\ref{lemma:uniform-hoeffding}. Then, plugging~\eqref{eq:uniform-hoeffding} and~\eqref{eq:uniform-var} into~\eqref{eq:differ-sample}, with probability at least $1-\delta$, $\forall x \in \Xscr_{L + 2}$, we have:
\begin{equation*}
    \begin{aligned}
     \hat{Z}^{P-DRO}(x) - \hat{Z}_m^{P-DRO}(x) &\leq \Delta_{\E} + \sqrt{2\varepsilon}\left(C_2 M \sqrt{\frac{\text{Comp}_m(\Hscr)}{m}} + \left(1-\sqrt{1-\frac{1}{m}}\right)\sqrt{\text{Var}_{\hat{\Q}}[h(x;\xi)]} + \frac{2M^2}{m}\right)\\
     &\leq \Delta_{\E} + C_2M \sqrt{\frac{2\varepsilon\text{Comp}_m(\Hscr)}{m}}+\frac{\sqrt{2\varepsilon}(2M^2 + \text{Var}_{\hat{\Q}}[h(x;\xi)])}{m}\\
     &\leq C_3 \Delta_{\E} + \frac{3\sqrt{2\varepsilon}M^2}{m} \leq C_3' \Delta_E,\\
    \end{aligned}
\end{equation*}
where $C_3, C_3'$ are numerical constants independent of the function complexity and sample size. The second inequality follows from the fact that TV distance is an IPM.

\textit{(2.2) $\forall x \in \Xscr_i, i \in \{1,\ldots, L + 1\}$}, we have $\text{Var}_{\hat{\Q}_m}[h(x;\xi)] \geq \max\left\{\frac{i - 2}{L}2\varepsilon M^2, 0\right\}$.

\textit{(2.2.1) If $i \geq 2$,} following~\eqref{eq:chi2-cauchy-theory} and~\eqref{eq:chi2-cauchy-lower} as well as the definition of $\Xscr_i$, we have:
\begin{equation}\label{eq:differ-sample-corner}
    \begin{aligned}
        \hat{Z}^{P-DRO}(x) -\hat{Z}_m^{P-DRO}(x) &\leq \left(\E_{\hat{\Q}}[h(x;\xi)] - \E_{\hat{\Q}_m}[h(x;\xi)]\right)+ \sqrt{\frac{i}{L}2\varepsilon\text{Var}_{\hat{\Q}}[h(x;\xi)]} - \sqrt{\frac{i - 2}{L}2\varepsilon\text{Var}_{\hat{\Q}_m}[h(x;\xi)]}\\
        &\leq \Delta_{\E} + \sqrt{\frac{2(i-2)\varepsilon}{L}}\left(\sqrt{\text{Var}_{\hat{\Q}}[h(x;\xi)]} - \sqrt{\text{Var}_{\hat{\Q}_m}[h(x;\xi)]}\right) + \sqrt{\frac{4\varepsilon}{L} \text{Var}_{\hat{\Q}}[h(x;\xi)]}\\
        &\leq \Delta_{\E} + \sqrt{\frac{2(i-2)\varepsilon}{L}}\frac{\text{Var}_{\hat{\Q}}[h(x;\xi)] - \text{Var}_{\hat{\Q}_m}[h(x;\xi)]}{\sqrt{\text{Var}_{\hat{\Q}}[h(x;\xi)]} + \sqrt{\text{Var}_{\hat{\Q}_m}[h(x;\xi)]}} + \sqrt{\frac{4\varepsilon}{L} \frac{2i}{L}M^2}\\
        &\leq \Delta_{\E} + \sqrt{\frac{2(i-2)\varepsilon}{L}} \frac{2\varepsilon M^2 /L}{2\sqrt{2(i-2)\varepsilon M^2/L}} + C_4'\sqrt{\frac{\varepsilon}{L}}M \leq \Delta_{\E} + C_4 \frac{\varepsilon}{L} M + C_4'\sqrt{\frac{\varepsilon}{L}}M.
    \end{aligned}
\end{equation}

\textit{(2.2.2) If $i = 1$,} following from the definition of $\Xscr_1$ and \eqref{eq:chi2-cauchy-theory}, we have:
\begin{equation}\label{eq:differ-sample-corner-2}
    \begin{aligned}
        \hat{Z}^{P-DRO}(x) -\hat{Z}_m^{P-DRO}(x) &\leq \left(\E_{\hat{\Q}}[h(x;\xi)] - \E_{\hat{\Q}_m}[h(x;\xi)]\right)+ \sqrt{\frac{2\varepsilon}{L}\text{Var}_{\hat{\Q}}[h(x;\xi)]}\\
        &\leq \Delta_{\E} + \frac{2\varepsilon M}{L}.
    \end{aligned}
\end{equation}

In general, combining cases in (2.2), with probability at least $1- \delta$, if $\frac{\varepsilon}{L}\leq \sqrt{\frac{\varepsilon}{L}}\leq 1$, then:
\begin{equation}\label{eq:differ-result-all}
\hat{Z}^{P-DRO}(x) - \hat{Z}_m^{P-DRO}(x) \leq C_0 \sqrt{\frac{\varepsilon}{L}} M, \forall x \in \Xscr \backslash \Xscr_{L + 2}.
\end{equation}
\noindent\textit{Step 3: Generalization Error Decomposition. }Plugging the solution $x^{\PDRO_m}$ into~\eqref{eq:differ-result-all}, we have:
\begin{equation}\label{eq:differ-result}
\begin{aligned}
    &~Z(x^{\PDRO_m})- \hat{Z}_m^{P-DRO}(x^{\PDRO_m}) \\
    &\leq \hat{Z}^{P-DRO}(x^{\PDRO_m}) - \hat{Z}_m^{P-DRO}(x^{\PDRO_m}) \leq \begin{cases}\Delta_{\E} + C_0 \sqrt{\frac{\varepsilon}{L}} M,~&\text{ if }x^{\PDRO_m}\not\in \Xscr_{L+2}\\C_0^{\prime}\Delta_{\E},~&\text{ otherwise}\end{cases}.
\end{aligned}
\end{equation}
for some constant $C_0$ when $L$ is large. The first inequality follows from Assumption~\ref{asp:oracle-param-est} and Theorem~\ref{thm:general-param-dro-ipm}, with probability at least $1-\delta$, $\P^* \in \hat{\Ascr}$ when $\varepsilon \geq \Delta(\delta, \Theta)$. Besides,
\begin{equation}\label{eq:diff-result}
    \hat{Z}_m^{P-DRO}(x^*) - Z(x^*) \leq (\E_{\hat{\Q}_m}[h(x^*;\xi)] - \E_{\P^*}[h(x^*;\xi)]) + \sqrt{2\varepsilon \text{Var}_{\hat{\Q}_m}[h(x^*;\xi)]}.
\end{equation}
The first term of RHS in~\eqref{eq:diff-result} can be further bounded by Bernstein inequality. That is, with probability at least $1-\delta$:
\begin{equation}\label{eq:differ-result-1}
    \begin{aligned}
       \E_{\hat{\Q}_m}[h(x^*;\xi)] - \E_{\P^*}[h(x^*;\xi)]
        &\leq (\E_{\hat{\Q}_m}[h(x^*;\xi)] - \E_{\hat{\Q}}[h(x^*;\xi)]) + (\E_{\hat{\Q}}[h(x^*;\xi)] - \E_{\P^*}[h(x^*;\xi)])\\
        &\leq \sqrt{\frac{2\text{Var}_{\hat{\Q}}[h(x^*;\xi)]\log(1/\delta)}{m}}+ \frac{\|h(x^*;\cdot)\|_{\infty}\log(1/\delta)}{3m} + \sqrt{2\varepsilon \text{Var}_{\P^*}[h(x^*;\xi)]}\\
        &\leq \|h(x^*;\cdot)\|_{\infty}\left(\sqrt{\frac{2\log(1/\delta)}{m}}+ \frac{\log(1/\delta)}{3m}\right) + \sqrt{2\varepsilon \text{Var}_{\P^*}[h(x^*;\xi)]}.\\
    \end{aligned}
\end{equation}

Following~\eqref{eq:var-concentration-upper}, with probability at least $1-\delta$, the second term of RHS in~\eqref{eq:diff-result} can be bounded by:
\begin{equation}\label{eq:differ-result-2}
\begin{aligned}
\sqrt{2\varepsilon \text{Var}_{\hat{\Q}_m}[h(x^*;\xi)]} &\leq \sqrt{2\varepsilon \text{Var}_{\hat{\Q}}[h(x^*;\xi)]} + 2\|h(x^*;\cdot)\|_{\infty}\sqrt{\frac{\varepsilon\log(1/\delta)}{m}}.
\end{aligned}
\end{equation}
Then following the same decomposition procedure,
\begin{equation*}
\begin{aligned}
    Z(x^{\PDRO_m}) - Z(x^*) &\leq (Z(x^{\PDRO_m}) - \hat{Z}_m^{P-DRO}(x^{\PDRO_m})) + (\hat{Z}_m^{P-DRO}(x^{*}) - Z(x^*))\\
    &\leq 2 \Escr_{\chi}^2 + C_1 \sqrt{\frac{\varepsilon}{L}}M\mathbb{I}_{\{x^{\PDRO_m} \in \Xscr_{L + 2}\}},
\end{aligned}
\end{equation*}
where the second inequality follows from the observation that the error $Z(x^{\PDRO_m}) - \hat{Z}_m(x^{\PDRO_m})$ in~\eqref{eq:differ-result} is bounded by $\hat{Z}_m(x^*) - Z(x^*)$ in~\eqref{eq:diff-result} when the required sample size $m$ is chosen as the way in Theorem~\ref{thm:p-dro-sample-chi2}.
Therefore, we can attain the given generalization bound.$\hfill \square$

\subsubsection{Proof of Theorem~\ref{thm:p-dro-sample-w1}} 
Before presenting the proof, we introduce the following lemma demonstrating the regularization effects of the 1-Wasserstein DRO model:
\begin{lemma}[Extracted from Theorem 6.3 of \cite{esfahani2018data}]\label{lemma:w1-dro-reg}
Suppose $h(x;\xi)$ is Lipschitz continuous and convex w.r.t. $\xi$ and $\Xi$ is unbounded. There exists $\xi_0 \in \Xi$ such that $\limsup_{\|\tilde{\xi} - \xi_0\| \to \infty}\frac{h(x;\xi) - h(x;\xi_0)}{\|\tilde{\xi} -\xi_0\|} = \|h(x;\cdot)\|_{\text{Lip}}$, 
Then for any $\hat{\P}$ we have:
$$\sup_{W_1(\P, \hat{\P}) \leq \varepsilon} \E_{\P}[h(x;\xi)] = \E_{\hat{\P}}[h(x;\xi)]+ \varepsilon \|h(x;\cdot)\|_{\text{Lip}}.$$
\end{lemma}
\textit{Proof of Theorem~\ref{thm:p-dro-sample-w1}.}~On one hand, with probability at least $1-\delta$, we have:
\begin{equation}\label{eq:w1-dro-sample-true}
\begin{aligned}
\E_{\P^*}[h(x^{\PDRO_m};\xi)]&\leq \sup_{\P: W_1(\P, \hat{\Q}) \leq \varepsilon} \E_{\P}[h(x^{\PDRO_m};\xi)]\\
&= \E_{\hat{\Q}}[h(x^{\PDRO_m};\xi)] + \varepsilon \|h(x^{\PDRO_m};\cdot)\|_{\text{Lip}}\\
&= \E_{\hat{\Q}}[h(x^{\PDRO_m};\xi)] + \left(\sup_{\P: W_1(\P, \hat{\Q}_m) \leq \varepsilon} \E_{\P}[h(x^{\PDRO_m};\xi)]  - \E_{\hat{\Q}_m}[h(x^{\PDRO_m};\xi)]\right),    
\end{aligned}
\end{equation}
where the first inequality follows from $W_1(\P^*,\hat{\Q}) \leq \varepsilon$ under the ambiguity size $\varepsilon \geq \Delta(\delta, \Theta)$. The first and second equalities follow from Lemma~\ref{lemma:w1-dro-reg} when taking the distribution center $\hat{\P}$ to be $\hat{\Q}$ and $\hat{\Q}_m$ respectively.

On the other hand, we have:
\begin{align*}
    \sup_{\P: W_1(\P, \hat{\Q}_m) \leq \varepsilon} \E_{{\P}}[h(x^*;\xi)] &\leq \E_{\hat{\Q}_m}[h(x^*;\xi)] + \varepsilon \|h(x^*;\cdot)\|_{\text{Lip}}\\
    &\leq \E_{\hat{\Q}_m}[h(x^*;\xi)] + (\E_{\P^*}[h(x^*;\xi)] - \E_{\hat{\Q}}[h(x^*;\xi)]) + 2\varepsilon \|h(x^*;\cdot)\|_{\text{Lip}}.
\end{align*}
Therefore, we have:
\begin{align*}
    \mathcal{E}(x^{\PDRO_m}) & \leq \left(\E_{\P^*}[h(x^{\PDRO_m};\xi)] - \sup_{\P: W_1(\P, \hat{\Q}_m) \leq \varepsilon} \E_{\P}[h(x^{\PDRO_m};\xi)]\right) + \left(\sup_{\P: W_1(\P, \hat{\Q}_m) \leq \varepsilon} \E_{{\P}}[h(x^*;\xi)] - \E_{\P^*}[h(x^*;\xi)]\right)\\
    &\leq 2\|h(x^*;\cdot)\|_{\text{Lip}}\varepsilon + 2\sup_{x\in \Xscr}\bigg|\E_{\hat{\Q}}[h(x;\xi)] - \E_{\hat{\Q}_m}[h(x;\xi)]\bigg|\\
    &\leq 2\|h(x^*;\cdot)\|_{\text{Lip}}\varepsilon + 2C_1 M\sqrt{\frac{\text{Comp}_m(\Hscr)}{m}}. 
\end{align*}
Then the required sample size $m$ is chosen such that the second term $C_1 M\sqrt{\frac{\text{Comp}_m(\Hscr)}{m}}$ above is smaller than $\|h(x^*;\cdot)\|_{\text{Lip}}\varepsilon$.
$\hfill \square$
\subsection{Generalization Results for $p$-Wasserstein Distance}\label{app:p-was-general}
In this part, we extend the result of bounding the Monte Carlo sampling error to the case when $d$ is taken as $p$-Wasserstein distance with $p \in (1,2]$. 
\begin{definition}[$p$-Wasserstein Distance]\label{def:p-was}
    $p$-Wasserstein distance $(1\leq p < \infty)$ between two distributions $\P$ and $\Q$ supported on $\Xi$ is defined as:
\begin{equation*}
        W_p(\P,\Q) = \inf_{\Pi \in \Mscr{(\Xi\times \Xi)}}\left\{\left(\int_{\Xi\times \Xi} \|x - y\|^p\Pi(dx,dy)\right)^{\frac{1}{p}}: \Pi_x=\P, \Pi_y=\Q\right\},
\end{equation*}
where $\Pi_x$ and $\Pi_y$ are the marginal distributions of $\Pi$.
\end{definition}
We present the following lemma establishing the regularization effect of $p$-Wasserstein distance.
\begin{lemma}[Extracted from Lemma 1 in \cite{gao2022wasserstein}]\label{lemma:wp-dro-reg}
Suppose some mild conditions hold for $\Hscr$ (i.e. the same conditions of Assumption 1 and 2 from \cite{gao2022wasserstein}, only replacing $\Fscr$ there with $\Hscr$). Consider $p$-Wasserstein distance with $p \in (1, 2]$, for any distribution $\hat{\Q}$, there exists $C > 0$ such that:
$$\bigg|\sup_{\P: W_p(\P, \hat{\Q}) \leq \varepsilon} \E_{\P}[h(x;\xi)] - \E_{\hat{\Q}}[h(x;\xi)] - \varepsilon \Vscr_{\hat{\Q}, q}(h(x;\cdot))\bigg| \leq C \varepsilon^{p}.$$
where $\Vscr_{\hat{\Q}, q}(h(x;\cdot))$ is the $L_{q}$ norm of the vectorized random variable $\nabla_{\xi} h(x;\xi)$ under the measure $\hat{\Q}$ with $\frac{1}{p} + \frac{1}{q} = 1$.
\end{lemma}
We abbreviate $h:=h(x;\cdot), h^*:= h(x^*;\cdot)$ in the following.  
\begin{corollary}\label{coro:p-dro-sample-wp}
Suppose Assumption~\ref{asp:oracle-param-est} in the main body and Assumption 1 and 2 in \cite{gao2022wasserstein} hold. The size of the ambiguity set $\varepsilon \geq \Delta(\delta, \Theta)$, when $d$ is $p$-Wasserstein distance with $p \in (1, 2]$, if the Monte Carlo size satisfies:
$$m \geq \max\left\{\left(\frac{C_1 + C_2 \tilde{M}\sqrt{\text{Comp}_m(\partial(\Hscr))}}{ \Vscr_{\P^*,q}(h^*)}\right)^q, C_0\left(\frac{M}{\varepsilon \Vscr_{\P^*, q}(h^*)}\right)^2\text{Comp}_m(\Hscr)\right\},$$
where $\tilde{M}:=\sup_{x\in \Xscr, \xi \in \Xi}\left\|\nabla_{\xi}h(x;\xi)\right\|_2$ and $\partial(\Hscr)= \left\{\left\|\nabla_{\xi}h(x;\xi)\right\|_2^q: x\in \Xscr\right\}$ for some constants $C_0, C_1, C_2$, Then with probability at least $1-\delta$, we have $\mathcal{E}(x^{\PDRO_m}) \leq 4 \varepsilon \Vscr_{\P^*, q}(h(x^*;\cdot)) + C\varepsilon^p$.
\end{corollary}
\proof{Proof of Corollary~\ref{coro:p-dro-sample-wp}.}Following from the result in Lemma~\ref{lemma:wp-dro-reg}, we have:
\begin{align*}
 \E_{\P^*}[h] - \sup_{\P: W_p(\P, \hat{\Q}_m) \leq \varepsilon}\E_{\P}[h]&\leq \sup_{\P: W_p(\P, \hat{\Q}) \leq \varepsilon} \E_{\P}[h] - \sup_{\P: W_p(\P, \hat{\Q}_m) \leq \varepsilon} \E_{\P}[h]\\
 &\leq (\E_{\hat{\Q}}[h] - \E_{\hat{\Q}_m}[h]) + \varepsilon (\Vscr_{\hat{\Q},q}(h) - \Vscr_{\hat{\Q}_m,q}(h)) + 2C \varepsilon^{p}.  
\end{align*}
Therefore, we obtain:
\begin{align*}
\mathcal{E}(x^{\PDRO_m}) &\leq 2\varepsilon \Vscr_{\P^*,q}(h^*) + C\varepsilon^{p} + \sup_{h \in \Hscr}\left|\E_{\P}[h] - \E_{\hat{\Q}_m}[h]\right| + \varepsilon\sup_{h \in \Hscr}\left|\Vscr_{\hat{\Q},q}(h) - \Vscr_{\hat{\Q}_m,q}(h)\right|\\
&\leq 2\varepsilon \Vscr_{\P^*,q}(h^*) + C\varepsilon^{p} + C_0 M \sqrt{\frac{\text{Comp}_m(\Hscr)}{m}}+ \varepsilon\left[C_1 + C_2 \tilde{M}\sqrt{\text{Comp}_m(\partial(\Hscr))}\right] m^{\frac{1}{p} - 1},
\end{align*}
And the last term of the second inequality follows from the uniform concentration inequality of $L_p$-norm for $\sup_{h \in \Hscr}|\Vscr_{\hat{\Q},q}(h) - \Vscr_{\hat{\Q}_m,q}(h)|$. Specifically, following from Theorem 6.10 in \cite{boucheron2013concentration} and Lemma 7 in \cite{duchi2021learning}, $\forall h \in \Hscr$, with probability at least $1-\delta$:
$$|\Vscr_{\hat{\Q}_m, q}(h) - \E[\Vscr_{\hat{\Q}_m, q}(h)]| \leq \tilde{M}m^{-\frac{1}{q}}\sqrt{\log(1/\delta)}.$$
Then following the same covering number argument to $\partial(\Hscr)$ as Lemma~\ref{lemma:uniform-var}, with probability at least $1-\delta$, we have:
\begin{equation}\label{eq:lp-norm-1}
|\Vscr_{\hat{\Q}_m, q}(h) - \E[\Vscr_{\hat{\Q}_m, q}(h)]| \leq C_2\tilde{M}m^{-\frac{1}{q}}\sqrt{\text{Comp}_m(\partial(\Hscr))}, \forall h\in \Hscr.
\end{equation}

Following from Lemma 9 in \cite{duchi2021learning} and the definition of $\Vscr_{\Q, q}(h)$, we have:
\begin{equation}\label{eq:lp-norm-2}
\Vscr_{\hat{\Q}, q}(h) - \frac{2}{p}\sqrt{C}n^{-\frac{1}{q}} \leq  \E[\Vscr_{\hat{\Q}_m, q}(h)] \leq \Vscr_{\hat{\Q}, q}(h).
\end{equation}

Combining~\eqref{eq:lp-norm-1} and~\eqref{eq:lp-norm-2}, we would obtain the bound for $\sup_{h \in \Hscr}|\Vscr_{\hat{\Q},q}(h) - \Vscr_{\hat{\Q}_m,q}(h)|$. 

Finally, the required Monte Carlo sample size $m$ is chosen such that:
\[\max\left\{C_0 M \sqrt{\frac{\text{Comp}_m(\Hscr)}{m}}, \varepsilon\left[C_1 + C_2 \tilde{M}\sqrt{\text{Comp}_m(\partial(\Hscr))}\right] m^{\frac{1}{p} - 1}\right\}\leq \varepsilon \Vscr_{\P^*, q}(h^*).\]
$\hfill \square$


\subsection{Proof of Theorem~\ref{thm:general-param-erm-sample}}
Denote $x^{\PERM_m} \in \argmin_{x \in \Xscr}\E_{\hat{\Q}_m}[h(x;\xi)]$, 
$\hat{Z}^{P-ERM}(x) := \E_{\hat{\Q}}[h(x;\xi)]$ and $\mathcal{E}_{P}$ being the upper bound of $Z(x^{\PERM}) - Z(x)$.
The result is directly from Lemma~\ref{lemma:nonparam-erm} and Theorem~\ref{thm:general-param-erm}. Specifically, we have:
\begin{equation}\label{eq:nonparam-erm-mc}
    \begin{aligned}
    Z(x^{\PERM_m}) - Z(x^*) &\leq \mathcal{E}_{P} + 2\sup_{x \in \Xscr}\bigg|\E_{\hat{\Q}}[h(x;\xi)] - \E_{\hat{\Q}_m}[h(x;\xi)]\bigg|\\
    &\leq \mathcal{E}_{P} +  2\left(\sqrt{\frac{\hat{Z}^{P-ERM}(x^{\PERM})\text{Comp}_m(\Hscr)M}{m}} + \frac{\text{Comp}_m(\Hscr)M}{m}\right)\\
    &\leq \mathcal{E}_{P} +  3\left(\sqrt{\frac{(Z(x^*) + \mathcal{E}_{P})\text{Comp}_m(\Hscr)M}{m}}\right) \leq 2\mathcal{E}_{P},\\
    \end{aligned}
\end{equation}
where the second inequality in~\eqref{eq:nonparam-erm-mc} follows from Lemma~\ref{lemma:nonparam-erm} since $x^{\PERM} \in \argmin_{x \in \Xscr}\E_{\hat{\Q}}[h(x;\xi)]$. And the third inequality in~\eqref{eq:nonparam-erm-mc} follows from the Monte Carlo size $m \geq \frac{M\text{Comp}_m(\Hscr)}{Z(x^*) + \mathcal{E}_{P}}$ and a result of inequalities based on Theorem~\ref{thm:general-param-erm}:
\begin{align*}
    \hat{Z}^{P-ERM}(x^{\PERM})&\leq \hat{Z}^{P-ERM}(x^*)=Z(x^*) + (\hat{Z}^{P-ERM}(x^*) - Z(x^*))\leq Z(x^*) + \mathcal{E}_{P}.
\end{align*}
The fourth inequality in~\eqref{eq:nonparam-erm-mc} follows from $m \geq \frac{(Z(x^*) + \mathcal{E}_{P}) M\text{Comp}_m(\Hscr)}{\mathcal{E}_{P}^2}$.
$\hfill \square$



\subsection{A Short Discussion of Stochastic Approximation Methods}\label{app:sa}
Besides the SAA approach mentioned in the main body, stochastic approximation (SA) comprises another common approach to solve stochastic optimization problems with underlying continuous distribution $\hat{\Q}$. In SA, we apply stochastic gradient descent to obtain a batch of samples from $\hat{\Q}$ at each step of each iteration. For example, in the \PERM\ case ($\min_{x \in \Xscr}\E_{\hat{\Q}}[h(x;\xi)]$), we can obtain a solution $\hat{x}$ with the expected generalization error after a polynomial number of iterations w.r.t. $\frac{1}{\gamma} (\gamma > 0)$ (such as \cite{nemirovski2009robust}):
\begin{equation}\label{eq:so-sa}
\E_{\hat{x}}[Z(\hat{x}) - Z(x^*)] \leq \mathcal{E}(x^{\PERM}) + \gamma.    
\end{equation}
Additionally, we can express our optimization problem as $\min_{x,y\in \Xscr\times \Yscr}\E_{\hat{\Q}}[G(x,y)]$ for some auxiliary variable $y$ and apply SA for some DRO formulations such as $f$-divergence ($y = (\lambda, \mu)$ in the dual problem). 


\section{Further Details and Proofs for Section~\ref{subsec:dist-shift}}\label{subsec:dist-shift-proof2}
Besides the notations in Section~\ref{subsec:dist-shift}, we denote $x^{tr} \in \argmin_{x\in\Xscr}\left\{Z^{tr}(x) := \E_{\P^{tr}}[h(x;\xi)]\right\}$ and $\hat{\P}^{tr}_n$ as the empirical distribution of the training set.  We show how each bound under distribution shifts is derived to explain the rationale of the curve shape in Figure~\ref{fig:concept}(b).
We use $\chi^2$-divergence between $\P^{te}$ and $\P^{tr}$ to evaluate the extent of distribution shifts. We also assume $\text{Var}_{\P^{tr}}[h(x^{tr};\xi)]\approx \text{Var}_{\P_{te}}[h(x^*;\xi)]$, i.e. ``variability" does not change across shifts.
For each data-driven solution $\hat{x}$ being a minimizer of $\hat Z(\cdot)$ where $\hat{Z}(\cdot)$ is some objective estimated from the training data, we have:
$$Z^{te}(\hat x)-Z^{te}(x^*)=[Z^{te}(\hat x)-\hat Z(\hat x)]+[\hat Z(\hat x)-\hat Z(x^*)]+[\hat Z(x^*)-Z^{te}(x^*)],$$
where the middle term $[\hat Z(\hat x)-\hat Z(x^*)]$ is at most 0. Then we only need to bound:
\begin{equation}\label{eq:decomposition-dist-shift}
[Z^{te}(\hat x)-\hat Z(\hat x)]+[\hat Z(x^*)-Z^{te}(x^*)],
\end{equation}
where the second term in~\eqref{eq:decomposition-dist-shift} can be further decomposed in the form of $\Vscr_d(x^*) d(\P^{tr}, \P^{te})$: 
$$(\hat Z(x^*)-Z^{tr}(x^*)) + (Z^{tr}(x^*) - Z^{te}(x^*)) \leq \hat Z(x^*)-Z^{tr}(x^*) + \Vscr_d(x^*) d(\P^{tr}, \P^{te}).$$
For example, using Lemma~\ref{lemma:chi2-cauchy}, we have:
$Z^{tr}(x^*) - Z^{te}(x^*) \leq \sqrt{2\chi^2(\P^{tr},\P^{te})\text{Var}_{\P^{te}}[h(x^*;\xi)]}.$

The term $\Delta^*:=Z^{tr}(x^*) - Z^{te}(x^*) \leq \Vscr_d(x^*) d(\P^{tr}, \P^{te})$ cannot be avoided in generalization error bounds across all methods. That is why we set the minimum value in the y-axis as $\Vscr_d(x^*) d(\P^{tr}, \P^{te})$ in Figure~\ref{fig:concept}(b). 
When comparing the generalization error bound in each method, we focus on the additional error term, in addition to $\Delta^*$ and the error without distribution shifts. 


\subsection{Generalization Errors of Existing ERM and DRO Approaches under Distribution Shifts}\label{app:erm-dro-dist-shift} 

\paragraph{ERM.} Consider $x^{\texttt{NP-ERM}} \in \argmin_{x \in \Xscr}\left(\hat{Z}^{ERM}(x): = \E_{\hat{\P}^{tr}_n}[h(x;\xi)]\right)$. Denote $\mathcal{E}^{ERM}(\P^{tr}, \Hscr)$ as the generalization error of \texttt{NP-ERM} in Lemma~\ref{lemma:nonparam-erm} while replacing $Z, x^*$ there with $Z^{tr}, x^{tr}$, and $\Escr^{ERM}(\P^{tr}, \P^{te}, \Hscr):=  d(\P^{te}, \P^{tr})M^{\frac{3}{4}}\left(\frac{\text{Comp}_n(\Hscr)}{n}\right)^{\frac{1}{4}}$. Then: 
\begin{align*}
    Z^{te}(x^{\texttt{NP-ERM}})-\hat Z^{tr}(x^{\texttt{NP-ERM}})&=[Z^{te}(x^{\texttt{NP-ERM}})-Z^{tr}(x^{\texttt{NP-ERM}})]+[Z^{tr}(x^{\texttt{NP-ERM}})-\hat Z^{ERM}(x^{\texttt{NP-ERM}})],\\
    &\leq [Z^{te}(x^{\texttt{NP-ERM}})-Z^{tr}(x^{\texttt{NP-ERM}})] + \mathcal{E}^{ERM}(\P^{tr}, \Hscr),
\end{align*}
where the first term above can be further bounded by:
\begin{align*}
Z^{te}(x^{\texttt{NP-ERM}})-Z^{tr}(x^{\texttt{NP-ERM}})&=\E_{\P^{tr}}\left[\left(\frac{d\P^{te}}{d\P^{tr}}-1\right)h(x^{\texttt{NP-ERM}};\xi)\right]\leq\sqrt{2\chi^2(\P^{te},\P^{tr})\text{Var}_{\P^{tr}}[h(x^{\texttt{NP-ERM}};\xi)]}\\
&\leq\sqrt{2\chi^2(\P^{te},\P^{tr}) \left(\E_{\P^{tr}}[h^2(x^{tr};\xi)]  + M\sqrt{\frac{\text{Comp}_n(\Hscr)M}{n}}\right)}\\
&\leq \chi^2(\P^{te}, \P^{tr}) \Vscr_d(x^*) + d(\P^{te}, \P^{tr})M^{\frac{3}{4}}\left(\frac{\text{Comp}_n(\Hscr)}{n}\right)^{\frac{1}{4}}\\
& = \Delta^* + \Escr^{ERM}(\P^{tr}, \P^{te}, \Hscr),
\end{align*}
where the first inequality follows from $\E_{\P^{tr}}[h^2(x^{\texttt{NP-ERM}};\xi)] - \E_{\P^{tr}}[h^2(x^{tr};\xi)] \leq 2M (Z^{tr}(x^{\texttt{NP-ERM}}) - Z^{tr}(x^{tr}))$. Then we use Lemma~\ref{lemma:nonparam-erm} to bound it further. Thus the generalization error is bounded by:
$$Z^{te}(x^{\texttt{NP-ERM}})-Z^{te}(x^*)\leq \mathcal{E}^{ERM}(\P^{tr},\P^{te},\Hscr) + \mathcal{E}^{ERM}(\P^{tr},\Hscr) + \Delta^*.$$
Therefore, we incur an additional term $\mathcal{E}^{ERM}(\P^{tr}, \P^{te}, \Hscr)$ in the error bound compared with the case without distribution shifts. See similar results in \cite{ben2010theory,lee2018minimax}.


Then, we consider the standard (nonparametric) DRO problem, i.e.:
\[x^{\texttt{NP-DRO}}\in \argmin_{x \in \Xscr}\left\{ \hat{Z}^{DRO}(x):=\max_{\Q: d(\Q,\hat{\P}^{tr}_n)\leq\varepsilon}\E_{\Q}[h(x;\xi)]\right\}.\] 

\paragraph{DRO from the regularization perspective.}  If we choose the ambiguity size $\varepsilon = O\left(\left(\frac{\text{Comp}_n(\Hscr)}{n}\right)^{\beta}\right)$, \eqref{eq:dro-reg} holds when replacing $Z, x^*$ there with $Z^{tr}, x^{tr}$. After the replacement, we denote $\mathcal{E}^{DRO_1}(\P^{tr}, \Hscr;\Ascr)$ 
as that generalization error bound. And denote $\mathcal{E}^{DRO}(\P^{tr}, \P^{te}, \Hscr) := d(\P^{te}, \P^{tr})\Vscr_{d}^{\frac{1}{2}}(x^*)M^{\frac{1}{2}}(\frac{\text{Comp}_n(\Hscr)}{n})^{\frac{1}{4}}$. Following the same error decomposition above, we obtain the following bound:
\begin{align*}
    Z^{te}(x^{\texttt{NP-DRO}})-Z^{te}(x^*)&\leq \sqrt{d(\P^{te}, \P^{tr})\left(\Vscr_d^2(x^*) + M \Vscr_d(x^*)\sqrt{\frac{\text{Comp}_n(\Hscr)}{n}}\right)} + \mathcal{E}^{DRO}(\P^{tr}, \Hscr;\Ascr)\\
     & \leq \mathcal{E}^{DRO}(\P^{tr}, \P^{te}, \Hscr) + \Delta^* + \mathcal{E}^{DRO_1}(\P^{tr}, \Hscr;\Ascr).
\end{align*}
Therefore, we incur an additional term $\mathcal{E}^{DRO}(\P^{tr}, \P^{te}, \Hscr)$ in the error bound compared with the case without distribution shifts. 
This term is smaller than $\mathcal{E}^{ERM}(\P^{tr}, \P^{te}, \Hscr)$ but still depends on $\text{Comp}(\Hscr)$.

We see that the generalization error in each method above incurs an additional term $\Escr^{\cdot}(\P^{tr}, \P^{te},\Hscr)$, i.e., $\frac{C}{n^{1/4}}$ for a large $C$ depending on $\text{Comp}_n(\Hscr)$. When the sample size is small, generalization error bounds under distribution shifts for these methods are larger than that of \PDRO\ (see Corollary~\ref{coro:param-dro-dist-shift}). That is why we draw the curve shape for these two methods in Figure~\ref{fig:concept}(b) when $n$ is small. We do not draw the curve of the generalization error bound for the nonparametric \textit{DRO method from the robustness perspective} but mention it here for completeness.

\paragraph{DRO from the robustness perspective.} When $d$ is an IPM,
if we choose the size $\varepsilon \geq d(\P^{te}, \P^{tr}) + O\left(n^{-1/g(D_{\xi})}\right)$, then~\eqref{eq:dro-pessim} holds when replacing $Z, x^*$ there with $Z^{tr}, x^{tr}$. After the replacement, we denote $\mathcal{E}^{DRO_2}(\P^{tr}, \Hscr;\Ascr)$ as that generalization error bound. We have:
\begin{equation}\label{eq:dro-pessim-dist-shift}
\begin{aligned}
Z^{te}(x^{\texttt{NP-DRO}})-Z^{te}(x^*)&\leq \Vscr_d(x^*)(n^{-1/g(D_{\xi})} + d(\P^{te},\P^{tr}))\\
& \leq \mathcal{E}^{DRO_2}(\P^{tr}, \Hscr;\Ascr) + \Delta^*.
\end{aligned}
\end{equation}
In the case of the 1-Wasserstein DRO model under distribution shifts, 
the bound would be typical of the order $\|h(x^*;\cdot)\|_{Lip}\left(n^{-1/D_{\xi}}+d(\P^{tr},\P^{te})\right)$ (see Theorem E.3 in~\cite{zeng2021generalization}). 



\subsection{Proofs of Corollaries~\ref{coro:param-dro-dist-shift} and~\ref{coro:param-erm-dist-shift}} 
Before presenting the results, we introduce the following lemma, which indicates \eqref{eq:dist-shift} in Assumption~\ref{asp:oracle-dist-shift} holds for some $c_1, c_2$.
\begin{lemma}[Pseudo triangle inequality for some $f$-divergence]\label{lemma:triangle-f-div}
Suppose the three distributions $\P^{te}, \P^{tr}, \hat{\Q}$ have the same support, we have:
$$\chi^2(\P^{te},\hat{\Q}) \leq 2 \left\|\frac{d\P^{tr}}{d\hat{\Q}}\right\|_{\infty} \chi^2(\P^{te},\P^{tr}) + 2 \chi^2(\P^{tr},\hat{\Q}).$$
$$KL(\P^{te},\hat{\Q}) \leq KL(\P^{te},\P^{tr}) + \left\|\frac{d\P^{te}}{d\P^{tr}}\right\|_{\infty} KL(\P^{tr},\hat{\Q}).$$
\end{lemma}
\proof{Proof of Lemma~\ref{lemma:triangle-f-div}.} Since we only consider the distribution class $\Pscr_{\Theta}$ with continuous distributions, we denote the density of $\P^{te}, \P^{tr},\hat{\Q}$ as $f, g, h$ under Lebesgue measure $\mu$ respectively. 

For $\chi^2$-divergence, we have:
$$\int \frac{(f-h)^2}{h} d\mu \leq \int \frac{2(f-g)^2 + 2(g-h)^2}{h} d\mu \leq 2 \left\|\frac{g}{h}\right\|_{\infty}\int \frac{(f-g)^2}{g}d\mu + 2\int \frac{(g-h)^2}{h}d\mu.$$

For KL-divergence, we have:
$$\int f \ln \frac{f}{h} d\mu  = \int f \left(\ln \frac{f}{g} + \ln \frac{g}{h}\right)d\mu \leq \int f \ln \frac{f}{g} d\mu + \bigg\|\frac{f}{g}\bigg\|_{\infty}\int g \ln \frac{g}{h} d\mu. ~\hfill \square$$

For the \texttt{P-DRO} problem (i.e. Corollary~\ref{coro:param-dro-dist-shift}), denote $\hat{Z}^{P-DRO}(x) := \sup_{d(\P, \hat{\Q})\leq \varepsilon}\E_{\P}[h(x;\xi)]$. If $\varepsilon \geq c_1 d(\P^{te}, \P^{tr}) + c_2 \Delta(\delta, \Theta) $, when $d$ is an IPM, by~\eqref{eq:decomposition-dist-shift}, with probability at least $1-\delta$, we have:
\begin{align*}
    \mathcal{E}(x^{\PDRO}) &\leq Z^{te}(x^{\PDRO}) - \hat{Z}^{P-DRO}(x^{\PDRO}) + \hat{Z}^{P-DRO}(x^*) - Z^{te}(x^*)\\
    &\leq 0 + \max_{\P: d(\P, \hat{\Q}) \leq \varepsilon}\E_{\P}[h(x^*;\xi)] - \E_{\P^{te}}[h(x^*;\xi)]\\
    &\leq 2\Vscr_d(x^*)\varepsilon,
\end{align*}
where the second inequality follows from $\P[d(\P^{te}, \hat{\Q})\leq \varepsilon] \geq 1-\delta$ such that $Z^{te}(\cdot)\leq \hat{Z}(\cdot)$ due to Assumption~\ref{asp:oracle-dist-shift} with probability at least $1-\delta$. And the third inequality is the same as in case $(a)$ in the proof of Theorem~\ref{thm:general-param-dro-ipm}. Other cases of metrics $d$ follow from the same proof argument as cases $(b)$ and $(c)$ in Theorem~\ref{thm:general-param-dro-ipm}. 

For the \texttt{P-ERM} problem (i.e. Corollary~\ref{coro:param-erm-dist-shift}), denote $\hat{Z}^{P-ERM}(x) := \E_{\hat{\Q}}[h(x;\xi)]$. When $d$ is an IPM, we have:
\begin{align*}
    \mathcal{E}(x^{\PERM}) &\leq 2\sup_{x \in \Xscr}|Z^{te}(x) - \hat{Z}^{P-ERM}(x)|\\
    &\leq 2(\sup_{x \in \Xscr}\Vscr_d(x) )d(\P^{te}, \hat{\Q}) \leq 2(\sup_{x \in \Xscr}\Vscr_d(x))(d(\P^{te}, \P^{tr}) + d(\P^{tr}, \hat{\Q})).
\end{align*}
Other cases of metrics $d$ in Corollary~\ref{coro:param-erm-dist-shift} follow from the same argument by replacing $\Delta(\delta, \Theta)$ in the proof of Theorem~\ref{thm:general-param-erm} with $c_1 d(\P^{te}, \P^{tr}) + c_2\Delta(\delta, \Theta) $ following Assumption~\ref{asp:oracle-dist-shift}. $\hfill \square$
%

\section{Further Details and Proofs for Section~\ref{sec:context}}\label{app:context-proof}
\subsection{Proof of Example~\ref{prop:conditional-oracle-estimator}}
We abbreviate $\hat{f}$ to be $f_{\hat{\theta}}$ and $f^*$ to be $f_{\theta^*}$. Before we prove the result, we introduce some technical assumptions for the model $\xi = f(y) + \eta$ and least square estimator $\hat{f}$:
\begin{assumption}\label{asp:new}

We assume:
\begin{itemize}
    \item Pointwise convergence, $\P(|f^*(y) - \hat{f}(y)| \geq \eta)\leq C\exp(-a_n\eta^2)$, $\forall y$, which is also the assumption in Theorem 8 in \cite{hu2022fast};
    \item $\Fscr$ is star-shaped;
    \item $\eta$ is bounded, $\|\eta\|_2^2 \leq C_{\eta}$. 
\end{itemize}
\end{assumption}
Here, the star-shaped condition of Assumption~\ref{asp:new} implies that the term $\|\hat{f}_i - f_i^*\|_n^2 := \sqrt{\frac{1}{n}\|\hat{f}_i(\hat{y}_i) - f_i^*(\hat{y}_i)\|_2^2} \leq C_0\sqrt{\frac{\log(1/\delta)}{n}}, \forall i \in [D_{\xi}]$ following from Lemma 9 in \cite{hu2022fast}.

By the triangle inequality, we have:
$$W_1(\P_{\xi|y}^*, \hat{\Q}_{\xi|y})\leq W_1(\P_{\xi|y}^*, \Q_{\xi|y}^*) + W_1(\Q_{\xi|y}^*, \tilde{\Q}_{\xi|y}) + W_1(\tilde{\Q}_{\xi|y}, \hat{\Q}_{\xi|y}),$$
where $\Q_{\xi|y}^* = \Nscr(f^*(y), \Sigma)$ and $\tilde{\Q}_{\xi|y} = \Nscr(\hat{f}(y), \Sigma)$.

\textit{The first term $W_1(\P_{\xi|y}^*, \Q_{\xi|y}^*)$} is the term $\Escr_{apx}(y)$.

\textit{The second term $W_1(\Q_{\xi|y}^*, \tilde{\Q}_{\xi|y})$} is the mean difference of two Gaussian distributions: $W_1(\Q_{\xi|y}^*, \tilde{\Q}_{\xi|y}) = \|f^*(y) - \hat{f}(y)\|_2$, then we obtain $W_1(\Q_{\xi|y}^*, \tilde{\Q}_{\xi|y}) \leq C_0\sqrt{\frac{\log(1/\delta)}{n}}$ following from Assumption~\ref{asp:new} and Theorem 8 of \cite{hu2022fast}.


\textit{The third term $W_1(\tilde{\Q}_{\xi|y}, \hat{\Q}_{\xi|y})$} is bounded by:
\begin{equation}\label{eq:context-cov-error}
W_1(\tilde{\Q}_{\xi|y}, \hat{\Q}_{\xi|y}) \leq \text{Tr}[(\sqrt{\Sigma} - \sqrt{\hat{\Sigma}})]^2\leq \text{Tr}[(\sqrt{\Sigma} - \sqrt{\hat{\Sigma}^*})]^2 + \text{Tr}[(\sqrt{\hat{\Sigma}} - \sqrt{\hat{\Sigma}^*})]^2    
\end{equation}
where $\hat{\Sigma}^* = \frac{1}{n}\sum_{i = 1}^n(\hat{\xi}_i - f^*(\hat{y}_i))(\hat{\xi}_i - f^*(\hat{y}_i))^{\top}=:\frac{1}{n}\sum_{i = 1}^n \hat{\eta}_i \hat{\eta}_i^{\top}$ is the sample covariance matrix of true noise. And the first term of RHS in~\eqref{eq:context-cov-error} can be bounded by the standard matrix concentration inequality. Specifically, following from Corollary 2 in \cite{delage2010distributionally}, with probability at least $1-\delta$, $(1-\alpha(n))\hat{\Sigma}^* \preceq \Sigma \preceq (1+\alpha(n))\hat{\Sigma}^*$ holds with $\alpha(n) = \frac{C_{\eta}^2(1 + \sqrt{\log(1/\delta)})}{\sqrt{n}}$. When $n$ is large enough such that $(\sqrt{1+\alpha(n)} - 1)^2 \leq \alpha(n)$, we have:
$$\text{Tr}[(\sqrt{\Sigma}-\sqrt{\hat{\Sigma}^*})]^2 \leq (\sqrt{1+\alpha(n)} - 1)^2 \text{Tr}(\Sigma) \leq \alpha(n)\text{Tr}(\Sigma).$$
The second term above in~\eqref{eq:context-cov-error} can be bounded by the matrix Frobenius norm (for a matrix $A = \{a_{ij}\}_{i \in [d_1], j \in [d_2]}$, $\|A\|_F = \sqrt{\sum_{i = 1}^{d_1} \sum_{j = 1}^{d_2} \|a_{ij}\|^2}$). Denote the sample noise to be $\hat{\eta}_i = (\hat{\eta}_{i,1}, \ldots, \hat{\eta}_{i,D_{\xi}})^{\top}, \forall i \in [n]$. Then with probability at least $1-\delta$, we have:
\begin{align*}
\|\hat{\Sigma} - \hat{\Sigma}^*\|_{F} &= \left\|\frac{1}{n}\sum_{i = 1}^n(\hat{\xi}_i - \hat{f}(\hat{y}_i))(\hat{\xi}_i - \hat{f}(\hat{y}_i))^{\top} - \frac{1}{n}\sum_{i = 1}^n(\hat{\xi}_i - f^*(\hat{y}_i))(\hat{\xi}_i - f^*(\hat{y}_i))^{\top}\right\|_{F}\\
&\leq\|\frac{1}{n}\sum_{i = 1}^n(\hat{f}(\hat{y}_i) - f^*(\hat{y}_i))\hat{\eta}_i^{\top}\|_{F} + \|\frac{1}{n}\sum_{i = 1}^n\hat{\eta}_i(\hat{f}(\hat{y}_i) - f^*(\hat{y}_i))^{\top}\|_{F}\\
&+ \left\|\frac{1}{n}\sum_{i = 1}^n (\hat{f}(\hat{y}_i) - f^*(\hat{y}_i)) (\hat{f}(\hat{y}_i) - f^*(\hat{y}_i))^{\top}\right\|_{F}\\
&\leq 2 C_{\eta} D_{\xi}\sum_{j \in [D_{\xi}} \|\hat{f}_j^* - f_j^*\|_n + \sum_{j \in [D_{\xi}]}\sum_{k \in [D_{\xi}]}\|\hat{f}_j - f_j^*\|_n \|\hat{f}_k - f_k^*\|_n\\
&\leq C_1\sqrt{\frac{\log(D_{\xi}/\delta)}{n}}:=\Delta(n),
\end{align*}
where the first inequality follows from the triangle inequality. And the second inequality follows from the decomposition:
\begin{align*}
\left\|\frac{1}{n}\sum_{i = 1}^n(\hat{f}(\hat{y}_i) - f^*(\hat{y}_i))\hat{\eta}_i^{\top}\right\|_{F} & =  \sqrt{\frac{\sum_{i = 1}^n\sum_{j \in [D_{\xi}]}\sum_{k \in [D_{\xi}]}[\hat{f}_j(\hat{y}_i) - f_j^*(\hat{y}_i)]^2\hat{\eta}_{i,k}^2}{n}}\\
&\leq \sum_{j \in [D_{\xi}]} \sum_{k \in [D_{\xi}]}\sqrt{\frac{\sum_{i = 1}^n [\hat{f}_j(\hat{y}_i) - f_j^*(\hat{y}_i)]^2\hat{\eta}_{i,k}^2}{n}}\\
&\leq \sum_{j \in [D_{\xi}]} \sum_{k \in [D_{\xi}]} \|\hat{f}_j - f_j^*\|_n C_{\eta}.
\end{align*}

Next we present another technical lemma to establish the inequality $(1-\beta(n))\hat{\Sigma}^* \preceq \hat{\Sigma} \preceq (1+\beta(n))\hat{\Sigma}^*$ for some $\beta(n)\to 0$ when $\|\hat{\Sigma} - \hat{\Sigma}^*\|_F \to 0$.
\begin{lemma}\label{lemma:matrix-ineq}
For two positive definite matrices $\Sigma_1, \Sigma_2 \in \R^{d\times d}$, $\|\Sigma_1 - \Sigma_2\|_{F} \leq \Delta(n)$ and $\beta(n) = \frac{\Delta(n)}{\lambda_{\min}(\Sigma_2)}$, then $ (1-\beta(n))\Sigma_2 \preceq \Sigma_1 \preceq (1+\beta(n))\Sigma_2$.
\end{lemma}
\proof{Proof of Lemma~\ref{lemma:matrix-ineq}.}
It is known that $\sqrt{\sum_{j = 1}^d |\lambda_j(\Sigma)|^2} \leq \|\Sigma\|_{F}$, where $\lambda_j(\Sigma)$ is the $j$-th largest eigenvalue of the matrix $\Sigma$. Then $\sqrt{\sum_{j = 1}^d |\lambda_j(\Sigma_1 - \Sigma_2)|^2} \leq \Delta(n)$. 

Consider any normalized vector $x \in \R^d, \|x\|_2 = 1$, we have: $x^{\top}(\Sigma_1 - \Sigma_2) x \leq \Delta(n) = \beta(n) \lambda_{\min}(\Sigma) \leq \beta(n) x^{\top}\Sigma_2 x$, which leads to $x^{\top}\Sigma_1 x \leq (1 + \beta) x^{\top}\Sigma_2 x, \forall x$. The other side follows the same argument if we consider $x^{\top}(\Sigma_2 - \Sigma_1)x$ in the beginning.
$\hfill \square$ 

Finally since $\Sigma_1, \Sigma_2$ are symmetric, $\beta(n) \leq \frac{\Delta(n)}{\lambda_{\min}(\Sigma_2)}.$ 
In terms of the second term in~\eqref{eq:context-cov-error}, we have the inequality $(1-\beta(n))\hat{\Sigma}^* \preceq \hat{\Sigma} \preceq (1+\beta(n))\hat{\Sigma}^*$ with $\beta(n) \leq \frac{\Delta(n)}{\lambda_{\min}(\Sigma_2)}$, which implies that: 
$$\text{Tr}\left[\left(\sqrt{\hat{\Sigma}} - \sqrt{\hat{\Sigma}^*}\right)\right]^2 \leq (\sqrt{1+ \beta(n)} - 1)^2 \text{Tr}(\hat{\Sigma}^*) \leq \alpha(n) (1 + \beta(n)) \text{Tr}(\Sigma).$$

Therefore, when $n$ is large enough, RHS of~\eqref{eq:context-cov-error} is bounded by $(\alpha(n) + \beta(n) + \alpha(n) \beta(n))\text{Tr}(\Sigma) \leq C \text{Tr}(\Sigma)\max\{C_{\eta}^2, D_{\xi}\}\sqrt{\frac{\log(D_{\xi}/\delta)}{n}}$ for some numerical constant $C$.

Combining the results together, we have:
$$W_1(\tilde{\Q}_{\xi|y}, \hat{\Q}_{\xi|y}) \leq \Escr_{apx}(y) + \frac{C_0 \sqrt{\log(1/\delta)}}{\sqrt{n}} + C \text{Tr}(\Sigma)\max\{C_{\eta}^2, D_{\xi}\}\sqrt{\frac{\log(D_{\xi}/\delta)}{n}},$$ 
which implies that Example~\ref{prop:conditional-oracle-estimator} holds. $\hfill \square$
\subsection{Proof of Corollary~\ref{coro:param-dro-erm-context}}
This result follows from the same argument as the proof of Theorem~\ref{thm:general-param-dro-ipm} and Theorem~\ref{thm:general-param-erm}. $\hfill \square$
\subsection{Proof of Example~\ref{prop:expect-conditional}}
Before presenting the proofs, we introduce the following technical assumptions and lemmas:
\begin{assumption}\label{asp:bound-cov}
We assume:
\begin{itemize}
    \item The true generating process of the $i$-th marginal distribution of $\xi$ is $(\xi)_i = (\theta_i^*)^{\top} y + g_i(y) + \eta_i(y)$, where the random $\eta_i(y)$ is sub-Gaussian with parameter $\sigma_{\eta}$. Conditions 1-4 in \cite{hsu2012random} hold for the determinstic approximation error $g_i(y)$ and random noise $\eta_i(y)$.
    \item The covariate $y$ is bounded, i.e. $\|y\|^2 \leq C, \forall y$.
\end{itemize}
\end{assumption}

\begin{lemma}\label{lemma:linear-context}
For some covariate $Y = y$, if $\hat{\Q}_{\xi|y}:=\Nscr(\hat{\theta}^{\top} y, \Sigma)$ with $\hat{\theta} = (\hat{\theta}_1,\ldots,\hat{\theta}_{D_{\xi}})^{\top}$, we have:
$$W_2(\P_{\xi|y}^*, \hat{\Q}_{\xi|y}) \leq \Escr_{apx}(y) + \sum_{i = 1}^{D_{\xi}}\|\hat{\theta}_i - \theta_i^*\|_{\Sigma_y} \|y\|_{\Sigma_y^{-1}},$$
where $\Escr_{apx}(y):=W_2(\P_{\xi|y}^*, \Q_{\xi|y}^*)$, $\Q_{\xi|y}^*:=\Nscr((\theta^*)^{\top}y , \Sigma)$, $\theta^* = (\theta_1^*,\ldots,\theta_{D_{\xi}}^*)^{\top}$ and $\Sigma_y = \E[y y ^{\top}]$.
\end{lemma}

\begin{lemma}[Extracted from Theorem 11 in \cite{hsu2012random}]\label{lemma:linear-estimate-error}
Suppose Assumption~\ref{asp:bound-cov} holds. Then with probability at least $1- 3 e^{-t}$, OLS estimator $\hat{\theta}$ satisfies:
\begin{equation}\label{eq:linear-estimate-error0}
    \|\hat{\theta}_i - \theta_i^*\|_{\Sigma_y}^2 \leq \frac{2\E[\|\Sigma_y^{-\frac{1}{2}}y g_i(y)\|^2]}{n}(1 + \sqrt{8t})^2 + \frac{\sigma_{\eta}^2 (D_{y} + 2\sqrt{D_{y} t} + 2 t)}{n} + \frac{C}{n^2},
\end{equation}
where $C$ is some constant independent with $n$.
\end{lemma}
By some algebraic manipulation, when $n$ is large we simplify~\eqref{eq:linear-estimate-error0} to be: with probability at least $1-\delta$:
\begin{equation}\label{eq:linear-estimate-error}
    \|\hat{\theta}_i - \theta_i^*\|_{\Sigma_y} \leq \sqrt{\frac{C_{M,i} + \sigma_{\eta}^2 D_{y}}{n}} + \sqrt{\frac{\log(1/\delta)}{C_0 n}}.
\end{equation}
where $C_{M,i} = 2\E[\|M^{-\frac{1}{2}}y g_i(y)\|^2]$ and $C_0$ is independent with $n$.
\begin{lemma}[Extracted from Lemma 30 in \cite{hsu2012random}]\label{lemma:quad-subgaussian}
Let $\xi$ be a random vector taking values in $\R^d$ such that for some $c \geq 0$, $\E[\exp(u^{\top} \xi)] \leq \exp(c\|u\|^2 / 2), \forall u \in \R^d$. Then for any symmetric positive semidefinite matrix $K \succeq 0$ and all $t > 0$, 
$$\P\left(\xi^{\top}K \xi > c(\text{tr}(K) + 2\sqrt{\text{tr}(K^2)t} + 2\|K\|t)\right)\leq \exp(-t).$$
\end{lemma}

\proof{Proof of Lemma~\ref{lemma:linear-context}.}
Following from the triangle inequality, we have:
\begin{align*}
    W_2(\P_{\xi|y}^*, \hat{\Q}_{\xi|y})&\leq  W_2(\P_{\xi|y}^*, \Q_{\xi|y}^*) + W_2(\Q_{\xi|y}^*, \hat{\Q}_{\xi|y}^*)\\
    &= \Escr_{apx}(y) + \|(\hat{\theta} - \theta^*)^{\top}y\|_2\\
    &= \Escr_{apx}(y) + \sqrt{\sum_{i = 1}^{D_{\xi}}\|(\hat{\theta}_i - \theta_i^*)^{\top}y\|_2}\\
    &\leq \Escr_{apx}(y) + \sqrt{\sum_{i = 1}^{D_{\xi}}\|\hat{\theta}_i - \theta_i^*\|_{\Sigma_y}^2}\|y\|_{\Sigma_y^{-1}},\\
    &\leq \Escr_{apx}(y) + \sum_{i = 1}^{D_{\xi}}\|\hat{\theta}_i - \theta_i^*\|_{\Sigma_y}\|y\|_{\Sigma_y^{-1}}
\end{align*}
where the second inequality follows from the Cauchy-Schwarz inequality and the third inequality follows from the fact that $\sqrt{\sum_{i \in [n]}x_i} \leq \sum_{i\in [n]}\sqrt{x_i}$. $\hfill \square$

\proof{Proof of Example~\ref{prop:expect-conditional}.} 

We denote $\hat{\E}[\|\hat{\theta}_i - \theta_i^*\|_{\Sigma_y}] = \sqrt{\frac{C_{M, i} + \sigma_{\eta}^2 D_{y}}{n}}, \forall I \in [D_{\xi}]$, then from~\eqref{eq:linear-estimate-error}, we have:
\begin{equation*}
    \P(\|\hat{\theta}_i - \theta_i^*\|_{\Sigma_y}- \hat{\E}[\|\hat{\theta}_i - \theta_i^*\|_{\Sigma_y}] > t)\leq \exp(-C_0 nt^2), \forall i \in [D_{\xi}].
\end{equation*}
This implies that:
\begin{equation}\label{eq:concentration-linear-param}
    \P\left(\sum_{i = 1}^{D_{\xi}}\|\hat{\theta}_i - \theta_i^*\|_{\Sigma_y}- \sum_{i = 1}^{D_{\xi}}\hat{\E}[\|\hat{\theta}_i - \theta_i^*\|_{\Sigma_y}] > t\right)\leq D_{\xi} \exp(-C_0 nt^2), \forall i \in [D_{\xi}].
\end{equation}
Denoting $\hat{d}:= \sum_{i = 1}^{D_{\xi}}\hat{\E}[\|\hat{\theta}_i - \theta_i^*\|_{\Sigma_y}]\E[\|y\|_{\Sigma_y^{-1}}] + \sup_{y \in \Yscr}\Escr_{apx}(y)$, we need to show \eqref{eq:context-joint-ineq} holds here, with $\text{Comp}(\Theta, y) = \hat{\E}[\|\hat{\theta}_i - \theta_i^*\|_{\Sigma_y}] \|y\|_{\Sigma_y^{-1}}  = O(D_{\xi}D_y \|y\|_{\Sigma_y^{-1}})$ and $\text{Comp}(\Theta) = \sum_{i = 1}^{D_{\xi}}\hat{\E}[\|\hat{\theta}_i - \theta_i^*\|_{\Sigma_y}]\E[\|y\|_{\Sigma_y^{-1}}] = O\left(\frac{D_{\xi}D_y}{\lambda_{\min}(\Sigma_y)}\right)$ corresponding to the term in Assumption~\ref{asp:expect-conditional-oracle-estimator}. 



Since $0 \leq \|y\|_{\Sigma_y^{-1}} \leq  \sqrt{C\lambda_{\min}(\Sigma_y)}$ by Assumption~\ref{asp:bound-cov}, applying Hoeffding inequality, then we obtain:
\begin{equation}\label{eq:hoeffding-covariate}
  \P(\|y\|_{\Sigma_y^{-1}}  - \E[\|y\|_{\Sigma_y^{-1}}]\geq t) \leq \exp\left(-\frac{2t^2}{C \lambda_{\min}(\Sigma_y)}\right).  
\end{equation}



Combining all the previous arguments, we have:
\small{
\begin{align*}
&\quad \P(W_2(\P^*_{\xi|y}, \hat{\Q}_{\xi|y}) - \hat{d} > 2t)\\
&\leq \P\left(\sum_{i = 1}^{D_{\xi}} \|\hat{\theta}_i - \theta_i^*\|_{\Sigma_y}\|y\|_{\Sigma_y^{-1}} - \sum_{i = 1}^{D_{\xi}}\hat{\E}[\|\hat{\theta}_i - \theta_i^*\|_{\Sigma_y}]\E[\|y\|_{\Sigma_y^{-1}}] > 2t\right) \\
&=\P\left(\sum_{i = 1}^{D_{\xi}} \|\hat{\theta}_i - \theta_i^*\|_{\Sigma_y}\|y\|_{\Sigma_y^{-1}} - \sum_{i = 1}^{D_{\xi}}\hat{\E}[\|\hat{\theta}_i - \theta_i^*\|_{\Sigma_y}]\|y\|_{\Sigma_y^{-1}} \right.\\
&\left.+ \sum_{i = 1}^{D_{\xi}}\hat{\E}[\|\hat{\theta}_i - \theta_i^*\|_{\Sigma_y}]\|y\|_{\Sigma_y^{-1}} - \sum_{i = 1}^{D_{\xi}}\hat{\E}[\|\hat{\theta}_i - \theta_i^*\|_{\Sigma_y}]\E[\|y\|_{\Sigma_y^{-1}}]\geq 2t\right).\\
&\leq \P\left(\sum_{i = 1}^{D_{\xi}}\|\hat{\theta}_i - \theta_i^*\|_{\Sigma_y}-\sum_{i = 1}^{D_{\xi}} \hat{\E}[\|\hat{\theta}_i - \theta_i^*\|_{\Sigma_y}] > \frac{t}{\max_{y \in \Yscr}\|y\|_{\Sigma_y^{-1}}}\right) + \P\left(\|y\|_{\Sigma_y^{-1}} - \E_y[\|y\|_{\Sigma_y^{-1}}] \geq \frac{t}{\sum_{i = 1}^{D_{\xi}}\hat{\E}[\|\hat{\theta}_i - \theta_i^*\|_{\Sigma_y}]}\right)\\
& \leq D_{\xi}\exp\left(-C_0 n \frac{t^2}{C \lambda_{\min}(\Sigma_y)}\right) + \exp\left(-\frac{2 n t^2}{C \lambda_{\min}(\Sigma_y)D_{\xi}^2(\max_{i \in [D_{\xi}]} C_{M,i} + \sigma_{\eta}^2 D_y)}\right) \leq c_1 \exp(-c_2 nt^2), 
\end{align*}}
\normalsize
where the first inequality follows from the definition of $\hat{d}$ and Lemma~\ref{lemma:linear-estimate-error}. The second inequality follows from partitioning the inequality into two parts. And the third inequality follows from the two previous concentration inequalities~\eqref{eq:concentration-linear-param} and~\eqref{eq:hoeffding-covariate}. Then we obtain the final inequality with $c_1, c_2$ independent with $n$. Since $W_1(\P, \Q) \leq W_2(\P, \Q)$, we can replace the results obtained here for $W_2$ with $W_1$.
$\hfill \square$
\subsection{Proof of Theorem~\ref{thm:expect-generalization-bd}}
Denote $Z(\hat{x}(y)) = \E_{\P_{\xi|y}^*}[h(\hat{x}(y);\xi)], Z(x^*(y)) = \E_{\P_{\xi|y}^*}[h(x^*(y);\xi)]$. Then we decompose the error as follows:
\begin{equation}\label{eq:context-decomp}
\begin{aligned}
   Z(\hat{x}(y)) - Z(x^*(y)) &= [Z(\hat{x}(y)) - \hat{Z}(\hat{x}(y))] + [\hat{Z}(\hat{x}(y)) - \hat{Z}(x^*(y))] + [\hat{Z}(x^*(y)) - Z(x^*(y))]\\
   & \leq [Z(\hat{x}(y)) - \hat{Z}(\hat{x}(y))] + [\hat{Z}(x^*(y)) - Z(x^*(y))].
\end{aligned}
\end{equation}

\paragraph{P-DRO.}In~\eqref{eq:context-decomp}, we set $\hat{Z}(\cdot) = \sup_{d(\P_{\xi|y}, \hat{\Q}_{\xi|y}) \leq \varepsilon}\E_{\P^*_{\xi|y}}[h(\cdot;\xi)]$. We partition the probability space $\Pscr = \P_{\Dscr^n}\otimes \P_y$ into the sets $\Ascr_1 = \{\{\Dscr_n, y): d(\P^*_{\xi | y}, \hat{\Q}_{\xi | y}) \leq \varepsilon\} \subseteq \P_{\Dscr^n}\otimes \P_y$ and $\Ascr_2 = ((\prod_{i = 1}^n \P_{(y,\xi)})\otimes \P_y)\backslash \Ascr_1$. Then we decompose the generalization error into two regions:
\begin{equation}\label{eq:p-dro-context-derive}
\begin{aligned}
\Escr_{\Yscr}(\hat{x}) &=\E_{\Dscr_n}\E_y\left(Z(\hat{x}(y)) - Z(x^*(y))\right) \\
&=\E_{\Dscr_n}\E_y\left([Z(\hat{x}(y)) - Z(x^*(y))]\mathbb{I}_{\{\Ascr_1\}} + [Z(\hat{x}(y)) - Z(x^*(y))]\mathbb{I}_{\{\Ascr_2\}}\right)\\
&\leq \E_{\Dscr_n}\E_{y}[2\varepsilon\Vscr_d(x^*(y))\mathbb{I}_{\{\Ascr_1\}}] + M \cdot \E[d(\P^*, \hat{\Q})\mathbb{I}_{\{\Ascr_2\}}] + \E_{\{(\hat{y}_i,\hat{\xi}_i)\}_{i = 1}^n}\E_{y}[2\varepsilon\Vscr_d(x^*(y))\mathbb{I}_{\{\Ascr_2\}}]\\
&= 2\varepsilon \E_y[\Vscr_d(x^*(y))] + M \cdot \E[d(\P^*, \hat{\Q})\mathbb{I}_{\{\Ascr_2\}}].
\end{aligned}
\end{equation}
where the first inequality follows from the error decomposition in~\eqref{eq:context-decomp}. The second term of~\eqref{eq:context-decomp} is non-positive since $\hat{x}(y)$ is the minimizer of $\hat{Z}(\cdot)$. We consider the error decomposition in the following two scenarios:

\textit{(1) Event $\Ascr_1$ holds.} The first term of RHS of~\eqref{eq:context-decomp} is non-positive from the definition of $\Ascr_1$. The third term follows from the same argument as in Theorem~\ref{thm:general-param-dro-ipm}, i.e.:
$\Vscr_d(x^*(y))\max_{d(\P_{\xi|y}, \hat{\Q}_{\xi|y})\leq \varepsilon}d(\P_{\xi|y}, \P^*_{\xi|y})\leq 2 \Vscr_d(x^*(y))\varepsilon.$

\textit{(2) Event $\Ascr_2$ holds.} The first term of RHS of~\eqref{eq:context-decomp} is bounded by $\E_{\P_{\xi|y}^*}[h(\hat{x}(y);\xi)] - \E_{\hat{\Q}_{\xi|y}}[h(\hat{x}(y);\xi)]\leq M d(\P^*_{\xi|y}, \hat{\Q}_{\xi|y})$. And the second term of RHS of~\eqref{eq:context-decomp} is bounded by $\E_{\hat{\Q}_{\xi|y}}[h(x^*(y);\xi)] + \varepsilon\Vscr_d(x^*(y)) - \E_{\P^*_{\xi|y}}[h(x^*(y);\xi)] \leq 2\varepsilon \Vscr_d(x^*(y))$ following from the definition of IPM and $\hat{Z}$.

Combining the two scenarios, we have \eqref{eq:p-dro-context-derive}. Denote $t:= \varepsilon - \hat{d}$ such that $\P(\Ascr_2) \leq c_1 \exp(-c_2 a_n t^2)$. and $\Delta:= d(\P_{\xi|y}^*, \hat{\Q}_{\xi|y}) - \hat{d}$ to be the random quantity deviating from $\hat{d}$. Therefore $\Ascr_2 = \{\Delta \geq t\}$. Expanding the second term in~\eqref{eq:p-dro-context-derive}, we have:
\begin{align*}
    \E_{\Dscr_n}\E_{y}[d(\P^*, \hat{\Q})\mathbb{I}_{\{\Ascr_2\}}] &= \E[(\Delta - t)\mathbb{I}_{\{\Delta \geq t\}}] + (t + \hat{d}) \P(\Delta \geq t)\\
    & = \int_0^{\infty} \P(\Delta \geq u + t)du + (t + \hat{d}) \P(\Delta \geq t)\\
    &= \int_0^{\infty}c_1 \exp(-c_2 a_n (u + t)^2)du + (t + \hat{d}) c_1 \exp(-c_2 a_n t^2)\\
    & \leq c_1 \exp(-c_2 a_n t^2)\left(t + \hat{d} + \frac{1}{\sqrt{c_2 a_n}}\right), 
\end{align*}
where the inequality above follows from $\exp(-(a+b)^2) \leq \exp(-a^2) \cdot \exp(-b^2)$ for $a, b \geq 0$ and $\int_0^{\infty}\exp(-ct^2) dt = \frac{1}{\sqrt{c}}$ for $c > 0$.


Plugging this $t$ into the generalization bound~\eqref{eq:p-dro-context-derive}, we have:
\begin{align*}
\Escr_{\Yscr}(\hat{x}) &\leq 2[\hat{d} + t]\E_y[\Vscr_d(x^*(y))] + M c_1 \exp(-c_2 a_n t^2)\left(t + \hat{d} + \frac{1}{\sqrt{c_2 a_n}}\right). 
\end{align*}
\paragraph{P-ERM.}
In~\eqref{eq:context-decomp}, we set $\hat{Z}(\cdot)=\E_{\hat{\Q}_{\xi|y}}[h(\cdot;\xi)]$. Then $Z(\hat{x}(y)) - Z(x^*(y))\leq 2M d(\P^*_{\xi|y}, \hat{\Q}_{\xi|y})$ following from the same analysis as in~\eqref{thm:general-param-erm}. Then the generalization error of \PERM\ is upper bounded by:
\begin{align*}
\Escr_{\Yscr}(\hat{x}) &\leq  2M \E_{\Dscr^n}\E_{y}[d(\P, \hat{\Q})]\\
&\leq 2M \int_{0}^{\infty}\P(d(\P, \hat{\Q}) \geq t) dt\\
&\leq 2M\left[\hat{d} + \int_{0}^{\infty} c_1 \exp(-c_2 a_n t^2)d t\right]\leq 2M \left(\hat{d} + \frac{c_1}{\sqrt{c_2 a_n}}\right).
\end{align*} $\hfill \square$
\section{Additional Experimental Details and Results}\label{app:experiment}
The optimization problems throughout this paper are solved in CVXPY and Gurobi implemented by Python 3.8.5. The computational environment is an Intel(R) Core(TM) i7-8650U CPU @1.90GHz personal computer. Specifically in this part, we denote $\xi_i, x_i$ (or $(\xi)_i$) as the $i$-th marginal component of the random variable and $\hat{\xi}^j, \hat{x}^j$ as the $j$-th sample in the dataset.

\subsection{Detailed Setups and Analysis in Section~\ref{subsec:synthetic}}\label{subsec:synthetic2}
In the base case, the unknown distribution $\P^*$ is fully parametrized such that $(\xi)_i:=2r \times Beta(\eta_i, 2) - r$ with $\{\eta_i\}_{\{i \in [D_{\xi}]\}}$ i.i.d. drawn from $[1.5, 3]$ for the $i$-th marginal distribution. First, we estimate $\text{Comp}(\Hscr)$ and $\text{Comp}(\Theta)$ here. 
\subsubsection{Comparison between $\text{Comp}(\Hscr)$ and $\text{Comp}(\Theta)$.}\label{subsec:comp-theta-h}
First, we present an upper bound for the covering number of $\Hscr$.
\begin{lemma}[Extracted from Theorem 5.4 in~\cite{matousek1999geometric}]\label{lemma:vc-poly}
If $\Hscr$ consists of polynomials up to degree $D$ with $d$ variables, i.e., each $h(\xi) \in \Hscr, \xi = (\xi_1,\ldots, \xi_d)^{\top} \in \R^d$ can be represented as $h(\xi) = \sum_{i_1 + \ldots + i_d \leq D}a_i \xi_1^{i_1}\ldots\xi_d^{i_d}$, then we have:
$$\text{VC}(\Hscr)\leq \tbinom{d+ D}{d} \sim (d+ D)^{\min\{d, D\}}.$$
\end{lemma}

And following from Theorem 2.6.7 in \cite{van1996weak}, we have:
$$N(\Hscr_n(\bm \xi), \varepsilon, \|\cdot\|_{\infty}) \leq \sup_{\Q} N\left(\Hscr, \frac{\varepsilon}{2n}, \|\cdot\|_{L^1(\Q)}\right) \leq c \text{VC}(\Hscr)\left(\frac{16M ne}{\varepsilon}\right)^{\text{VC}(\Hscr) - 1},$$
for some numerical constant $c$.

Then combining it with Lemma~\ref{lemma:vc-poly}, we have the following upper bound for the covering number of $\Hscr$ in \eqref{eq:high-order-simulation}:
\begin{equation}\label{eq:covering-number2}
\begin{aligned}
N_{\infty}(\Hscr, \varepsilon, n) &\leq C (D_{\xi}+\gamma)^{\min\{D_{\xi},\gamma\}}\left(\frac{n M}{\varepsilon}\right)^{(D_{\xi}+\gamma)^{\min\{D_{\xi},\gamma\}}}\\
& = C (D_{\xi}+\gamma)^{\gamma}\left(\frac{n M}{\varepsilon}\right)^{(D_{\xi}+\gamma)^{\gamma}},
\end{aligned}
\end{equation}
where $M:=\sup_{x \in \Xscr}\|h(x;\cdot)\|_{\infty} \leq (D_{\xi} \tau r + \mu)^{\gamma}$. In these generalization error bounds, we approximate the variance term in $\chi^2$-divergence by $\text{Var}_{\P^*}[h(x^*;\cdot)] \leq \|h(x^*;\cdot)\|_{\infty}^2$. Setting $\varepsilon = \frac{1}{n}$ in~\eqref{eq:covering-number2}, and ignoring the term $\log(1/\delta)$ and $\log n$, $\text{Comp}(\Hscr) = O(D_{\xi}^{\gamma} \log M)$ in this setup.


In the base case, $\P^*$ is set to be a variant of the Beta distribution. 
We obtain the following results for $\text{Comp}(\Theta)$ as follows:
\begin{proposition}\label{prop:chi2-oracle-estimator}
When $d$ is $\chi^2$-divergence, if $\P^* \in \Pscr_{\Theta}(:= \{\P: \xi
= (\xi_1,\ldots,\xi_{D_{\xi}})^{\top}\sim \P, (\xi)_i:=
2r\times Beta(\eta_i, 2) - r, \eta_i \in [k, 2k],\forall i \in
[D_{\xi}]\})$ for some constant $k$, then Assumption~\ref{asp:oracle-param-est} holds for $\hat{\Q}$ with the $i$-th marginal distribution being $2r \times Beta(\hat{\eta}_i, 2) - r$, where $\hat{\eta}_i$ is computed from the moment method below and $\mathcal{E}_{apx} = 0, \text{Comp}(\Theta) = CD_{\xi}, \alpha = 1$ with some constant $C$ when $n$ is large.
\end{proposition}

Finally, combing the results of $\text{Comp}(\Hscr)$ and $\text{Comp}(\Theta)$ above, we show the main terms of generalization errors in the four methods in Table~\ref{tab:general-synthetic-bound}.


\textit{Proof of Proposition~\ref{prop:chi2-oracle-estimator}.}
The formula of $\chi^2$-divergence under Beta distribution is given as follows:

\begin{example}
    For $\P_1 \sim Beta(\eta_1, \beta_1), \P_2 \sim Beta(\eta_2, \beta_2)$, and $\eta_1, \eta_2, \beta_1, \beta_2 > 0$, we have:
    $$\chi^2(\P_1, \P_2) = \frac{B(\eta_1, \beta_1)B(2\eta_2 - \eta_1, 2\beta_2 - \beta_1)}{B(\eta_2, \beta_2)} - 1,$$
    where $B(\eta, \beta) = \int_0^1 x^{\eta - 1}(1-x)^{\beta - 1}dx.$ 
    If the true distribution $\xi \sim \prod_{i = 1}^{D_{\xi}} \P_i^*(:=Beta(\eta_i, \beta_i))$ and the estimated distribution $\hat{\xi}\sim \prod_{i = 1}^{D_{\xi}}\Q_i^*(:=Beta(\hat{\eta}_i, \hat{\beta}_i))$, then by the product rule:
    $$\chi^2(\prod_{i = 1}^{D_{\xi}} \P_i^*, \prod_{i = 1}^{D_{\xi}}\hat{\Q}_i) = \prod_{i = 1}^{D_{\xi}}\frac{B(\eta_i, \beta_i)B(2\hat{\eta}_i - \eta_i, 2\hat{\beta}_i - \beta_i)}{(B(\hat{\eta}_i, \hat{\beta}_i))^2} - 1.$$  
\end{example}
Changing the support of Beta distributions to $\Pscr_{\Theta}$ does not change the value of the $\chi^2$-divergence. 
Since $\beta_i = \hat{\beta}_i = 2$ in our problem setup, the divergence is:
\begin{equation}\label{eq:beta-simple}
 \chi^2(\prod_{i = 1}^{D_{\xi}} \P_i^*, \prod_{i = 1}^{D_{\xi}}\hat{\Q}_i) = \prod_{i = 1}^{D_{\xi}}\frac{\hat{\eta}_i}{\eta_i}\cdot\frac{\hat{\eta}_i + 1}{\eta_i + 1}\cdot\frac{\hat{\eta}_i}{2\hat{\eta}_i - \eta_i}\cdot\frac{\hat{\eta}_i + 1}{2\hat{\eta}_i - \eta_i + 1}- 1.   
\end{equation}
If we can control the estimation error such that $\forall i \in [D_{\xi}]$, with probability at least $1-\delta$, we have:
\begin{equation}\label{eq:beta-simple-reduce1}
    1 - u\sqrt{\Delta} \leq \frac{\hat{\eta_i}}{\eta_i} \leq 1 + u \sqrt{\Delta},
\end{equation}
\begin{equation}\label{eq:beta-simple-reduce2}
    1 - v\sqrt{\Delta} \leq \frac{\hat{\eta}_i + 1}{\eta_i + 1} \leq 1 + v \sqrt{\Delta},
\end{equation}
where $\Delta: = \frac{\log(1/\delta)}{n}$ in~\eqref{eq:beta-simple-reduce1} and~\eqref{eq:beta-simple-reduce2} and $u, v$ are independent with the sample size $n$. Then based on the formula in~\eqref{eq:beta-simple}, following from~\eqref{eq:beta-simple-reduce1} and~\eqref{eq:beta-simple-reduce2}, we have:
\begin{align*}
\chi^2\left(\prod_{i = 1}^{D_{\xi}} \P_i^*, \prod_{i = 1}^{D_{\xi}}\hat{\Q}_i\right) &\leq \left((1+u^2\Delta)(1+v^2\Delta)\right)^{D_{\xi}} -1   \\
&= \left[1 + 2(u^2+v^2)\Delta + o(\Delta)\right]^{D_{\xi}} - 1\\
&\leq 4D_{\xi}(u^2+v^2)\Delta + o(\Delta),
\end{align*}
where the first inequality follows from $\frac{\hat{\eta}_i}{\eta_i} \cdot\frac{\hat{\eta}_i}{2\hat{\eta}_i - \eta_i} = \frac{(\hat{\eta}_i/\eta_i)^2}{2(\hat{\eta}_i/\eta_i) - 1} \leq 1 + \frac{(u\sqrt{\Delta})^2}{1+2u\sqrt{\Delta}} \leq 1+ (u\sqrt{\Delta})^2$ (as long as $u\sqrt{\Delta} \leq \frac{1}{2}$). And the second inequality holds when $n$ is large. Thus we have that $\chi^2(\prod_{i = 1}^{D_{\xi}} \P_i^*, \prod_{i = 1}^{D_{\xi}}\hat{\Q}_i) = O\left(\frac{D_{\xi}}{n}\right)$. In the following, we argue that the moment method leads to a controllable estimation error in~\eqref{eq:beta-simple-reduce1} and~\eqref{eq:beta-simple-reduce2}.

In the moment method, we estimate each $\hat{\eta}_i = \frac{2}{\frac{1}{2}- \frac{\sum_{j = 1}^n \hat{\xi}_i^j}{2nr}} - 2=:\frac{2\hat{\E}[\gamma_i]}{1-\hat{\E}[\gamma_i]}$ through the first-order moment equation $\E[\frac{\xi_i}{2r} + \frac{1}{2}] = \frac{\eta_i}{\eta_i + 2}$, where $\hat{\E}[{\xi}_i] := \frac{1}{n}\sum_{j = 1}^n \hat{\xi}_i^j$ is the empirical average of the $i$-th marginal distribution of $\{\hat{\xi}^j\}_{j \in [n]}$, $\hat{\E}[\gamma_i] := \frac{\hat{\E}[{\xi}_i] + r}{2r}$ and $ \gamma_i := \frac{\xi_i + r}{2r}$. Then we have:
$$\frac{\hat{\eta}_i}{\eta_i}  = \frac{\hat{\E}[\gamma_i]}{\E[\gamma_i]}\cdot \frac{1-\E[\gamma_i]}{1-\hat{\E}[\gamma_i]}.$$
Then we apply the Hoeffding-type concentration argument to bound $|\hat{\E}[\gamma_i] - \E[\gamma_i]|$, which is the same case as for $\frac{\hat{\eta}_i + 1}{\eta_i + 1}$.
$\hfill \square$

\subsubsection{Detailed Experimental Results}
Besides the Beta model class, we fit the model with a normal class $\Pscr_{\Theta}$, i.e. the parametric model in Section~\ref{subsec:synthetic1}.  Intuitively, for $\chi^2$-divergence, $\text{Comp}(\Theta) = C D_{\xi}^2$ with $\mathcal{E}_{apx} > 0$ when $\Pscr_{\Theta}$ is taken as the normal class. Even with distribution misspecification, the Normal-DRO model still outperforms the nonparametric counterpart. 

\begin{figure*}[h] 
    \centering    
    \subfloat[$(\gamma, n) = (1, 25)$] 
    {
        \begin{minipage}[t]{0.25\textwidth}
            \centering          
            \includegraphics[width=0.9\textwidth, trim = 20 10 45 45]{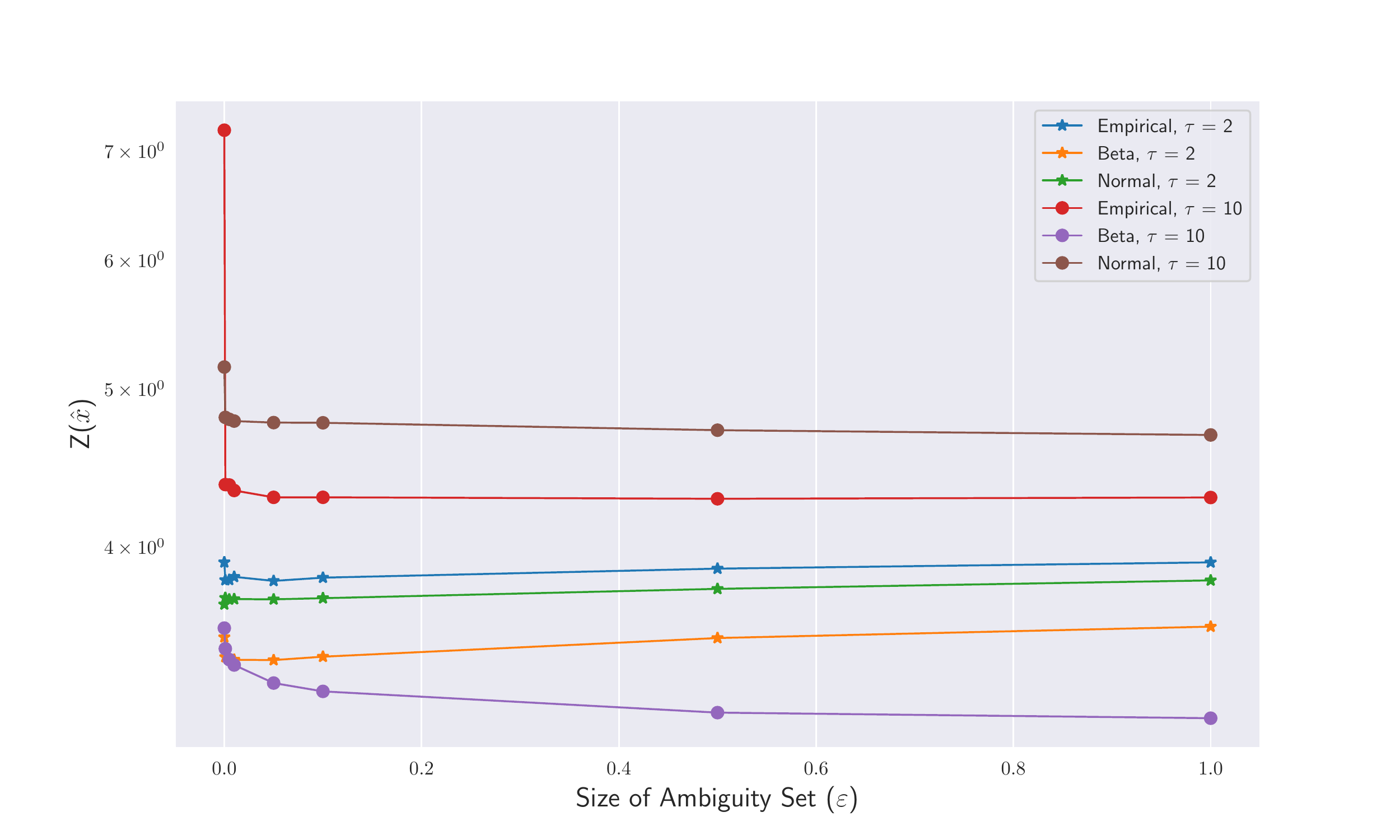}  
        \end{minipage}
    }
    \subfloat[$(\gamma, n) = (1, 50)$] 
    {
        \begin{minipage}[t]{0.25\textwidth}
            \centering  
            \includegraphics[width=0.9\textwidth, trim = 20 10 45 45]{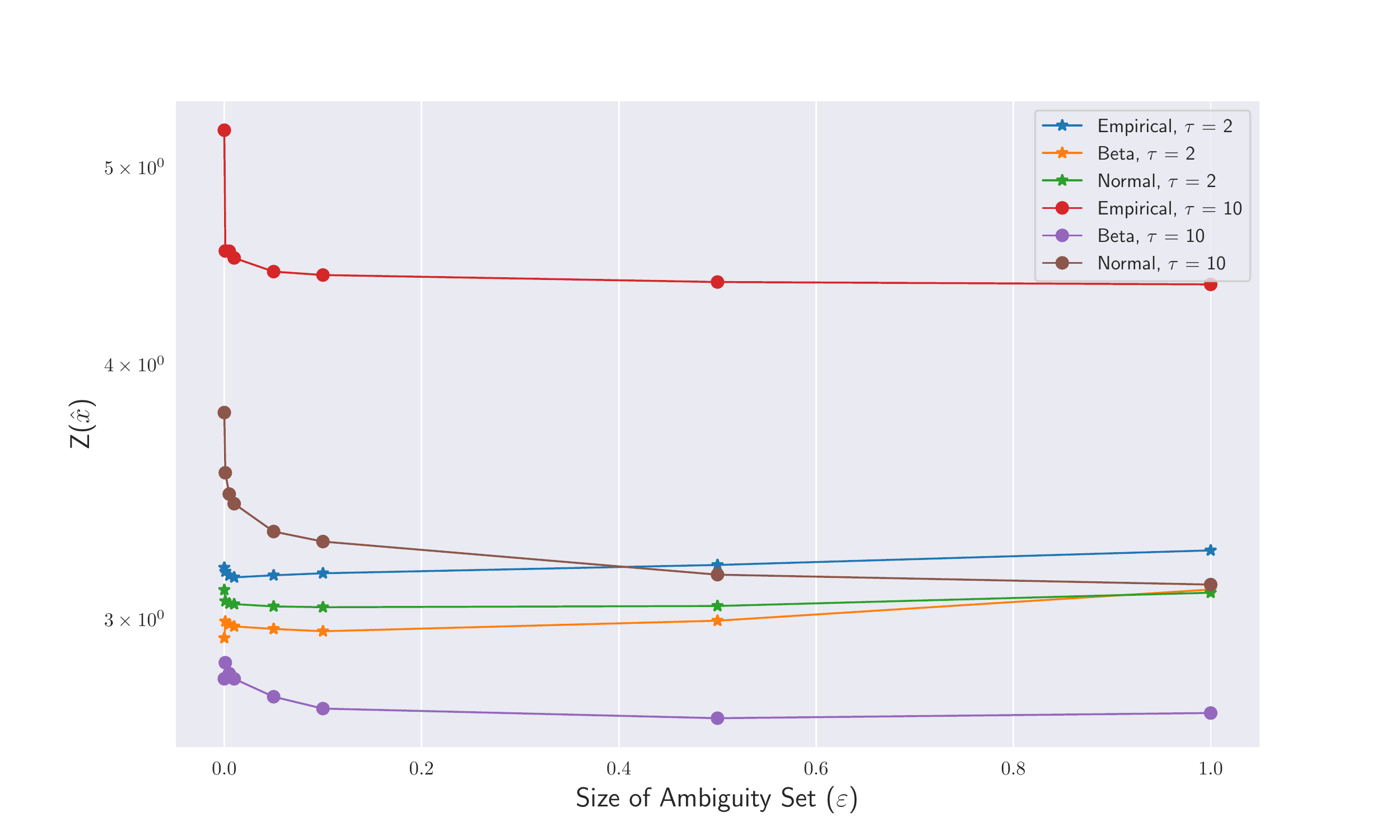} 
        \end{minipage}
    }%
    \subfloat[$(\gamma, n) = (1, 100)$] 
    {
        \begin{minipage}[t]{0.25\textwidth}
            \centering          
            \includegraphics[width = 0.9\textwidth, trim = 20 10 45 25 ]{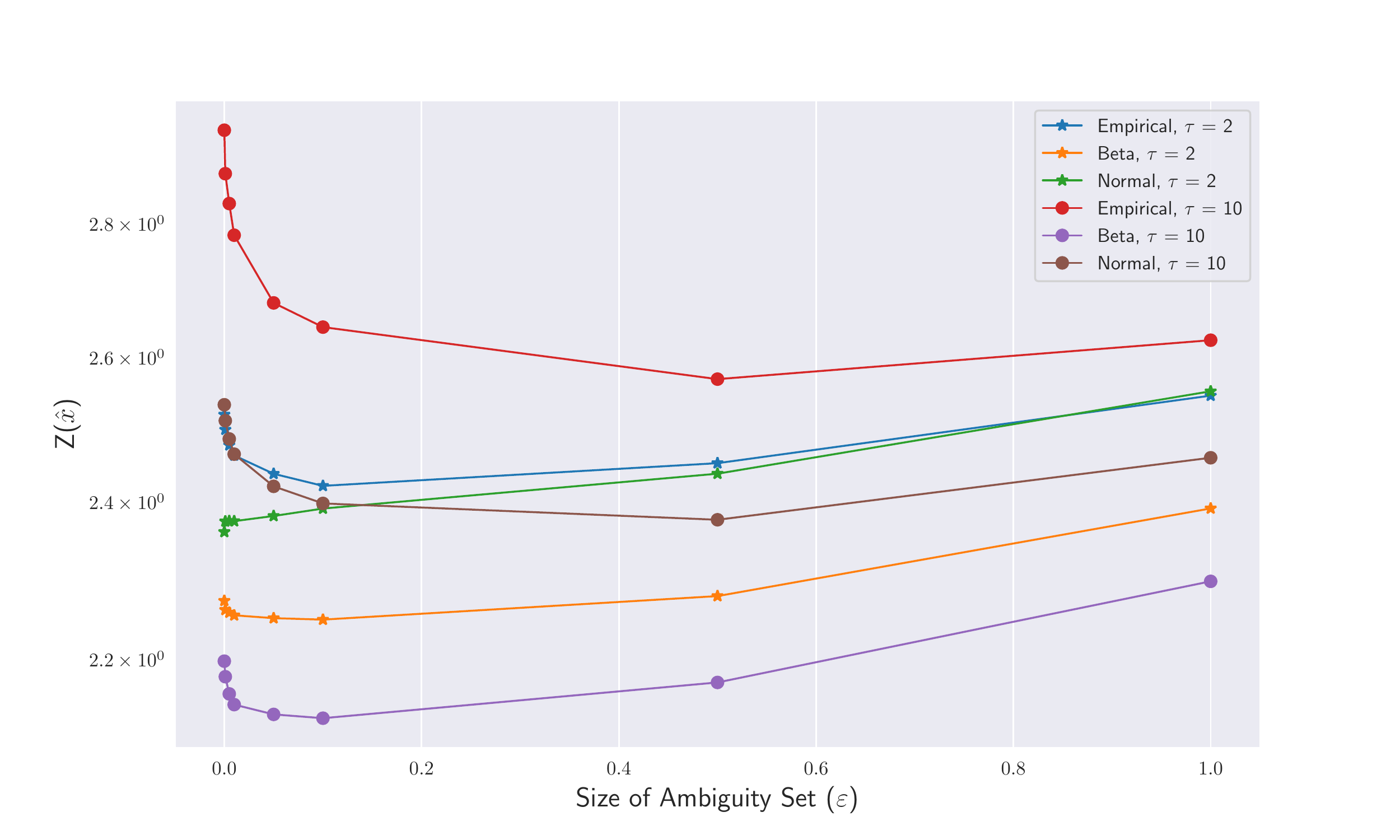} 
        \end{minipage}%
    }
    \subfloat[$(\gamma, n) = (1, 200)$]
    {
        \begin{minipage}[t]{0.25\textwidth}
            \centering      
            \includegraphics[width = 0.9\textwidth, trim = 20 10 45 25]{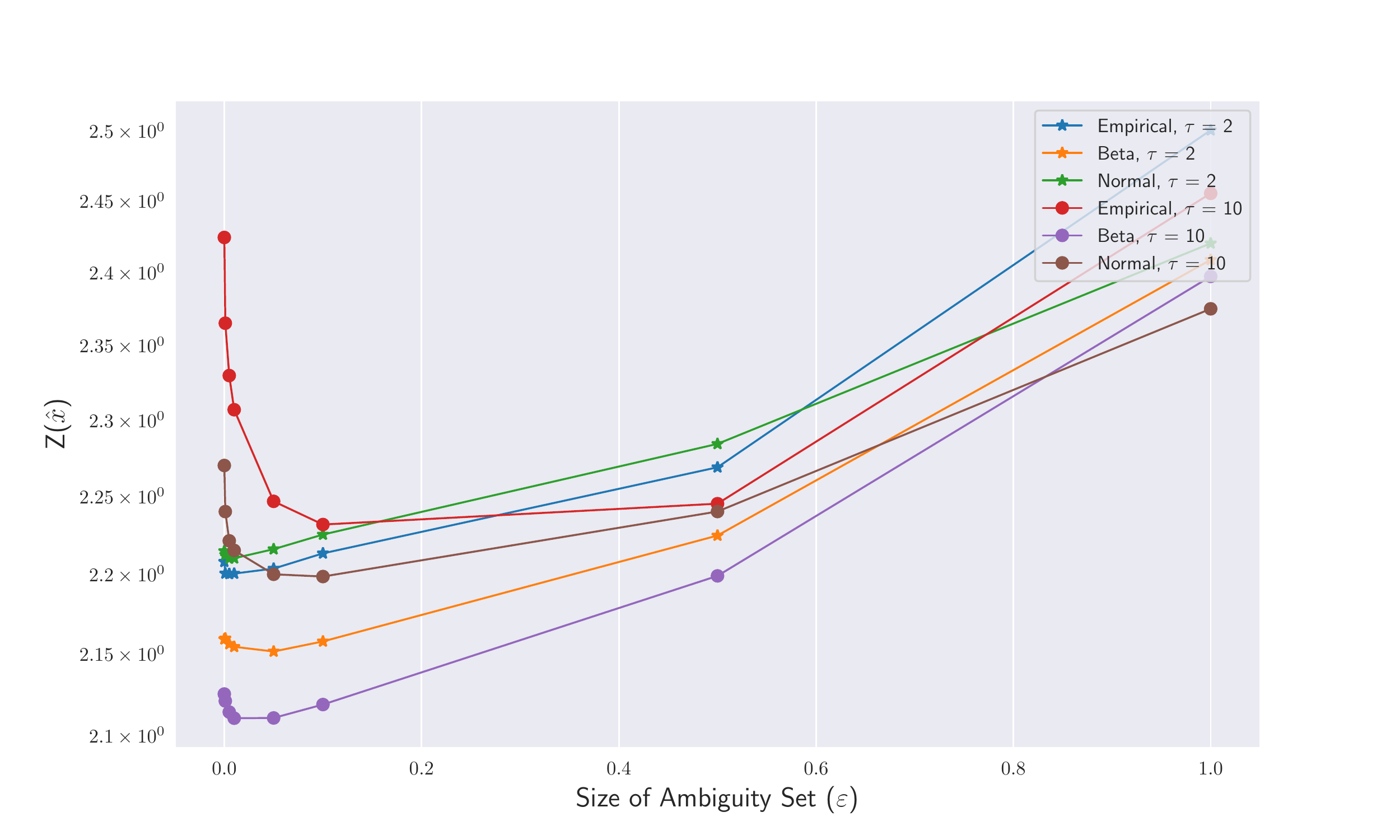}  
        \end{minipage}
    }%
    
    \subfloat[$(\gamma, n) = (2, 25)$] 
    {
        \begin{minipage}[t]{0.25\textwidth}
            \centering          
            \includegraphics[width=0.9\textwidth, trim = 20 10 45 45]{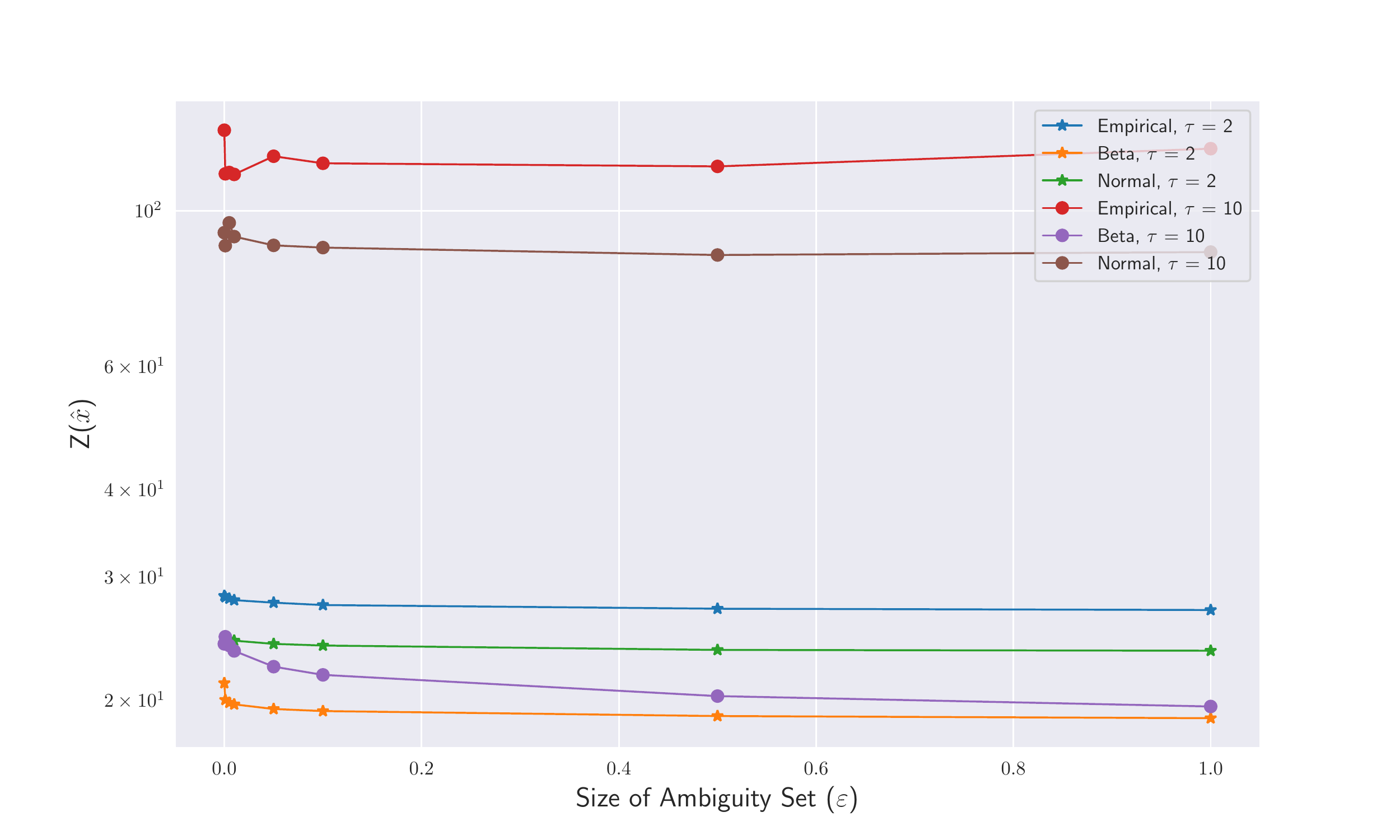}  
        \end{minipage}
    }
    \subfloat[$(\gamma, n) = (2, 50)$] 
    {
        \begin{minipage}[t]{0.25\textwidth}
            \centering  
            \includegraphics[width=0.9\textwidth, trim = 20 10 45 45]{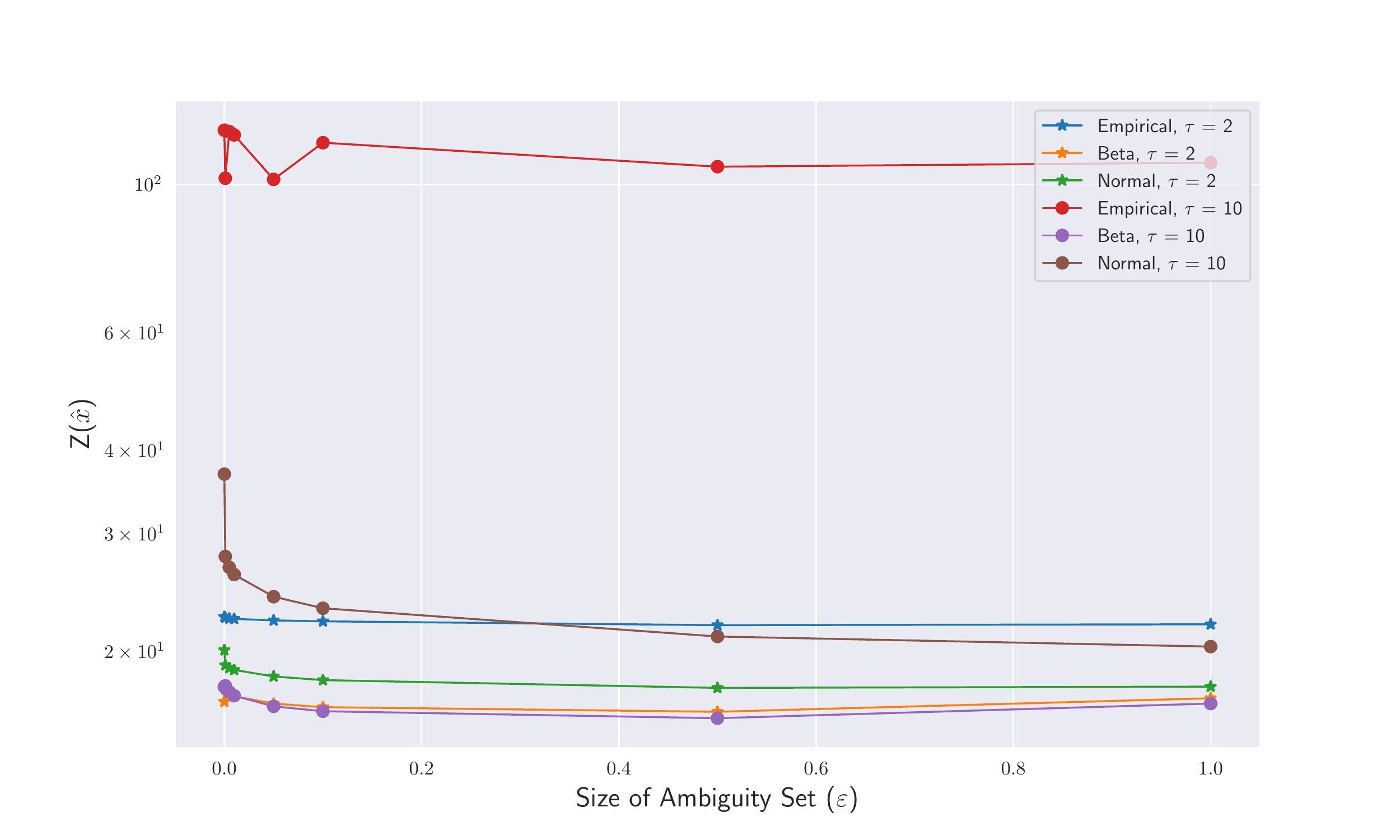} 
        \end{minipage}
    }%
    \subfloat[$(\gamma, n) = (2, 100)$] 
    {
        \begin{minipage}[t]{0.25\textwidth}
            \centering          
            \includegraphics[width = 0.9\textwidth, trim = 20 10 45 25 ]{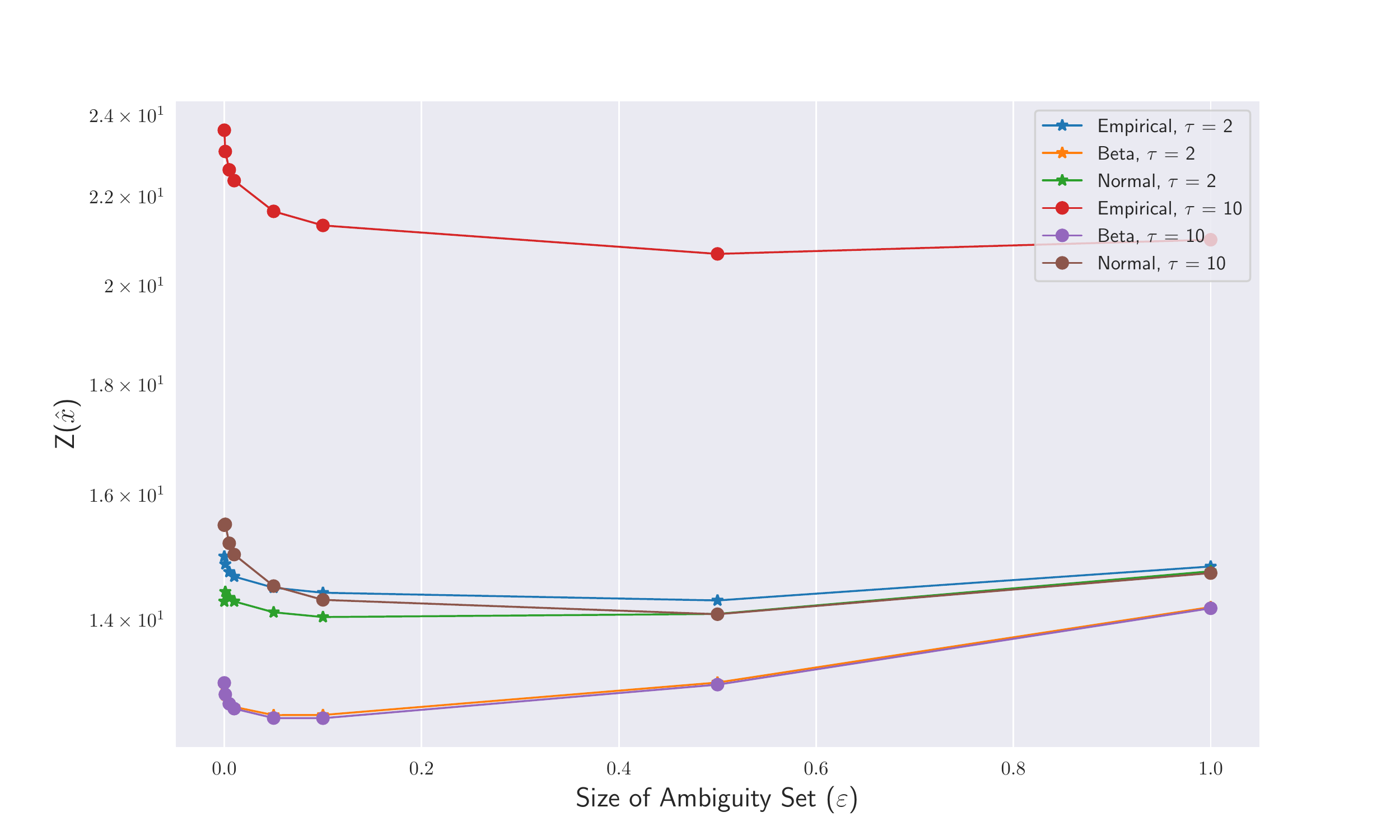} 
        \end{minipage}%
    }
    \subfloat[$(\gamma, n) = (2, 200)$]
    {
        \begin{minipage}[t]{0.25\textwidth}
            \centering      
            \includegraphics[width = 0.9\textwidth, trim = 20 10 45 25]{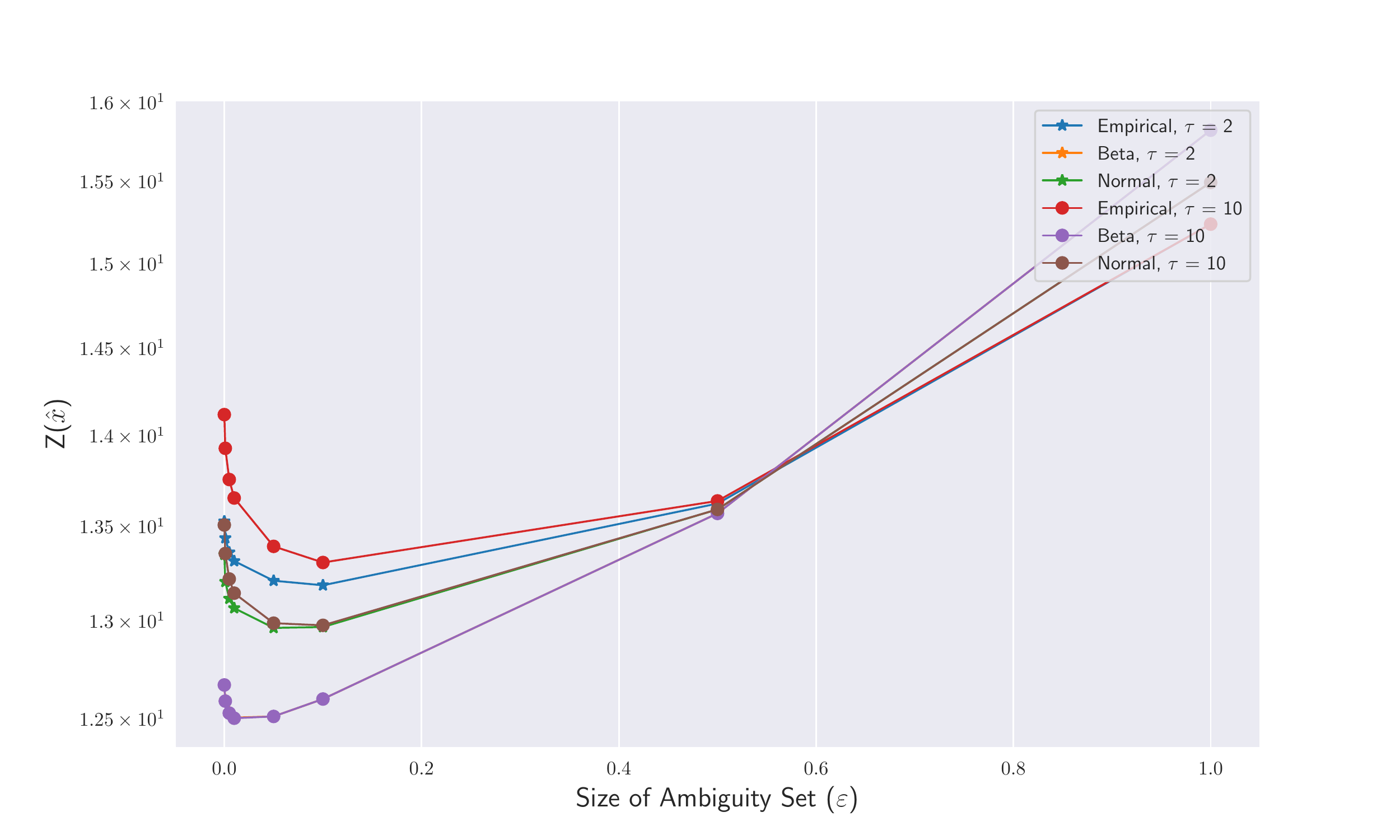}  
        \end{minipage}
    }%
    
    \subfloat[$(\gamma, n) = (4, 25)$] 
    {
        \begin{minipage}[t]{0.25\textwidth}
            \centering          
            \includegraphics[width=0.9\textwidth, trim = 20 10 45 45]{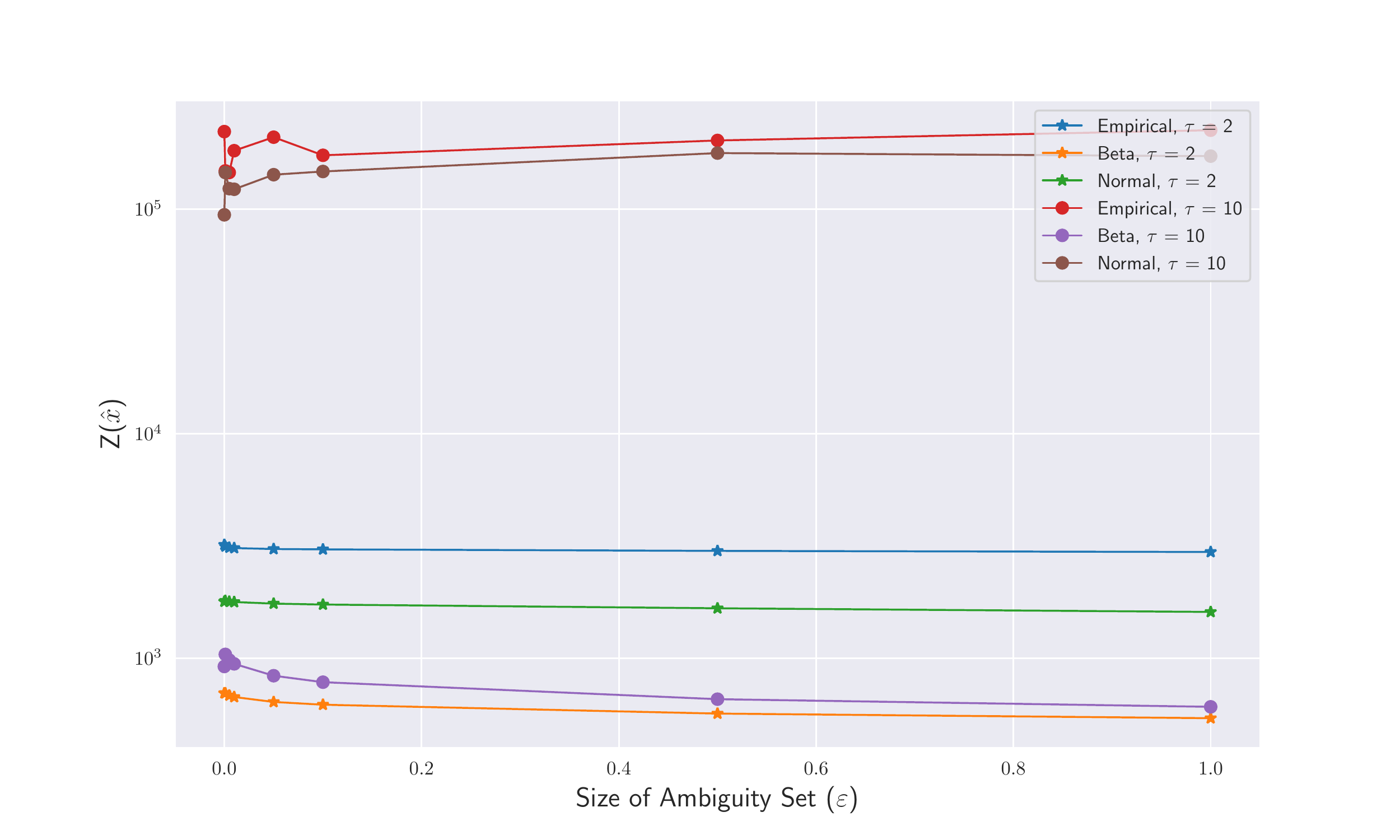}  
        \end{minipage}
    }
    \subfloat[$(\gamma, n) = (4, 50)$] 
    {
        \begin{minipage}[t]{0.25\textwidth}
            \centering  
            \includegraphics[width=0.9\textwidth, trim = 20 10 45 45]{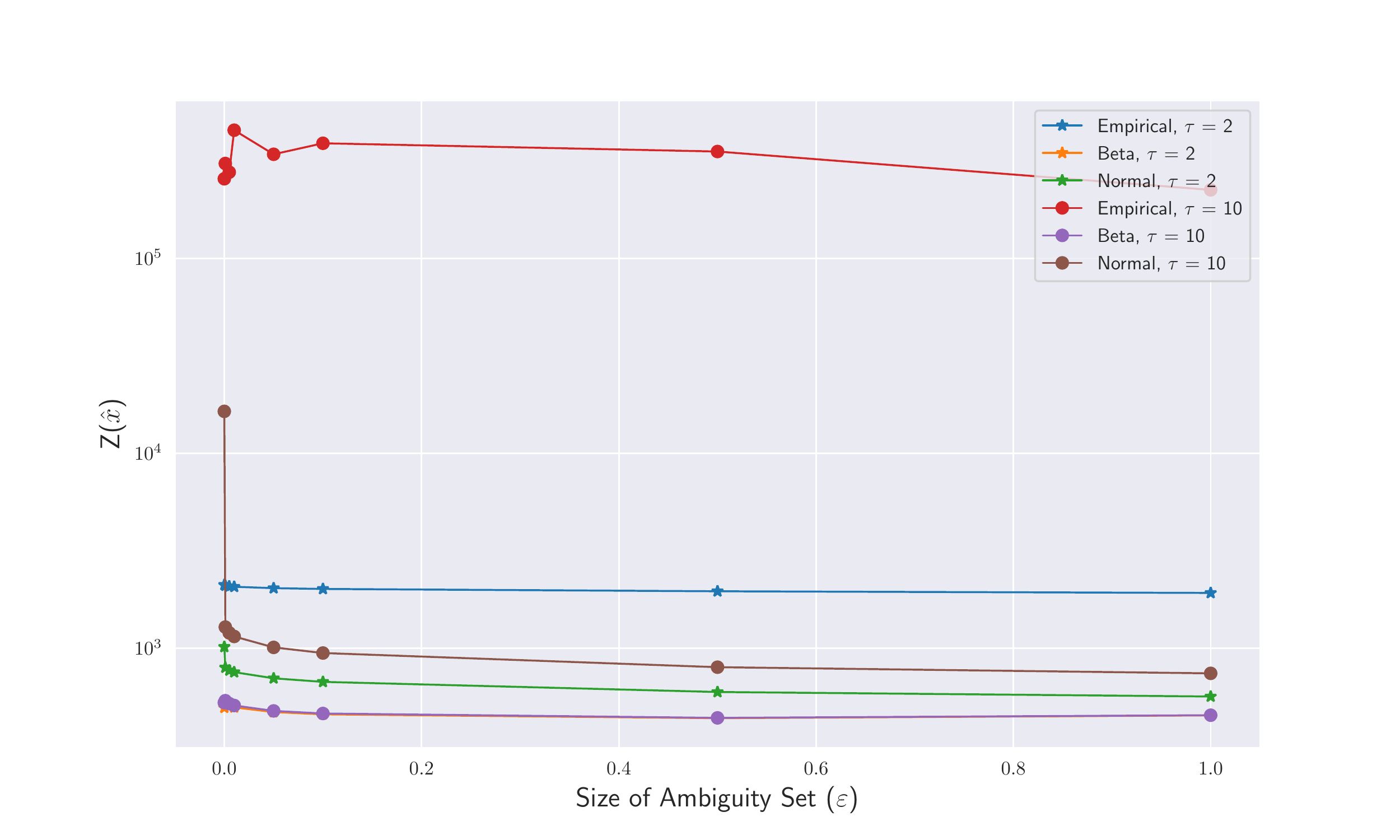} 
        \end{minipage}
    }%
    \subfloat[$(\gamma, n) = (4, 100)$] 
    {
        \begin{minipage}[t]{0.25\textwidth}
            \centering          
            \includegraphics[width = 0.9\textwidth, trim = 20 10 45 25 ]{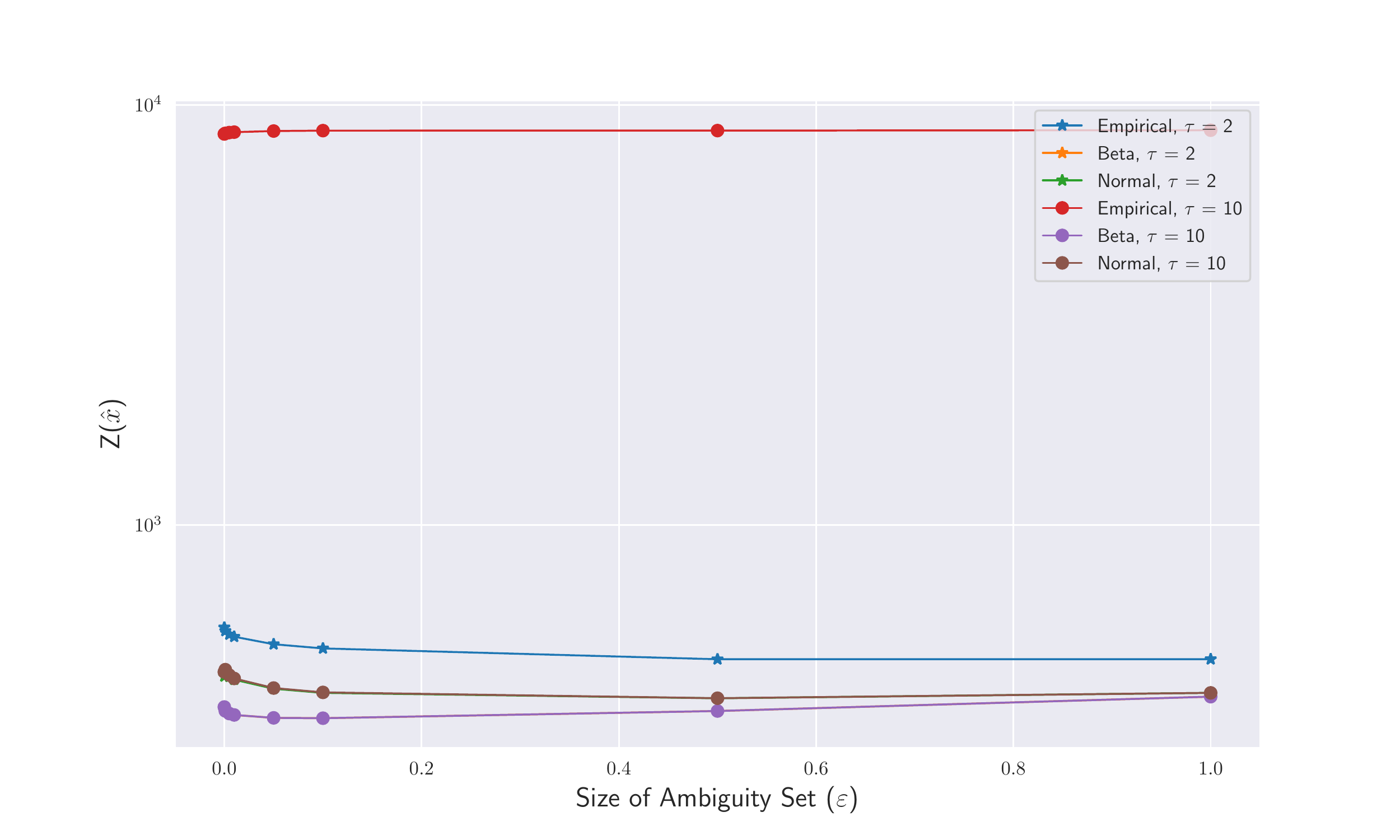} 
        \end{minipage}%
    }
    \subfloat[$(\gamma, n) = (4, 200)$]
    {
        \begin{minipage}[t]{0.25\textwidth}
            \centering      
            \includegraphics[width = 0.9\textwidth, trim = 20 10 45 25]{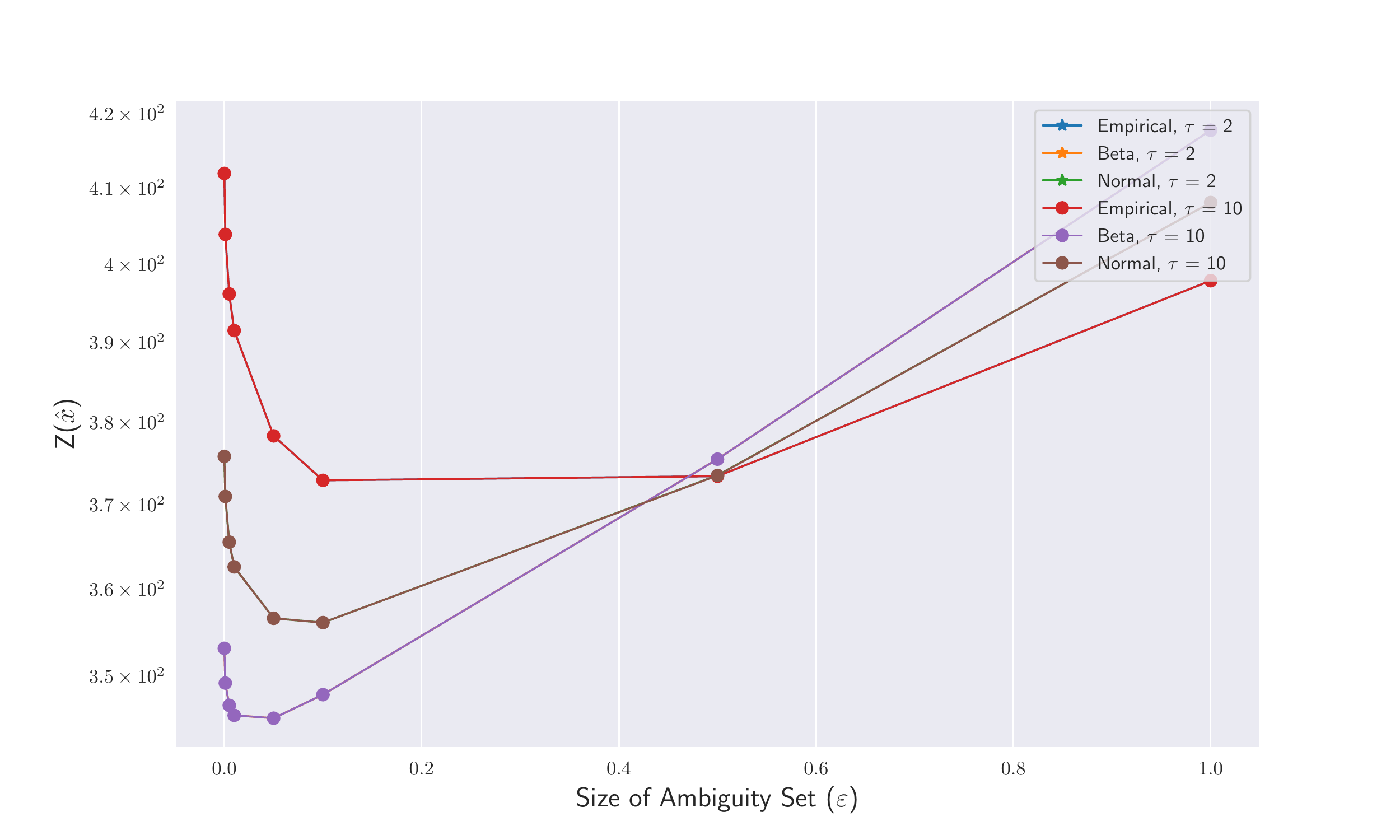}  
        \end{minipage}
    }%
    \caption{Value of cost function across different ERM-DRO models varying sample size $n$ and $\eta$.}
    \label{fig:simu2-param}   
\end{figure*}

We demonstrate results in three different cases depending on the generating process of $\xi$ to illustrate the effects of the ambiguity size and complexity term under a fixed ambiguity size $\varepsilon$ as follows:

(1) \textit{Fully Parametrized Case (Base Case)}. The results are shown in Figure~\ref{fig:simu2-param}. 
When $(\gamma, \tau)$ is large, $\text{Comp}(\Hscr)$ is much larger than $\text{Comp}(\Theta)$. Since parametric approaches do not depend on the complexity term, \texttt{P-DRO} enjoys relatively better performance, especially under small samples. When $(\gamma, \tau)$ is small, the gaps would be small. 

\begin{figure*}[h] 
    \centering    
    \subfloat[$(\gamma, n) = (1, 25)$] 
    {
        \begin{minipage}[t]{0.25\textwidth}
            \centering          
            \includegraphics[width=0.9\textwidth, trim = 20 10 45 45]{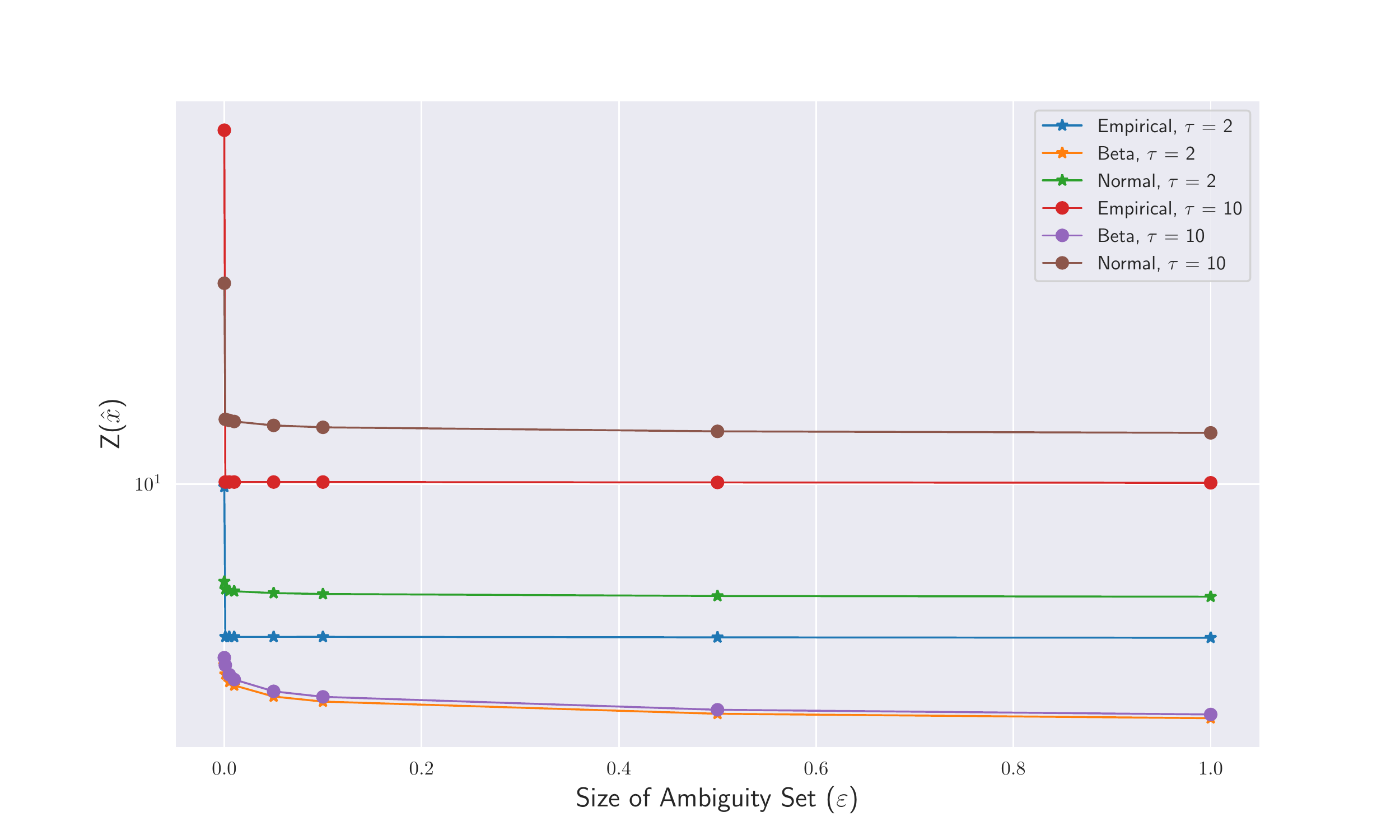}  
        \end{minipage}
    }
    \subfloat[$(\gamma, n) = (1, 50)$] 
    {
        \begin{minipage}[t]{0.25\textwidth}
            \centering  
            \includegraphics[width=0.9\textwidth, trim = 20 10 45 45]{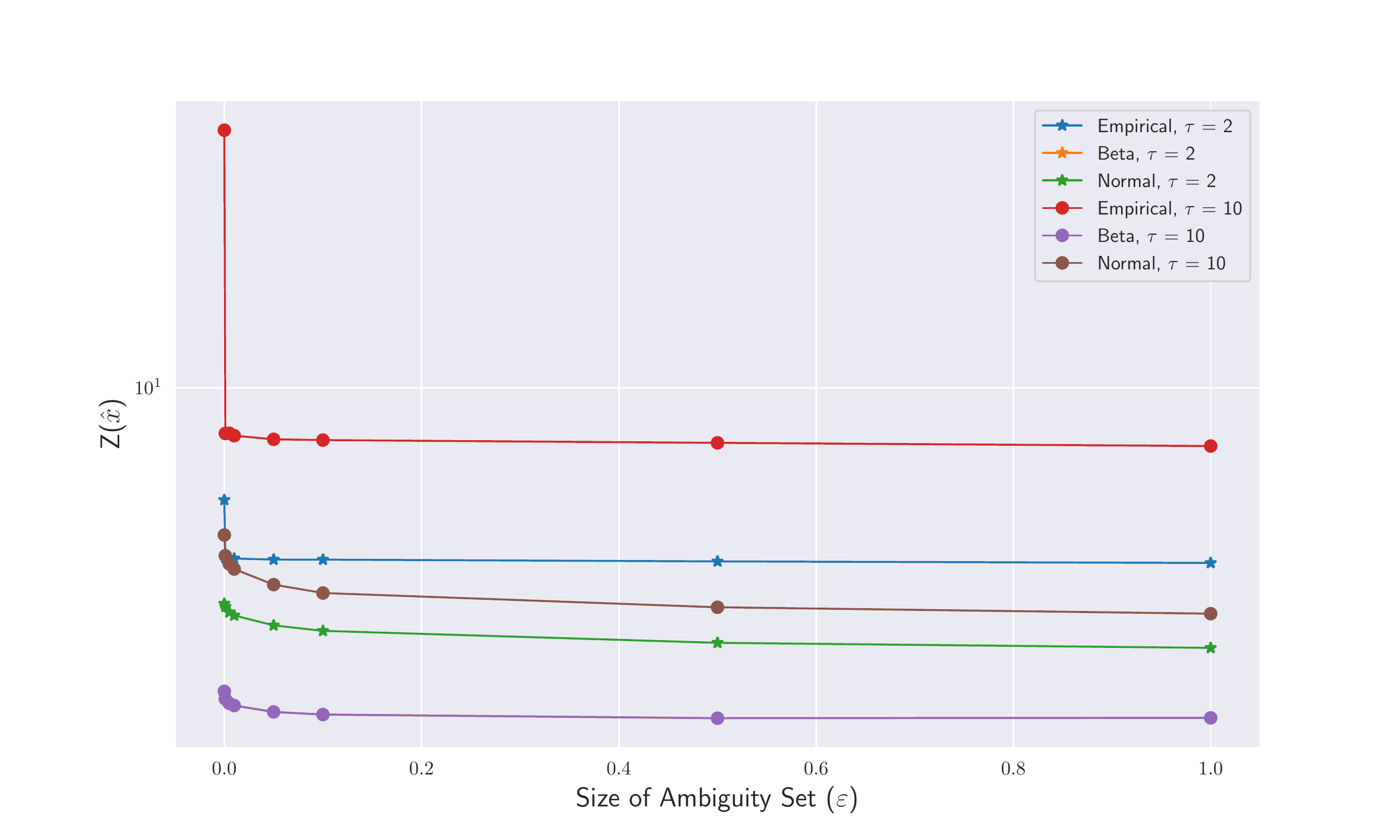} 
        \end{minipage}
    }%
    \subfloat[$(\gamma, n) = (1, 100)$] 
    {
        \begin{minipage}[t]{0.25\textwidth}
            \centering          
            \includegraphics[width = 0.9\textwidth, trim = 20 10 45 25 ]{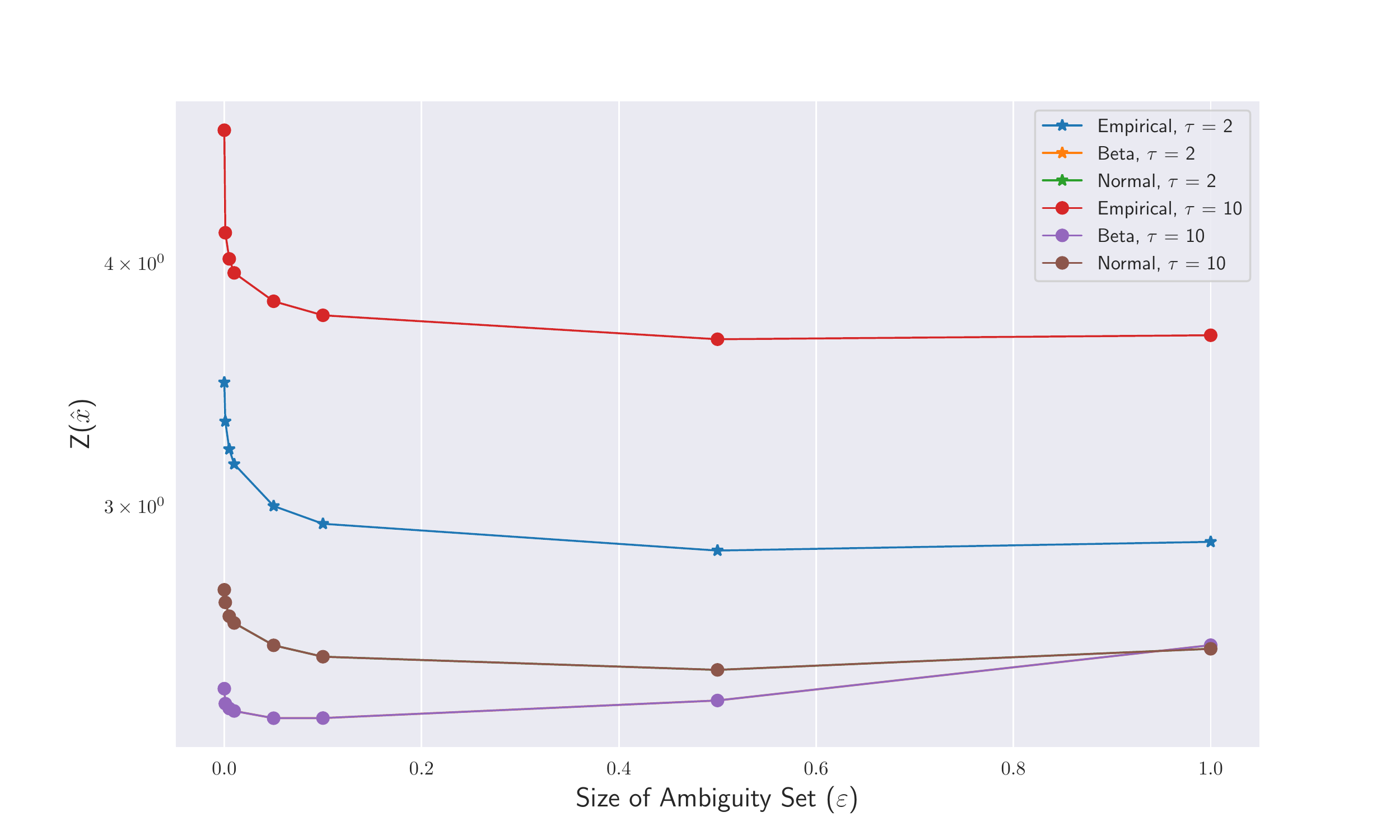} 
        \end{minipage}%
    }
    \subfloat[$(\gamma, n) = (1, 200)$]
    {
        \begin{minipage}[t]{0.25\textwidth}
            \centering      
            \includegraphics[width = 0.9\textwidth, trim = 20 10 45 25]{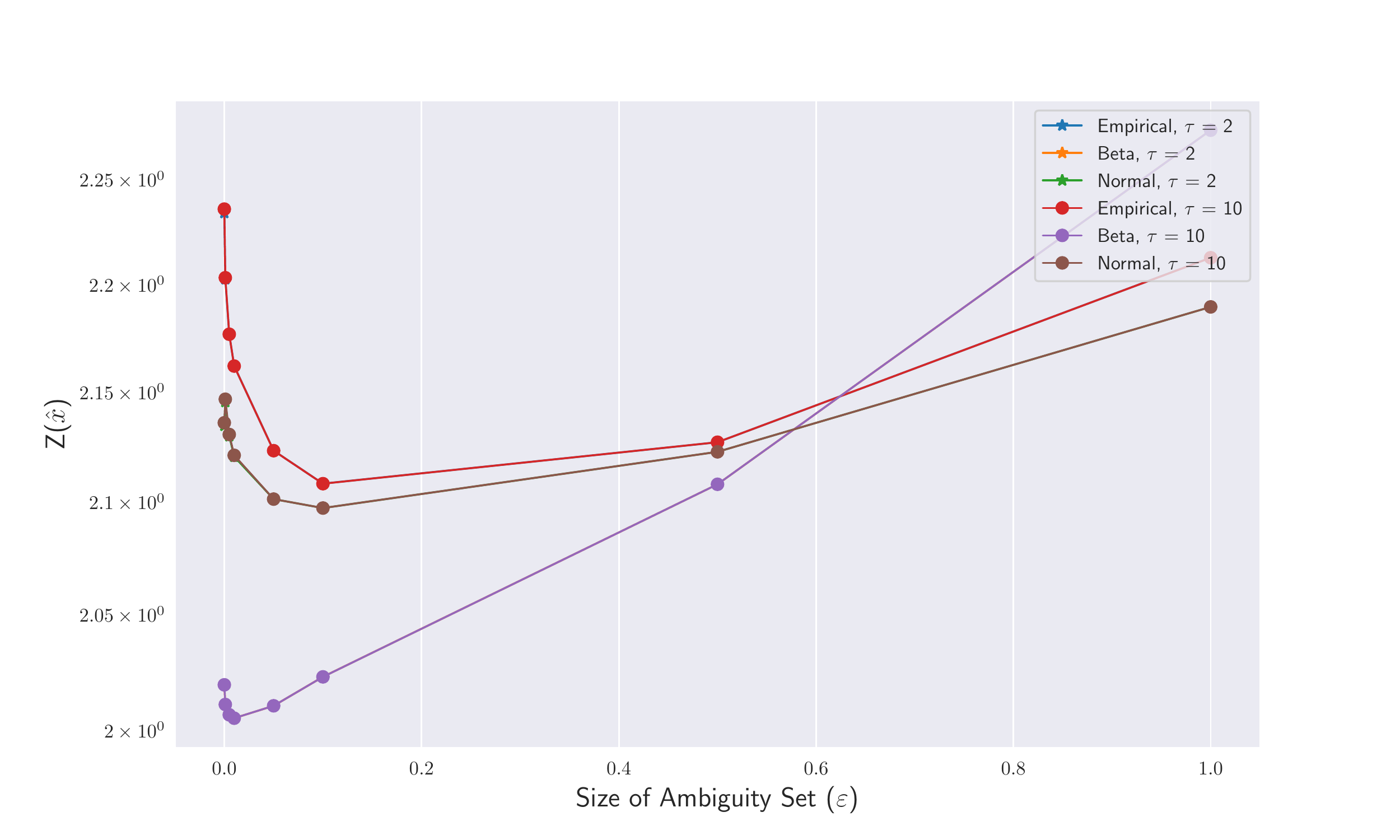}  
        \end{minipage}
    }%
    
    \subfloat[$(\gamma, n) = (2, 25)$] 
    {
        \begin{minipage}[t]{0.25\textwidth}
            \centering          
            \includegraphics[width=0.9\textwidth, trim = 20 10 45 45]{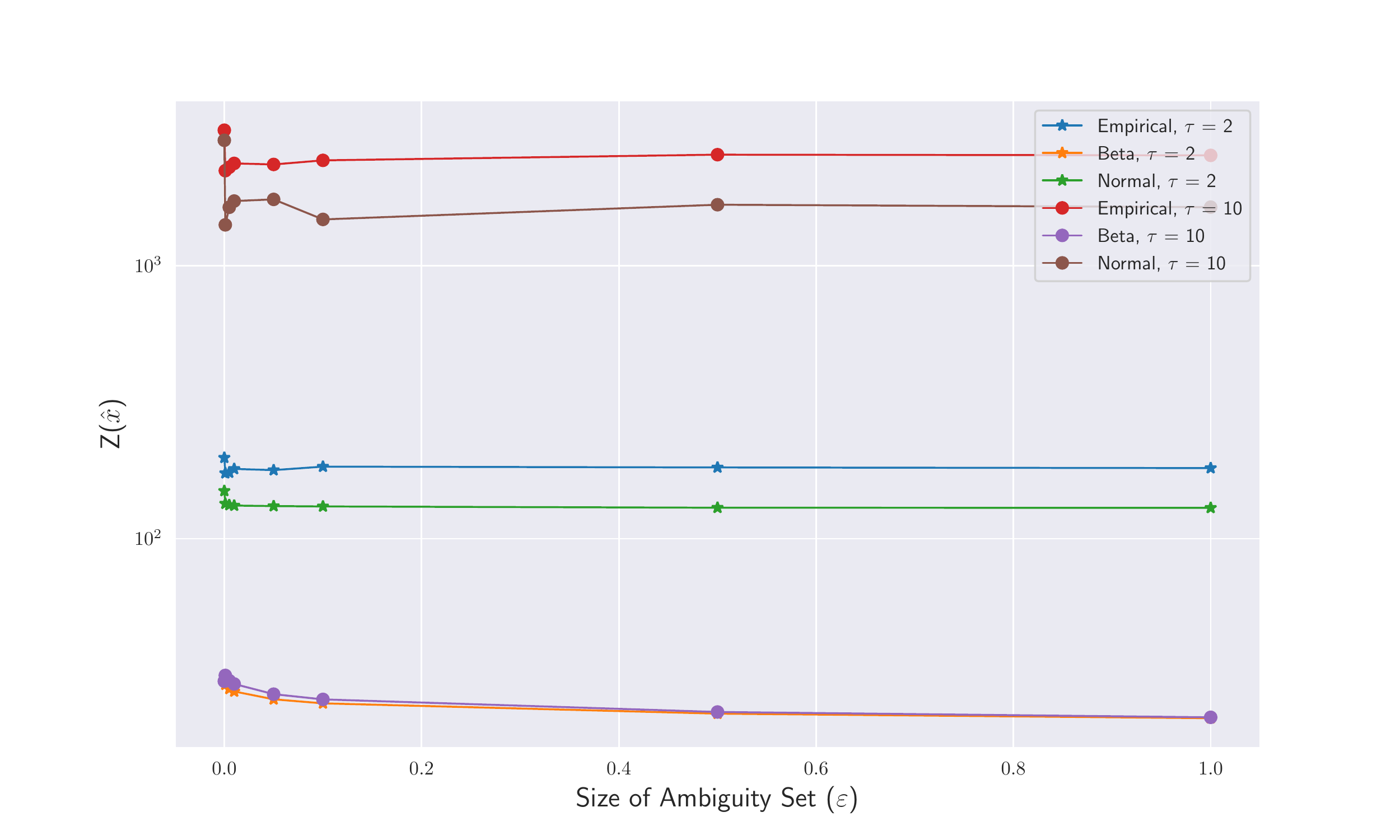}  
        \end{minipage}
    }
    \subfloat[$(\gamma, n) = (2, 50)$] 
    {
        \begin{minipage}[t]{0.25\textwidth}
            \centering  
            \includegraphics[width=0.9\textwidth, trim = 20 10 45 45]{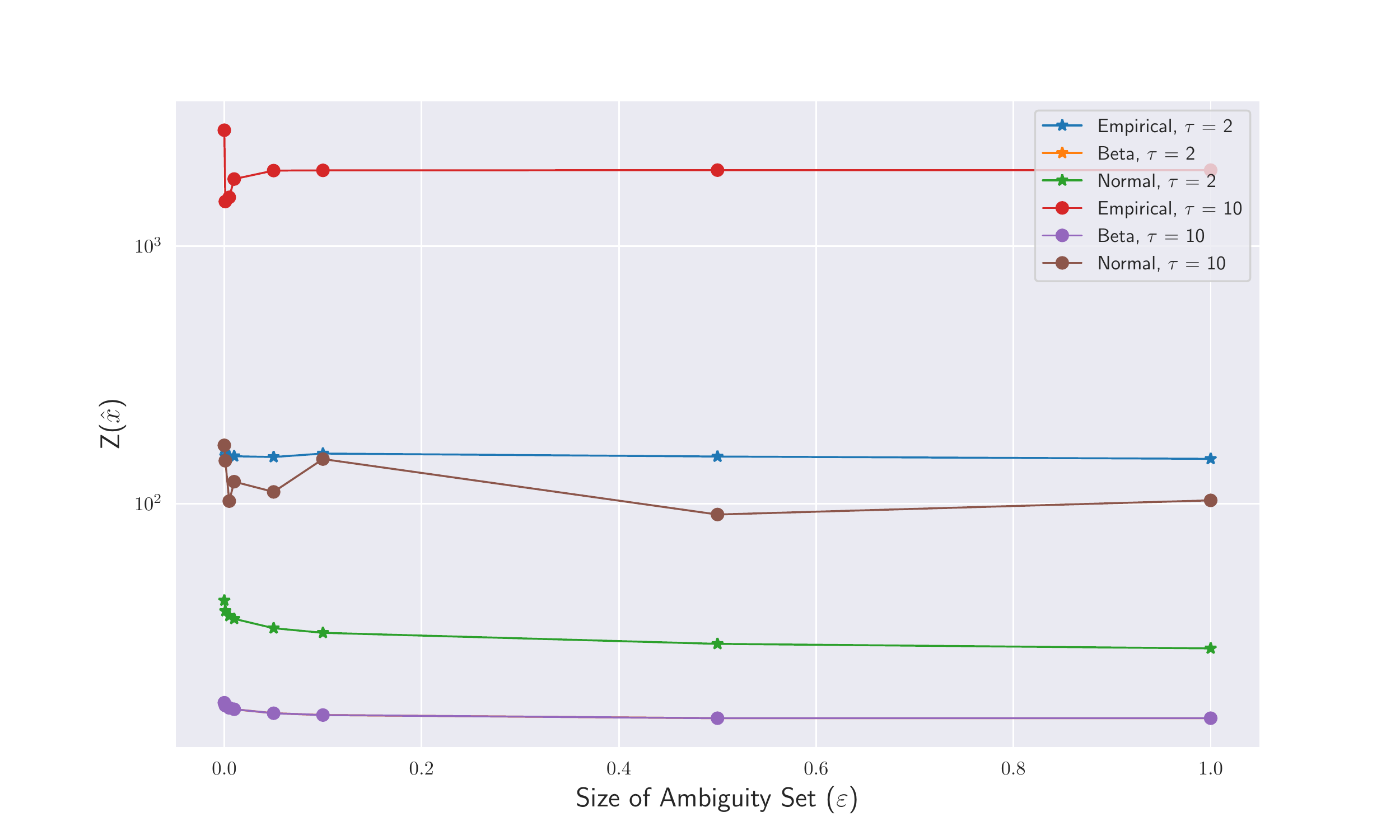} 
        \end{minipage}
    }%
    \subfloat[$(\gamma, n) = (2, 100)$] 
    {
        \begin{minipage}[t]{0.25\textwidth}
            \centering          
            \includegraphics[width = 0.9\textwidth, trim = 20 10 45 25 ]{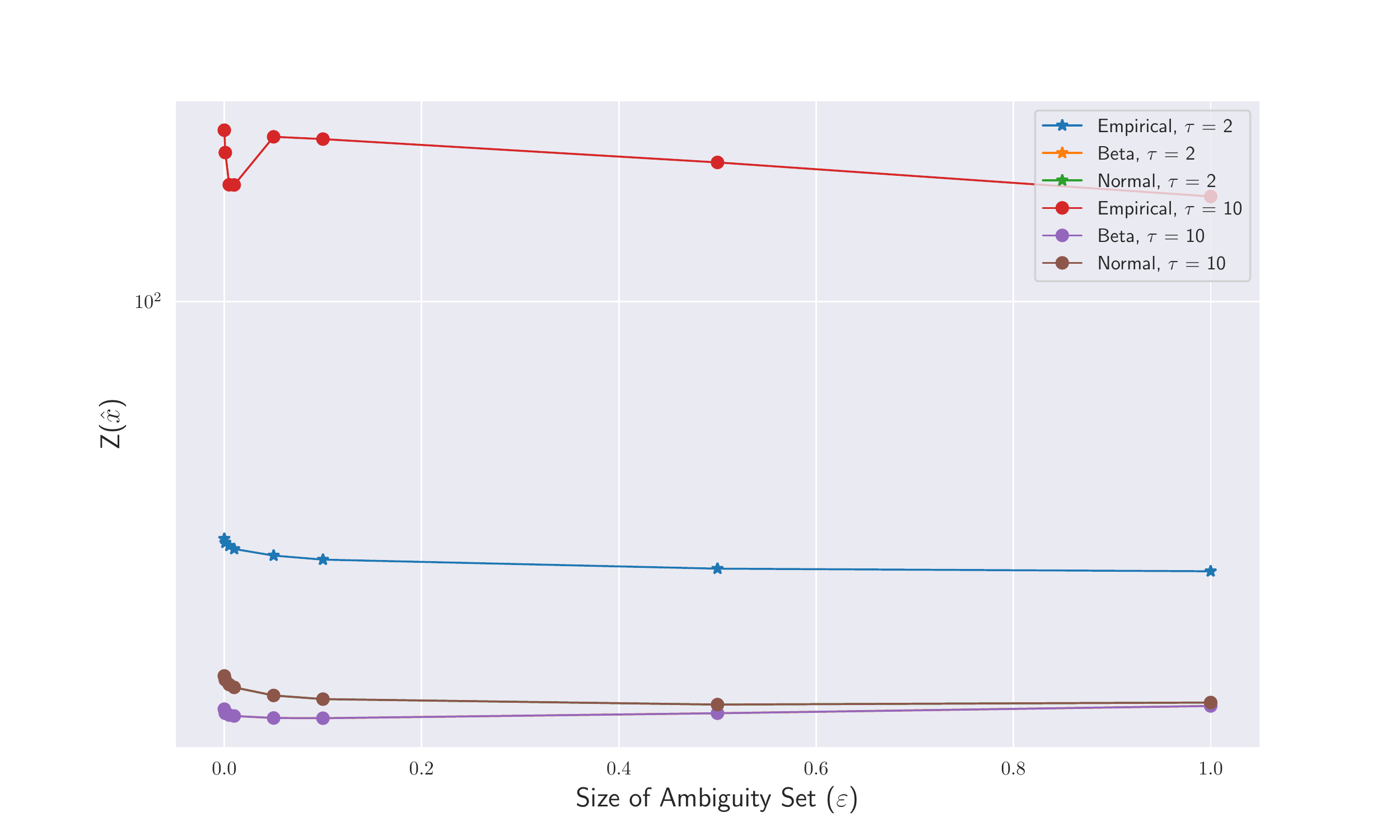} 
        \end{minipage}%
    }
    \subfloat[$(\gamma, n) = (2, 200)$]
    {
        \begin{minipage}[t]{0.25\textwidth}
            \centering      
            \includegraphics[width = 0.9\textwidth, trim = 20 10 45 25]{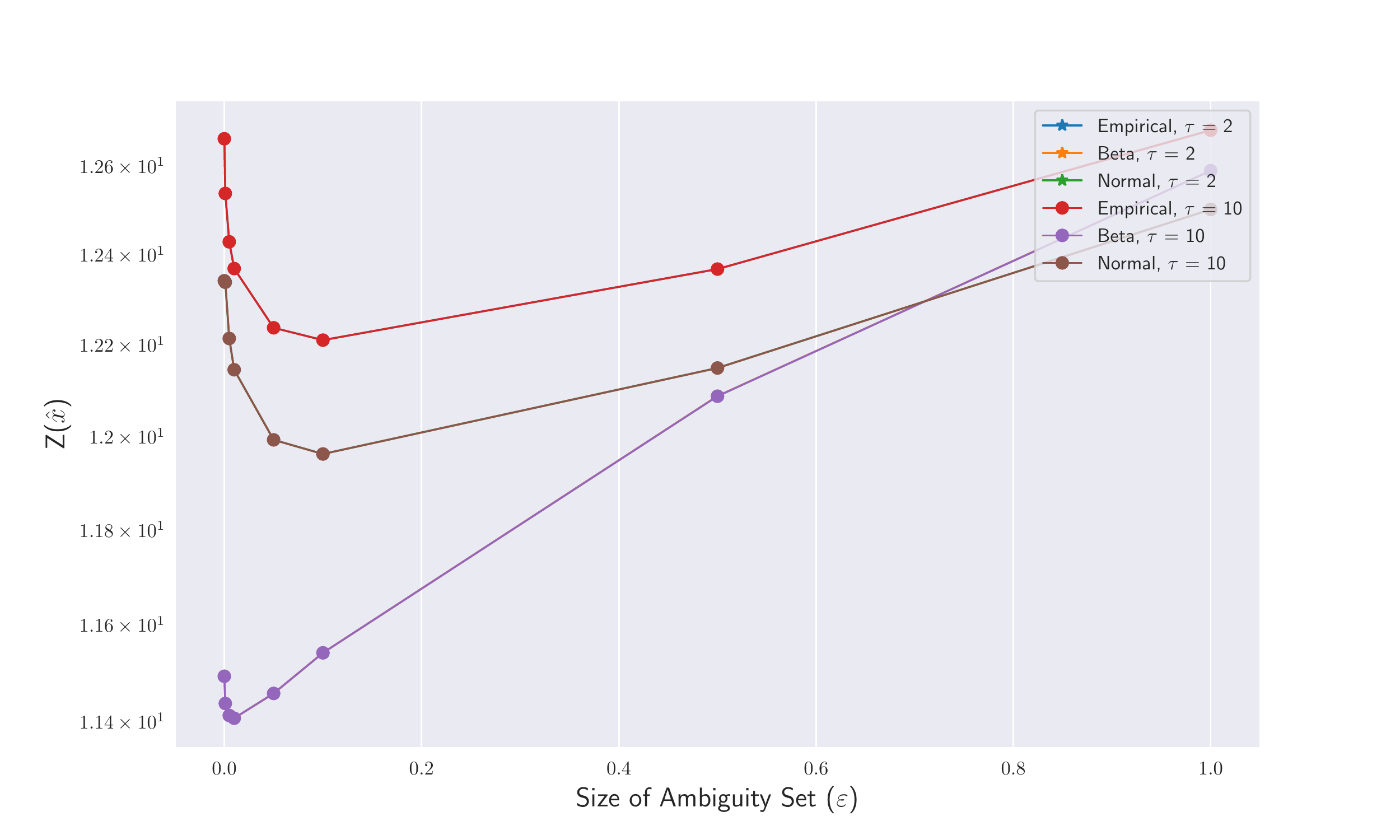}  
        \end{minipage}
    }%
    
    \subfloat[$(\gamma, n) = (4, 25)$] 
    {
        \begin{minipage}[t]{0.25\textwidth}
            \centering          
            \includegraphics[width=0.9\textwidth, trim = 20 10 45 45]{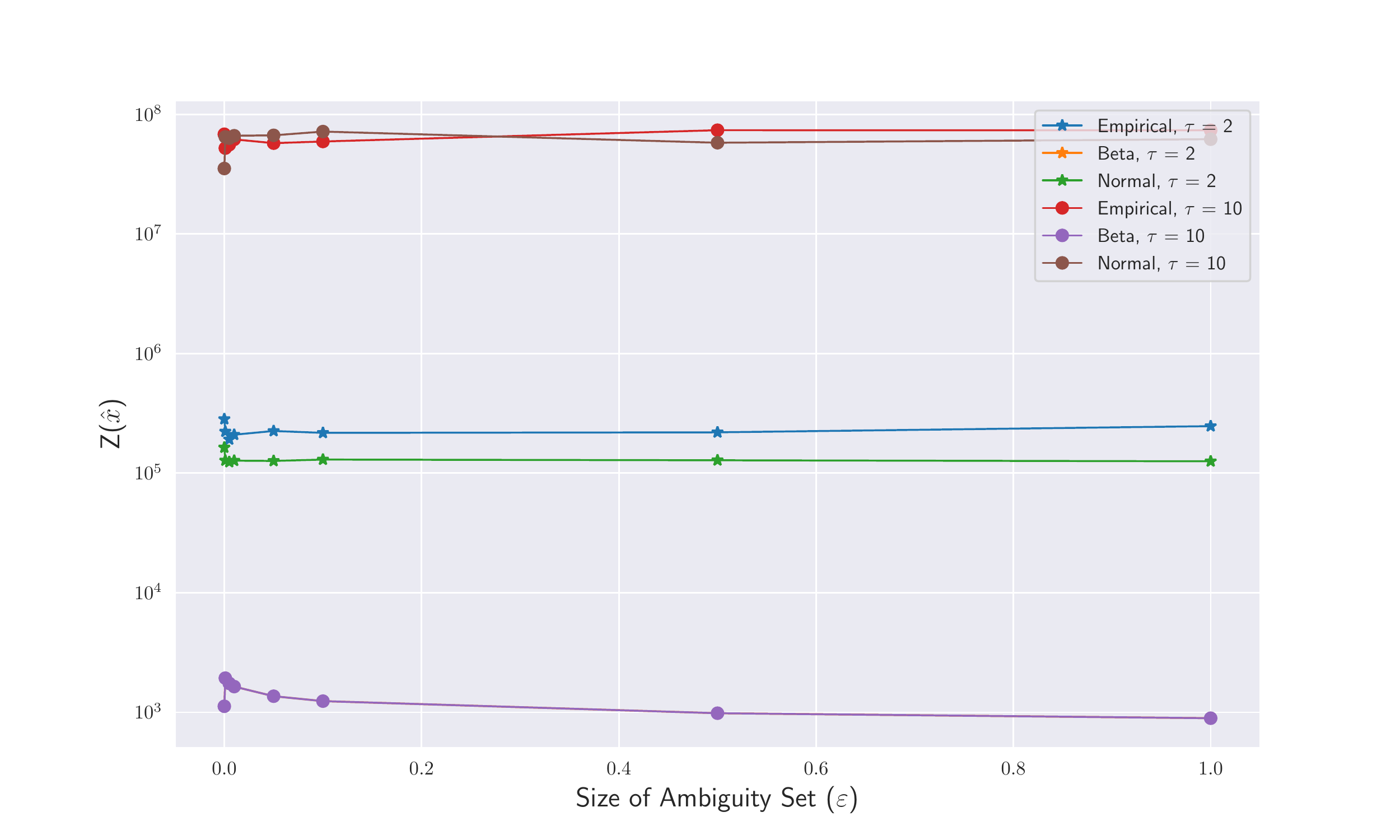}  
        \end{minipage}
    }
    \subfloat[$(\gamma, n) = (4, 50)$] 
    {
        \begin{minipage}[t]{0.25\textwidth}
            \centering  
            \includegraphics[width=0.9\textwidth, trim = 20 10 45 45]{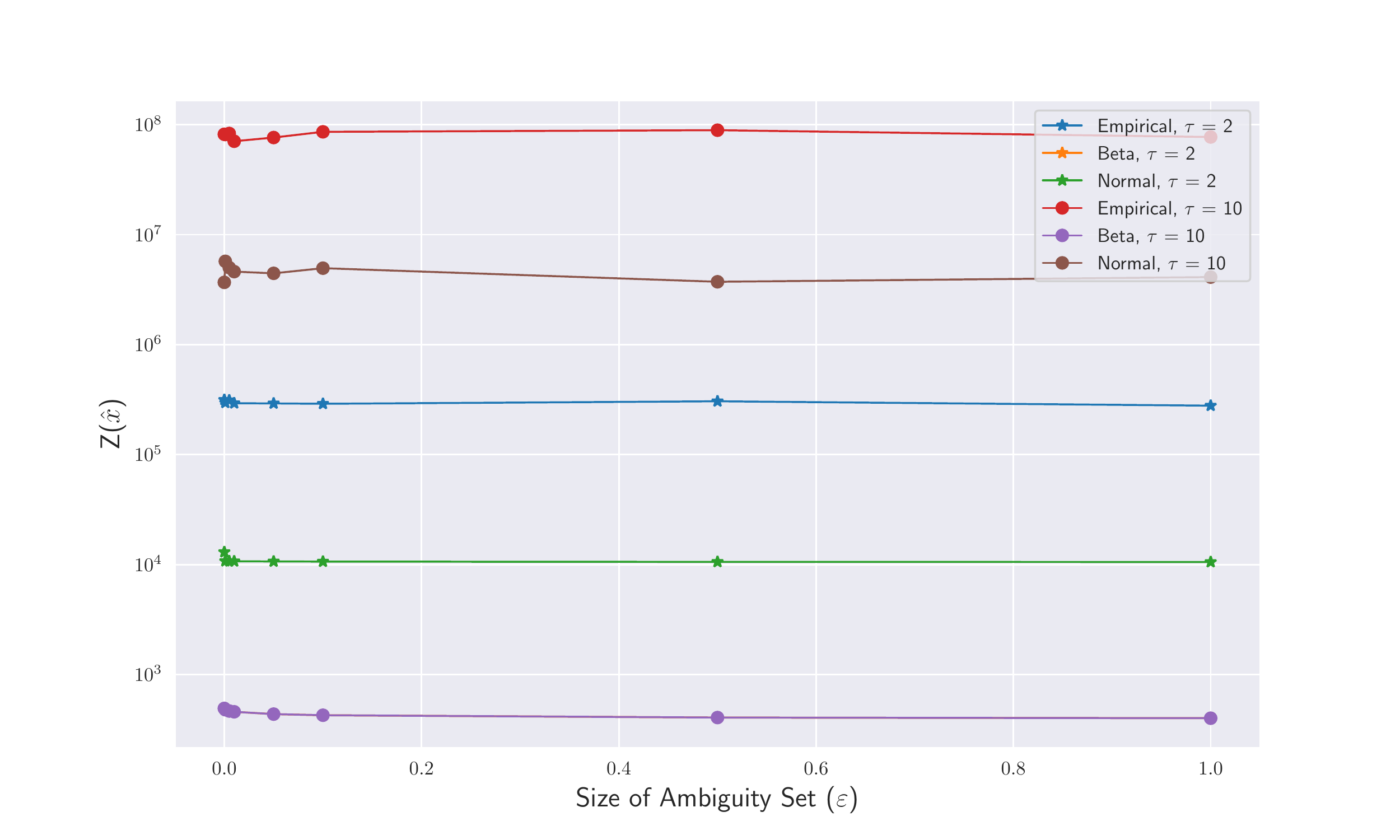} 
        \end{minipage}
    }%
    \subfloat[$(\gamma, n) = (4, 100)$] 
    {
        \begin{minipage}[t]{0.25\textwidth}
            \centering          
            \includegraphics[width = 0.9\textwidth, trim = 20 10 45 25 ]{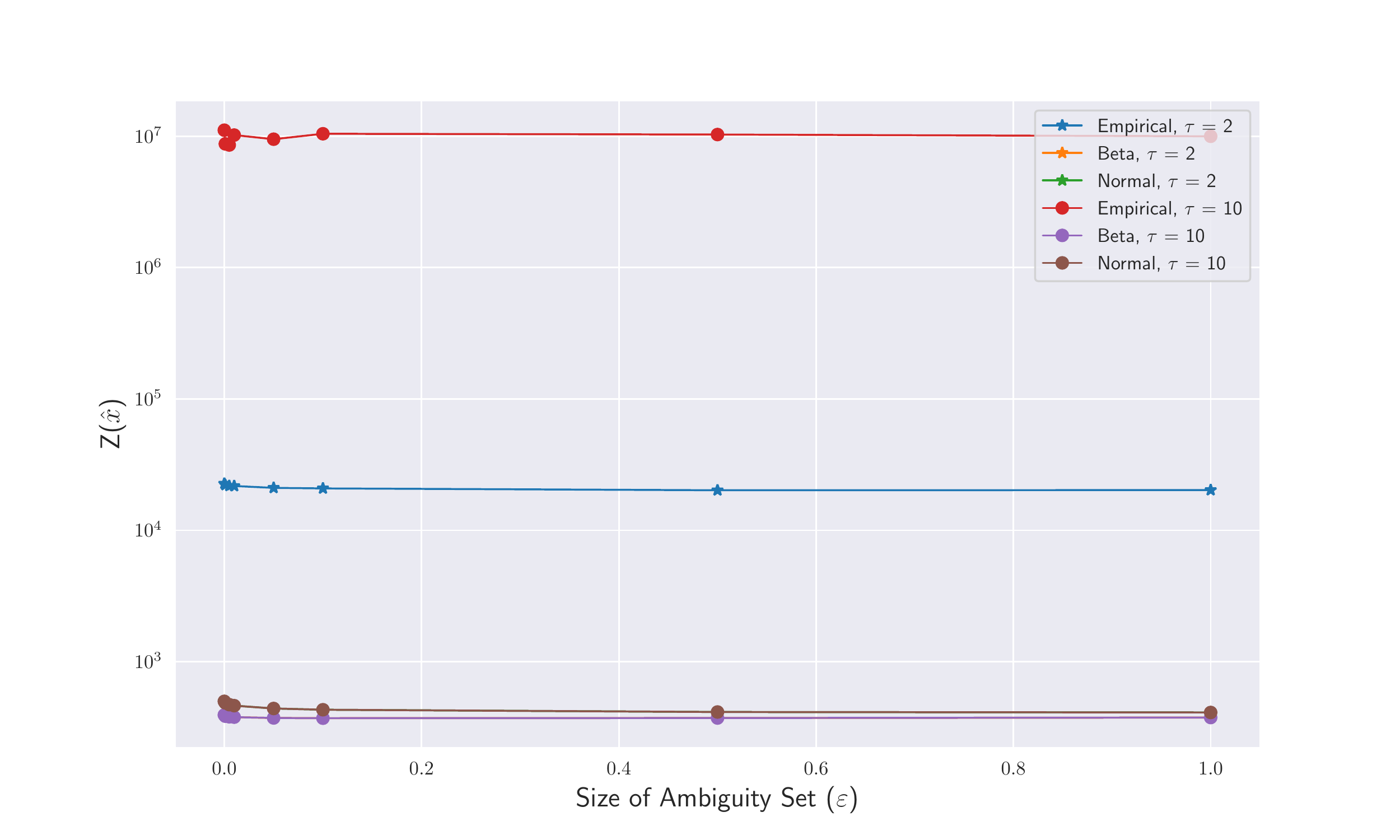} 
        \end{minipage}%
    }
    \subfloat[$(\gamma, n) = (4, 200)$]
    {
        \begin{minipage}[t]{0.25\textwidth}
            \centering      
            \includegraphics[width = 0.9\textwidth, trim = 20 10 45 25]{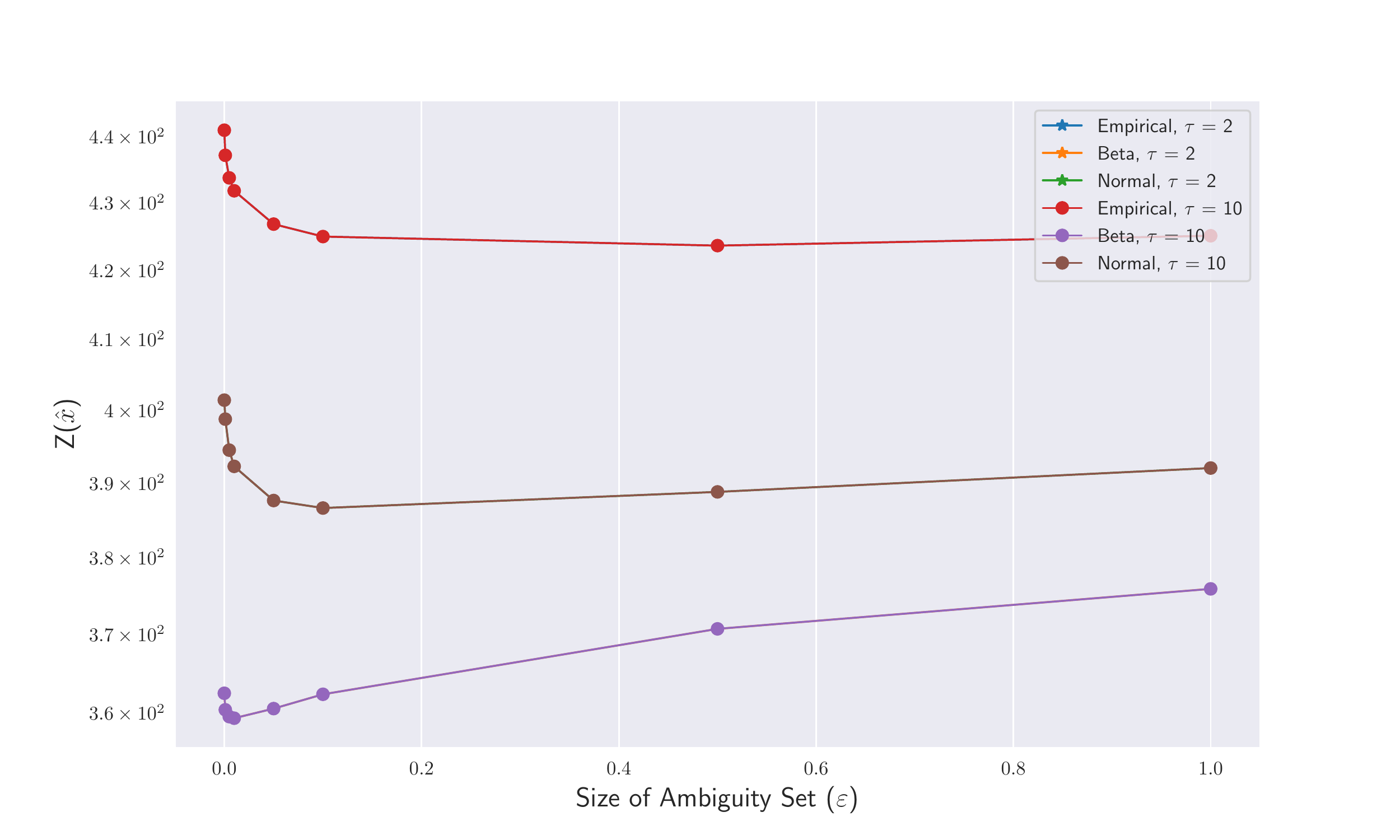}  
        \end{minipage}
    }%
    \caption{Value of cost function across different ERM-DRO models varying sample size $n$ and $\eta$ under misspecification.}
    \label{fig:simu2-mis}   
\end{figure*}

(2) \textit{Distribution Misspecification.} We perturb the $i$-th marginal distribution of the random variable $\xi_i$ to $\xi_i + \zeta_i$, where each $\zeta_i \sim U(-2,2)$ and is independent with $\xi_i$. The results are shown in Figure~\ref{fig:simu2-mis} with more noticeable performance advantages for \texttt{P-DRO} models.

\begin{figure*}[h] 
    \centering    
        \subfloat[$(\gamma, n) = (1, 25)$] 
    {
        \begin{minipage}[t]{0.25\textwidth}
            \centering          
            \includegraphics[width=0.9\textwidth, trim = 20 10 45 45]{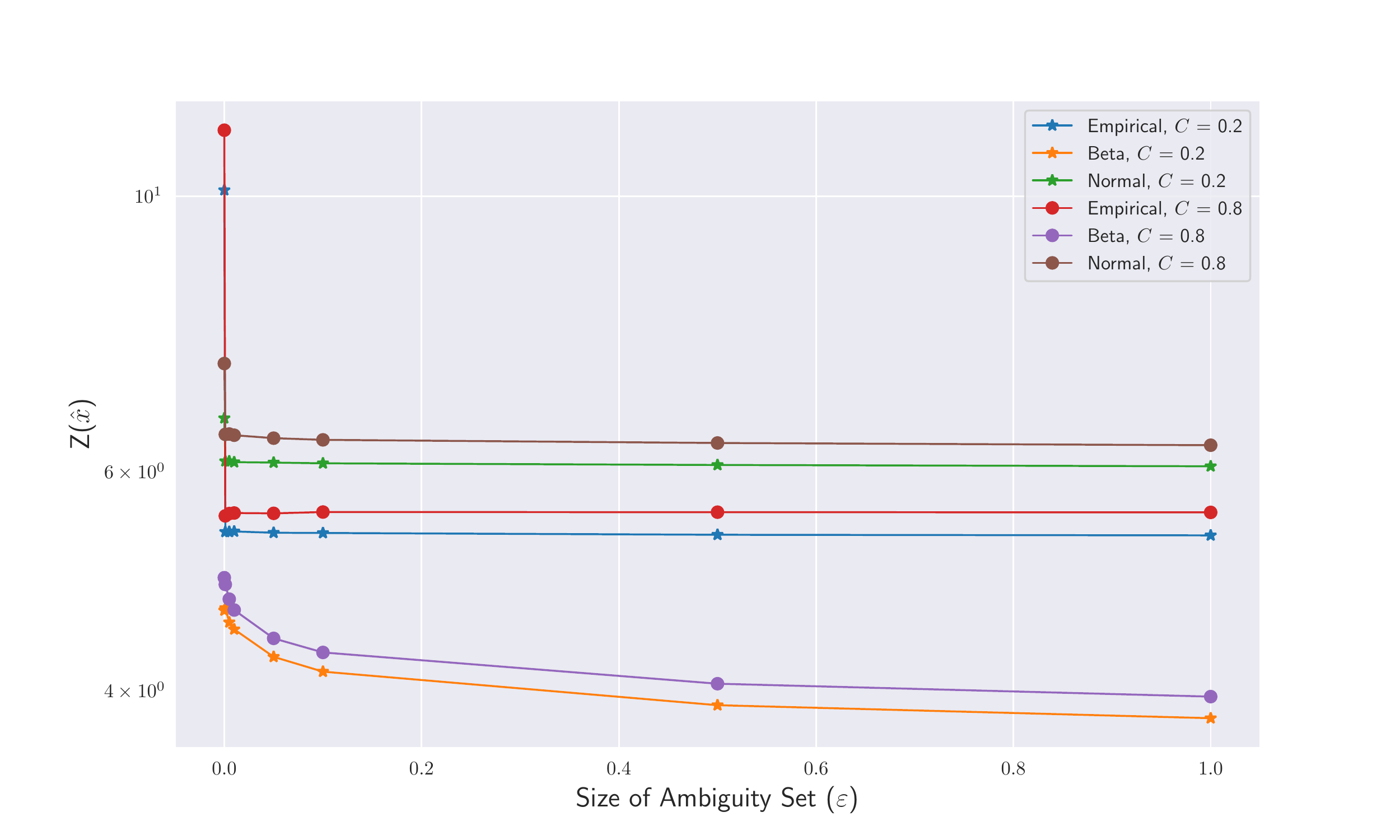}  
        \end{minipage}
    }
    \subfloat[$(\gamma, n) = (1, 50)$] 
    {
        \begin{minipage}[t]{0.25\textwidth}
            \centering  
            \includegraphics[width=0.9\textwidth, trim = 20 10 45 45]{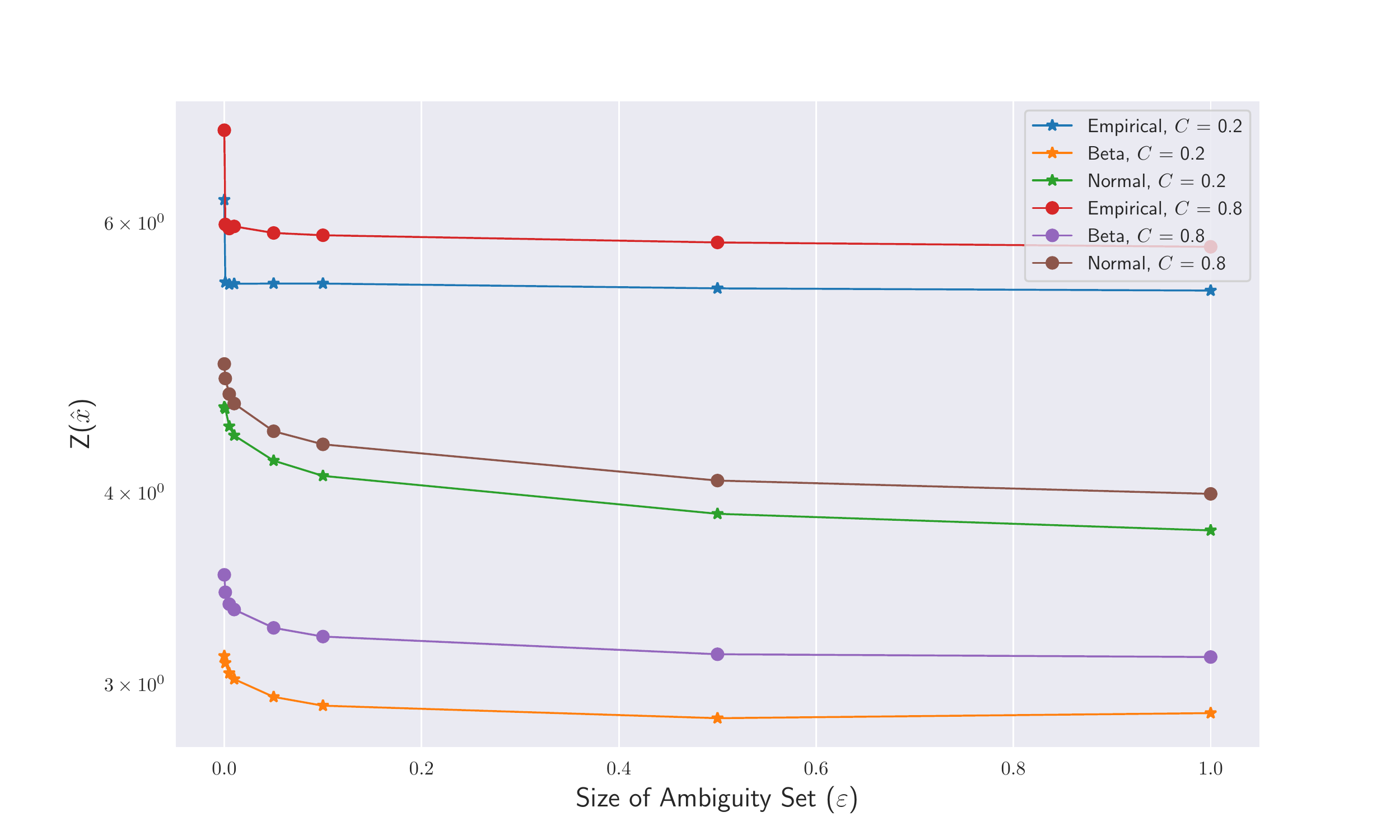} 
        \end{minipage}
    }%
    \subfloat[$(\gamma, n) = (1, 100)$] 
    {
        \begin{minipage}[t]{0.25\textwidth}
            \centering          
            \includegraphics[width = 0.9\textwidth, trim = 20 10 45 25 ]{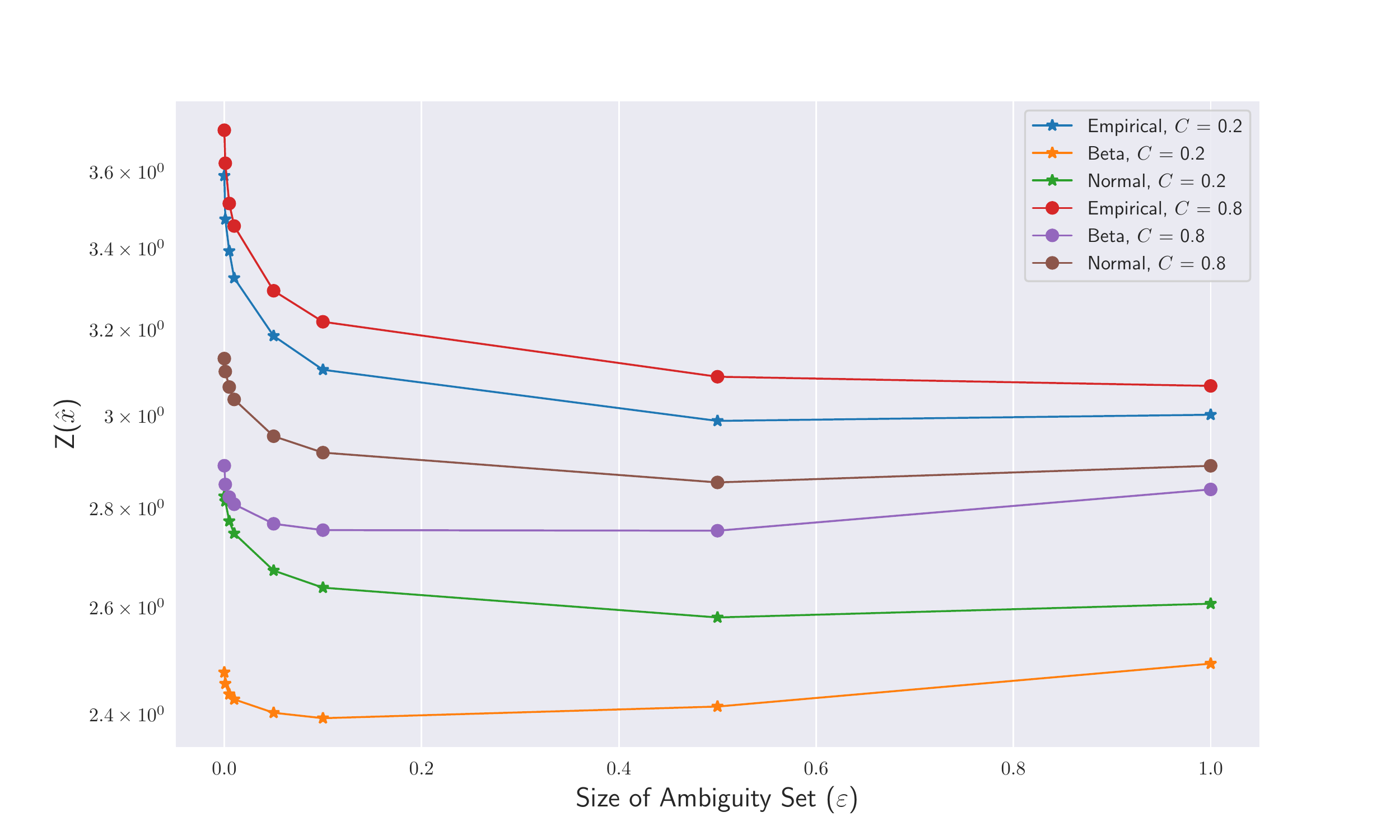} 
        \end{minipage}%
    }
    \subfloat[$(\gamma, n) = (1, 200)$]
    {
        \begin{minipage}[t]{0.25\textwidth}
            \centering      
            \includegraphics[width = 0.9\textwidth, trim = 20 10 45 25]{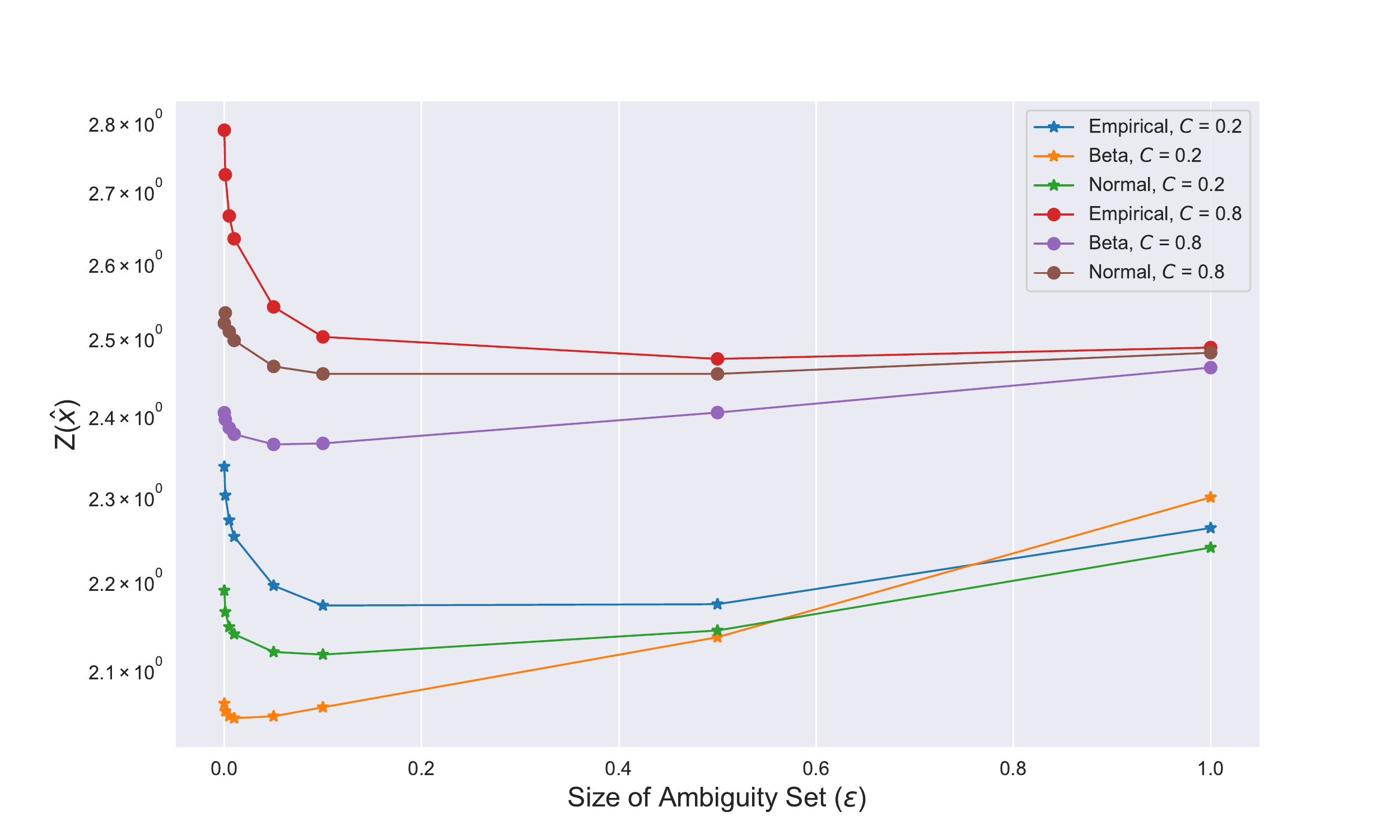}  
        \end{minipage}
    }%
    
    \subfloat[$(\gamma, n) = (2, 25)$] 
    {
        \begin{minipage}[t]{0.25\textwidth}
            \centering          
            \includegraphics[width=0.9\textwidth, trim = 20 10 45 45]{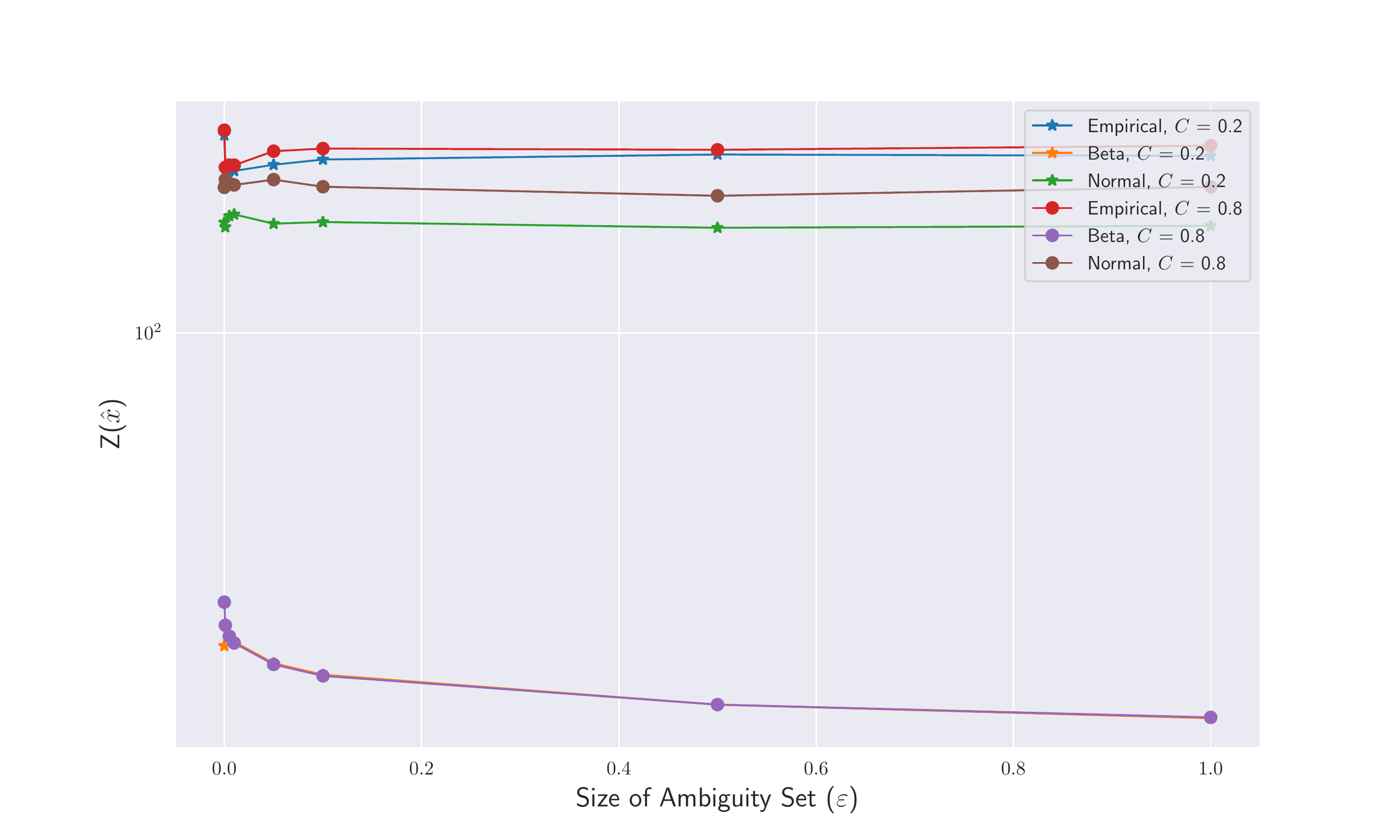}  
        \end{minipage}
    }
    \subfloat[$(\gamma, n) = (2, 50)$] 
    {
        \begin{minipage}[t]{0.25\textwidth}
            \centering  
            \includegraphics[width=0.9\textwidth, trim = 20 10 45 45]{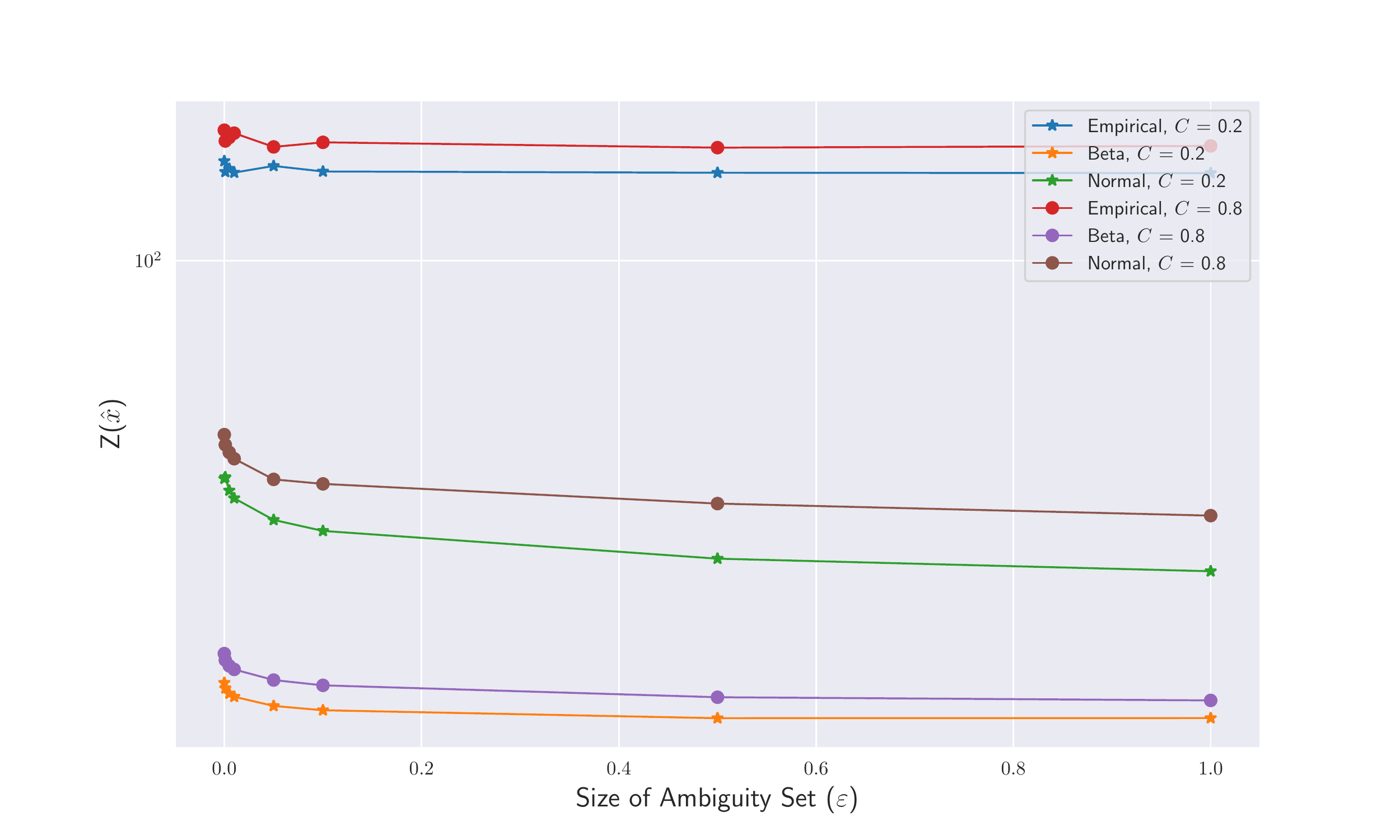} 
        \end{minipage}
    }%
    \subfloat[$(\gamma, n) = (2, 100)$] 
    {
        \begin{minipage}[t]{0.25\textwidth}
            \centering          
            \includegraphics[width = 0.9\textwidth, trim = 20 10 45 25 ]{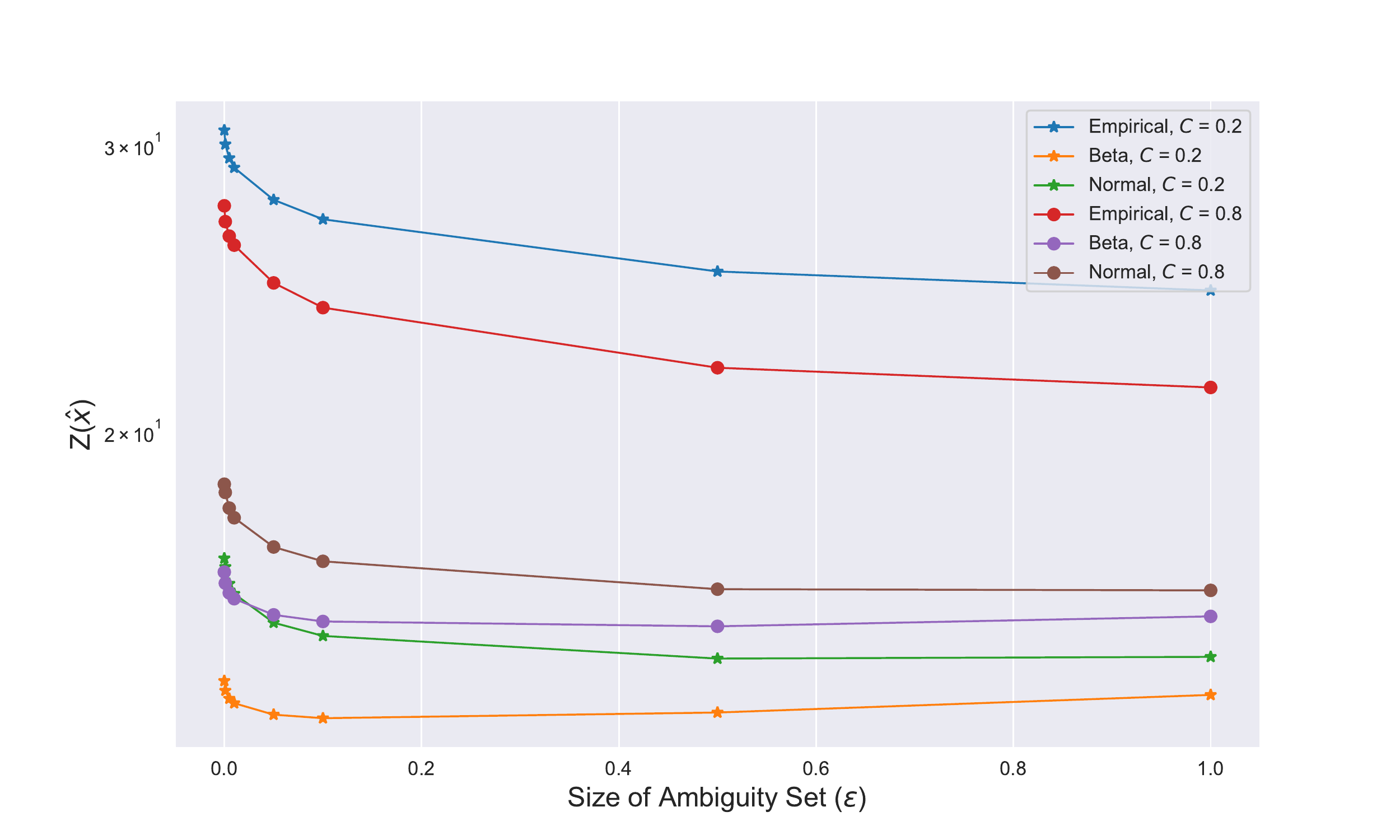} 
        \end{minipage}%
    }
    \subfloat[$(\gamma, n) = (2, 200)$]
    {
        \begin{minipage}[t]{0.25\textwidth}
            \centering      
            \includegraphics[width = 0.9\textwidth, trim = 20 10 45 25]{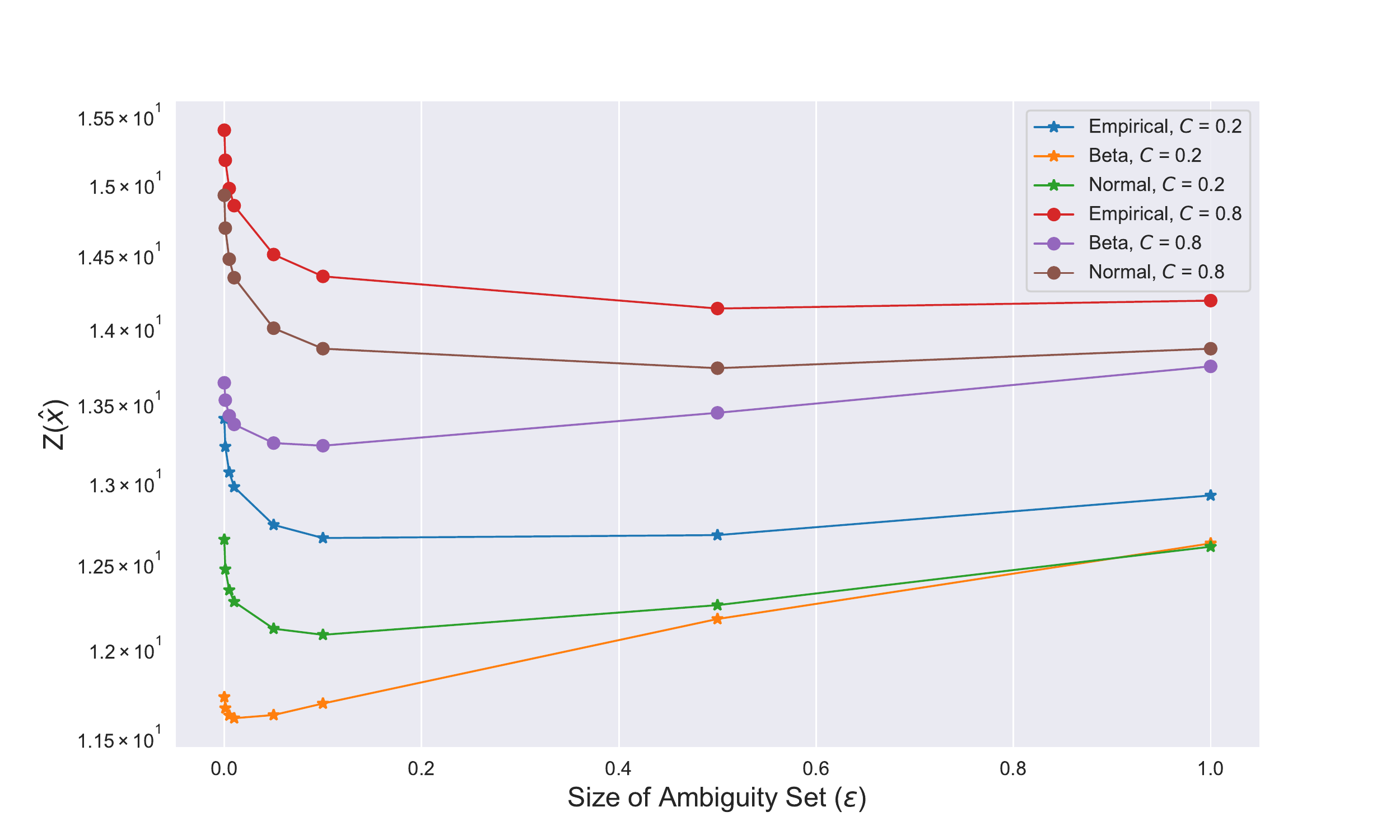}  
        \end{minipage}
    }%
    
    \subfloat[$(\gamma, n) = (4, 25)$] 
    {
        \begin{minipage}[t]{0.25\textwidth}
            \centering          
            \includegraphics[width=0.9\textwidth, trim = 20 10 45 45]{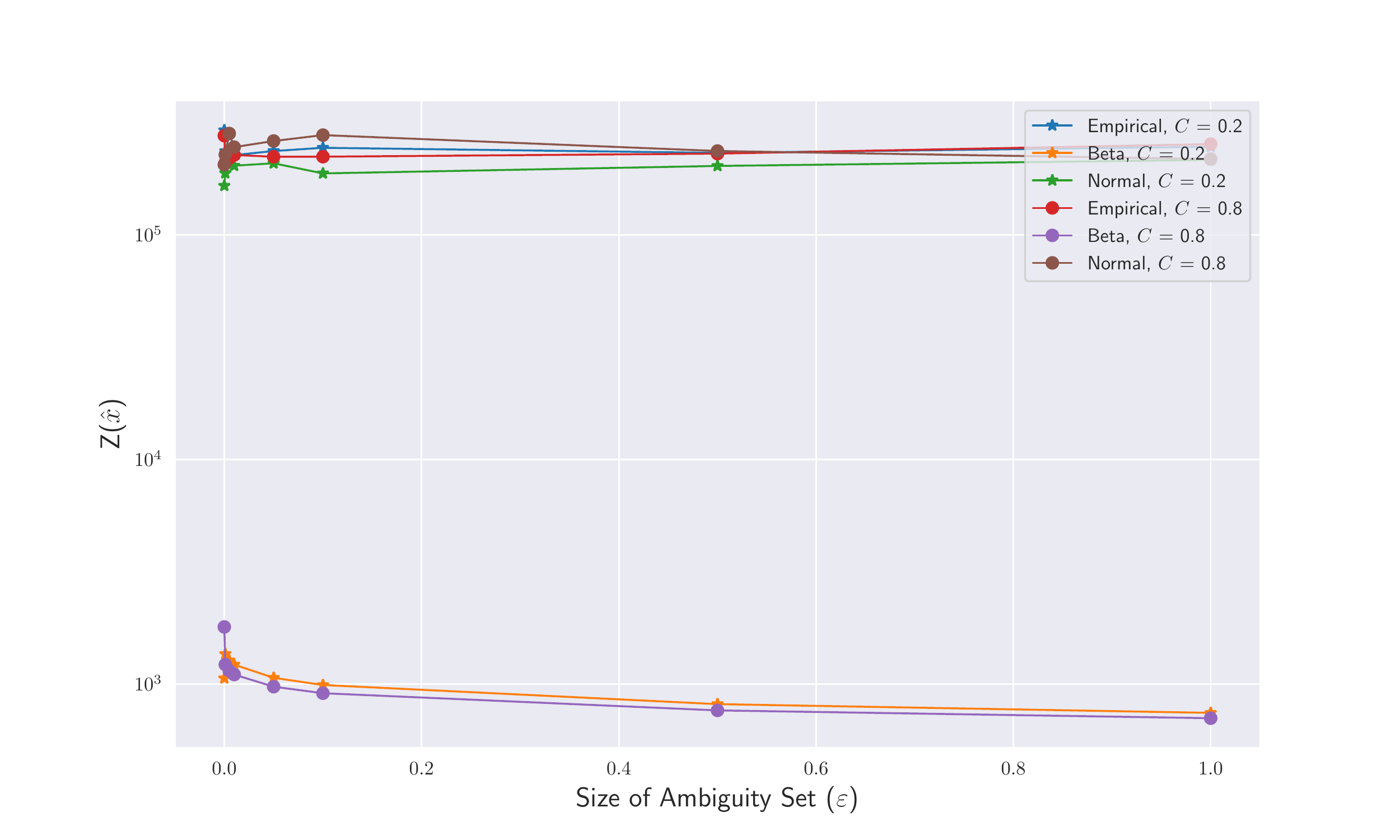}  
        \end{minipage}
    }
    \subfloat[$(\gamma, n) = (4, 50)$] 
    {
        \begin{minipage}[t]{0.25\textwidth}
            \centering  
            \includegraphics[width=0.9\textwidth, trim = 20 10 45 45]{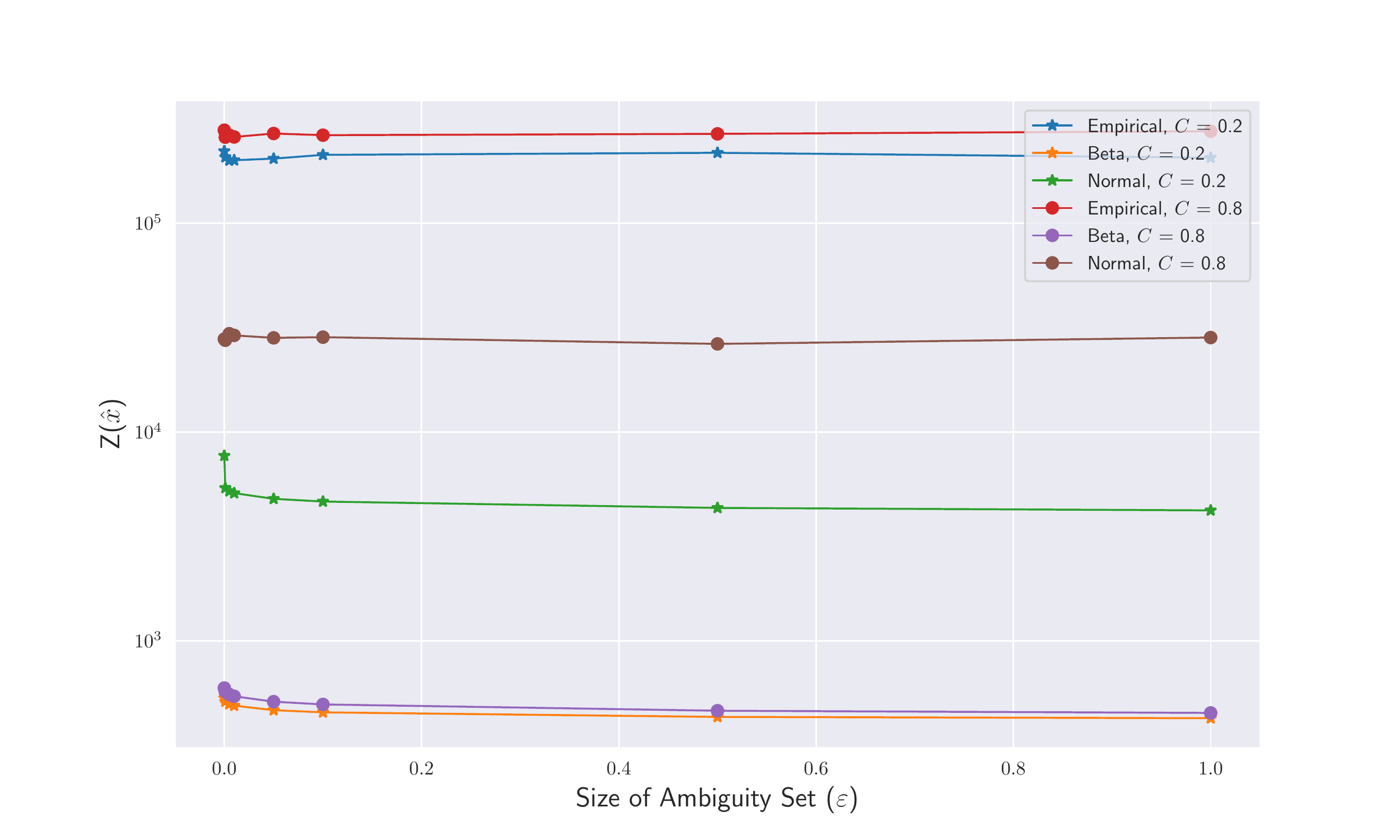} 
        \end{minipage}
    }%
    \subfloat[$(\gamma, n) = (4, 100)$] 
    {
        \begin{minipage}[t]{0.25\textwidth}
            \centering          
            \includegraphics[width = 0.9\textwidth, trim = 20 10 45 25 ]{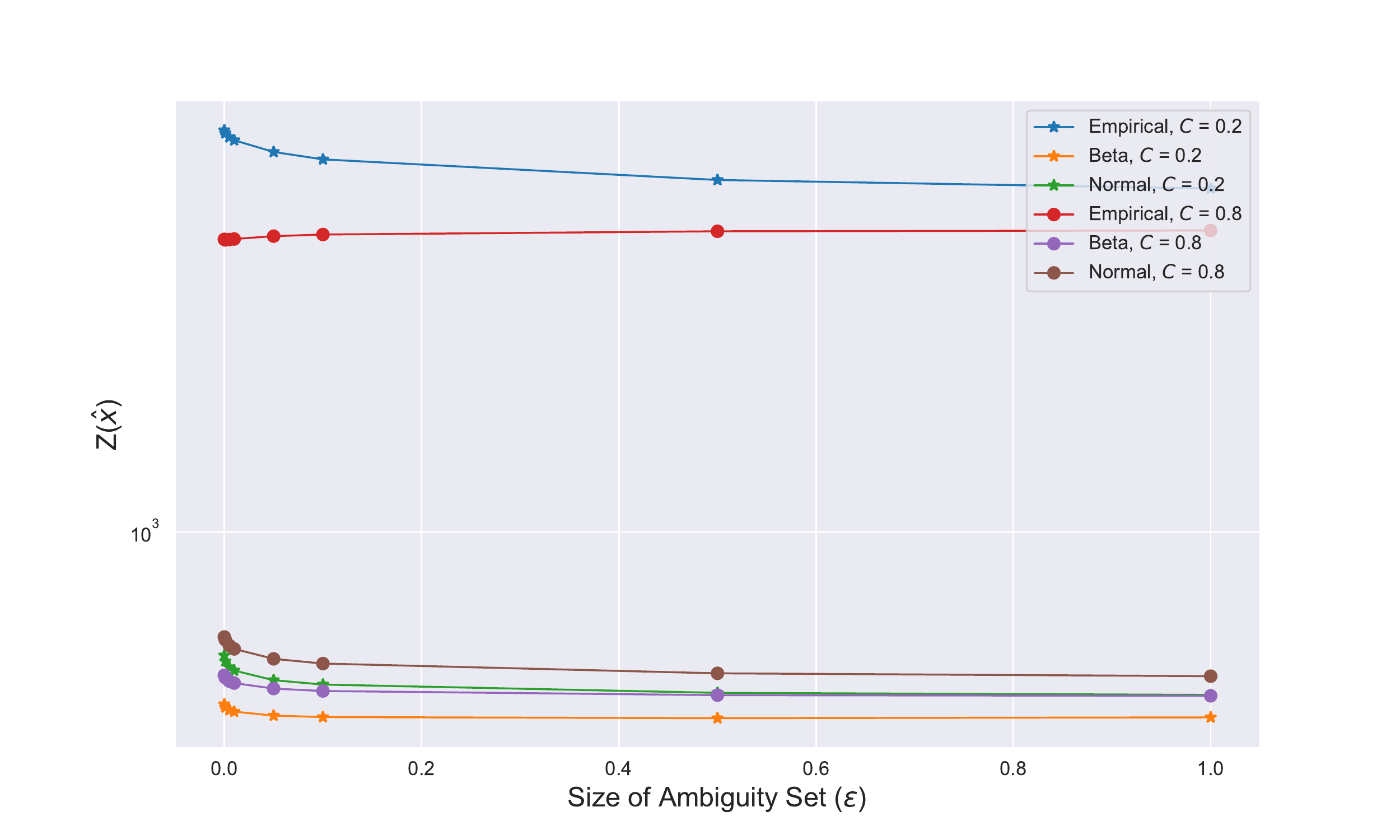} 
        \end{minipage}%
    }
    \subfloat[$(\gamma, n) = (4, 200)$]
    {
        \begin{minipage}[t]{0.25\textwidth}
            \centering      
            \includegraphics[width = 0.9\textwidth, trim = 20 10 45 25]{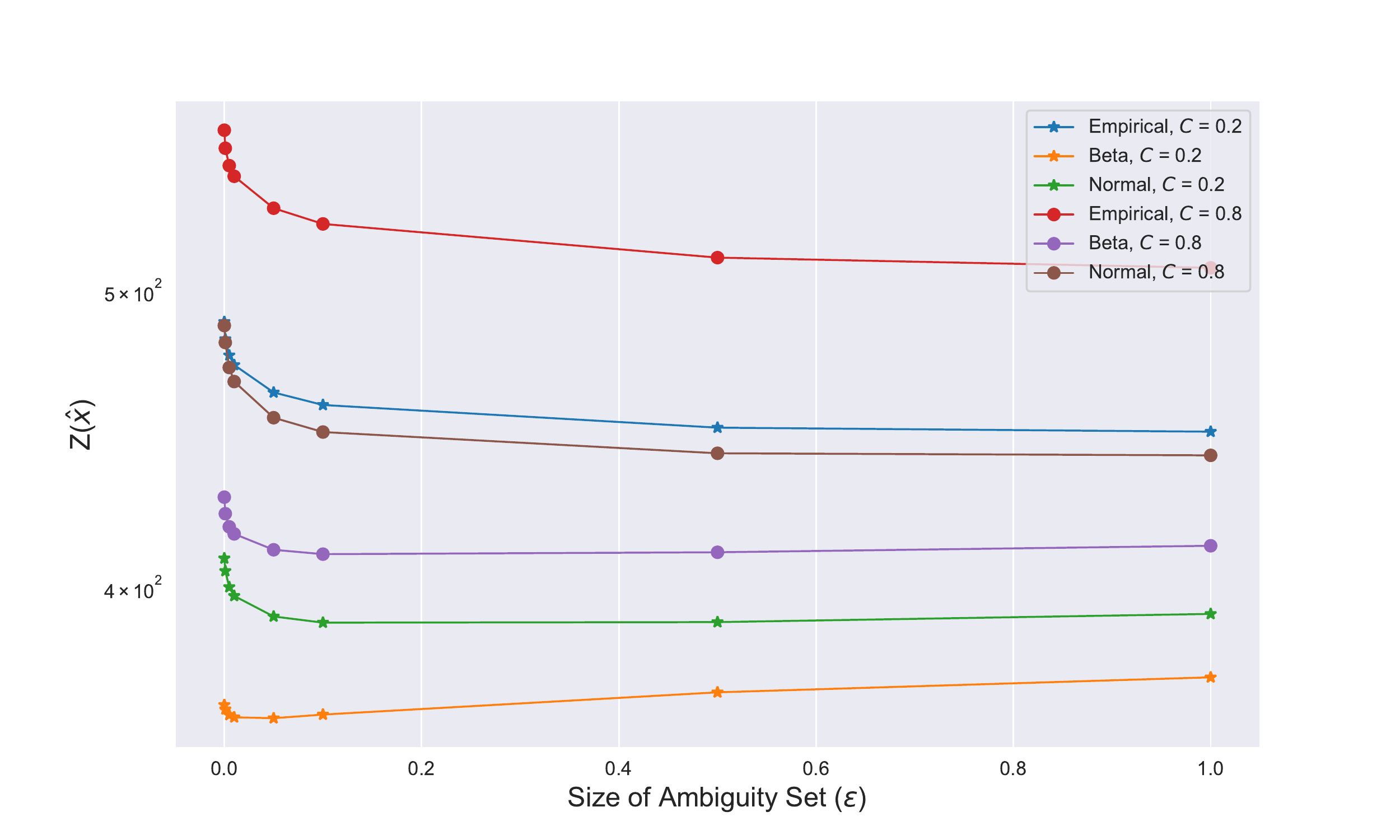}  
        \end{minipage}
    }%

    \caption{Value of cost function across different ERM-DRO models varying sample size $n$ and $\eta$ under distribution shiftS.}
    \label{fig:simu2-shift}   
\end{figure*}

\begin{figure*}
    \centering    
    \subfloat[Fully Parametrized] 
    {
        \begin{minipage}[t]{0.33\textwidth}
            \centering  
            \includegraphics[width=0.9\textwidth, trim = 20 10 45 45]{figs/concept/simu2_tempX_True.pdf} 
        \end{minipage}
    }%
    \subfloat[Distribution Misspecification] 
    {
        \begin{minipage}[t]{0.33\textwidth}
            \centering          
            \includegraphics[width = 0.9\textwidth, trim = 20 10 45 25 ]{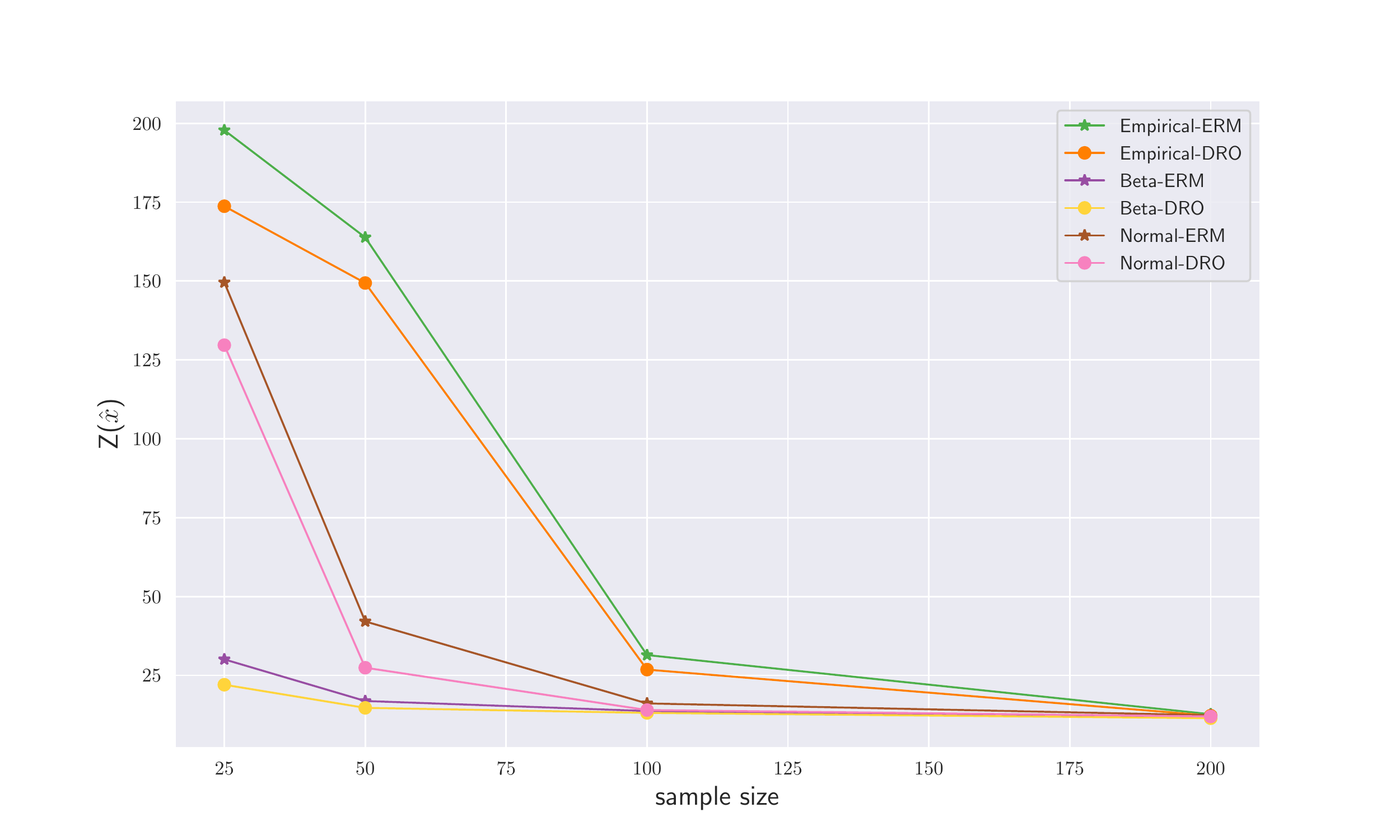} 
        \end{minipage}%
    }
    \subfloat[Distribution Shifts]
    {
        \begin{minipage}[t]{0.33\textwidth}
            \centering      
            \includegraphics[width = 0.9\textwidth, trim = 20 10 45 25]{figs/concept/simu0_tempX_Shift.pdf}  
        \end{minipage}
    }%
    \caption{Value of Cost function across different ERM-DRO models varying sample size $n$ with $(\tau, \eta) = (2,2)$.}
    \label{fig:simu2-concept}  
\end{figure*}

(3) \textit{Distribution Shifts.} We take the $i$-th marginal distribution of the train distribution $\xi_{i}^{tr}:=2r \times Beta(\eta_{i,1}, 2) - r$ and the test distribution
$\xi_{i}^{te}:=2r \times Beta(\eta_{i,2}, 2) - r,$
where $\eta_{i, 2} = \eta_{i,1} + C\min\{3 - \eta_{i,1}, \eta_{i,1} - 1.5\}$ with a shift parameter $C\in [-1, 1]$ which is randomly generated. 
We also perturb the $i$-th marginal distribution with uniform noise $\zeta_i \sim U(-2,2)$ to the test distribution as before. Figure~\ref{fig:simu2-shift} shows that \texttt{P-DRO} models have better performance than their \texttt{P-ERM} counterparts.
We summarize the results in Figure~\ref{fig:simu2-concept}. For each DRO method, we tune the best hyperparameter $\varepsilon \in \{0.001, 0.005, 0.01, 0.05, 0.1, 0.5, 1\}$ and average results over 50 independent runs. In the base case $(a)$ without distribution misspecification, the Beta-DRO model performs significantly better, especially under a small sample size. Although the performance gap between \texttt{P-ERM} and \texttt{P-DRO} is small under large sample size, the difference in the values of the cost function is still statistically significant with $p<0.001$. Similar results can be found in $(b)(c)$ of Figure~\ref{fig:simu2-concept} under the cases of distribution misspecification and distribution shifts.

\subsection{Detailed Setups in Section~\ref{subsec:real1}}\label{app:experiment-real-portfolio}
To obtain the out-of-sample performance from the empirical data, we apply the \textit{rolling-sample} approach from \cite{demiguel2009optimal} on the monthly data from July 1963 to December 2018 ($T = 666$) with an estimated window size of $M = 60$ months. The procedure of the \textit{rolling-sample} approach is as follows: to construct portfolios in month $t + M$ (from $t = 1$), we use the data from months $t$ to $t + M - 1$ as observed samples and solve the corresponding problem $\min_{x\in \Xscr}\hat{Z}(x)$ to obtain $\hat{x}$. Then we obtain the returns $\hat{r}_t = \xi_{t + M}^{\top}\hat{x}$ in month $t + M$. We repeat this procedure to construct portfolios in the following months by adding the next and dropping the earliest month until $t = T - M$. This gives us $T - M$ monthly out-of-sample returns $\{\hat{r}_{i}\}_{i = 1}^{T-M}$. Finally, we report the experimental results of different models by the empirical performance $\hat{h}:=\frac{1}{T-M}\sum_{i = 1}^{T-M}(\mu - \hat{r}_i)_+^2$.

We use $\chi^2$-divergence in our DRO models with cross-validation of the hyperparameter $\varepsilon \in \{0.2, 0.4, \ldots, 1.6, 1.8\}$ in each period. And we fit the observed samples with (1) Normal families (the same in Section~\ref{subsec:synthetic1} and~\ref{subsec:synthetic2}); (2) variants of Beta distributions, where we still fix $\eta_j = 2$, $\alpha_j$ using the formula in Section~\ref{subsec:synthetic2}, and choose the boundary parameter $r_j = \max{|(\xi)_j|}$ for each asset $j$.
\subsection{Detailed Setups in Section~\ref{subsec:real2}}\label{app:experiment-real-reg}
We vary the ambiguity size $\varepsilon \in \{0.01, 0.05, 0.1, 0.5, 1, 5, 10, 50, 100\}$ as the hyperparameter candidate set and tune $\varepsilon$ through cross-validation for DRO models. We set the Monte Carlo size $M = 10n$ and $\hat{\Q}$ is constructed as follows: 
\paragraph{Construction of the Gaussian Mixture Distribution Class.} Denote the feature vector $[x_1,\ldots,x_D] \in \R^D$ and among them, $x_1,\ldots,x_K$ are $K$ binary features with $x_i \in \{0, 1\}, \forall i \in [K]$. Then we consider the following Gaussian Mixture Distribution Class with $2^K$ groups, where $\overline{x_1x_2\ldots x_k}$ represents the decimal number of these binary digits $x_i$, e.g. $\overline{101} = 5$; $\Delta_n$ represents the $n$-dimension probability simplex: 
\begin{align*}
\Pscr_{\Theta} = &\left\{\P: (x_1,\ldots, x_D, y)^{\top}\sim \P: \P(\overline{x_1 x_2\ldots x_K} = k - 1) = p_s, \forall k \in [2^K], \right.\\
& \left.(x_{K + 1}, \ldots, x_D, y)|(x_1,\ldots, x_K)\sim \Nscr(\mu_k, \Sigma_k)\right.\\
&\left.\bigg| (p_1,\ldots, p_{2^K}) \in \Delta_{2^K}, \mu_k \in \R^{D-K+1}, \Sigma \in \mathbb{S}_{++}^{D-K+1}, \forall k \in [2^K]\right\}.    
\end{align*}
 For a dataset with $\{(\hat{\bm x}^i, \hat{y}^i)\}_{i \in [n]}$ with $\hat{x}_{j}^i$ denoting the $j$-th feature of the $i$-th sample, we output $\hat{\Q}$ parametrized by $\{(\hat{p}_k, \hat{\mu}_k, \hat{\Sigma}_k)\}_{k = 1}^{2^K}$:
\begin{align*}
    \hat{p}_k &= \frac{1}{n}\sum_{i = 1}^n\mathbb{I}_{\{\overline{\hat{x}_{1}^i\ldots \hat{x}_{K}^i} = k-1\}}\\
    \hat{\mu}_k &= \frac{1}{n \hat{p}_k}\sum_{i = 1}^n (\hat{x}_{ K+1}^{i},\ldots,\hat{x}_{D}^i, \hat{y}^{i})\mathbb{I}_{\{\overline{\hat{x}_{1}^i,\ldots \hat{x}_{K}^i} = k-1\}}\\
    \hat{\Sigma}_k &= \frac{1}{n \hat{p}_k}\sum_{i = 1}^n [(\hat{x}_{K+1}^i,\ldots,\hat{x}_{D}^i,\hat{y}^{i}) - \hat{\mu}_k][(\hat{x}_{K+1}^i,\ldots \hat{x}_{D}^i, \hat{y}^{i}) - \hat{\mu}_k]^{\top}\mathbb{I}_{\{\overline{x_{1}^i,\ldots x_{K}^i} = k-1\}}
\end{align*}
Since we have $K = 4$ binary features \textsf{black, Hispanic, married, non-degree}, we obtain a Gaussian mixture model with 16 subgroups with $\hat{p}_k$ being the empirical frequency. Then we estimate $\hat{\mu}_k, \hat{\Sigma}_k$ from each group for the other continuous variables \textsf{age, education, RE74, RE75, RE78} above. After that, in our Monte Carlo sampling, if there are few raw samples (<10) in this group, we directly use the original raw data within that group and copy each sample 10 times as the Monte Carlo sampling output for that group. Besides, we project the values of some features of the Monte Carlo data from $\hat{\Q}$ onto their value boundary if these values violate some rules: 1) If the value of the earnings (\textsf{RE74, RE75, RE78}) for one year in one sample is negative, we project that value to 0; 2) If the value of the \textsf{age} in one sample is too low or too high, we project that value onto the interval $[18, 60]$. 


\paragraph{Model Selection.}To illustrate that \PDRO\ can eliminate the model misspecification error and show consistently good performance over \PERM\ and the nonparametric models, we replace (A) Gaussian Mixture in $\Pscr_{\Theta}$ (i.e. the model in Section~\ref{subsec:real2}) to (B) joint Gaussian of all variables; (C) joint Gaussian fixing categorical variables zero correlation following the same setup. All models still have $\Escr_{apx} > 0$. However, Table~\ref{tab:model-selection} shows that \PDRO\ with different parametric models still enjoys superior performance, which indicates the robustness of our method.
\begin{center}

\begin{table}[htb]
\caption{Average model performance with different $\Pscr_{\Theta}$ without distribution shifts}
\label{tab:model-selection}
\begin{tabular}{c|c|c|c|c|c|c|c|c}
    $n = 200$ & \texttt{NP-ERM}&\texttt{NP-DRO}&\texttt{P-ERM}-(A)&\texttt{P-DRO}-(A)& \texttt{P-ERM}-(B)& \texttt{P-DRO}-(B) & \texttt{P-ERM}-(C)& \texttt{P-DRO}-(C) \\
    \hline
    Avg-$R^2$&0.1589&0.4433&$<0$&0.5050&$<0$&0.5266&0.4926&0.5033\\
\end{tabular}
\end{table}
\end{center}

\subsection{Another Numerical Example}\label{subsec:synthetic1}
We use a quadratic cost function with linear perturbations following from Section 5.2 in \cite{duchi2019variance}: $h(x;\xi) = \frac{1}{2}\|x-v\|^2 + \xi^{\top}(x- v).$ We set $D_{\xi} = 50$, the decision space $\Xscr = \{x\in \R^{D_{\xi}}: \|x\|_2 \leq B\}$ and $v = \frac{B}{2\sqrt{D_{\xi}}}\mathbf{1}$. 

We construct the $i$-th marginal distribution of the random variable $(\xi)_i = (\xi_{\theta})_i + (\tilde{\xi})_i,\forall i$, where $\xi_{\theta} \sim \Nscr(0, \Sigma)$, $(\tilde{\xi})_i\overset{d}{\sim} \text{Exp}(\lambda) - \frac{1}{\lambda}, \forall i \in [D_{\xi}]$. Intuitively, the smaller $\lambda$ is, the larger difference the distribution $(\xi)_i$ is compared to the normal $(\xi_{\theta})_i$. $\E_{\P^*}[(\xi)_i] = 0,\forall i \in [D_{\xi}]$ and $\Vscr_d(x^*) = 0$. 


We vary the decision boundary $B$ from $\{2, 10\}$, noise ratio $\lambda$ from $\{\frac{1}{5},\frac{1}{2}\}$. We consider DRO models with $d$ taken as $\chi^2$-divergence and 1-Wasserstein distance. We choose the parametric class as $\Pscr_{\Theta} = \left\{\Nscr(\mu, \Sigma):\mu \in \R^{D_{\xi}}, \Sigma \in \mathbb{S}_{++}^{D_{\xi}}\right\}$ with unknown $\mu$ and $\Sigma$. Then, we have $\mathcal{E}_{apx}(\P^*, \Theta) > 0$. We fit the distribution with $\hat{\Q}\sim \Nscr(\hat{\mu}, \hat{\Sigma})$, where $\hat{\mu} = \frac{1}{n}\sum_{i = 1}^{n}\hat{\xi}_i, \hat{\Sigma} = \frac{1}{n}\sum_{i = 1}^{n} (\hat{\xi}_i -\hat{\mu})(\hat{\xi}_i -\hat{\mu})^{\top}.$ We set each DRO model with ambiguity set sizes ranging from $\{0.1, 0.2, 0.5, 1, 2.5, 5, 10\}$ to show the trend.

\begin{figure*}[h] 
    \centering    
    \subfloat[$(n, B, \lambda) = (50, 2, \frac{1}{2})$] 
    {
        \begin{minipage}[t]{0.25\textwidth}
            \centering          
            \includegraphics[width=0.9\textwidth, trim = 20 10 45 45]{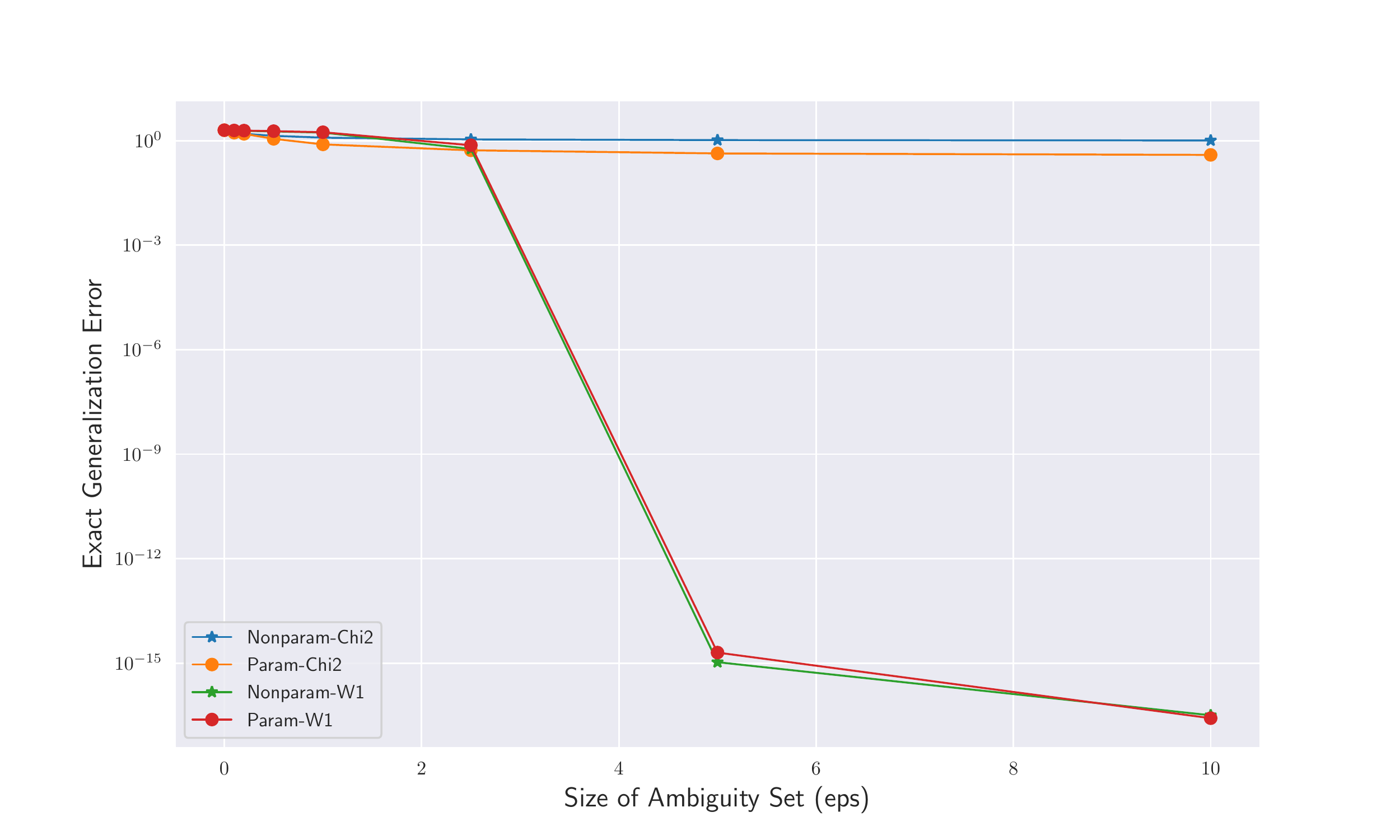}  
        \end{minipage}
    }
    \subfloat[$(n, B, \lambda) = (100, 2, \frac{1}{2})$] 
    {
        \begin{minipage}[t]{0.25\textwidth}
            \centering  
            \includegraphics[width=0.9\textwidth, trim = 20 10 45 45]{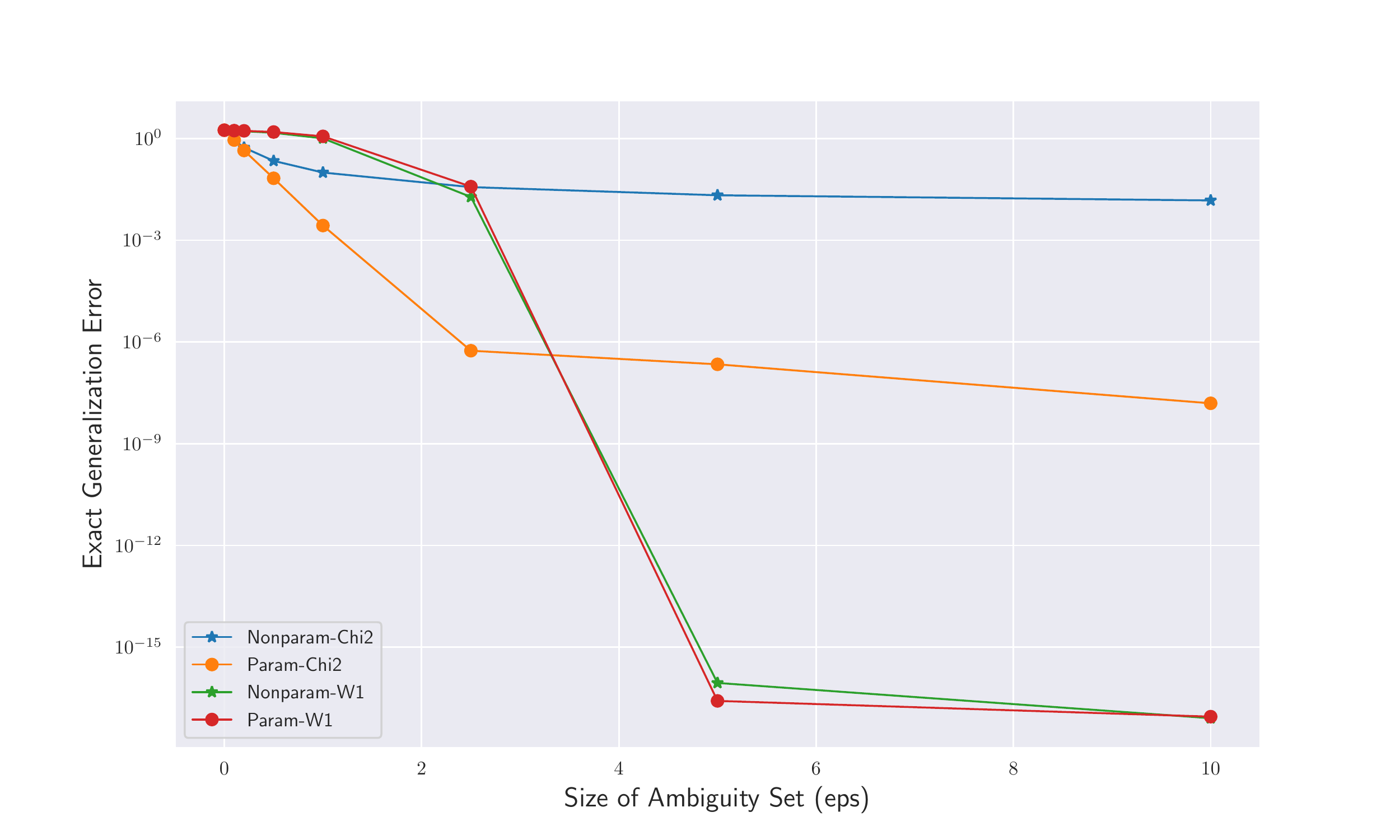} 
        \end{minipage}
    }%
    \subfloat[$(n, B, \lambda) = (150, 2, \frac{1}{2})$] 
    {
        \begin{minipage}[t]{0.25\textwidth}
            \centering          
            \includegraphics[width = 0.9\textwidth, trim = 20 10 45 25 ]{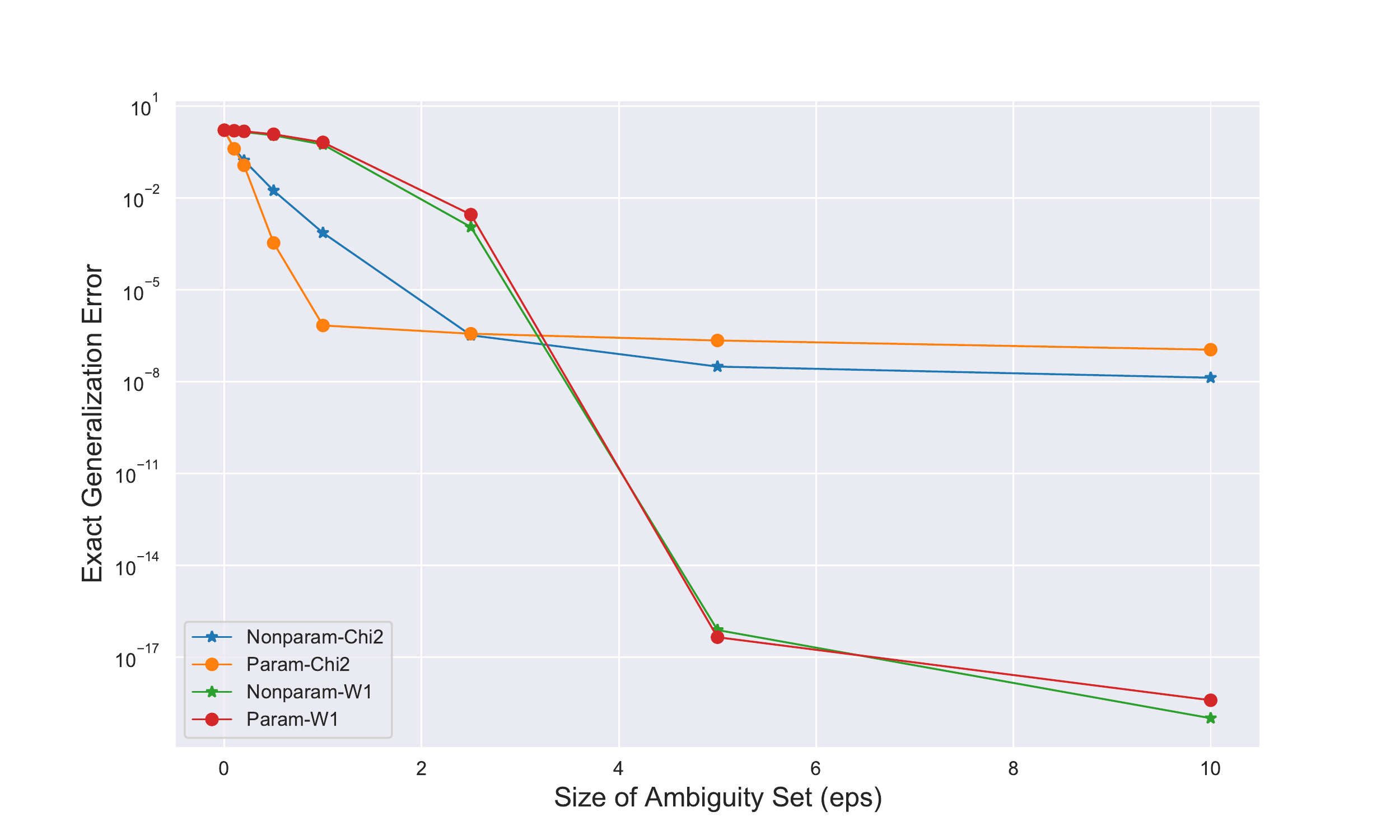} 
        \end{minipage}%
    }
    \subfloat[$(n, B, \lambda) = (200, 2, \frac{1}{2})$]
    {
        \begin{minipage}[t]{0.25\textwidth}
            \centering      
            \includegraphics[width = 0.9\textwidth, trim = 20 10 45 25]{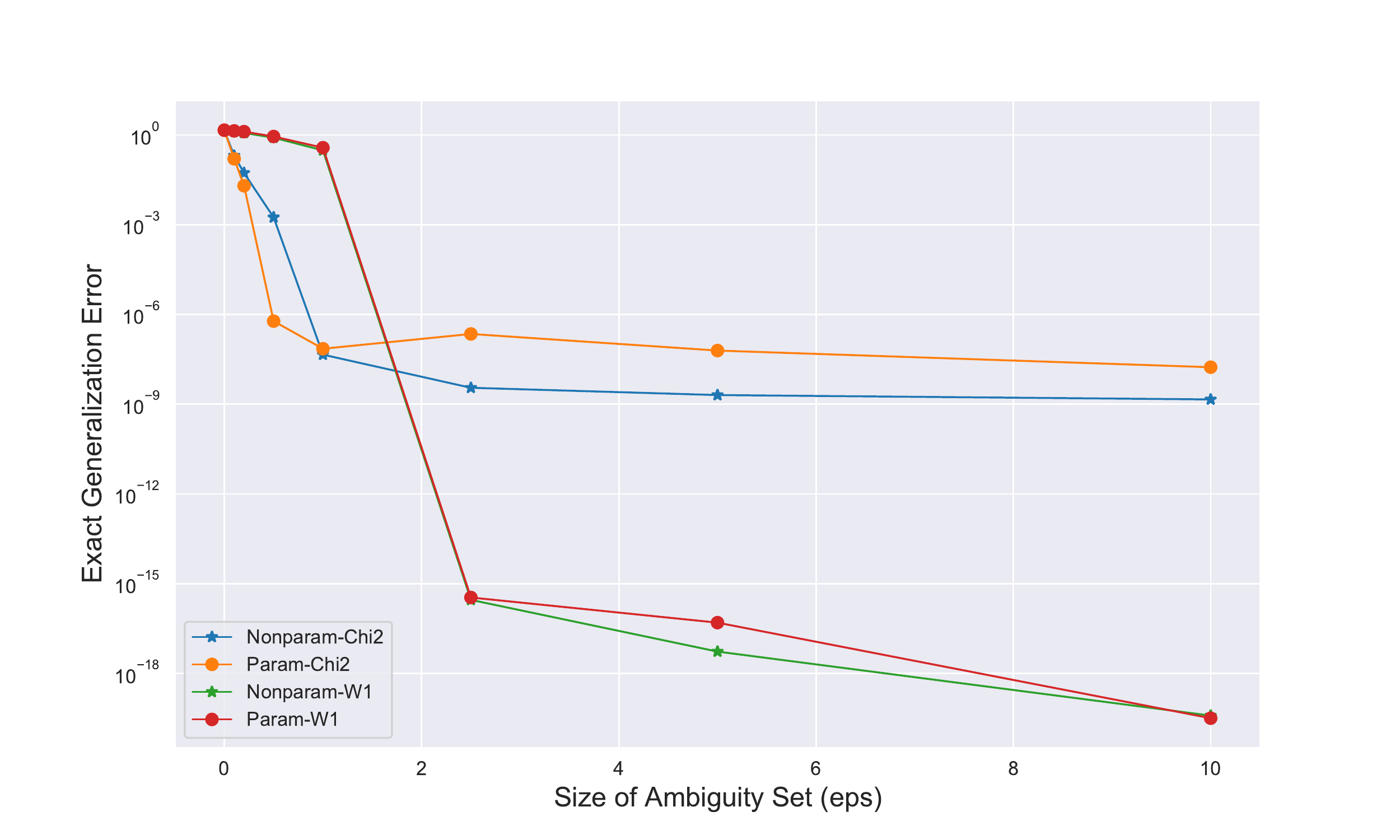}   
        \end{minipage}
    }%
    
        \subfloat[$(n, B, \lambda) = (50, 2, \frac{1}{5})$] 
    {
        \begin{minipage}[t]{0.25\textwidth}
            \centering          
            \includegraphics[width=0.9\textwidth, trim = 20 10 45 45]{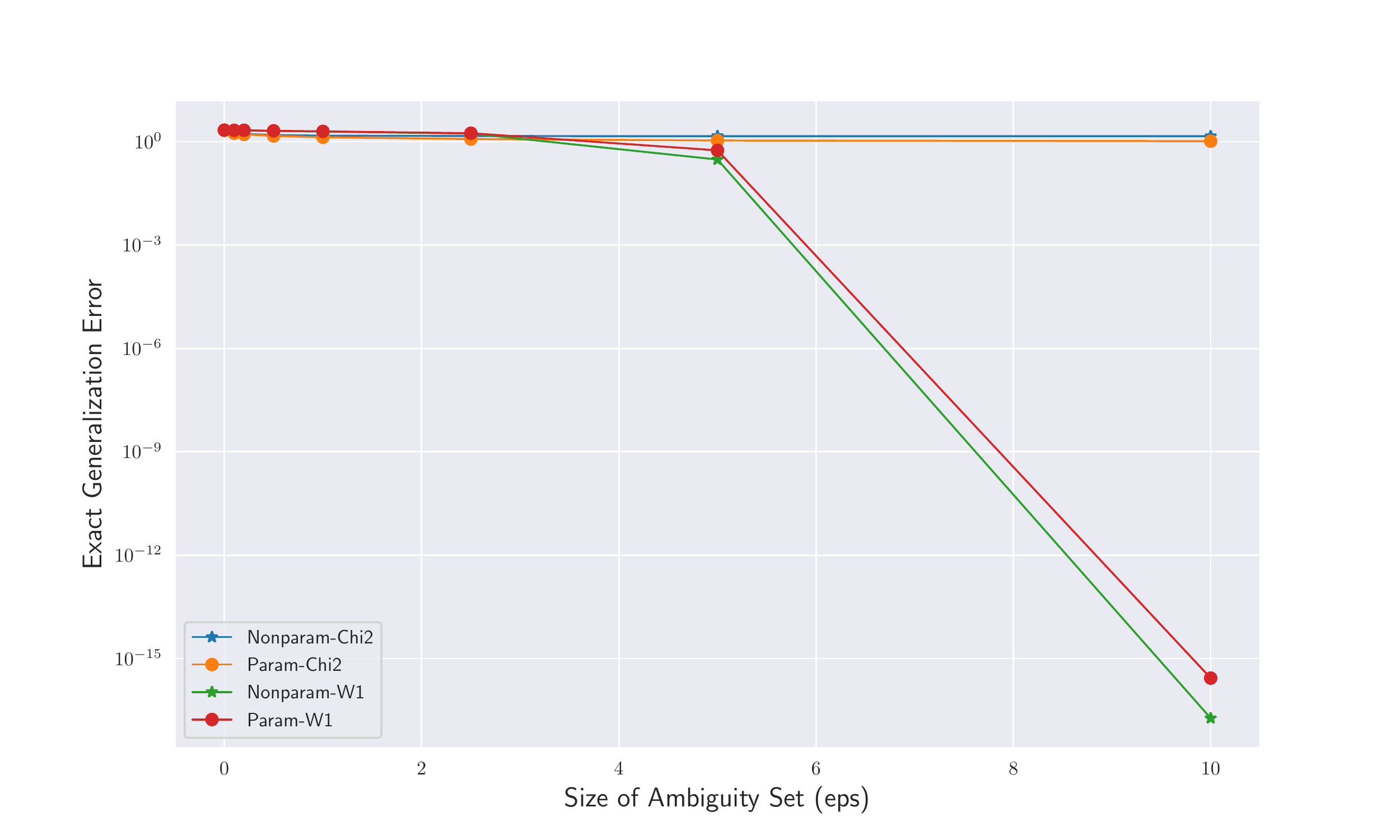}  
        \end{minipage}
    }
    \subfloat[$(n, B, \lambda) = (100, 2, \frac{1}{5})$] 
    {
        \begin{minipage}[t]{0.25\textwidth}
            \centering  
            \includegraphics[width=0.9\textwidth, trim = 20 10 45 45]{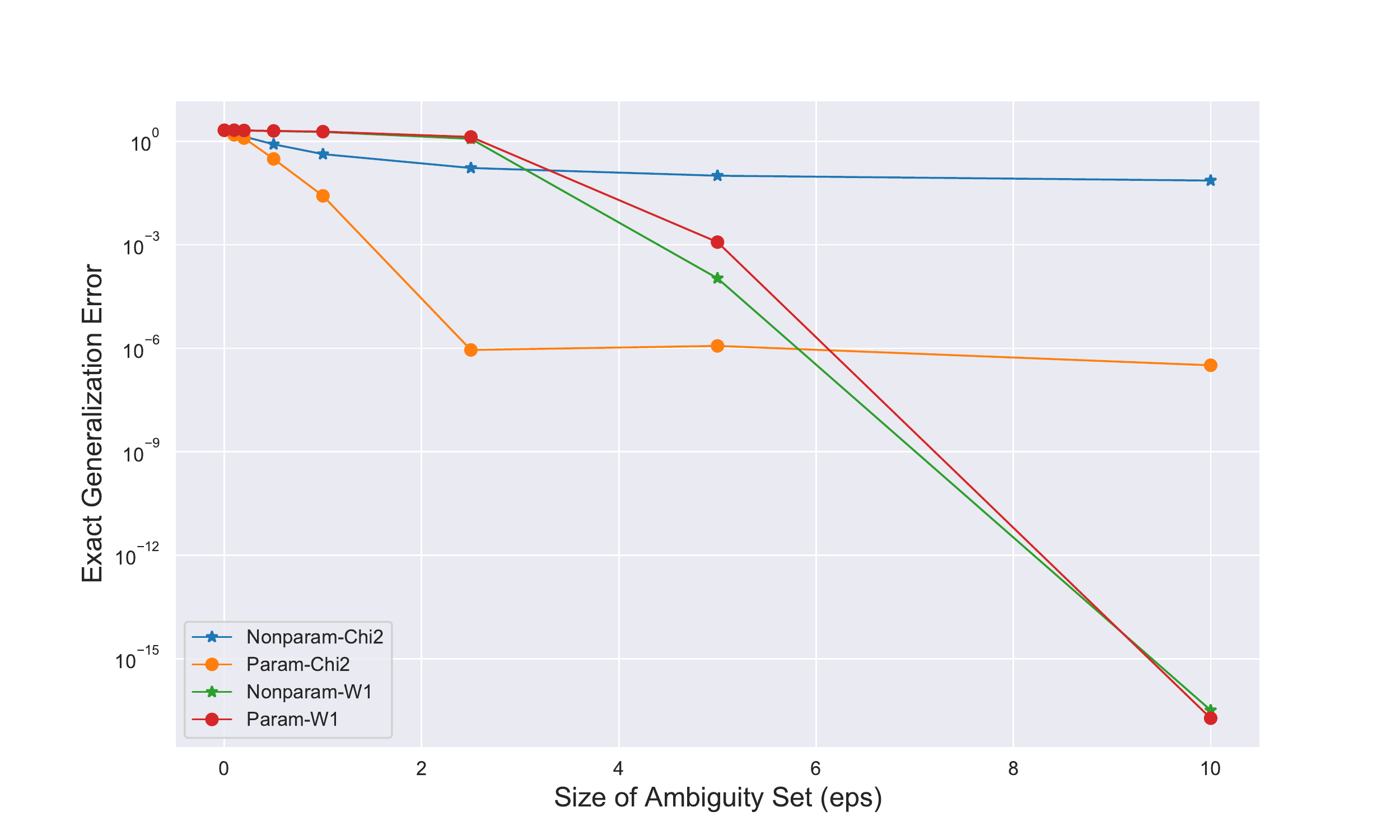} 
        \end{minipage}
    }%
    \subfloat[$(n, B, \lambda) = (150, 2, \frac{1}{5})$] 
    {
        \begin{minipage}[t]{0.25\textwidth}
            \centering          
            \includegraphics[width = 0.9\textwidth, trim = 20 10 45 25 ]{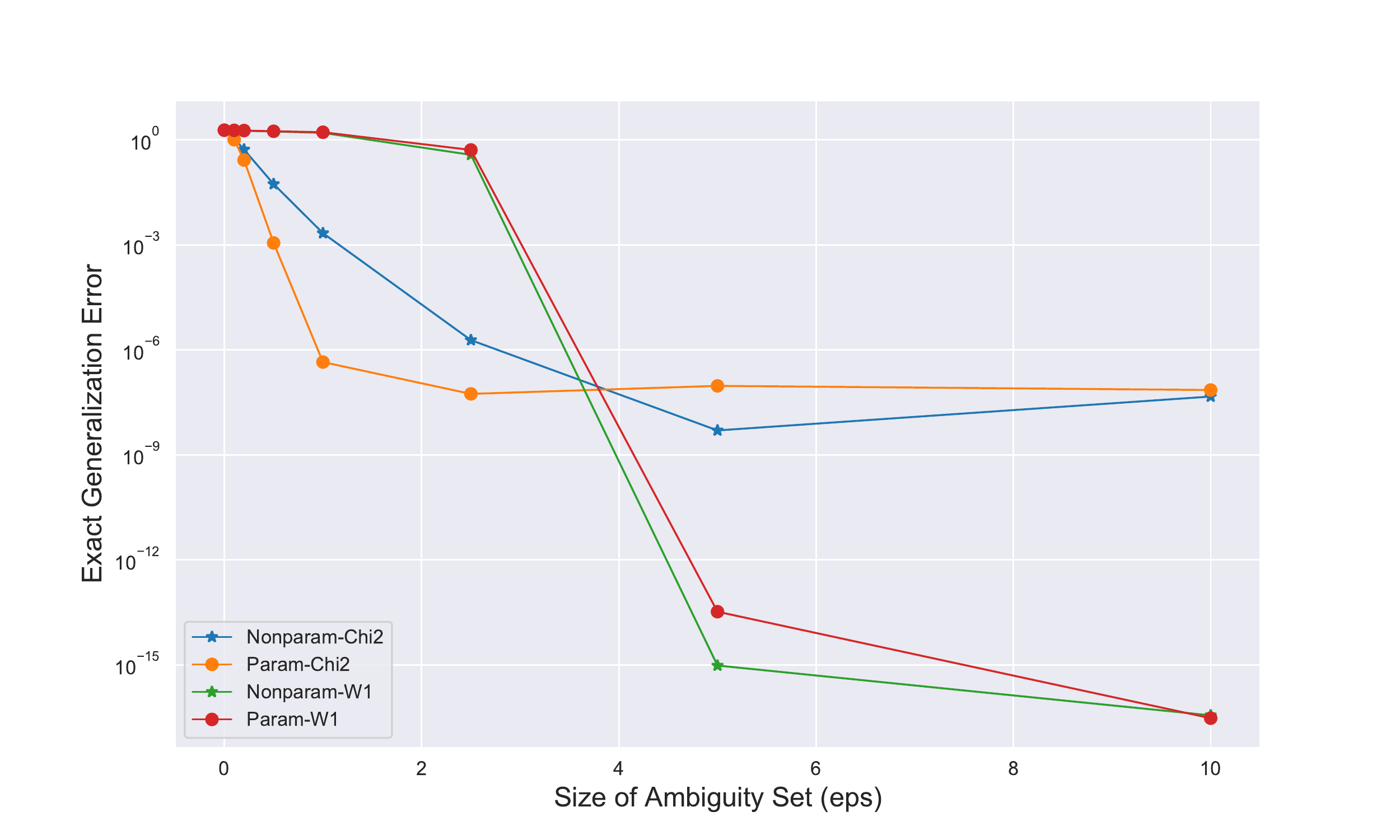} 
        \end{minipage}%
    }
    \subfloat[$(n, B, \lambda) = (200, 2, \frac{1}{5})$]
    {
        \begin{minipage}[t]{0.25\textwidth}
            \centering      
            \includegraphics[width = 0.9\textwidth, trim = 20 10 45 25]{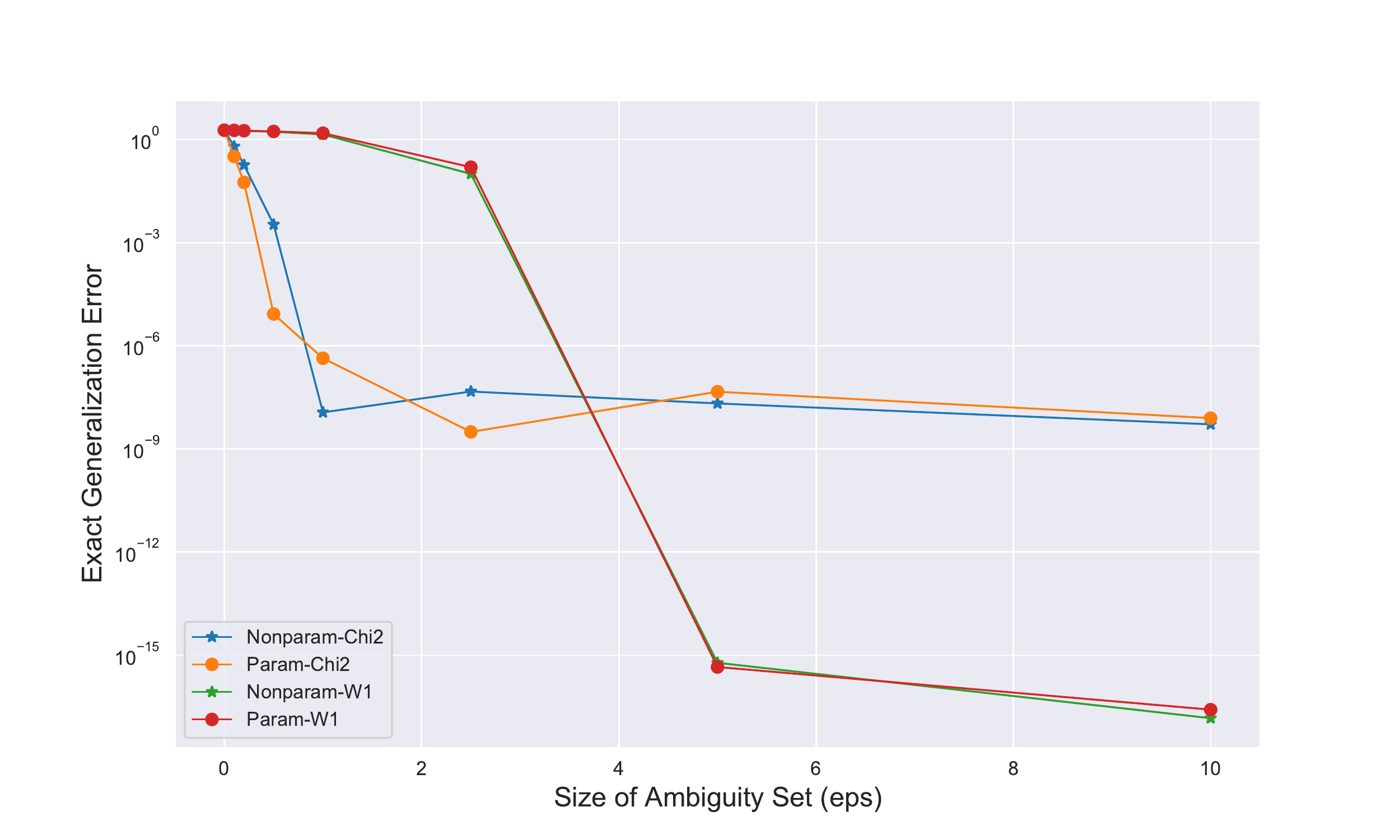}   
        \end{minipage}
    }%
    
    \subfloat[$(n, B, \lambda) = (50, 10, \frac{1}{2})$] 
    {
        \begin{minipage}[t]{0.25\textwidth}
            \centering          
            \includegraphics[width=0.9\textwidth, trim = 20 10 45 45]{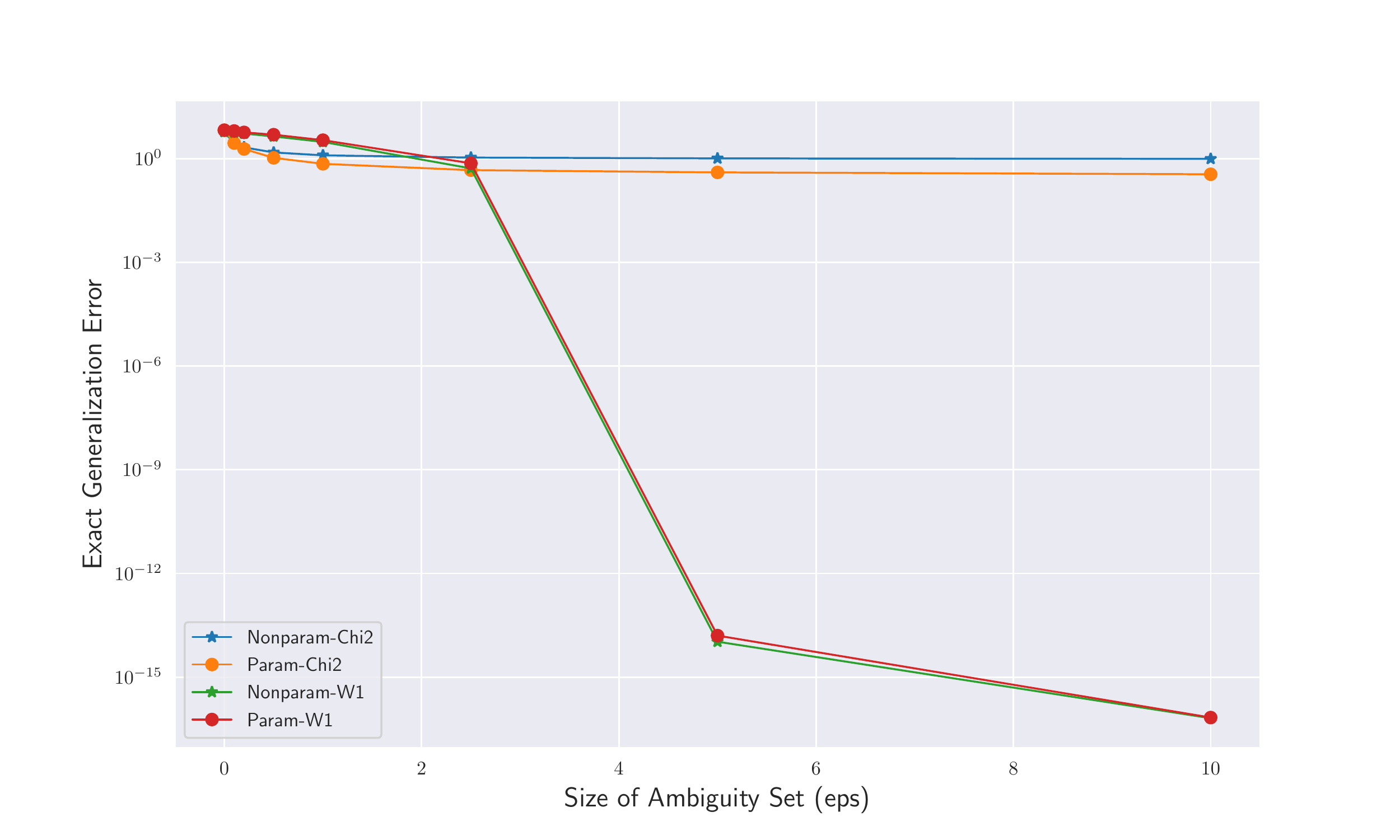}  
        \end{minipage}
    }
    \subfloat[$(n, B, \lambda) = (100, 10, \frac{1}{2})$] 
    {
        \begin{minipage}[t]{0.25\textwidth}
            \centering  
            \includegraphics[width=0.9\textwidth, trim = 20 10 45 45]{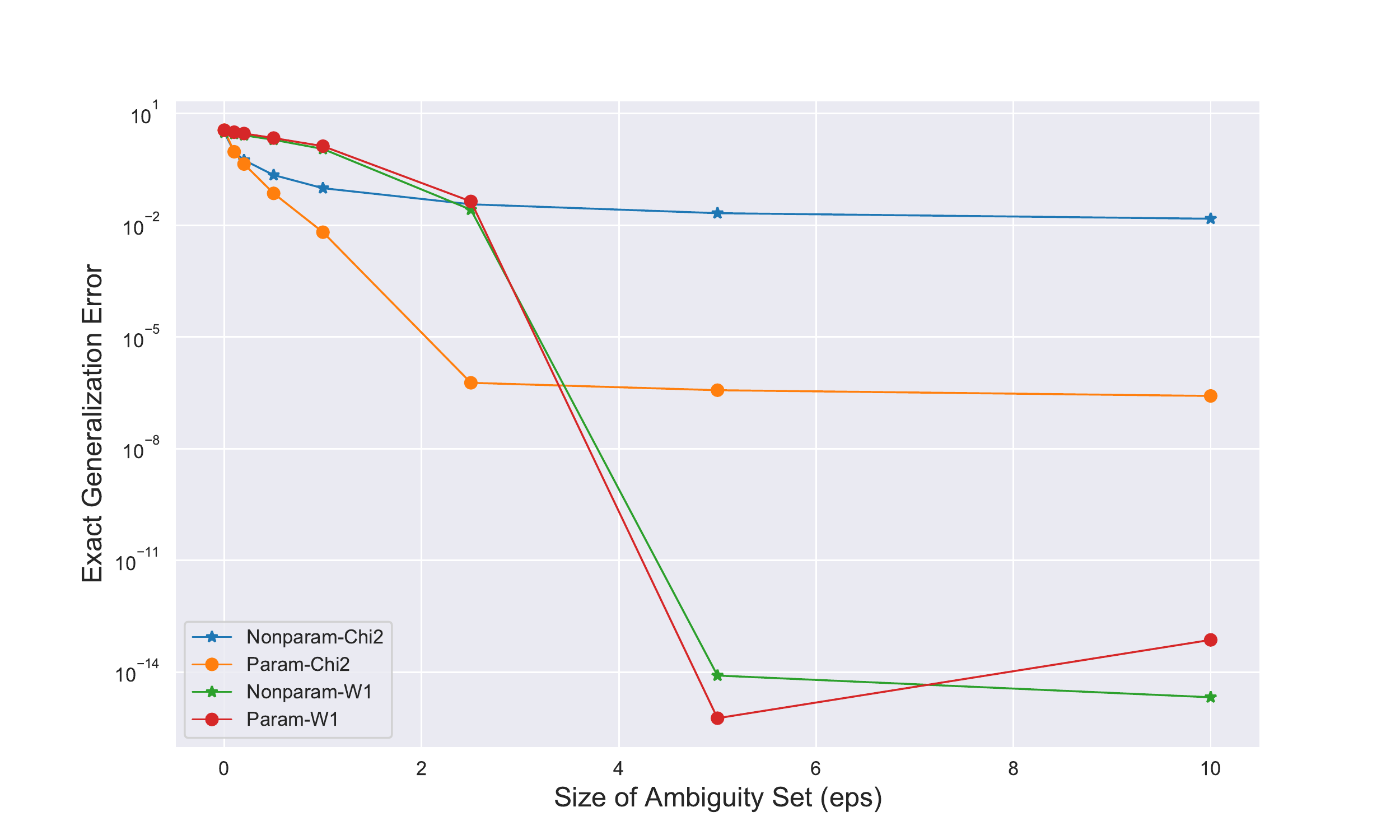} 
        \end{minipage}
    }%
    \subfloat[$(n, B, \lambda) = (150, 10, \frac{1}{2})$] 
    {
        \begin{minipage}[t]{0.25\textwidth}
            \centering          
            \includegraphics[width = 0.9\textwidth, trim = 20 10 45 25 ]{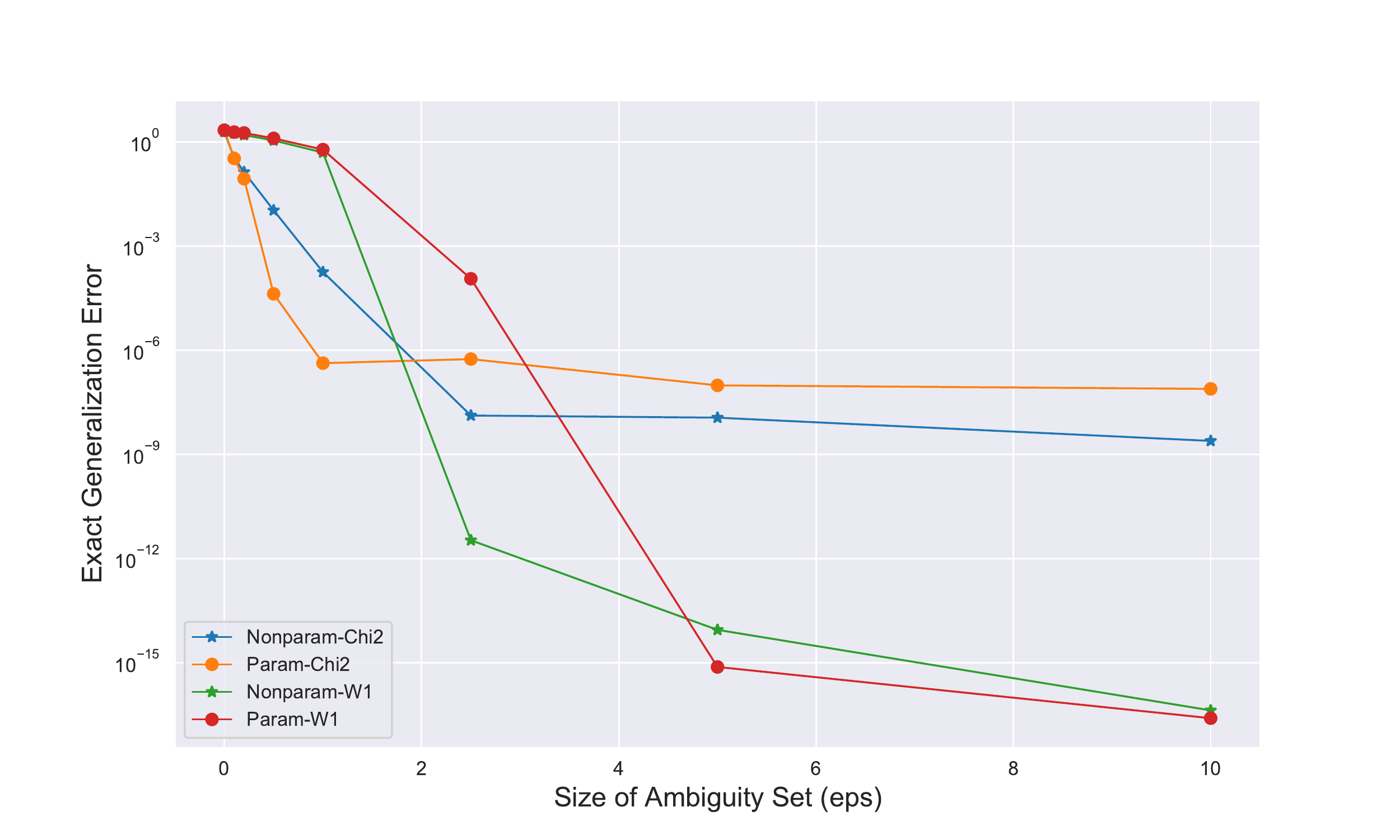} 
        \end{minipage}%
    }
    \subfloat[$(n, B, \lambda) = (200, 10, \frac{1}{2})$]
    {
        \begin{minipage}[t]{0.25\textwidth}
            \centering      
            \includegraphics[width = 0.9\textwidth, trim = 20 10 45 25]{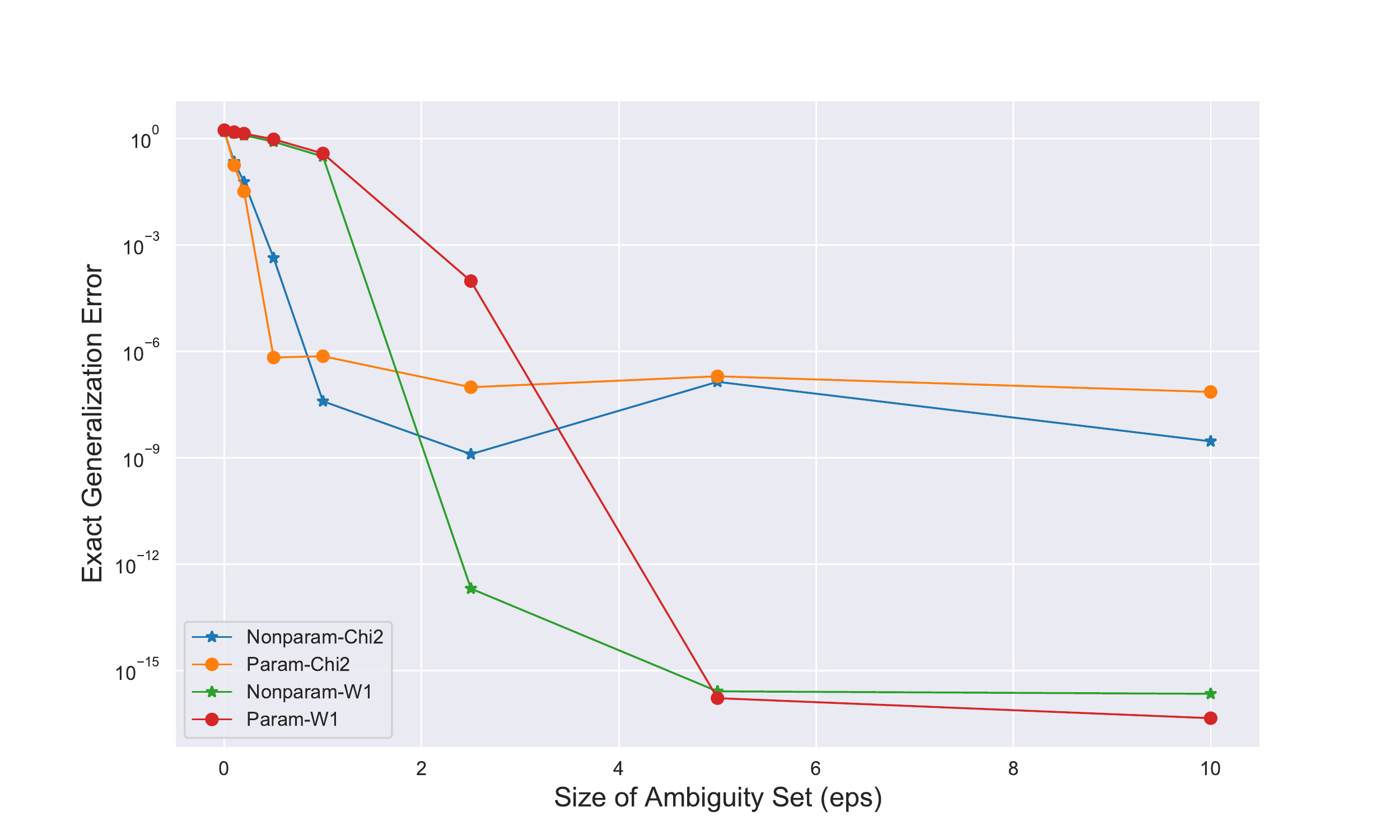}   
        \end{minipage}
    }%
    
        \subfloat[$(n, B, \lambda) = (50, 10, \frac{1}{5})$] 
    {
        \begin{minipage}[t]{0.25\textwidth}
            \centering          
            \includegraphics[width=0.9\textwidth, trim = 20 10 45 45]{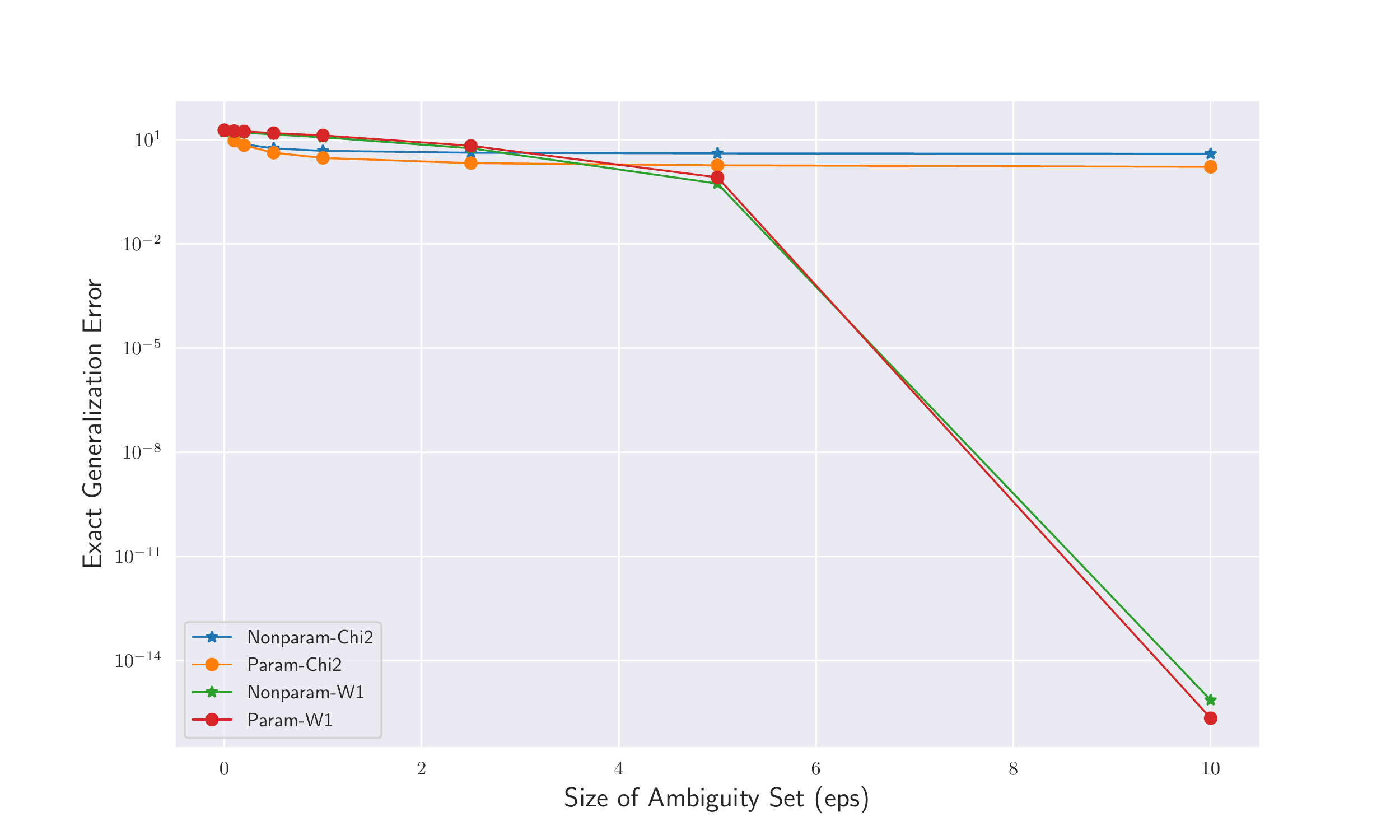}  
        \end{minipage}
    }
    \subfloat[$(n, B, \lambda) = (100, 10, \frac{1}{5})$] 
    {
        \begin{minipage}[t]{0.25\textwidth}
            \centering  
            \includegraphics[width=0.9\textwidth, trim = 20 10 45 45]{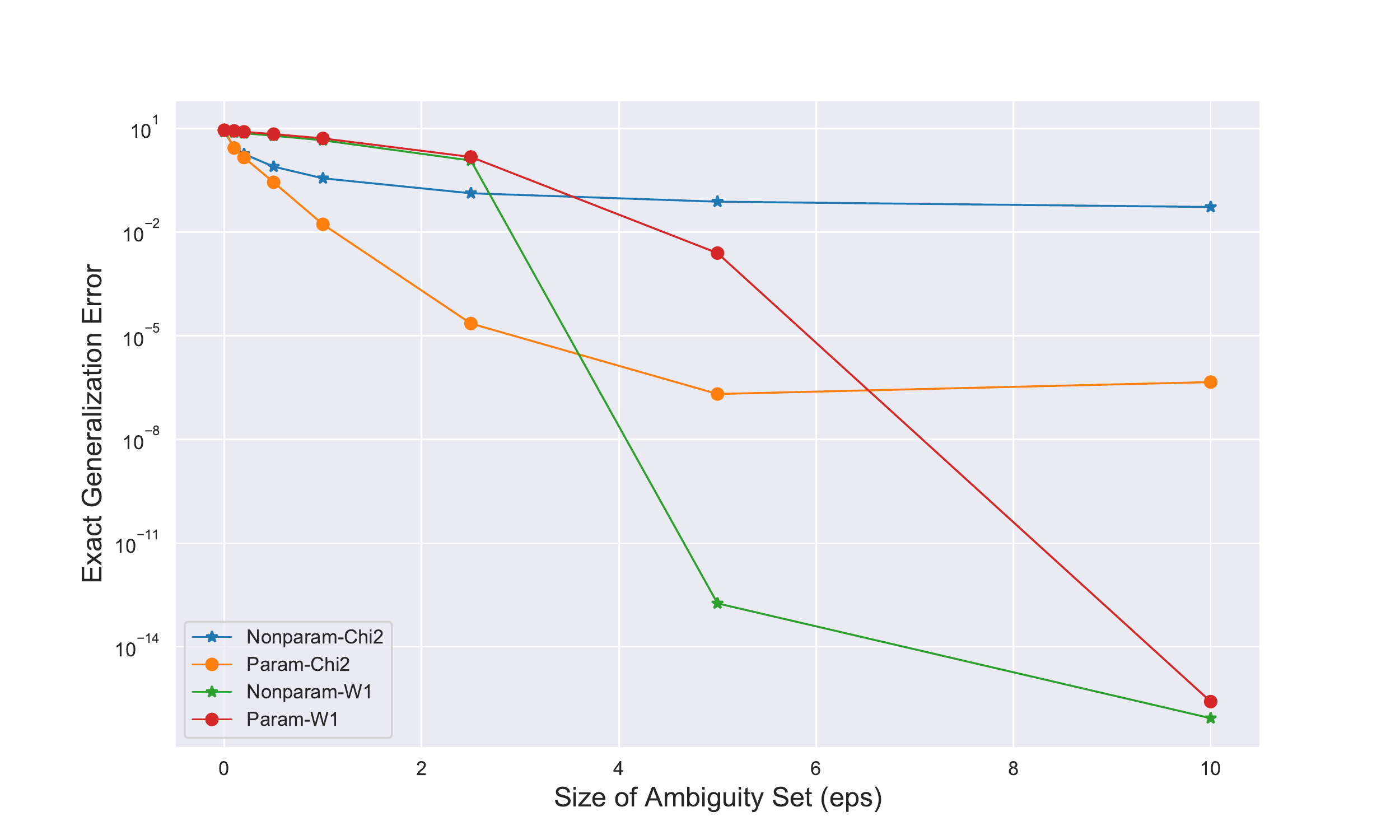} 
        \end{minipage}
    }%
    \subfloat[$(n, B, \lambda) = (150, 10, \frac{1}{5})$] 
    {
        \begin{minipage}[t]{0.25\textwidth}
            \centering          
            \includegraphics[width = 0.9\textwidth, trim = 20 10 45 25 ]{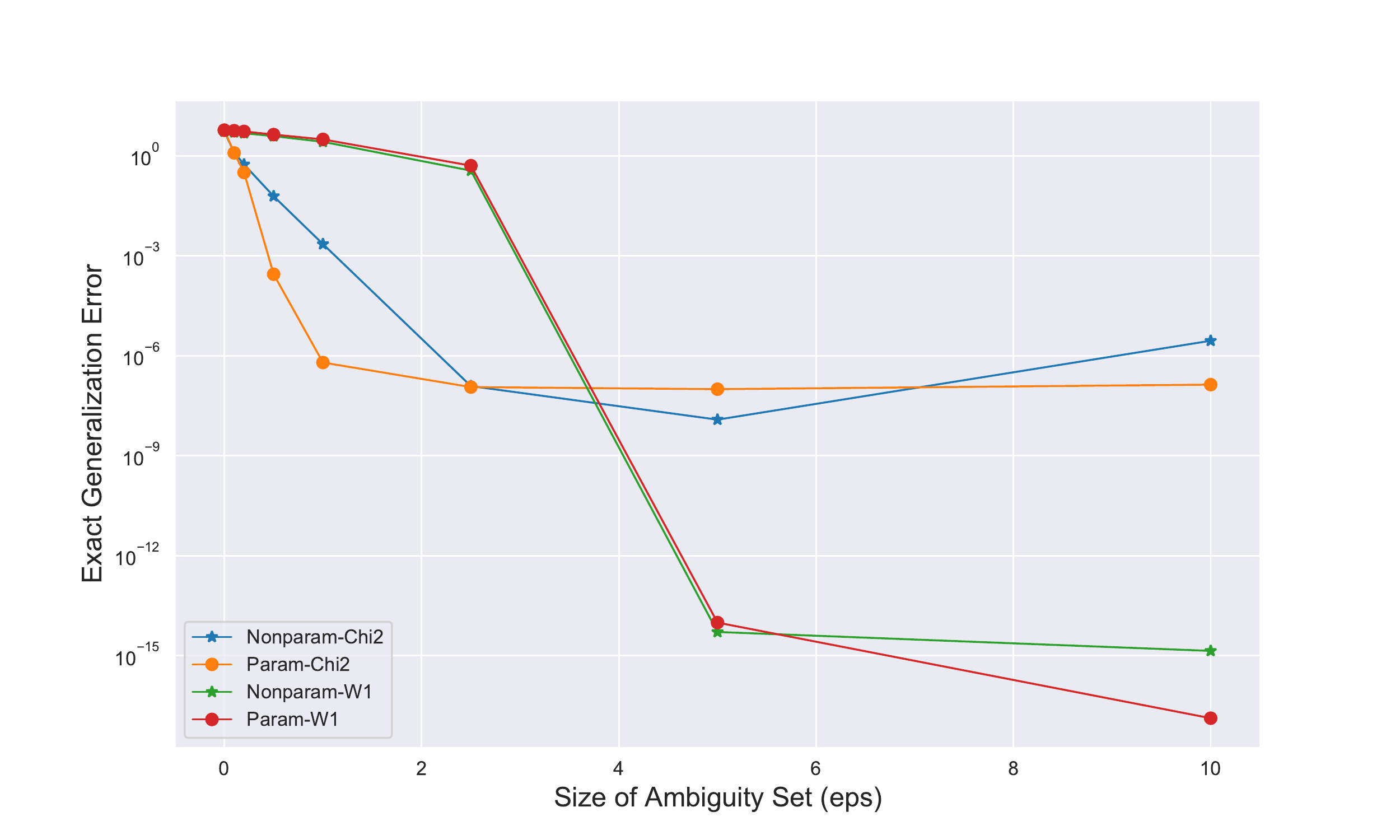} 
        \end{minipage}%
    }
    \subfloat[$(n, B, \lambda) = (200, 10, \frac{1}{5})$]
    {
        \begin{minipage}[t]{0.25\textwidth}
            \centering      
            \includegraphics[width = 0.9\textwidth, trim = 20 10 45 25]{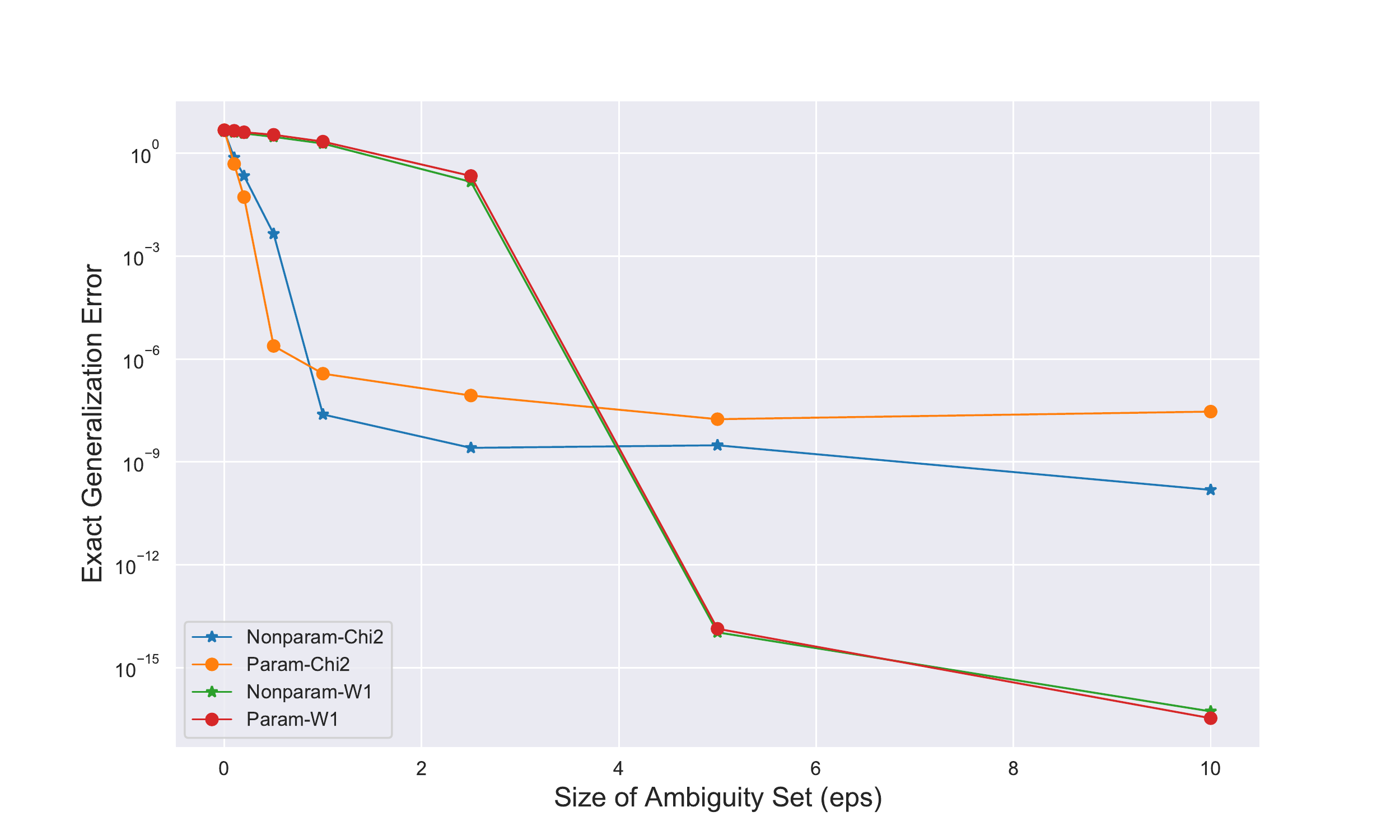}   
        \end{minipage}
    }%
    \caption{Exact generalization errors of DRO models varying sample size $n$, decision boundary $B$ and noise ratio $\lambda$}
    \label{fig:simu-quad-2-2}  
\end{figure*}

Figure~\ref{fig:simu-quad-2-2} shows different subcases varying $(n, B, \lambda)$ across 50 independent runs.
Across all these setups, since the optimal solution (i.e. $x^* = v$) does not depend on the decision boundary chosen here, increasing $B$ does not greatly affect the decision quality. \texttt{P-DRO} performs competitively against \texttt{NP-DRO} in all setups under the same ambiguity size $\varepsilon$. Both 1-Wasserstein and $\chi^2$-DRO outperform ERM models to a great extent with a large ambiguity size. This can be explained by the error bound that replaces the term $\sup_{x\in \Xscr}\Vscr_d(x)$ with $\Vscr_d(x^*)$.
In this case, \texttt{P-DRO} still keeps the property of \texttt{NP-DRO} in eliminating the dependence of the complexity term only to $\Vscr_d(x^*)$. \PDRO\ outperforms \PERM\ significantly and achieves almost zero generalization errors under a large ambiguity size $\varepsilon$ across 1-Wasserstein distance and $\chi^2$-divergence, which demonstrates the correctness of Theorem~\ref{thm:general-param-dro-ipm}.